\documentclass[11pt]{article}

\usepackage{latexsym}
\usepackage{amssymb}
\usepackage{amsmath}
\usepackage{amsfonts}
\usepackage{mathrsfs}

\usepackage{epsfig}
\usepackage{subfigure}

\usepackage{color}

\usepackage{theorem}
\newtheorem{theorem}{Theorem}[section]
\newtheorem{dfn}[theorem]{Definition}
{\theorembodyfont{\rmfamily} 
\newtheorem{rmk}[theorem]{Remark}}
\newtheorem{corollary}[theorem]{Corollary}
\newtheorem{proposition}[theorem]{Proposition}
\newtheorem{lemma}[theorem]{Lemma}

\newcommand{\C}{\mathbb{C}}
\newcommand{\B}{\mathscr{B}}
\newcommand{\Bn}{\mathscr{B}}
\newcommand{\G}{\mathscr{G}}
\newcommand{\F}{\mathscr{F}}
\newcommand{\Rc}{\mathscr{R}_1}
\newcommand{\Sc}{\mathscr{S}_1}
\newcommand{\R}{\mathbb{R}}
\newcommand{\N}{\mathbb{N}}
\newcommand{\Z}{\mathbb{Z}}
\newcommand{\Q}{\mathbb{Q}}
\renewcommand{\k}{k}
\renewcommand{\l}{\mathbb{L}}

\newcommand{\n}{\mathbf{n}}
\newcommand{\bu}{\mathbf{u}}
\newcommand{\bv}{\mathbf{v}}
\newcommand{\bw}{\mathbf{w}}
\newcommand{\z}{\mathbf{z}}
\newcommand{\p}{\mathbf{p}}
\newcommand{\q}{\mathbf{q}}
\newcommand{\vv}{\mathbf{v}}
\newcommand{\ww}{\mathbf{w}}

\def\ch#1{{{\bf H}^{#1}_{\C}}}
\def\chb#1{{{\bf \overline{H}}^{#1}_{\C}}}
\def\cp#1{{{\bf P}^{#1}_{\C}}}
\def\rh#1{{{\bf H}^{#1}_{\R}}}

\newcommand{\Pf}{{\sc Proof}. }
\newcommand{\EPf}{\hbox{}\hfill$\Box$\vspace{.5cm}}

\title{New non-arithmetic complex hyperbolic lattices}
\author{Martin Deraux, John R.~Parker, Julien Paupert}

\date{May 15, 2015}

\usepackage{a4wide}

\begin{document}

\maketitle
\begin{abstract}
 We produce a family of new, non-arithmetic lattices in ${\rm
   PU}(2,1)$. All previously known examples were commensurable with
 lattices constructed by Picard, Mostow, and Deligne--Mostow, and fell
 into 9 commensurability classes. Our groups produce 5 new distinct
 commensurability classes. Most of the techniques are completely
 general, and provide efficient geometric and computational tools for
 constructing fundamental domains for discrete groups acting on the
 complex hyperbolic plane.
\end{abstract}

\section{Introduction}

The general context of this paper is the study of lattices in
semisimple Lie groups and their classification. In what follows, $X$
denotes an irreducible symmetric space of non-compact type, and $G$
denotes its isometry group. Recall that $X$ is a homogeneous space
$G/K$ with $K$ a maximal compact subgroup of $G$, and up to scale, it carries
a unique $G$-invariant Riemannian metric.

A subgroup $\Gamma\subset G$ is called a lattice if $\Gamma\setminus
G$ has finite Haar measure (equivalently, if $\Gamma\setminus X$ has
finite Riemannian volume for the invariant metric on
$X$). Mostow-Prasad rigidity says that, for most symmetric spaces, a
given lattice $\Gamma$ in $G$ admits a unique discrete faithful
representation into $G$ (unique in the sense that all such
representation are conjugate in $G$). In order for this to hold, the
only case to exclude is $X=\rh 2$ (equivalently $G={\rm PSL}_2(\R)$),
where lattices are known to admit pairwise non-conjugate continuous
deformations.

The discreteness assumption in the statement of the Mostow-Prasad
rigidity theorem cannot be removed, since Galois automorphisms
sometimes produce non-discrete representations that are clearly
faithful and type-preserving (see the non-standard homomorphisms
constructed in~\cite{mostowpacific}, and also
Section~\ref{sec:arithmetic} of the present paper). This is related to
the notion of arithmeticity of lattices, which we now briefly recall.

There is a general construction of lattices, obtained by taking
integer matrices in the set of real points of linear algebraic groups
defined over $\mathbb{Q}$; this is a special case of the more general
notion of arithmetic group, where one considers lattices commensurable
with the image of an integral group under surjective homomorphisms
with compact kernel (see~\cite{margulisBook} for instance). Recall for
further reference that subgroups $\Gamma_1$ and $\Gamma_2$ of a group
$G$ are commensurable if there is a $g\in G$ such that
$\Gamma_1\cap g\Gamma_2 g^{-1}$ has finite index in both
$\Gamma_1$ and $g\Gamma_2 g^{-1}$.

Arithmetic groups are classified, and they can be understood (up to
commensurability) from purely number-theoretical data
(see~\cite{weilinvolutions},~\cite{tits}). Indeed, it follows from
rigidity that any lattice in a simple real Lie group not isomorphic to
${\rm PSL}_2(\R)$ is defined over a number field, and semisimple
groups over arbitrary number fields can be classified in a way similar
to the well known classification of semisimple groups over
$\mathbb{Q}$ (one uses Dynkin diagrams with a symmetry corresponding
to the action of the Galois group of the number field). The
arithmeticity condition says that the lattice is given up to
commensurability by taking the corresponding $\Z$-points in the group,
hence there are only countably many arithmetic lattices (since there
are only countably many number fields).

It follows from work of Margulis~\cite{margulis}, together with work
of Corlette~\cite{corlette} and Gromov-Schoen~\cite{gromovschoen},
that most lattices are arithmetic; given the discussion in the
previous paragraph, this means that there is a satisfactory
classification of lattices in \emph{most} semisimple Lie groups.
More specifically, if the isometry group $G$ of a symmetric space
$X=G/K$ of non-compact type contains a non-arithmetic lattice, then
$X$ can only be $\rh n$ or $\ch n$ for some $n$; in other words, up to
index 2, $G$ can only be ${\rm PO}(n,1)$ or ${\rm PU}(n,1)$.

Real hyperbolic space $\rh n$ is the model $n$-dimensional negatively
curved space form, i.e. it is the unique complete simply connected
Riemannian manifold of constant negative curvature (we may assume that
constant is $-1$). Similary, $\ch n$ is the model negatively curved
\emph{complex} space form, in the sense it is the unique complete,
simply connected K\"ahler manifold whose holomorphic sectional
curvature is a negative constant, which we may assume is $-1$.  With
that normalization, the real sectional curvatures of $\ch n$ vary
between $-1$ and $-1/4$. This makes complex hyperbolic spaces very
different from real hyperbolic ones (apart from the coincidence that
$\ch 1$ is isometric to $\rh 2$).

In the real hyperbolic case, the hybridization procedure described by
Gromov and Piatetski-Shapiro~\cite{gps} allows the construction of
(infinitely many commensurability classes of) non-arithmetic lattices
in ${\rm PO}(n,1)$ for any $n\geqslant 2$. In ${\rm PU}(n,1)$, it is
not clear how to make sense of the construction of such hybrids, and
the construction of non-arithmetic lattices in the complex hyperbolic
setting is a longstanding challenge.

In fact, for $n\geqslant 2$, only a handful of examples are known. The first ones
were constructed by Mostow~\cite{mostowpacific}, who showed that some
complex reflection groups are non-arithmetic lattices in ${\rm
  Isom}(\ch 2)$ (the analogous statement was well known for real
reflection groups in ${\rm Isom}(\rh n)$ for low values of $n$, see~\cite{vinbergcriterion}).

 Complex hyperbolic manifolds are of special interest in complex
 geometry, because of the extremal properties of their ratios of Chern
 numbers. Indeed, compact complex hyperbolic manifolds are complex
 manifolds of general type (i.e. their Kodaira dimension is maximal),
 and they have the same ratios of Chern numbers as complex projective
 space $\cp n$, which is the compact symmetric space dual to $\ch n$
 (this is explained by the Hirzebruch proportionality principle). In
 fact, any compact complex manifold of general type that realizes
 equality in the Miyaoka-Yau inequality $(-1)^n\,c_1^n \leq (-1)^n\,
 \frac{2n+2}{n}\,c_1^{n-2}\,c_2$ is biholomorphic to a quotient of $\ch n$
 (this is a corollary of the Calabi conjecture, proved by Aubin and
 Yau). In the special case $n=2$, the corresponding equality reads
 $c_1^2=3\,c_2$, and it plays an important role in so-called surface
 geography.

  In principle, this characterization should give a simple way to
  produce many complex hyperbolic manifolds, but so far this has been
  successful in very few cases, only in dimension 2. A spectacular
  example was obtained by Mumford, using sophisticated techniques from
  algebraic geometry to produce a fake projective plane, i.e. a
  compact surface with the same Betti numbers as $\cp 2$ but not
  biholomorphic to it. The fundamental group of any such surface is
  infinite, and admits a faithful and discrete representation into
  ${\rm PU}(2,1)$, even though it took a long time for this
  representation to be made somewhat explicit for Mumford's example
  (see~\cite{kato}). It is now well known that fake projective planes
  must be arithmetic (see~\cite{klingler} and~\cite{yeung}), and this
  was used to classify them (see~\cite{prasadYeung}
  and~\cite{cartwrightsteger}).
  
  Hirzebruch systematically explored coverings of $\cp 2$ branched along
  configurations of lines, and found all the ball
  quotients that could be produced in this way,
  see~\cite{hirzebrucharrangements},~\cite{bhh}.  The arithmetic
  properties of Hirzebruch ball quotients were studied later
  (see~\cite{ossip}). It turns out that all the non-arithmetic examples can
  also be obtained from the Deligne-Mostow
  construction~\cite{delignemostow}.

  Our construction is closer in spirit to Mostow's approach
  in~\cite{mostowpacific}. We start with some well-chosen complex
  reflection groups in ${\rm PU}(2,1)$ and show that they are lattices
  by constructing explicit fundamental domains (the domains are shown
  to be fundamental domains by applying the Poincar\'e polyhedron
  theorem).

  The application of this general strategy in the context of complex
  hyperbolic space is difficult and subtle, as can be seen in the
  difficulties in Mostow's proof that were analyzed
  in~\cite{derauxexpmath}. It has been successfully carried out in
  several places, but always for groups closely related to Mostow's
  (see~\cite{dfp},~\cite{parkerlivne},~\cite{falbelparker},~\cite{zhaoPJM}).
  For the groups we consider in this paper, the previously used
  techniques seem difficult to implement (for instance, the Dirichlet
  domains we studied in~\cite{dpp1} have extremely complicated
  combinatorial structure).

Mostow's lattices turn out to be closely related to monodromy groups
of hypergeometric functions studied a century earlier by
Picard~\cite{picard}. The hypergeometric interpretation was extended
by Deligne and Mostow~\cite{delignemostow}, who showed that these
monodromy groups produce a handful of non-arithmetic lattices in ${\rm
  PU}(2,1)$, and a single one in ${\rm PU}(3,1)$ (which is currenctly
the only known non arithmetic lattice in ${\rm PU}(n,1)$ wih
$n\geqslant 3$). 

Note that the above hypergeometric monodromy groups can also be
interpreted as modular groups for moduli spaces of
weighted points on $\cp 1$, or equivalently moduli
spaces of flat metrics on the sphere with prescribed cone angle
singularities (see~\cite{thurstonshapes}).  The few other examples of
moduli spaces that are known to admit a complex hyperbolic
uniformization (see~\cite{act},~\cite{kondo} for instance) produce
arithmetic lattices.

The commensurability relations between Deligne-Mostow lattices were
studied in detail in~\cite{sauter} and~\cite{delignemostowbook}, then
more recently in~\cite{kappesmoeller} and~\cite{mcmullen}. In
particular, it turns out that non-arithmetic Deligne-Mostow lattices
in ${\rm PU}(2,1)$ fall into exactly nine commensurability classes.

The main result of this paper produces five new commensurability
classes of non-arithmetic lattices in ${\rm PU}(2,1)$. They are the
first such groups to be constructed since the work of Deligne and Mostow.

 We consider groups generated by two
  isometries, namely a complex reflection and a regular elliptic
  element of order 3. For historical reasons (see
    \cite{mostowpacific}), these two isometries will be denoted
  respectively by $R_1$ and $J$; we will also denote $R_2=JR_1J^{-1}$ and
  $R_3=JR_2J^{-1}$.  It turns out that conjugacy classes of such
  groups are parametrized by the pair $(\psi,\tau)$, where $\psi$ is
  the rotation angle of $R_1$ and $\tau$ is the trace of
  $R_1J$. There are restrictions on the pair $(\psi,\tau)$ for a group
  with a given angle/trace pair to exist (see the discussion
  in~\cite{parkerpaupert}); when it exists, we denote the
  corresponding group by $\Gamma(\psi,\tau)$.

 Following previous work on complex reflection groups
  (see~\cite{mostowpacific}, \cite{schwartzICM}, \cite{deraux4445},
  \cite{thompsonthesis}), one expects such groups to be lattices only
  when well-chosen short words in the generators are elliptic.  For
  instance, for (4,4,4)-triangle groups, Schwartz conjectured that the
  discreteness was controlled by the word 1232, in the sense that the
  group should be discrete if and only if $R_1R_2R_3R_2$ is
  loxodromic, parabolic, or elliptic of finite order (here the $R_j$
  denote the generating complex reflections).  The (4,4,4)-triangle
  group where $R_1R_2R_3R_2$ is elliptic of order $n$ (more
  specifically has eigenvalues $1$, $e^{\pm 2\pi i/n}$) is referred to
  as the $(4,4,4;n)$-triangle group. This seems to have infinite
  covolume for $n>5$, see for instance the group
  in~\cite{schwartz4447}, which is the $(4,4,4;7)$-triangle group. On
  the other hand, it was observed in~\cite{deraux4445} that the
  $(4,4,4;5)$-triangle group is a lattice. When considering
  $(4,4,4;n)$-triangle groups, an important distinction between $n=5$
  and $n\geqslant 6$ is the type of the isometry $R_1R_2R_3=(R_1J)^3$, which is
  loxodromic for $n\geqslant 6$, but elliptic of order $10$ when $n=5$. It
  is then natural to expect the ellipticity of the words 123 and 1232
  to control the property of triangle groups to be a lattice. This is
  confirmed for some other triangle groups by the results
  in~\cite{thompsonthesis}.

This motivated us to study groups $\Gamma(\psi,\tau)$ where $R_1R_2$
and $R_1J$ are both elliptic of finite order (in fact we allow $R_1R_2$ to
be parabolic as well, but not $R_1J$). The groups $\Gamma(\psi,\tau)$ satisfying these
  conditions were listed by the last two authors
  (see~\cite{parkerunfaithful} and~\cite{parkerpaupert}); the key
  ingredient in their arguments is a result of Conway and Jones~\cite{conwayjones} that
  classifies rational relations between roots of unity.

We proved in~\cite{dpp1} that only finitely many of these groups can
be lattices. We also gave a conjectural list of ten groups that we
strongly suspected to be lattices, corresponding to three special
choices of $\tau$, namely
  $$ 
  \sigma_1=  -1+i\sqrt{2}, \quad
  \overline{\sigma}_4=\frac{-(1+i\sqrt{7})}{2}, \quad \sigma_5= e^{-\pi
    i/9}\, \frac{i\sqrt{3}+\sqrt{5}}{2},
  $$
  for some well-chosen rotation angles $\psi$.
The main result of this paper is the following:

  \begin{theorem}\label{thm:main}
    Let $\tau=-(1+i\sqrt{7})/2$, and let $p\in \Z$. Then
    $\Gamma(2\pi/p,\tau)$ is a lattice in ${\rm PU}(2,1)$ if and only if
    $p=3,4,5,6,8$ or $12$, and that lattice is
    \begin{enumerate}
    \item cocompact if and only if $p=3,5,8$ or $12$;
    \item arithmetic if and only if $p=3$;
    \end{enumerate}
    Moreover, the six lattices lie in distinct commensurability
    classes, and they are not commensurable to any Deligne-Mostow
    lattice.
  \end{theorem}
The values of $p$ in the theorem may seem mysterious, but one can
rephrase the theorem as an integrality condition, in the vein of the
Picard integrality condition (see~\cite{delignemostow} for
instance). Indeed (see section~\ref{sec:braiding}), $(R_1R_2)^2$ is a complex reflection with angle
$2\pi/c$, and $(R_1R_2R_3R_2^{-1})^3$ is a reflection with angle
$2\pi/d$, where:
  $$
  c=2p/(p-4); \quad d=2p/(p-6).
  $$ 
The values of $p$ in the theorem are precisely those for which both
$c$ and $d$ are integers (or infinite).
  
Non-arithmeticity of these groups was proved in~\cite{paupert}.  For
completeness, we prove that the five commensurability classes are new
\emph{and distinct} in Section~\ref{sec:arithmeticgalois}, and recall
the proof of non-arithmeticity in Section~\ref{sec:arithmetic}; see
also~\cite{paupert}.  The fact that the groups are non-discrete for
all other integer values of $p$ was proved in~\cite{dpp1}.

We prove that the six groups that appear in Theorem~\ref{thm:main} are
lattices by constructing explicit fundamental polyhedra for their
respective actions on $\ch 2$, using Theorem~\ref{thm:poin}.  By this
method we obtain presentations of the lattices:

\begin{theorem}\label{thm:pres}
  Let $\tau=-(1+i\sqrt{7})/2$. For
  $p=3,4,5,6,8,12$ the group $\Gamma(2\pi/p,\tau)$ has the
  presentation:
  $$
    \left\langle\ R_1,\,R_2,\,R_3,\,J\ |\ 
    \begin{array}{c} R_1^p=J^3=(R_1J)^7=id,\ R_2=JR_1J^{-1},\ 
      R_3=J^{-1}R_1J, \\
      (R_1R_2)^2=(R_2R_1)^2,\ (R_1R_2)^{4p/(p-4)}=(R_1R_2R_3R_2^{-1})^{6p/(p-6)}=id
    \end{array}\right\rangle.
  $$
  If an exponent is infinite or negative then the relation should be omitted. 
\end{theorem}

Moreover, from the detailed analysis of the polyhedra and the groups,
we obtain the orbifold Euler characteristic $\chi$ of the quotient
orbifolds (hence their complex hyperbolic volume which is
$8\pi^2\chi/3$):

\begin{theorem}\label{thm:euler}
  Let $\tau=-(1+i\sqrt{7})/2$. Let
  $\chi(2\pi/p,\tau)$ be the orbifold Euler characteristic of the
  orbifold $\Gamma(2\pi/p,\tau)\backslash \ch 2$. Then
  $\chi(2\pi/p,\tau)$ has the following values
  $$
    \begin{array}{c|cccccc}
      p & 3 & 4 & 5 & 6 & 8 & 12 \\ 
      \hline
      \chi(2\pi/p,\tau) &
      2/63 & 25/224 & 47/280 & 25/126 & 99/448 & 221/1008 \\
    \end{array}
  $$
\end{theorem}

These statements prove our conjecture from~\cite{dpp1} for five groups out
of ten.  The corresponding statement about lattices with
$\tau=\sigma_1$ and $\tau=\sigma_5$ would follow from similar
techniques, but some of the verifications are too tedious to be
written up here. 

For a fixed value of $\tau$ the number of sides of the
  domain is the same, and the
  combinatorial structure is very similar. For each of the groups
  presented in the paper, the domain has 28 sides, but there are only
  two isometry types of sides (the domain is invariant under a regular
  elliptic element of order 7, and sides are paired
  isometrically). This decreases the complexity enough to make it
  possible to present most of the proof in detail on paper.

  For $\sigma_1$ (resp. $\sigma_5$), our construction would produce a
  polyhedron with 64 (resp. 130) sides, and there would be four
  (resp. three) isometry types of sides. This makes some of the
  verifications longer, but not more difficult in any essential
  way. One important place where the arguments would be more painful
  is the verification that the bounding pyramids glue to build a copy
  of $S^3$, see section~\ref{sec:sphere} of the present paper. In
  these more complicated cases, one could resort to the solution
  of the Poincar\'e conjecture, and check only that the bounding 3-complex is a manifold with trivial fundamental group; that argument can easily be automatized, but the computation would be difficult to write down on
  paper.

Our main results are proved using a standard application of the
Poincar\'e polyhedron theorem, but we also provide new, useful
techniques for studying fundamental domains for lattices in ${\rm
  PU}(2,1)$. 

We apply systematic techniques for producing a fundamental domain with
fairly simple combinatorics and with sides contained in bisectors (as
other nice features, all vertices are fixed by specific isometries in
the group, and all $1$-faces of our polyhedron are chosen to be geodesic
arcs). These techniques, which were in part inspired by the construction in \cite{schwartz4447}, will be briefly summarized in section~\ref{sec:algo}.

Another novel aspect is the provision of efficient computational tools
in Section~\ref{sec:realization} for certification of the combinatorics
of a polyhedron bounded by (finitely many) bisectors. We summarize these
techniques, since they have not systematically been used in the
literature.

The determination of the combinatorics amounts to solving a system of
(finitely many) quadratic inequalities in the real and imaginary parts
of ball coordinates, for which no standard computational techniques
are known. The main difficulty is to determine which pairs of
bisectors contain a non-empty $2$-face, and to give a precise list of
the vertices and edges adjacent to that $2$-face.

Using geometric properties of bisectors and bisector intersections
(see Section~\ref{sec:bisectors}), we can give parameters for any
intersection $\B_1\cap \B_2$ of bisectors, and write explicit
equations for the intersection $\B_1\cap\B_2\cap\B_3$ for any third
bisector $\B_3$.  For well-chosen parametrizations, these equations
are given by polynomials in two variables, quadratic in each
variable. This is a consequence of the fact that the (hyperbolic
cosine of the) distance between two points in $\ch 2$ is given by a
quadratic polynomial in the real and imaginary part of affine
coordinates for the points.  The coefficients of these polynomial
equations are algebraic numbers (in fact, they lie in the adjoint
trace field), which allows for arbitrary precision calculations.

The determination of the combinatorics is then reduced to the
minimization of finitely many polynomials in two variables on finitely
many polygonal regions, where the polygons are bounded by plane curves
of degree at most two. We perform this minimization by computing the
critical points of the above polynomials, which can be done by exact
computations since the coefficients are known exactly.
  Even though these verifications can in principle be performed by
  hand, it is most reasonable to use a computer, given the large
  number of computations. Open source software that performs the
  necessary checks is available on the first author's web page.
 
Most of these techniques are completely general, and can be used to
reduce the certification of fundamental domains for discrete subgroups
of ${\rm PU}(2,1)$ to a finite number of verifications. Note that the
determination of fundamental domains is clearly important in the
context of the search of non-arithmetic lattices, but it is also
important for arithmetic ones, where discreteness comes for free. For
instance, the determination of explicit group presentations for small
covolume arithmetic lattices in ${\rm PU}(2,1)$ allowed Cartwright and
Steger to complete the classification of fake projective planes
(see~\cite{prasadYeung} and~\cite{cartwrightsteger}).

The paper is organized as follows. In Sections~\ref{sec:chg}
and~\ref{sec:discrete} we give well-known background material.
Section~\ref{sec:domain} starts by establishing basic notation for the
specific groups we study in the paper.  We then sketch the general
procedure that was used in order to produce the combinatorial model
for our fundamental polyhedron (Sections~\ref{sec:braiding}
and~\ref{sec:algo}).

The details of the combinatorial model $\widehat{E}$ for the groups that appear in
this paper (incidence properties of facets of various dimensions) are
discussed in Section~\ref{sec:combinatorics}, and the
strategy to show that it is homeomorphic to its geometric realization $E$
is outlined in Section~\ref{sec:model}.
To that end, an important point is to show that the $3$-skeleton of
$\widehat{E}$ is homeomorphic to $S^3$, which is proved in
section~\ref{sec:sphere}.
We then show that the geometric realization of $\widehat{E}$ is
embedded, see Section~\ref{sec:realization}; 
apply the Poincar\'e polyhedron theorem to $E$ in order to show that
our groups are lattices, Section~\ref{sec:applyingpoincare}; and
finally we show that these lattices are non-arithmetic and not
commensurable to each other or to any Deligne-Mostow lattice,
Section~\ref{sec:arithmeticgalois}.

\smallskip
\noindent {\bf Acknowledgements:} The authors would like to thank the
following institutions for their support during the preparation of
this paper, in chronological order: the University of Utah,
Universit\'e de Fribourg, Durham University, Universit\'e de Grenoble,
Arizona State University. The authors acknowledge support from the ANR
through the program ``Structures G\'eom\'etriques et Triangulations'',
NSF grants DMS 1107452, 1107263, 1107367 (the GEAR Network) and ICERM
at Brown University. The third author was also partially supported by
SNF grant 200020-121506/1 and NSF grant DMS 1007340/1249147. The
authors would also like to thank Bernard Parisse and Fabrice Rouillier
for useful assistance on the computational aspects of this project,
as well as the referee for several suggestions which improved the
exposition of the paper.

\section{Complex hyperbolic geometry}\label{sec:chg}

We define complex hyperbolic space, introduce totally geodesic
subspaces and bisectors, and review properties of bisector
intersections. Most of this material may be found in~\cite{goldman}.

\subsection{Complex hyperbolic space} \label{sec:chn}

We define $\ch n$ to be the subset of $\cp n$ consisting of negative
complex lines in $\C^{n,1}$. Here $\C^{n,1}$ denotes $\C^{n+1}$
equipped with a Hermitian form of signature $(n,1)$, and a negative
line is one spanned by a vector $\bu\in \C^{n,1}$ with $\langle
\bu,\bu\rangle<0$. There is a natural action on $\ch n$ of the unitary
group of the Hermitian form, which is denoted by ${\rm U}(n,1)$. We
shall often work with ${\rm SU}(n,1)$, which is an $(n+1)$-fold cover
of the projective group ${\rm PU}(n,1)$.

Up to scaling, $\ch n$ carries a unique ${\rm U}(n,1)$-invariant
Riemannian metric, which makes it a symmetric space.
It is well known that $\ch n$ is a complex space form, in fact the
holomorphic sectional curvature of the above metric is equal to
$-1$. This implies that real sectional curvatures are pinched between
$-1$ and $-1/4$ (when $n=1$ the curvature is in fact constant).

It is a standard fact that the group of holomorphic isometries of $\ch
n$ is precisely ${\rm PU}(n,1)$. The full group ${\rm Isom}(\ch n)$
contains ${\rm PU}(n,1)$ with index 2, the other component consisting
of all anti-holomorphic isometries (an example of which is complex
conjugation in affine coordinates, provided that the Hermitian form
has real entries).

The only metric information we will need in this paper is the
following distance formula:
\begin{equation}\label{eq:distance}
  \cosh \left(\frac{d(u,v)}{2}\right)=\frac{ |\langle
    \bu,\bv\rangle|}{\sqrt{\langle \bu,\bu\rangle\langle \bv,\bv\rangle}},
\end{equation}
where $\bu$, $\bv$ denote lifts of $u,v$ to $\C^{n+1}$.  The standard
choice of a Hermitian form of signature $(n,1)$ is the Lorentzian form
\begin{equation}\label{eq:stdform}
  \langle \bu,\bv\rangle = -u_0\overline{v}_0 + u_1 \overline{v}_1 +
  \dots + u_n\overline{v}_n.
\end{equation} 
In the affine coordinates $(z_1,\dots,z_n)=(u_1/u_0,\dots,u_n/u_0)$
for the chart $\{u_0\neq 0\}$ of $\cp n$, complex hyperbolic space
corresponds to the unit ball $|z_1|^2+\dots+|z_n|^2 < 1$. We
will denote $\chb n = \ch n\cup \partial \ch n$ the
corresponding closed ball.

We will use the standard classification of isometries of negatively
curved spaces into elliptic, parabolic and hyperbolic isometries, with
the following refinements. An elliptic element of ${\rm PU}(n,1)$ is
called {\bf regular elliptic} if any of its lifts to ${\rm U}(n,1)$
has distinct eigenvalues. In particular, regular elliptic isometries
have isolated fixed points in $\ch n$ (given by the single negative
eigenspace). When $n=2$, non-regular elliptic isometries of $\ch 2$
have two distinct eigenvalues, one simple and one double
eigenvalue. These elements are called {\bf complex reflections} in
points or in lines, depending on the sign of the simple eigenspace.

\subsection{Totally geodesic subspaces} \label{sec:totallygeodesic}

If $V\subset \C^{n,1}$ is any complex linear subspace where the
restriction of the Hermitian form has signature $(k,1)$, then the set
of negative complex lines in $V$ gives a copy of $\ch k$, which can
easily be checked to be totally geodesic. In terms of the ball model,
they correspond to the intersection with the unit ball of complex
affine subspaces in $\C^n$. When $k=n-1$, the subspaces $V$ as above
are simply obtained by taking the orthogonal complement of a positive
vector ${\bf n}$, which we refer to as a {\bf polar vector} to the
corresponding complex hyperplane. When $k=1$, the complex totally
geodesic submanifolds described above are called {\bf complex lines},
or {\bf complex geodesics}. Note that there is a unique complex line
joining any two distinct points in $\ch n$.

Similarly, some $\R$-linear subspaces can be used to produce totally
geodesic submanifolds, namely those where the Hermitian form restricts
to a quadratic form of hyperbolic signature. These yield copies of
$\rh k$, $k\leqslant n$ (with real sectional curvature $-1/4$). For $k=1$ these
are precisely the real geodesics, and for $k=2$, they are called {\bf
  $\R$-planes}.

It is a standard fact that every complete totally geodesic submanifold
is in one of the two families described above (see~\cite{goldman},
Section~3.1.11). Note that in particular, there are no totally
geodesic real hypersurfaces in complex hyperbolic spaces.

\subsection{Bisectors} \label{sec:bisectors}

The basic building blocks for our fundamental domains are {\bf bisectors},
which are hypersurfaces equidistant from two given points. Their
geometric structure has been analyzed in great detail
in~\cite{mostowpacific},~\cite{goldman}, we only recall a few facts
that will be needed in the paper.

Given two distinct points $p_0,p_1\in \ch n$, write 
$$
  \B(p_0,p_1)=\{u\in \ch n:d(u,p_0)=d(u,p_1)\}
$$ 
for the bisector equidistant from $p_0$ and $p_1$. 

There is some freedom in choosing the pair of points $p_0$,
$p_1$, but perhaps surprisingly, the possible choices are much more
constrained than in constant curvature geometries. Indeed, the points
$q_0$ for which there exists a $q_1$ with $\B(q_0,q_1)=\B(p_0,p_1)$ all
lie on the complex line through $p_0$ and $p_1$; thus the complex line
through the two points is canonically attached to the bisector, and is
called its {\bf complex spine}.

The {\bf real spine} is the intersection of the complex spine with the
bisector itself, which is a (real) geodesic; it is the locus of
points inside the complex spine which are equidistant from $p_0$ and
$p_1$. It consists of points $z\in \ch n$ associated to
negative vectors ${\bf u}\in \C^{n,1}$ that are in the
complex span of $\p_0$ and $\p_1$ and that satisfy
\begin{equation}\label{eq:realspine}
|\langle {\bf u},\p_0\rangle|=|\langle {\bf u},\p_1\rangle|,
\end{equation}
where $\p_0$ and $\p_1$ denote lifts of $p_0$ and $p_1$ to $\C^{n,1}$
with the same norm (e.g. one could take
$\langle\p_0,\p_0\rangle=\langle\p_1,\p_1\rangle=-1$).  The nonzero
vectors in ${\rm Span}_\C(\p_0,\p_1)$ satisfying~(\ref{eq:realspine}),
but that are not necessarily negative, span a real projective line in
${\bf P}^n_{\C}$, which we call the {\bf extended real spine} of the
bisector.

We will sometimes describe bisectors by giving two points of their
extended real spines; as mentioned in
Section~\ref{sec:totallygeodesic}, we can think of a real geodesic as
a real $1$-dimensional hyperbolic space. Hence, we can describe it as
the projectivization of a totally real $2$-dimensional subspace of
$\C^{n,1}$, i.e. we take two vectors $\mathbf v$ and $\mathbf w$ in
$\C^{n+1}$ with $\langle \mathbf v, \mathbf w\rangle\in\R$, and
consider their real span. The simplest way to guarantee that the span
really yields a geodesic in $\ch n$ is to require moreover that
$\mathbf v$ and $\mathbf w$ form a Lorentz basis, i.e. $\langle
\mathbf v, \mathbf v\rangle=-1$, $\langle \mathbf w,\mathbf
w\rangle=1$ and $\langle \mathbf v,\mathbf w\rangle=0$.

Recall that bisectors are \emph{not} fixed point sets of any isometric
involution, since they are not totally geodesic.  According to the
discussion in Section~\ref{sec:totallygeodesic}, there are two kinds
of maximal totally geodesic submanifolds contained in a given
bisector.  The complex ones are called {\bf complex slices} of $\B$;
they are the preimages of points of the real spine under orthogonal
projection onto the complex spine. They can also be described as
hyperplanes with polar vectors ${\bf n}\in \C^{n,1}$ satisfying
$\langle {\bf n},{\bf n}\rangle>0$ and $|\langle {\bf
  n},\p_0\rangle|=|\langle {\bf n},\p_1\rangle|$, see
equation~(\ref{eq:realspine}).  The totally real submanifolds that are
contained in $\B$ are precisely those containing the real spine and
are called {\bf real slices} of $\B$.

\subsection{Bisectors and geodesics} \label{sec:bisectorsandgeodsics}

Several times in the paper, we will need to determine when a geodesic
arc is contained in a bisector. Recall that bisectors in $\ch 2$ are
not convex: given two distinct points $p$, $q$ in a bisector $\B$, the
geodesic through $p$ and $q$ is in general not contained in $\B$. This
is stated more precisely in the following (see~Theorem~5.5.1
of~\cite{goldman}).
\begin{lemma}\label{lem:nonconvex}
  Let $\alpha\subset\ch 2$ be a real geodesic, and let $\B$ be a
  bisector. Then $\alpha$ is contained in $\B$ if and only if it is
  contained either in a real slice or in a complex slice of $\B$.
\end{lemma}

In particular, we have:
\begin{lemma} \label{lem:realspine}
  Let $u$ and $v$ be distinct points of a bisector $\B$, such that $u$
  is on the real spine of $\B$. Then the geodesic through $u$ and $v$
  is contained in $\B$.
\end{lemma}
Indeed, if $v$ is in the complex slice through $u$ this is
obvious. Otherwise, taking lifts ${\bf u}$, ${\bf v}$ of $u$, $v$, the
point $v$ is in a slice ${\bf w}^\perp$ with $\langle {\bf u},{\bf
  w}\rangle \neq 0$ in which case ${\bf u},{\bf v},{\bf w}$ are
linearly independent. Since $\langle {\bf v}, {\bf w} \rangle = 0$, we
can normalize ${\bf u}$, ${\bf v}$, ${\bf w}$ to have pairwise real
inner products, so they span a copy of $\rh 2$ which contains the real
spine of $\B$, i.e. they span a real slice of $\B$.

The following result will also be useful (see also Lemma~2.2
of~\cite{dfp} for a related statement).
\begin{lemma} \label{lem:orthslice}
  Let $\B$ be a bisector and let $C$ be a
  complex line orthogonal to a complex slice of $\B$. Then $C\cap \B$ is a
  real geodesic contained in a real slice of $\B$.
\end{lemma}
\Pf 
In the unit ball, we may normalize the real spine $\sigma$ of $\B$
to be the set of points of the form $(x_1,0)$, with $x_1\in\R$, the
slice of $\B$ to be $\{z_1=0\}$, and $C$ to be of the form
$\{z_2=\beta\}$ for some $\beta$ (with $|\beta|<1$). Since orthogonal
projection onto $\{z_2=0\}$ is just the usual projection onto the
first coordinate axis, and a bisector is the preimage of its real
spine under orthogonal projection onto its complex spine, $C\cap \B$
is given by $(x_1,\beta)$ with $x_1\in\R$ (and
$|x_1|<\sqrt{1-|\beta|^2}$). The latter curve is a geodesic; in fact
it is the intersection $C\cap L$ where $L$ is obtained from the real
ball $B^2_\R\subset B^2_\C$ by applying the isometry $(z_1,z_2)\mapsto
(z_1,\frac{\beta}{|\beta|}z_2)$.  
\EPf

\subsection{Coequidistant bisector intersections} \label{sec:coequidistant}

It is well-known that bisector intersections can be somewhat
complicated, see the detailed analysis in~\cite{goldman}.  In what
follows, for simplicity of notation, we restrict ourselves to the case
of complex dimension $n=2$.

A simple way to get these intersections to be somewhat reasonable is
to restrict to {\bf coequidistant} pairs, i.e. intersections $\B_1\cap
\B_2$ where $\B_1$, $\B_2$ are equidistant from a common point $p_0\in
\ch n$. We write $\Sigma_1$, $\Sigma_2$ for the complex spines of
$\B_1$, $\B_2$ respectively, and $\sigma_1$, $\sigma_2$ for their real
spines.

In view of the discussion in Section~\ref{sec:bisectors}, it should be
clear that this is a restrictive condition, which implies that
$\Sigma_1$, $\Sigma_2$ intersect inside $\ch n$. Since complex lines
with two points in common coincide, there are two possibilities for
coequidistant pairs $\B_1,\B_2$; either the complex spines coincide or
they intersect in a single point.

When the complex spines coincide, the bisectors are called {\bf
  cospinal}, and the intersection is easily understood; it is
non-empty if and only if the real spines $\sigma_1$ and $\sigma_2$
intersect, and in that case the intersection $\B_1\cap \B_2$ consists
of a complex line (namely the complex line orthogonal to
$\Sigma_1=\Sigma_2$ through $\sigma_1\cap \sigma_2$).

Now suppose $\Sigma_1$ and $\Sigma_2$ intersect in a point which lies
outside of the real spines $\sigma_1$ and $\sigma_2$, so that
$\B_1\cap \B_2$ can be written as the equidistant locus from three
points $p_0$, $p_1$, $p_2$ which are not contained in a common complex
line.  The following important result is due to
G. Giraud~\cite{giraud} (see also Theorem 8.3.3 of \cite{goldman}).
\begin{proposition} \label{prop:giraud}
  Let $p_0,p_1,p_2$ be distinct points in $\ch 2$, not all contained in a
  complex line. When it is non empty, the intersection
  $\B(p_0,p_1)\cap \B(p_0,p_2)$ is a (non-totally
  geodesic) smooth disk. Moreover, it is contained in precisely three
  bisectors, namely $\B(p_0,p_1)$, $\B(p_0,p_2)$ and $\B(p_1,p_2)$. 
\end{proposition}
\begin{dfn} \label{dfn:girauddisk}
  We define a {\bf Giraud disk} to be a bisector intersection as in
  Proposition~\ref{prop:giraud} and we denote it by $\G(p_0,p_1,p_2)$.
\end{dfn}
Given a Giraud disk $\G=\B_1\cap \B_2\cap \B_3$, the complex slices of
the three bisectors $\B_j$ intersect $\G$ in hypercycles, so that $\G$
has three pairwise transverse foliations by arcs. 

\smallskip
\noindent {\bf Parametrizing Giraud disks:} Let $\p_j$ denote a lift
of $p_j$ to $\C^3$. By rescaling the lifts, we may assume that the
three square norms $\langle \p_j,\p_j\rangle$ are equal, and also that
$\langle \p_0,\p_1\rangle$ and $\langle \p_0,\p_2\rangle$ are real and
positive.

Now for $j=1,2$, consider $\tilde{\bv}_j =\p_0-\p_j$ and 
$\tilde{\bw}_j = i(\p_0+\p_j)$ (note that $\tilde{\bv}_j$ corresponds to the midpoint
of the geodesic segment between $p_0$ and $p_j$), and normalize these
to unit vectors $\bv_j=\tilde{\bv}_j/\sqrt{-\langle
  \tilde{\bv}_j,\tilde{\bv}_j\rangle}$ and
$\mathbf{w}_j=\tilde{\bw}_j/\sqrt{\langle\tilde{\bw}_j,\tilde{\bw}_j\rangle}$.

Here recall that $\mathbf{u}\boxtimes\mathbf{v}$ is the Hermitian
cross product of $\mathbf{u}$ and $\mathbf{v}$ associated to the
Hermitian form $H$ (see p.~43 of~\cite{goldman}). By definition,
$\mathbf{u}\boxtimes\mathbf{v}$ is $\mathbf{0}$ if $\mathbf{u}$ and
$\mathbf{v}$ are collinear, and spans their $H$-orthogonal complement
otherwise. In other words, it is the Euclidean cross product of the
vectors $\mathbf{u}^*H$ and $\mathbf{v}^*H$, where $H$ is the matrix defining the Hermitian form (the cross product makes
sense because the vector space of homogeneous coordinates has
dimension three).

Then the extended real spine of $\B(p_0,p_j)$ is given by real linear
combinations of $\mathbf{v}_j$ and $\mathbf{w}_j$, so (lifts of)
points in $\B_1\cap \B_2$ are given by negative vectors of the form
\begin{equation}\label{paramgiraud}
  V(t_1,t_2)=(\mathbf{w}_1 + t_1 \mathbf{v}_1)\boxtimes(\mathbf{w}_2 + t_2
  \mathbf{v}_2),
\end{equation}
with $t_1,t_2\in\R$. The only linear combinations we are missing with
this parametrization of the extended real spines are $\mathbf{v}_1$
and $\mathbf{v}_2$, but these are negative vectors so the
projectivization of their orthogonal complement does not intersect
$\ch 2$. We will call $t_1,\,t_2$ {\bf spinal coordinates} for $\G$.

We will use this parametrization repeatedly; for now, note that, given
three points $p_0$, $p_1$ and $p_2$, it is easy to determine whether
the intersection $\B(p_0,p_1)\cap \B(p_0,p_2)$ is empty or
not. Indeed, negative vectors in $\C^{2,1}$ are characterized by the
fact that their orthogonal complement is a positive definite complex
$2$-plane, so the condition that $(\mathbf{w}_1 + t_1
\mathbf{v}_1)\boxtimes(\mathbf{w}_2 + t_2 \mathbf{v}_2)$ is negative
is equivalent to requiring that $g(t_1,t_2)>0$ where
\begin{equation} \label{eq:det}
g(t_1,t_2)=\det \left[ \begin{matrix} \langle \mathbf{w}_1+ t_1 \mathbf{v}_1,
    \mathbf{w}_1+ t_1 \mathbf{v}_1\rangle & \langle \mathbf{w}_1+ t_1
    \mathbf{v}_1, \mathbf{w}_2+ t_2 \mathbf{v}_2\rangle\\ 
    \langle \mathbf{w}_2+ t_2 \mathbf{v}_2, \mathbf{w}_1+ t_1
    \mathbf{v}_1\rangle &  \langle \mathbf{w}_2+ t_2 \mathbf{v}_2, \mathbf{w}_2+
    t_2 \mathbf{v}_2\rangle\\
\end{matrix} \right].
\end{equation}

\subsection{Cotranchal bisectors} \label{sec:cotranchal}

We will encounter more complicated intersections than coequidistant
ones. Two bisectors are called {\bf cotranchal} if they have a common
complex slice; the intersection of cotranchal bisectors is completely
understood; see Thm~9.2.7 in~\cite{goldman}.

\begin{theorem}[Goldman]\label{thm:cotranchal} Let $\B_1$, $\B_2$ be distinct
 bisectors, with real spines $\sigma_1$, $\sigma_2$. Assume that $\B_1$
 and $\B_2$ share a complex slice $C_0$.
\begin{itemize}
\item[(1)] If $\sigma_1$ and $\sigma_2$ intersect, then $\B_1 \cap \B_2$
  is equal to $C_0$;
\item[(2)] If $\sigma_1$ and $\sigma_2$ lie in a common $\R$-plane $L$
  and are ultraparallel, then $\B_1 \cap \B_2=C_0 \cup
  L$;
\item[(3)] If $\sigma_1$ and $\sigma_2$ do not lie in a totally
  geodesic plane, then either $\B_1 \cap \B_2 = C_0$, or $\B_1 \cap \B_2 =
  C_0 \cup R$ where $R$ is diffeomorphic to a disk and $R \cap C_0$ is
  a hypercycle.
\end{itemize}
\end{theorem}

Goldman's proof also gives a way to distinguish between the two very
different possibilities that are described in case~(3), working with 
Heisenberg coordinates for the vertices of the bisectors
(which are the endpoints in $\partial\ch 2$ of their real spines).
We briefly review this in terms of our notation.

Let $C$ be a complex line with polar vector ${\bf n}$. Note that this
implies $\langle {\bf n},{\bf n}\rangle>0$. Let ${\bf v}_1$
and ${\bf v}_2$ be two vectors in ${\bf n}^\perp$ with
$\langle {\bf v}_j,{\bf v}_j\rangle<0$. Hence ${\bf v}_j$
projects to a point $v_j$ in $\ch 2$, which of course is in $C$. 
Consider bisectors $\B_1$ and $\B_2$ so that the extended real spine of 
$\B_j$ contains ${\bf n}$ and ${\bf v}_j$. (Note that we have to be careful 
when choosing the lift ${\bf v}_j$ in $\C^{2,1}$ of the corresponding points 
$v_j$ in $C$ since we want the extended spine to be the real span of
${\bf n}$ and ${\bf v}_j$; choosing a different lift of ${\bf v}_j$
rotates $\B_j$ around $C$.)
We want to decide when the intersection of $\B_1$ and $\B_2$ is precisely
$C$; see \cite{goldman}, \cite{hsieh}.

\begin{proposition}\label{prop:cotranchal}
  Let $\B_1$, $\B_2$ be cotranchal bisectors with common slice $C$.
  Let ${\bf n}$ be a polar vector for $C$ and let ${\bf v}_1$, ${\bf
    v}_2$ be negative vectors with $\langle{\bf v}_j,{\bf n}\rangle=0$
  so that the extended real spine of $\B_j$ is the real span of ${\bf
    n}$ and ${\bf v}_j$. Then the intersection of $\B_1$ and $\B_2$ is
  precisely $C$ if and only if
  $$
    \langle{\bf v}_1,{\bf v}_1\rangle\langle{\bf v}_2,{\bf v}_2\rangle
    \ge \bigl(\Re\langle{\bf v}_1,{\bf v}_2\rangle\bigr)^2.
  $$
\end{proposition}

\Pf 
By assumption $C$ is a slice of $\B_1$ and $\B_2$. The other slices are parametrized
by the polar vectors $t_j{\bf n}+{\bf v}_j$ subject to
$\langle t_j{\bf n}+{\bf v}_j, t_j{\bf n}+{\bf v}_j\rangle>0$. In
other words,
\begin{equation}\label{eq:x-range}
t_j^2>\frac{-\langle {\bf v}_j,{\bf v}_j\rangle}
{\langle {\bf n},{\bf n}\rangle}.
\end{equation}
Hence we want to find conditions under which the corresponding slices
are disjoint, for all $t_1$ and $t_2$ satisfying this inequality.

Two complex lines ${\bf a}^\perp$ and ${\bf b}^\perp$ are disjoint if
and only if the restriction of the Hermitian form to the span of ${\bf
  a}$ and ${\bf b}$ is indefinite, so the above condition is
equivalent to the following inequality
being satisfied for all $t_1$ and $t_2$ in
the range \eqref{eq:x-range}:
\begin{eqnarray}
0 & \le &  \bigl|\langle t_1{\bf n}+{\bf v}_1,t_2{\bf n}+{\bf v}_2\rangle\bigr|^2
-\langle t_1{\bf n}+{\bf v}_1,t_1{\bf n}+{\bf v}_1\rangle
\langle t_2{\bf n}+{\bf v}_2,t_2{\bf n}+{\bf v}_2\rangle \notag \\
& = & -t_1^2\langle{\bf n},{\bf n}\rangle\langle{\bf v}_2,{\bf v}_2\rangle
+t_1t_2\langle{\bf n},{\bf n}\rangle2\Re\langle{\bf v}_1,{\bf v}_2\rangle
-t_2^2\langle{\bf n},{\bf n}\rangle\langle{\bf v}_1,{\bf v}_1\rangle 
\label{eq:cotranchal}\\
&& \quad +\bigl|\langle{\bf v}_1,{\bf v}_2\rangle\bigr|^2
-\langle{\bf v}_1,{\bf v}_1\rangle\langle{\bf v}_2,{\bf v}_2\rangle.\notag
\end{eqnarray}
First, note that the constant term is non-negative since 
$$
1\le \cosh^2\left(\frac{d(v_1,v_2)}{2}\right)
=\frac{\bigl|\langle{\bf v}_1,{\bf v}_2\rangle\bigr|^2}
{\langle{\bf v}_1,{\bf v}_1\rangle\langle{\bf v}_2,{\bf v}_2\rangle}.
$$ 
Therefore, the expression in \eqref{eq:cotranchal} is always
non-negative if and only if its quadratic part is positive
semi-definite. (Recall $\langle{\bf n},{\bf n}\rangle>0$ and
$\langle{\bf v}_j,{\bf v}_j\rangle<0$.)  This is equivalent to
$$
\langle{\bf v}_1,{\bf v}_1\rangle\langle{\bf v}_2,{\bf v}_2\rangle
\ge \bigl(\Re\langle{\bf v}_1,{\bf v}_2\rangle\bigr)^2.
$$
Finally, note that the range \eqref{eq:x-range} is irrelevant in the
argument as we may rescale $t_1$ and $t_2$ without affecting the sign
of the quadratic form. 
\EPf

\subsection{Computational issues} \label{sec:numberfields}
As described in the previous sections, bisector intersections can be
studied by checking the sign of various functions, sometimes defined
in a somewhat complicated compact region.

This will be facilitated by the fact that all geometric objects we
introduce are defined over specific number fields. For general
computational tools in number fields, see~\cite{cohen}. Throughout
this section, we denote by $\k$ a number field, and we assume complex
hyperbolic space is given as the set of lines in projective space that
are negative with respect to a Hermitian inner product \emph{defined
  over $\k$}.

\begin{dfn}
  A point $p\in\cp 2$ will be called $\k$-rational if it can be
  represented by a vector in $\k^3$. A real geodesic will be called
  $\k$-rational if it can be parametrized by vectors of the form
  $\vv+t\ww$ $(t\in\R)$ with $\vv,\ww\in\k^3$. A bisector will be called
  $\k$-rational if its real spine is $\k$-rational.
\end{dfn}

Note that a real geodesic through two $\k$-rational points is
automatically $\k$-rational. Indeed, if $v\neq w$ are two
$\k$-rational points in $\ch 2$, then they have lifts to $\vv,\ww\in
\k^3$ such that $\langle \vv,\ww\rangle$ is real and positive (given
any two lifts, multiply $\ww$ by $\langle \vv,\ww\rangle$, which is in
$\k$ because the Hermitian form is defined over $\k$). The real
geodesic segment from $v$ to $w$ is then parametrized by vectors of
the form $\vv+t(\ww-\vv)$, $t\in [0,1]$, where $\vv$ and $\ww$ are
$\k$-rational. In these circumstances, we will also say that the
geodesic segment $[v,w]$ is $\k$-rational.

Note also that a bisector is $\k$-rational if and only if it can be written
as $\B(p,q)$ for some $\k$-rational points $p$, $q$.

\begin{proposition}
  Let $[p,q]$ be a geodesic segment defined by $\k$-rational
  endpoints. Let $\B$ be a $\k$-rational bisector. The following
  statements can be checked by performing arithmetic in $\k$: 
  \begin{itemize}
    \item $p$ is in $\B$ (resp. $p$ is not in $\B$);
    \item $[p,q]$ is contained in $\B$;
    \item $[p,q]$ is in a specific component of $\ch 2\setminus \B$.
    \item $p$ is in $\B$ and $[p,q]$ is tangent to $\B$.
  \end{itemize}
\end{proposition}

\Pf 
  We parametrize the geodesic segment by vectors
  $\z_t=\vv+t(\ww-\vv)$, $t\in[0,1]$. If $\B=\B(q_0,q_1)$ and $\q_0$,
  $\q_1\in\k^3$ are corresponding lifts, the intersection of the
  extended real geodesic with (the extension to $\cp 2$ of) $\B$ is described by
  the following equation,
  \begin{equation}\label{eq:tvals}
  \frac{|\langle \z_t,\q_0\rangle|^2}{\langle \q_0,\q_0\rangle} =
    \frac{|\langle \z_t,\q_1\rangle|^2}{\langle \q_1,\q_1\rangle},
  \end{equation}
  which has degree at most two in $t$.

  The geodesic is contained in $\B$ if and only if this equation is
  identically zero; checking whether this holds amounts to checking
  whether $a=b=c=0$ where 
  $$
  \begin{array}{c}
    a=|\langle \ww-\vv,\q_0\rangle|^2/\langle\q_0,\q_0\rangle-|\langle
    \ww-\vv,\q_1\rangle|^2/\langle\q_1,\q_1\rangle\\ 
    b=\Re\left(\ \langle \vv,\q_0\rangle \langle \q_0,\ww-\vv\rangle/\langle\q_0,\q_0\rangle -
    \langle \vv,\q_1\rangle \langle \q_1,\ww-\vv\rangle/\langle\q_1,\q_1\rangle\ \right)\\
    c=|\langle \vv,\q_0\rangle|^2/\langle\q_0,\q_0\rangle-|\langle \vv,\q_1\rangle|^2/\langle\q_1,\q_1\rangle.
  \end{array}
  $$ 
  This can be decided with arithmetic in $\k$.

  In order to check that the segment is on a specific side of $\B$, we
  need to check that a strict inequality holds between the two sides
  of~\eqref{eq:tvals}. A natural way to do this is to find the
  solutions of~\eqref{eq:tvals}, and to check that they are all
  negative (resp. positive); the latter requires extracting
  square-roots, which is not quite arithmetic in $\k$. Alternatively,
  one can simply check that there is no sign change between the
  endpoints, in conjunction with a computation of the value of the
  polynomial at the point where its derivative vanishes.

  Suppose $p\in\B$ and $q\notin\B$. Then the geodesic is tangent to
  $\B$ if and only if 
  $$
  \Re\{\langle \ww,\q_0\rangle\langle \q_0,\vv\rangle\}/\langle\q_0,\q_0\rangle=\Re\{\langle
  \ww,\q_1\rangle\langle \q_1,\vv\rangle)\}/\langle\q_1,\q_1\rangle.
  $$ 
  This condition is
  checked by arithmetic in $\k$, since the vectors $\vv,\ww,\q_0,\q_1$
  are rational, and the Hermitian form is defined over $\k$. 
\EPf 

\section{Discrete groups and lattices}\label{sec:discrete}

\subsection{Lattices and arithmeticity}

A given lattice $\Gamma\subset {\rm SU}(n,1)$ can be conjugated in
complicated ways, but it is a classical fact that one can always
represent it as a subgroup of some ${\rm GL}(N,\k)$ where $\k$ is a
number field (of course one needs to require $n>1$ to exclude lattices
in ${\rm SU}(1,1)\simeq {\rm SL}(2,\R)$). For a general result along
these lines, see Theorem~7.67 in~\cite{raghunathan}. The smallest
field one can use is given by the field
  $\mathbb{Q}(\rm{tr}\,\rm{Ad}\,\Gamma)$ generated by traces
in the adjoint representation, see~Proposition~(12.2.1)
of~\cite{delignemostow} and also Section~\ref{sec:fields}.

In this manner, $\Gamma$ determines a $\k$-form of ${\rm SU}(n,1)$,
and the $\k$-forms of ${\rm SU}(n,1)$ are known to be obtained as
${\rm SU}(H)$ for some Hermitian form $H$ on $F^r$, where $F$ is a
division algebra with involution over 
a quadratic imaginary extension $\l$ of the totally real field $\k$
(see~\cite{weilinvolutions},~\cite{tits}). For dimension reasons, $r$
must divide $n+1$; in fact $n+1=rd$, where $d$ is the degree of
$F$. In particular, if $n=2$ then $r$ can only be $1$ or $3$. This
gives two types of lattices in ${\rm PU}(2,1)$, those related to
Hermitian forms over number fields ($d=1$, hence $r=3$) and those
related to division algebras ($d=3$, $r=1$). These are often referred
to as arithmetic lattices of the first and second type, respectively.

Note that the groups we construct in the present paper are clearly not
arithmetic of the second type, because they contain Fuchsian
subgroups.  For the relationship between Fuchsian subgroups and groups
of the second type, see~\cite{reznikov},~\cite{stover} (and also
\cite{mcreynoldsfinite}).  For lattices preserving a Hermitian form
over a number field, there is a fairly simple arithmeticity criterion,
which we now state for future reference. The following result is
essentially Mostow's Lemma~4.1 of \cite{mostowpacific} (following
Vinberg \cite{vinbergcriterion}, see also Corollary~12.2.8
of~\cite{delignemostow}). We refer to this statement as the {\bf
  Mostow--Vinberg arithmeticity criterion}.

\begin{proposition}\label{prop:crit}
  Let $\l$ be a purely imaginary quadratic extension of a totally real
  field $\k$, and $H$ a Hermitian form of signature $(n,1)$ defined
  over $\l$. Suppose $\Gamma\subset {\rm SU}(H;\mathcal{O}_\l)$ is a
  lattice. Then $\Gamma$ is arithmetic if and only if for all $\varphi
  \in {\rm Gal}(\l)$ not inducing the identity on
      $\mathbb{Q}(\rm{tr}\,\rm{Ad}\,\Gamma)$, the
      form $^\varphi H$ is definite.
\end{proposition}  
Note that when $\Gamma$ is as in Proposition~\ref{prop:crit} and it is not
arithmetic, the whole group of integral matrices ${\rm SU}(H;\mathcal{O}_\l)$ 
is non discrete in ${\rm SU}(H)$, and in particular
$\Gamma$ must have infinite index in ${\rm SU}(H;\mathcal{O}_\l)$.

\subsection{The Poincar\'e polyhedron theorem}\label{sec:poin}

Various versions of the Poincar\'e polyhedron theorem for complex
hyperbolic space have been given (see for example
\cite{mostowpacific}, \cite{dfp}, \cite{falbelparker},
\cite{parkerlivne}). For the purpose of the present paper, just as
in~\cite{mostowpacific}, we need to consider fundamental polyhedra for
coset decompositions, where the polyhedron is invariant under a
non-trivial subgroup.

Since the hypotheses as well as the conclusions of the theorem require
quite a bit of notation, we now give a detailed statement of the
Poincar\'e polyhedron theorem. A detailed proof can be found
in~\cite{JRP-book}. For simplicity, we only state this theorem for
finite-sided polyhedra, as this is sufficient when considering
lattices (for a more general statement that applies to locally finite
polyhedra, see~\cite{JRP-book}).

\smallskip
\noindent {\bf Polyhedra:} For the purpose of the present paper, we
only consider finite CW complexes which are regular, in the sense that
every attaching map is an embedding. In particular, the closure of each cell is
homeomorphic to a closed ball of the appropriate dimension, with
embedded boundary sphere. A {\bf (finite-sided) polyhedron} is the
geometric realization in $\chb n$ of such a complex with a single
top-dimensional cell of dimension $2n$. 

We will refer to closed cells of a polyhedron as facets, and denote by
$\F_k(E)$ the set of codimension $k$ facets of $E$. We give special
names to facets of each codimension: facets of codimension one, two,
three and four will be called {\bf sides}, {\bf ridges}, {\bf edges}
and {\bf vertices}, respectively. Because the focus of the present
paper is mainly about lattices, we will assume $E\cap \partial \ch n$
consists of finitely many vertices, which we call {\bf ideal
  vertices}.

Given a facet $f$, we denote by $f^\circ$ the relative interior of $f$
(equivalently, $f^\circ$ is the set of points that are on $f$ but not
in any other facet of the same codimension).
It follows from the regularity assumption that for each $k$-cell $f$,
the $(k-2)$-cells in $\partial f$ are contained in precisely two
$(k-1)$-cells of $\partial f$. In particular, in the above
terminology, a ridge of a polyhedron is on precisely two sides. 

In the context of complex hyperbolic space, there is no canonical
choice of hypersurfaces that can be used to bound polyhedra. Indeed,
as mentioned in Section~\ref{sec:chg}, there are no totally geodesic
real hypersurfaces in $\ch n$, $n\geqslant 2$.  In fact the polyhedra that
appear in this paper are bounded exclusively by bisectors, and they
actually have piecewise smooth boundary (one can relax this to allow
for more general faces, see~\cite{mostowpacific} and~\cite{JRP-book}
for a specific set of hypotheses).

\smallskip
\noindent {\bf Side pairings:} A map 
$\sigma:\F_1(E)\longrightarrow {\rm Isom}(\ch n)$ is called a 
{\bf side pairing} for $E$ if it satisfies the following conditions:
\begin{enumerate}
\item For each side $s\in\F_1(E)$ there is another side $s^-$ in $\F_1(E)$
so that $S=\sigma(s)$ maps $s$ onto $s^-$ preserving the cell structure.
Moreover, $\sigma(s^-)=S^{-1}$. In particular, if $s=s^-$ then $S=S^{-1}$ 
and so $S$ is an involution. We call $S^2=id$ a {\bf reflection relation}.
(In fact $s \neq s^-$ in all the cases we consider in this paper.)
\item If $s\in\F_1(E)$ and $\sigma(s)=S$ then $S^{-1}(E)\cap E=s$ and
$S^{-1}(E^\circ)\cap E^\circ=\emptyset$.
\item If $w\in s^\circ$ then there is an open neighborhood $U(w)$ of $w$ 
contained in $E\cup S^{-1}(E)$.
\end{enumerate}
We say that $S=\sigma(s)$ is the {\bf side-pairing map} associated to
the side $s\in\F_1(E)$.

We suppose the polyhedron $E$ is preserved by a finite group
$\Upsilon< {\rm Isom}(\ch n)$, which acts by cell-preserving
automorphisms. We assume moreover that we have a presentation of
$\Upsilon$ in terms of generators and relations (in the examples
treated in this paper, $\Upsilon$ will simply be a finite cyclic
group).

Let $\Gamma< {\rm Isom}(\ch n)$ be the group generated by $\Upsilon$
and the side-pairing maps.  In what follows, we will need to consider
$\Upsilon$-orbits of sides and ridges.  We suppose that the side
pairing $\sigma$ is {\bf compatible with} $\Upsilon$ in the following
sense: for all $s\in{\mathcal F}_1(E)$ and all $P\in\Upsilon$ we have
$\sigma(Ps)=P\sigma(s)P^{-1}$.

\smallskip
\noindent {\bf Ridge cycles:}
Consider a ridge $\rho_1\in\F_2(E)$. We know that $\rho_1$ is contained in 
exactly two sides of $E$, so we can write $\rho_1=s_1\cap s'$ where 
$s_1,\,s'\in\F_1(E)$. Let $S_1=\sigma(s_1)$ be the side pairing map
associated to $s_1$. Then we know that $S_1(s_1)$ is another side
$s_1^-$ of $E$. Since the side pairing maps are bijections preserving the 
cell structure, it must be the case that $\rho_2=S_1(\rho_1)$ is a ridge 
contained in $s_1^-$. 
In other words, there is a side $s_2\in\F_1(E)$ so that 
$\rho_2=s_2\cap s_1^-$. Let $S_2=\sigma(s_2)$ be the side pairing map 
associated to $s_2$. We repeat the above process and thereby obtain a 
sequence of ridges $\rho_i$, sides $s_i$ and side pairing maps 
$S_i$ so that $\rho_i=s_i\cap s_{i-1}^-$ and $S_i=\sigma(s_i)$. 

If $\rho_2=\rho_1$, we set $m=0$; otherwise let $m>0$ be the smallest
positive integer such that $\rho_{m+1}$ is in the same
$\Upsilon$-orbit as $\rho_1$ (there must exist such an $m$, since the
set of ridges is finite). 

Note that there is then a unique $P\in\Upsilon$ such that
$P\rho_{m+1}=\rho_1$ (or else there would be a non-trivial element of
$\Upsilon$ fixing $\rho_1$ pointwise and preserving the pair of sides
containing it). We then have $\rho_1=s_1\cap Ps_m^-$, and the other
side $s'$ containing $\rho_1$ must be $Ps_m^-$.

We say that $(\rho_1,\,\ldots,\,\rho_{m+1})$ is the {\bf ridge cycle}
of $\rho_1$ and we define the {\bf cycle transformation} $T=T(\rho_1)$
of $\rho_1$ to be $T=P\circ S_m\circ\cdots\circ S_1$. Then $T$ maps
$\rho_1$ to itself:
$$
  T(\rho_1)=P\circ S_m\circ\cdots\circ S_2\circ S_1(\rho_1)
  =\cdots=P\circ S_m(\rho_m)=P(\rho_{m+1})=\rho_1.
$$
Note that $T$ may not act as the identity on $\rho_1$ and, even 
if it does, then it may not be the identity on the whole of $\ch n$. We 
assume that $T$ has finite order $\ell$.
The relation $T^\ell=id$ is called the {\bf cycle relation} associated to 
$\rho_1$.

Note that the sequence $\rho_1,\,\ldots,\,\rho_m$ is entirely determined by 
$\rho_1$ and the choice of the first side $s_1$.
Observe that if we had started at another ridge in the cycle, say $\rho_i$
then, using our hypothesis that the side pairings are compatible with
$\Upsilon$, we would obtain the ridge cycle
$(\rho_i,\,\rho_{i+1},\,\ldots,\,\rho_m,\,P^{-1}\rho_1,\,P^{-1}\rho_2,\,
\ldots,\,P^{-1}\rho_i)$. Because $\Upsilon$ is compatible with the 
side pairings, the cycle transformation of $\rho_i$ is 
\begin{eqnarray*}
  T_i & = & P\circ(P^{-1}S_{i-1}P)\circ\cdots \circ(P^{-1}S_1P)\circ S_m\circ\cdots\circ S_1 \\
  & = & S_{i-1}\circ \cdots\circ S_1\circ P\circ S_m\circ\cdots \circ S_i.
\end{eqnarray*}
This is just a cyclic permutation of $T$.
Likewise, if we had started at $\rho_1$ but had chosen the first side pairing
to be the one associated to the other side $s'$ instead, we would have 
obtained (up to elements of $\Upsilon$) the same sequence of sides 
in the reverse order, with the side pairing maps inverted.
Furthermore, if $Q\in\Upsilon$ then the ridge cycle of $Q\rho_1$ is
$(Q\rho_1,\,\ldots,\,Q\rho_{m+1})$ and the corresponding cycle transformation 
is $QTQ^{-1}$. Hence it suffices to consider the cycle associated to a single 
ridge in each $\Upsilon$-orbit of ridge cycles.

Writing out $T$ in terms of $P$ and the $S_i$, we let 
${\mathcal C}={\mathcal C}(\rho_1)$ be the collection of suffix subwords of
$T^\ell$, that is
$$
  {\mathcal C}(\rho_1)=\Bigl\{S_i\circ\cdots\circ S_1\circ T^j\ :\  
  0\leqslant i \leqslant m-1,\ 0\leqslant j\leqslant   \ell-1
  \Bigr\}
$$
where $i=0$ means we write none of the $S_i$ (this includes the case where
$m=0$) and so $i=j=0$ corresponds to the identity map.

For $i\in\{1,\,\ldots,\, m\}$ we have $\rho_i=s_i\cap s_{i-1}^-$, where
$s_0^-=Ps_m^-$. From the side pairing conditions, 
$s_i=E\cap S_i^{-1}(E)$ and $s_i^-=E\cap S_{i-1}(E)$ and so
$\rho_i\subset E\cap S_i^{-1}(E)\cap S_{i-1}(E)$. Furthermore, for $2\leqslant i\leqslant m$ we have 
$\rho_i=S_{i-1}\circ\cdots\circ S_1(\rho_1)$. Therefore
\begin{eqnarray*}
  \rho_1 & = & S_1^{-1}\circ\cdots\circ S_{i-1}^{-1}(\rho_i) \\
  & \subset  & S_1^{-1}\circ\cdots\circ S_{i-1}^{-1}
  \bigl(S_i^{-1}(E)\cap E\cap S_{i-1}(E)\bigr) \\
  & = &  \Bigl(S_1^{-1}\circ\cdots\circ S_i^{-1}(E)\Bigr)
  \cap \Bigl(S_1^{-1}\circ\cdots\circ S_{i-1}^{-1}(E)\Bigr)
  \cap \Bigl(S_1^{-1}\circ\cdots\circ S_{i-2}^{-1}(E)\Bigr).
\end{eqnarray*}
Hence, $\rho_1$ is contained in $S_1^{-1}\circ\cdots\circ
S_{i-1}^{-1}(E)$ for all $i$ between $1$ and $m$ where, as before,
$i=1$ corresponds to the identity map.  Note that, since $P(E)=E$ we
have
$$
   T^{-1}(E)
   =S_1^{-1}\circ\cdots\circ S_{m}^{-1}\circ P^{-1}(E)
   =S_1^{-1}\circ\cdots\circ S_{m}^{-1}(E).
$$
Hence we have shown that
$$
   \rho_1\subset \bigcap _{C\in{\mathcal C}(\rho_1)}C^{-1}(E).
$$
We say that $E$ and $\Gamma$ satisfy the {\bf cycle condition} at $\rho_1$ 
if this intersection is precisely $\rho_1$ and
all these copies of $E$ tessellate around $\rho_1$. That is: 
\begin{enumerate}
\item $$
  \rho_1= \bigcap _{C\in{\mathcal C}(\rho_1)}C^{-1}(E).
  $$
\item If $C_1,\,C_2\in {\mathcal C}(\rho_1)$ with $C_1\neq C_2$ then
$C_1^{-1}(E^\circ)\cap C_2^{-1}(E^\circ)=\emptyset$.
\item For all $w\in\rho_1^\circ$ there exists an open neighborhood $U(w)$ of $w$
so that 
$$
  U(w)\subset \bigcup_{C\in{\mathcal C}(\rho_1)} C^{-1}(E).
$$
\end{enumerate}

\smallskip
\noindent {\bf Consistent system of horoballs:} When $E$ has cusps, we
need an extra hypothesis, related to metric completeness of the
quotient of the polyhedron under the side-pairing maps. Let
$\xi_1,\dots,\xi_m$ be the cusps of $E$.  We will assume the existence
of a {\bf consistent system of horoballs}, which is a collection
$\{U_1,...,U_m\}$, where each $U_j$ is a horoball based at $\xi_j$
which is preserved by the stabilizer of $\xi_j$ in $\Gamma$.
By shrinking if necessary, we may suppose that a consistent system is
made up of pairwise disjoint horoballs.  The existence of a consistent
system of horoballs can be checked by verifying that all cycle
transformations fixing a given cusp are non-loxodromic (since in that
case they automatically preserve every horoball based at that ideal
vertex).

\smallskip
Then the statement of the complex hyperbolic Poincar\'e polyhedron theorem
(see~\cite{mostowpacific} or~\cite{JRP-book}) is the following.

\begin{theorem}\label{thm:poin}
  Suppose $E$ is a smoothly embedded finite-sided polyhedron in $\ch
  n$, together with a side pairing $\sigma:\F_1(E)\longrightarrow {\rm
    Isom}(\ch n)$. Let $\Upsilon<{\rm Isom}(\ch n)$ be a group of
  automorphisms of $E$.  Let $\Gamma$ be the group generated by
  $\Upsilon$ and the side-pairing maps. Suppose the cycle condition is
  satisfied for each ridge in ${\mathcal F}_2(E)$, and that there is a
  consistent system of horoballs at the cusps of $E$ (if it has any).

  Then the images of $E$ under the cosets of $\Upsilon$ in $\Gamma$
  tessellate $\ch n$. That is
\begin{enumerate}
\item 
$$
  \bigcup _{A\in\Gamma} A(E)=\ch n.
$$
\item If $A\in\Gamma-\Upsilon$ then $E^\circ\cap A(E^\circ)=\emptyset$.
\end{enumerate}
Moreover, $\Gamma$ is discrete and a fundamental domain for its action
on $\ch n$ is obtained by intersecting $E$ with a fundamental domain
for $\Upsilon$.  

Finally, one obtains a presentation for $\Gamma$ in terms of the
generators given by the side pairing maps together with a generating
set for $\Upsilon$; the relations are given by
the reflection relations, the cycle relations and the relations in a
presentation for $\Upsilon$.
\end{theorem}

\smallskip
\noindent {\bf Finite volume:} Note that the statement of the
Poincar\'e polyhedron theorem says nothing about the volume of the
quotient $\Gamma\setminus \ch n$, which is also the volume of
$E$. This volume is of course finite when $E$ is entirely contained in
$\ch n$, in which case $\Gamma$ is a cocompact lattice in ${\rm PU}(n,1)$.

Some of the groups we study in this paper are not cocompact (namely,
when $p=4,6$). In that case, $E$ has some ideal vertices, but one
easily sees that the volume of $E$ is finite by studying the structure
of the stabilizer of the ideal vertices. Let $x$ be an ideal vertex
with stabilizer $\Gamma_x$.  By the existence of a consistent system
of horoballs, the ideal point corresponding to $x$ in the quotient
$\Gamma\setminus \ch n$ has a neighborhood diffeomorphic to
$(0,+\infty)\times H_x$, where $H_x$ is the quotient of any horosphere
based at $x$ by the action of $\Gamma_x$. It is well known that the
fact that $\Gamma_x$ acts cocompactly on $\partial \ch 2 \setminus \{
x \}$ (or equivalently, on horospheres based at $x$) implies that cusp
neighborhoods in the quotient have finite volume (see for example
Lemma~5.2 in~\cite{hersonskypaulin}). This is clear for our polyhedra,
whose cusp cross-sections are compact (see
Section~\ref{sec:stabilisers}, Figure~\ref{fig:cusplinks} for their
explicit description).

\begin{figure}[htbp]
  \epsfig{figure=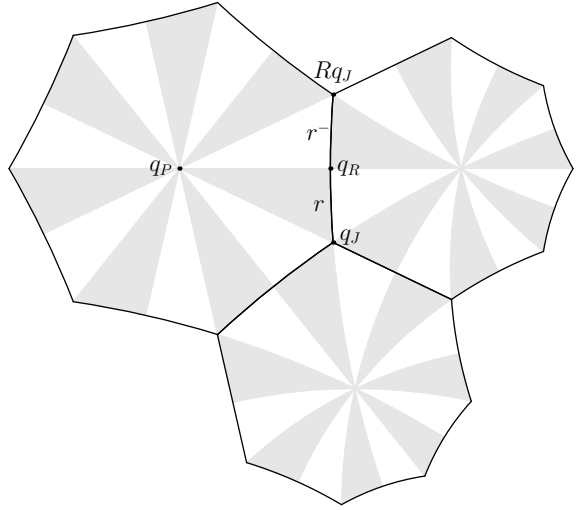, height=0.3\textheight}
  \centering
  \caption{Toy model for fundamental domains for coset
    decompositions.}\label{fig:toymodel}
\end{figure}

\smallskip
\noindent {\bf Toy model:} To illustrate this theorem and to help
understand our family of examples, it is instructive to work through a
toy example with $n=1$. Consider a triangle in the Poincar\'e disk
$\ch 1$ with vertices $q_R$ with angle $\pi/2$, $q_J$ with angle
$\pi/3$ and $q_P$ with angle $\pi/7$ (see
Figure~\ref{fig:toymodel}). Let $R$, $J$ and $P$ be elements of ${\rm
  Isom}(\ch 1)$, with orders $2$, $3$ and $7$, fixing the respective
points and satisfying $P=RJ$. Let $\Gamma$ be the group generated by
$R$, $J$ and $P$. The usual fundamental domain $D$ for this group
consists of this triangle together with its image under reflection in
the side joining $q_P$ and $q_R$. There is a natural $P$-invariant
hyperbolic heptagon $E$, obtained as the union of the seven images of
$D$ under powers of $P$.  This heptagon is a fundamental domain for
the cosets of $\Upsilon=\langle P\rangle$ in $\Gamma$.

The vertices of the heptagon $E$ are of the form $P^kq_J$, where
$k=0,\,\ldots,\,6$, and the midpoints of its sides are given by the
points $P^kq_R$. Now let $r$ be the edge from $q_J$ to $q_R$ and let
$r^-$ be the edge from $q_R$ to $Pq_J=Rq_J$. Then we define the side
pairing map $R:r\longrightarrow r^-$, and extend this to the other
sides in a way that is compatible with $\Upsilon$. Namely, we have
side pairings $P^kRP^{-k}:P^kr\longrightarrow P^kr^-$. We obtain two
vertex cycles, which we write in terms of the edges and the vertices.
$$
  \begin{array}{|r|l|l|}
    \hline
    R & r\cap r^- & q_R \\
    & r\cap r^- & q_R \\
    \hline
    R & r\cap P^{-1}r^- & q_J \\
    P^{-1} & Pr\cap r^- & Rq_J \\
    & r\cap P^{-1}r^- & q_J \\
    \hline
  \end{array}
$$
It is easy to check local tessellation around the vertices and
so Poincar\'e's theorem shows that the heptagons tile the
Poincar\'e disk. 

For $r\cap r^-$, the cycle transformation is $R$, and the cycle
relation is $R^2=id$. For $r\cap P^{-1}r^-$, the cycle transformation
is $P^{-1}R$ and the cycle relation is $(P^{-1}R)^3=id$.  Adding the
generator $P$ of $\Upsilon$ and its relation $P^7=id$ gives the well
known presentation
$$
  \Gamma=\langle R,\, P\ | \ R^2=(P^{-1}R)^3=P^7=id\rangle.
$$
Note that a fundamental domain for $\Gamma$ is obtained by
intersecting $E$ with a fundamental domain for $\Upsilon$ (one example
of such a fundamental domain is of course given simply by $D$).

We can also calculate the orbifold Euler characteristic $\chi$ of
$\Gamma\backslash\ch 1$ as follows. We consider 
each $\Gamma$-orbit of facets of $E$ and we weight
them by the reciprocal of the order of their stabilizer:
$$
  \begin{array}{|l|l|l|}
    \hline
    \hbox{Facet} & \hbox{Stabilizer} & \hbox{Order} \\
    \hline
    q_R & \langle R\rangle & 2 \\
    q_J & \langle P^{-1}R\rangle & 3 \\
    \hline
    r & id & 1 \\
    \hline
    E & \langle P\rangle & 7 \\
    \hline
  \end{array}
$$
Thus 
$$
  \chi=\left(\frac{1}{2}+\frac{1}{3}\right)-1+\frac{1}{7}=\frac{-1}{42}.
$$

We now briefly review standard techniques used to verify the
hypotheses on ridge cycles.

\smallskip
\noindent {\bf Tessellating around Giraud disks. } For cycles around
Giraud ridges the conditions of the Poincar\'e polyhedron theorem
are easily checked, as we now recall.
 
Let $\G=\G(p_0,p_1,p_2)$ be a Giraud disk. Then $\G$ defines three
regions $Y_0$, $Y_1$ and $Y_2$ given, for indices $j=0,\,1,\,2$ 
taken mod 3, by:
\begin{equation}\label{eq:y-piece}
  Y_j=\bigl\{u\in\ch 2\ :\ d(u,p_j)\leqslant d(u,p_{j+1}),\ 
  d(u,p_j)\leqslant d(u, p_{j-1})\bigr\}
\end{equation}
Any point in $\ch 2$ is contained in (at least) one of these three
regions according to the minimum of $d(u,p_j)$ for $j=0,\,1,\,2$. The
interior $Y^\circ_j$ of $Y_j$ is the set of points where $d(u,y_j)$ is
strictly smaller than $d(u,p_{j+1})$ and $d(u,p_{j-1})$. Clearly the
interiors of these three regions are disjoint.

\begin{lemma}\label{lem:giraud-tess}
  Let $E$ be a polyhedron bounded by bisectors and let $\rho=s_1\cap
  s_2$ be a Giraud ridge of $E$ contained in the sides $s_1$ and $s_2$
  of $E$. Let $S_1=\sigma(s_1)$ and $S_2=\sigma(s_2)$. Let $Y_0$,
  $Y_1$ and $Y_2$ be the three regions given in \eqref{eq:y-piece},
  defined by the Giraud disk containing $\rho$. Suppose the $Y_j$ each
  contain exactly one of $E$, $S_1^{-1}(E)$ and $S_2^{-1}(E)$ and that
  $S_1^{-1}(E)\cap S_2^{-1}(E)$ is contained in the third bisector
  containing $\rho$. Then $E$, $S_1^{-1}(E)$ and $S_2^{-1}(E)$
  tessellate a neighborhod of the interior of $\rho$.
\end{lemma}

\Pf 
Without loss of generality, suppose that $E\subset Y_0$,
$S_1^{-1}(E)\subset Y_1$ and $S_2^{-1}(E)\subset Y_2$. Since the $Y_j$
have disjoint interiors, it is clear that $E^\circ$,
$S_1^{-1}(E^\circ)$ and $S_2^{-1}(E^\circ)$ are disjoint.

We must show that the three copies of $E$ cover a neighborhood of the
interior of $\rho$. Suppose that $w\in\rho^\circ$ and $U(w)$ is a
neighborhood of $w$. Then any point in $U(w)$ is in (at least) one of
$Y_0$, $Y_1$, $Y_2$.  If the point is sufficiently close to $\rho$
then it is in $E$, $S_1^{-1}(E)$ or $S_2^{-1}(E)$ respectively.  
\EPf

\smallskip
\noindent
{\bf Tessellating around complex lines. } Cycles around ridges
contained in complex lines can be more complicated than ridges
contained in Giraud disks. However, by looking at complex lines
orthogonal to the intersection, tessellation is reduced to an angle
condition resembling the classical version of the Poincar\'e
polyhedron theorem for the hyperbolic plane.

Let $\B$ and $\B'$ be two bisectors that intersect in a complex line
$C$. Suppose $E$ is a polyhedron contained in the intersection of two
half-spaces defined by $\B$ and $\B'$ and that these bisectors contain
sides $s$ and $s'$ of $E$ intersecting in a simply connected ridge
$\rho$ in $C$. Using Lemma~\ref{lem:orthslice} any complex line
$C^\perp$ orthogonal to $C$ intersects $\B$ and $\B'$ along geodesics,
and $C^\perp\cap E$ is contained in a wedge bounded by these
geodesics, as $\B$ and $\B'$ both bound $E$.  We must keep track of
the angle subtended by this wedge. It is important to note that this
angle will depend on the point where $C$ and $C^\perp$ intersect, but
that the total angle subtended over a ridge cycle will remain the
same.

In particular, suppose $T$ is the cycle transformation of $\rho$.  By
construction, $T$ maps $C$ to itself. There are two possibilities:
either $T$ fixes $C$ pointwise, and so is a complex reflection in $C$,
or $T$ has a unique fixed point $o$ in $C$. In the latter case, the
point $o$ must lie in $\rho$ as $T$ is a symmetry of $\rho$, which is
simply connected.  In either case, if $o$ is a point of $\rho$ fixed
by $T$ then the complex line $C^\perp_o$ through $o$ orthogonal to $C$
is mapped to itself by $T$.

Suppose the cycle transformation is $T=P\circ S_m\circ\cdots\circ
S_1$.  The definition of a ridge cycle leads to ridges $\rho_i=s_i\cap
s_{i-1}^-$ where $S_i=\sigma(s_i)$. The intersection of $C^\perp_o$
with $S_1^{-1}\circ\cdots\circ S_{i-1}^{-1}(E)$ for $1\leqslant
i\leqslant m$ is a wedge, say with angle $\beta_i$, bounded by the
intersection of $C^\perp_o$ with the sides $S_1^{-1}\circ\cdots\circ
S_{i-1}^{-1}(s_i)$ and $S_1^{-1}\circ\cdots\circ
S_{i-1}^{-1}(s_{i-1}^-)$. These wedges fit together to give a larger
wedge bounded by the intersection of $C^\perp_o$ with a bisector and
its image under $T$. In this larger wedge, the total angle at $o$
subtended by copies of $E$ under the cycle is
$\beta_1+\cdots+\beta_m$.

\begin{lemma}\label{lem:cx-tess}
  Let $E$ be a polyhedron bounded by bisectors and let $\rho$ be a
  ridge of $E$ contained in a complex line $C$. Let $T$ be the cycle
  transformation of $\rho$ and let $o$ be a fixed point of $T$ in
  $\rho$.  Suppose that $T$ has order $\ell>0$ and that $T$ acts on
  $C^\perp_o$ as a rotation by angle $2\pi/\ell$. Suppose that the
  total angle at $o$ in $C^\perp_o$ subtended by copies of $E$ under
  the cycle is $2\pi/\ell$.  Then any point $w\in\rho^\circ$ has an
  open neighborhood tessellated by images of $E$.
\end{lemma}

\Pf 
For simplicity, we begin by supposing the ridge cycle has length one and
$T=P\circ S$ where $S$ is the side-pairing map of $s$ and
$P\in\Upsilon$. This means we only have to show that
$E,\,T^{-1}(E),\,\ldots,\,T^{-(\ell-1)}(E)$ tessellate around $\rho$.
(This is the only case we need in the applications in this paper.)

First consider the action of $T$ on $C^\perp_o$. Let $\B$ and $\B'=T(\B)$ 
be the bisectors containing the sides $s$ and $s'=T(s)$ of $E$. 
Write $b_o=\B\cap C^\perp_o$ and $b'_o=\B'\cap C^\perp_o$. Since $T$ acts on 
$C^\perp_o$ as a rotation by $2\pi/\ell$ then $b'_o=T(b_o)$ and the arcs 
$b_o$, $b'_0$ bound a wedge with angle $2\pi/\ell$ containing $C^\perp\cap E$. 
Applying powers of $T$ we obtain $\ell$ wedges containing 
images of $E$ that tessellate a neighborhood of $o$ in $C^\perp_o$. 

When $T$ is a complex reflection in $C$, this argument applies
to all points of $\rho$ and the result follows.

We now consider the case where the fixed point $o$ of $T$ is unique.
For a general point $u\in C\cap E$ consider $C^\perp_u$, the orthogonal complex 
line to $C$ at $u$. Since $T$ does not fix $u$, we see that $T$
sends $C^\perp_u$ to $C^\perp_{Tu}$. However, by continuity, the 
intersection of $C^\perp_u$ with $E$ is contained in some wedge with apex $u$.
The angle may change as $u$ varies but will always be positive
since $\B\cap\B'$ is $C$. Let $k$ (which divides $\ell$)
be the smallest positive integer so that $T^k$ is the identity
on $C$. Then the angles of these wedges at all the images of
$C^\perp_u,\,T(C^\perp_u),\,\ldots,\, T^k(C^\perp_u)$ add up to
$2k\pi/\ell$. Since $T^k$ is a complex reflection in $C$ with
angle $2k\pi/\ell$ we see that the intersection of 
$E,\,T(E),\,\ldots,\,T^{\ell-1}(E)$ with $C^\perp_u$ tessellate a
neighborhood of $u$ in $C^\perp_u$. 

By carrying out this process for all points of $C\cap E$, we see that
for any point $u$ in the (relative) interior of $C\cap E$ there
is an open neighborhood of $u$ tessellated by copies of $E$.

Now consider the general case where the length of the ridge cycle is $m\geqslant 1$
and the cycle transformation is $T=P\circ S_m\circ\cdots\circ S_1$. 
In this case, we need to keep track of wedges in
$T^{-j}\circ S_1^{-1}\circ\cdots \circ S_i^{-1}(C^\perp_u)$ for $0\leqslant i\leqslant m-1$
and $0\leqslant j\leqslant \ell-1$. Nevertheless, the same argument gives the result.
\EPf

\section{Construction of the fundamental polyhedron} \label{sec:domain}

In section~\ref{sec:notation} and~\ref{sec:wordnotation}, we define
the six groups that appear in Theorem~\ref{thm:main}, and establish
basic notation that will be used throughout the paper. 

In section~\ref{sec:braiding}, we study braiding properties of
some complex reflections in our groups, which are used in an essential
way in order to build our fundamental polyhedra, as explained in
section~\ref{sec:algo}. These braiding pairs of reflections are also
important because they will correspond to the stabilizers of various
vertices and complex ridges of our fundamental polyhedron, see
section~\ref{sec:stabilisers}; this will allow us to compute the
orbifold Euler characteristic, see section~\ref{sec:euler}.

We define our polyhedra in section~\ref{sec:defE} by building the
$k$-skeleton for increasing values of $k$. The general strategy for
this definition is described in section~\ref{sec:realization}; an
important point for this inductive definition to make sense is that
the boundary of every cell is a topological (piecewise smooth)
sphere. The only difficult part of this verification is the one that
concerns the boundary of the whole polyhedron. We show that the
boundary of $E$ is homeomorphic to $S^3$ in section~\ref{sec:sphere}.

We then start from the combinatorial model of the polyhedron (this is
described in section~\ref{sec:model}), and prove that the geometric
realization is well-defined, and that it gives an embedding of the
combinatorial model (sections~\ref{sec:vertices}
through~\ref{sec:3-cells}).

\subsection{Generators} \label{sec:notation}

In this section, we give explicit matrices that generate our
groups. For more details on sporadic triangle groups,
see~\cite{parkerpaupert} and \cite{dpp1}.

Recall from the introduction that our groups are generated by a
complex reflection $R_1$ and a regular elliptic element $J$ of order
three. We write
$$
R_2=JR_1J^{-1},\quad R_3=JR_2J^{-1}.
$$ 

We write ${\mathbf n}_j$ for polar vector to the mirror of the complex
reflection $R_j$ ($j=1,2,3$). Matrix representatives for $R_j$ in
${\rm SU}(2,1)$ will have eigenvalues $a^2,\overline{a},\overline{a}$
where
$$
  a = e^{2\pi i/3p}.
$$
Since $R_{j+1}=JR_jJ^{-1}$ we see that ${\mathbf n}_{j+1}=J{\mathbf
  n}_j$ (with $j$ mod 3). In the cases that interest us,
  $({\bf n}_1,{\bf n}_2,{\bf n}_3)$ forms a basis for $\C^3$, which we will use throughout the paper.
The matrix for $J$ in this basis is then simply the permutation matrix:
$$
  J = \left[
    \begin{matrix} 
      0 & 0 & 1 \\ 
      1 & 0 & 0 \\ 
      0 & 1 & 0
    \end{matrix}
    \right]
$$
and the Hermitian form is given by a matrix of the form 
$$
  H = \left[\begin{matrix}
      \alpha & \beta & \overline{\beta}\\
      \overline{\beta} & \alpha & \beta\\
      \beta & \overline{\beta} & \alpha
    \end{matrix}\right].
$$ 
It is convenient to choose
$$
  \alpha=4\sin^2(\pi/p)=2-a^3-\overline{a}^3,
$$ 
in which case the condition ${\rm tr}(R_1J)=\tau$ imposes
$$
  \beta=(\overline{a}^2-a)\tau.
$$
The matrix of the reflection $R_1$ (adjusted to have determinant $1$,
since we want to work with ${\rm SU}(2,1)$), is easily seen to be
$$
  R_1 =
  \left[\begin{matrix} 
      a^2 & \tau & -a\overline{\tau} \\
      0 & \overline{a} & 0 \\ 
      0 & 0 & \overline{a} 
    \end{matrix}\right]
$$
and $R_2=JR_1J^{-1}$, $R_3=J^{-1}R_1J$.

\begin{figure}
\centering
\epsfig{figure=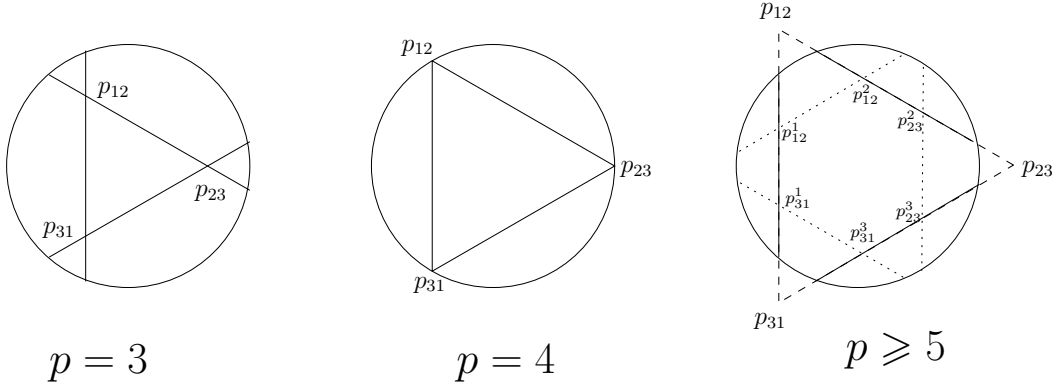, width=0.9\textwidth}
\caption{Relative positions of the mirrors of $R_1$, $R_2$, $R_3$. For
  $p\geqslant 5$, the dotted lines correspond to the lines polar to
  $p_{jk}$, or in other words the common perpendicular between the
  mirrors of $R_j$ and $R_k$.}\label{fig:mirrors}
\end{figure}

The configuration of three complex lines given by the mirrors of
$R_1$, $R_2$ and $R_3$ varies with $p$. Note that the relative
position of the mirrors of $R_1$ and $R_2$ is controlled by the
restriction of the Hermitian form to ${\rm Span}(\n_1,\n_2)$, i.e by the
sign of the determinant of the $2\times 2$ matrix in the upper left
corner of $H$.
Using the specific value $\tau=-(1+i\sqrt{7})/2$ (in fact all we need is that
$|\tau|^2=2$), we have
$$
  \alpha^2-|\beta|^2=8\sin^2(\pi/p)\bigl(2\sin^2(\pi/p)-1\bigr)=
  -8\sin^2(\pi/p)\cos(2\pi/p).
$$ 
Hence the mirrors intersect inside complex hyperbolic space only
for $p=3$, they intersect at infinity when $p=4$, and they are
ultraparallel when $p\geqslant 5$ (see Figure~\ref{fig:mirrors}).

For future reference, we define 
$$
P = R_1J
$$
and 
\begin{equation}\label{eq:S1}
S_1 = P^2 R_1 P^{-2} R_1 P^2
\end{equation}
which will turn out to be key isometries when constructing our
fundamental domains. Since $P=R_1J$ has trace 
$\tau=-(1+i\sqrt{7})/2=e^{-2\pi i/7}+e^{-4\pi i/7}+e^{-8\pi i/7}$,
we immediately see that $P$ is a regular elliptic map of order $7$,
with isolated fixed point
\begin{equation}
  {\bf p} = \left[\begin{matrix} 
     ae^{-2\pi i/7} \\
     1 \\
     \overline{a}e^{2\pi i/7}
  \end{matrix}\right].
\end{equation}
For completeness we give an explicit matrix
for $S_1$, which is also regular elliptic of order $2p$:
$$
  S_1 = \left[\begin{matrix}
      -a & \overline{a}\,\overline{\tau} & -1 \\
      0 & -a\tau+\overline{a}^2\overline{\tau} 
      & a^2\overline{\tau}-a^2-\overline{a} \\
      0 & -1-\overline{a}^3+\overline{a}^3\tau 
      & a\tau-\overline{a}^2\overline{\tau}
    \end{matrix}\right].
$$

\subsection{Word notation} \label{sec:wordnotation}

We shall often use word notation and write $j$ for $R_j$ and $\bar j$
for $R_j^{-1}$, so that a word in $R_1$, $R_2$, $R_3$ and their
inverses is described by a sequence of (possibly overlined)
integers. Note that there is very little confusion possible between
$R_j^{-1}$ and the complex conjugate matrix $\overline{R}_j$ (which we
shall never use in this paper). For example, $23\bar 2$ denotes 
$R_2 R_3 R_2^{-1}$, and $1\bar 3 23\bar{1}$
denotes $R_1R_3^{-1}R_2R_3R_1^{-1}$.

When an isometry given by a word $w$ in 1, 2, 3  has an
isolated fixed point, we denote this point by $p_w$ or $q_w$. In
particular, consider two reflections with distinct mirrors that can be
expressed as words $u$, $v$. If $w=uv$ is conjugate in $\Gamma$ to $12$
(this is a word, not a number!) then we denote by $p_w$ the
intersection point of their mirrors of reflections (or the
intersection of the extension of their mirrors to projective lines,
i.e. these points may be in $\cp n$ rather than $\ch n$). Similarly,
if $w=uv$ is conjugate to $123\bar{2}$ then we denote the fixed point of
$w$ by $q_w$. The reason for using different letters of the
alphabet is that points the of the form $p_w$ and those of the form
$q_w$ are in different group orbits. We will refer to
vertices of our polyhedra that have the form $p_w$ for some word $w$ as
$p_*$-vertices (and similarly for $q_*$-vertices).

Let $\n_{w}$ denote a lift to $\C^3$ of the corresponding point
$p_w$. Of course, by definition, these lifts are only determined up to
multiplication by a scalar; so we give some explicit formulae for
future reference.
$$
  {\bf n}_{23\bar{2}}=R_2\n_3=\left[\begin{matrix} 
      0 \\ \tau \\ \overline{a}
    \end{matrix}\right],
  \quad\n_{123\bar 1}=a R_1\n_{23}=\left[\begin{matrix} 
      a^6+a^3-a^3\tau+\overline{\tau} \\
      a^2\overline{\tau}-\overline{a} \\
      \overline{a}^2\tau-a 
    \end{matrix}\right],
$$
$$
  \n_{123\bar 2}=\left[\begin{matrix} 
      \overline{a}^2+a\tau \\
      \overline{a}^3-a^3-\overline{\tau} \\
      \overline{a}^4\overline{\tau}+ \overline{a} 
    \end{matrix}\right],\quad
  \n_{1\bar 3 23}=\left[\begin{matrix}
      a^2+\overline{a}\overline{\tau} \\
      a^4\tau+a \\
      a^3-\overline{a}^3-\tau
    \end{matrix}\right], \quad
  \n_{1\bar 3 23\bar 1}=a R_1R_3^{-1}\ \n_2=\left[\begin{matrix}
      a^2 \\
      a \\
      \overline{\tau}
    \end{matrix}\right].
$$

We will also denote by $m_w$ the mirror of the complex reflection
corresponding to a word $w$, so that $\n_w$ is a polar vector to $m_w$
(in other words, $m_w$ corresponds to $\n_w^\perp$). We will also
extend this notation to group elements that have a complex reflection
as a power, so that (when $p\geqslant 5$) $m_{23}$ denotes the mirror
of $(R_2R_3)^2$, and (when $p\geqslant 8$) $m_{123\bar2}$ denotes the
mirror of $(R_1R_2R_3R_2^{-1})^3$.

\subsection{Higher braiding}\label{sec:braiding}

The key to the general procedure to build the fundamental domains for
our groups (see section~\ref{sec:algo}) is the following result from
\cite{parkerpaupert} (see also section~2.2 of~\cite{mostowpacific}).

\begin{proposition}[Proposition 4.4 of \cite{parkerpaupert}]\label{prop:pencil}
Suppose that $R_x$ and $R_y$ are complex reflections in ${\rm
  SU}(2,1)$ both with eigenvalues $e^{2i\psi/3}$, $e^{-i\psi/3}$,
$e^{-i\psi/3}$. Let $L_x$ and $L_y$ be the $e^{-i\psi/3}$ eigenspaces
of $R_x$ and $R_y$, respectively. Then $L_x\cap L_y$ is an
$e^{-2i\psi/3}$ eigenspace of $R_xR_y$. Suppose $R_xR_y$ is
non-loxodromic and that its other eigenvalues are $-e^{i\psi/3\pm
  i\phi}$. Then the group $\Gamma_{xy}=\langle
  R_x,R_y\rangle$ acts on the orthogonal complement of $L_x\cap L_y$
  as (the orientation preserving subgroup of) a $(\psi/2,\psi/2,\phi)$
  triangle group $\Delta_{xy}$.
\end{proposition}

Consider the action of $\langle R_x,R_y\rangle$ on the
projectivization of $(L_x\cap L_y)^\perp$, which we denote by $m_{xy}$.  The space
$m_{xy}$ is a sphere, Euclidean plane or hyperbolic plane
respectively, depending on whether $m_x$ and $m_y$ intersect, are
asymptotic or are ultraparallel. Note that there is a close connection
between the values of $\psi$ and $\phi$ and these three cases (see
Figure~\ref{fig:mirrors}).  In what follows, we will suppose
$\psi=2\pi/p$ and $\phi=2\pi/q$ (which will be the case in our
examples). This means that the triangle associated to $\Delta_{xy}$
has internal angles $\pi/p$, $\pi/p$, $2\pi/q$.  In particular,
$R_xR_y$ acts on $m_{xy}$ as a rotation through angle
$2\phi=4\pi/q$. If $q$ is even then $\Delta_{xy}$ is a $(p,p,q/2)$
triangle group. Therefore $(R_xR_y)^{q/2}$ acts as the identity on
$m_{xy}$.  On the other hand, if $q$ is odd this triangle can be
divided into two $(2,p,q)$ triangles, and $\Delta_{xy}$ is a $(2,p,q)$
triangle group.  In this case a simple geometric argument in $m_{xy}$
shows that $(R_xR_y)^{(q-1)/2}(m_x\cap m_{xy})=(m_y\cap m_{xy})$.
Since $m_x$ and $m_y$ are orthogonal to $m_{xy}$ we immediately see
that $R_x$ and $R_y$ satisfy a (generalized) braid relation.  When $q$
is even, this braid relation is 
\begin{equation}\label{eq:braidingqeven}
(R_xR_y)^{q/2}=(R_yR_x)^{q/2}
 \end{equation}
and when $q$ is odd, it is 
\begin{equation}\label{eq:braidingqodd}
(R_xR_y)^{(q-1)/2}R_x=R_y(R_xR_y)^{(q-1)/2}.
\end{equation}
Furthermore, by examining the eigenvalues we can see that
$(R_xR_y)^{q/2}$ (respectively $(R_xR_y)^q$) has a repeated
eigenvalue. Therefore it is either is a complex reflection (in a line
or a point), the angle depending on $p$ and $q$ or it is possibly
parabolic. The latter case arises when $p_{xy}$ lies on $\partial{\bf
  H}^2_{\mathbb C}$. In either case, a power of $R_xR_y$ is in the
center of $\langle R_x,R_y\rangle$.

The following result is clear from the analysis
in~\cite{parkerpaupert}, but it can also be checked by explicit
computations (see also Section~9.2.1 of~\cite{dpp1}, and Section~2.2
of~\cite{mostowpacific}). Geometrically, it corresponds to the
relative positions of the mirrors of $R_1$ and $R_2$ described in
Section~\ref{sec:notation} and illustrated in
Figure~\ref{fig:mirrors}. When we use the Poincar\'e polyhedron
theorem in Section~\ref{sec:pres}, we will derive these equations from
the cycle relations, see Lemma~\ref{lem:poin-braid}.

\begin{proposition}\label{prop:4-braiding}
  $(R_1R_2)^2$ is a complex reflection with angle $(p-4)\pi/p$ when
  $p\neq 4$, and it is parabolic when $p=4$. For $p=3$, $(R_1R_2)^2$
  is a complex reflection in the point $p_{12}$ where the mirrors of
  $R_1$ and $R_2$ intersect. For $p\geqslant 5$, it is a
  complex reflection in the common perpendicular complex line $m_{12}$ to the
  mirrors of $R_1$ and $R_2$. Moreover, for any $p$, we have
  the higher braid relation
  \begin{equation}\label{eq:braiding}
    (R_1R_2)^2=(R_2R_1)^2
  \end{equation}
\end{proposition}
Note that it follows from~(\ref{eq:braiding}) that $R_1$ and $R_2$
both commute with $(R_1R_2)^2$, which implies the orthogonality
statement about their mirrors. More specifically, we have 

\begin{proposition}\label{prop:stab-12}
  The group $\langle R_1,R_2\rangle$ is a central extension of a
  $(2,p,p)$ orientation preserving triangle group with center
  $\langle(R_1R_2)^2\rangle$ of order $2p/|p-4|$ (which is infinite
  for $p=4$).  In particular, $\langle R_1,R_2\rangle$ has order
  $8p^2/(4-p)^2$ when $p\leqslant 3$ and infinite order when
  $p\geqslant 4$.
\end{proposition}

\Pf 
The first statement follows from the fact that $R_1$ and $R_2$
have order $p$ and $(R_1R_2)^2$ is central with order $2p/|p-4|$. The
second statement follows from the fact that, for $p\leqslant 3$, a
$(2,p,p)$ triangle group has order $4p/(4-p)$ and for $p\geqslant 4$
it is infinite.  \EPf

Geometrically, when $p=3$ the spherical $(2,3,3)$ triangle group acts
on the sphere of complex lines through $p_{12}$; when $p=4$ the vertex
$p_{12}$ is on the ideal boundary and the Euclidean $(2,4,4)$ triangle
group acts on the horizontal factor of the Heisenberg group based at
that vertex; when $p\geqslant 5$ the hyperbolic triangle group
$(2,5,5)$ acts on the complex line $m_{12}$.

It is useful to have a formula for the mirror of $(R_1R_2)^2$ (when
$p\geqslant 5$), which is polar to
$$
  \n_{12}=\left[\begin{matrix}
      a^2\overline{\tau}-\overline{a}\\
      \overline{a}^2\tau-a\\
      a^3+\overline{a}^3
    \end{matrix}\right].
$$ 
The notation is set up to indicate that it is the intersection point
of the mirrors of $R_1$ and $R_2$ (see also
Section~\ref{sec:braiding}). When $p=3,4$ the vector ${\bf n}_{12}$
projects to $p_{12}$, which is in $\ch 2$ or on its boundary
respectively. The obvious extension of this notation allows us to
describe the $J$-orbit of this vector, namely
$$
  \n_{23}=J\n_{12}, \quad \n_{13}=J\n_{23}.
$$
Note that the higher braid relation holds of course between any pair of
reflections $R_j$ and $R_k$ with $j\neq k$, since
$JR_kJ^{-1}=R_{k+1}$, so that
$$
  (R_2R_3)^2=(R_3R_2)^2,\quad (R_3R_1)^2=(R_1R_3)^2.
$$

The following observation will be useful later.
\begin{proposition}\label{prop:betterS1}
  $P^2S_1$ is a complex reflection, in fact
  \begin{equation}\label{eq:betterS1}
    P^2S_1=R_2R_3^{-1}R_2^{-1}.
  \end{equation}
\end{proposition}
\Pf
Using the identity~(\ref{eq:S1}), we see that
$
  P^2S_1=R_3^{-1}R_2^{-1}R_3^{-1}R_2R_3,
$
by repeatedly using $P=R_1J$ and $JR_kJ^{-1}=R_{k+1}$.  The
equality~(\ref{eq:betterS1}) then follows from the braid relation
$(R_2R_3)^2=(R_3R_2)^2$.  
\EPf

A similar analysis holds for $R_1$ and $R_2R_3R_2^{-1}$. In fact, in
this case we recover the braid relation (compare with the calculations
in \cite{dfp}, or in Section 2.2 of~\cite{mostowpacific}, in or
Section 6 of~\cite{parkerlivne}).
\begin{proposition}\label{prop:3-braiding}
  $(R_1R_2R_3R_2^{-1})^3$ is a complex reflection with angle
  $(p-6)\pi/p$ when $p\neq 6$, and it is parabolic when $p=6$. For
  $p\leqslant 5$, $(R_1R_2R_3R_2^{-1})^3$ is a complex reflection in
  the point $q_{123\bar{2}}$. For $p\geqslant 7$, it is a complex
  reflection in the complex line $m_{123\bar{2}}$ perpendicular to the
  mirrors of $R_1$ and $R_2R_3R_2^{-1}$. Moreover, for any $p\in
  \N^*$, we have the braid relation
  \begin{equation}\label{eq:3braiding}
    R_1(R_2R_3R_2^{-1})R_1=(R_2R_3R_2^{-1})R_1(R_2R_3R_2^{-1}).
  \end{equation}
\end{proposition}

\begin{proposition}\label{prop:stab-1232i}
  The group $\langle R_1,R_2R_3R_2^{-1}\rangle$ is a central extension
  of a $(2,3,p)$ orientation preserving triangle group with center
  $\langle(R_1R_2R_3R_2^{-1})^3\rangle$ of order $2p/|p-6|$ (which is
  infinite for $p=6$).  In particular, $\langle
  R_1,R_2R_3R_2^{-1}\rangle$ has order $24p^2/(6-p)^2$ when
  $p\leqslant 5$ and infinite order when $p\geqslant 6$.
\end{proposition}

\Pf
Since $R_1$ has order $p$ and
$(R_1R_2R_3R_2^{-1}R_1)^2=(R_1R_2R_3R_2^{-1})^3$ is central with order
$2p/|p-6|$, we obtain the first statement.  The second statement
follows since a $(2,3,p)$ triangle group has order $12p/(6-p)$ for
$p\leqslant 5$ and is infinite for $p\geqslant 6$.  
\EPf

\subsection{General procedure for building fundamental domains}\label{sec:algo}

In this section, we give a rough idea of the procedure that allowed
us to produce the fundamental domains that appear in this paper.

The general context is that of triangle groups, i.e. groups generated
by three complex reflections. An important example of such a group was
studied in~\cite{schwartz4447}, where Schwartz gave a fundamental
domain for a group generated by complex reflections $I_1$, $I_2$,
$I_3$ of order 2 with pairwise products of order 4 and such that $I_1I_2I_3I_2$ has order 7.
The combinatorial structure of Schwartz's fundamental domain inspired
several subsequent constructions, see~\cite{parkerlivne}
and~\cite{thompsonthesis}, which in turn evolved into a fairly general
procedure that applies to triangle groups. We now outline that
procedure and relate it to Schwartz's construction.

In Section~\ref{sec:braiding}, we discussed the fact that some
pairs of complex reflections in the group satisfy the (generalized) braid relations (\ref{eq:braidingqeven}) and (\ref{eq:braidingqodd}).
In the present paper, we consider only $q=3$ or $q=4$
(for groups with $\tau=\sigma_1$ or $\sigma_5$ in the notation
of~\cite{parkerpaupert}, we would need to consider $q=5$ and $q=6$ as
well).
\medskip

We will now outline how to go from a pair of braiding complex reflections $R_x$ and
$R_y$ to (the combinatorial model of) a side in our polyhedron.  We
first explain the construction in the simplest case, namely when
$L_x\cap L_y$ is negative for the Hermitian inner product,
or equivalently, $m_x\cap m_y=p_{xy}\in{\bf H}^2_{\mathbb C}$.

\subsubsection{Basic construction of pyramids using braiding}

The mirrors $m_x$ and $m_y$ of $R_x$ and $R_y$ pass through $p_{xy}$.
Now consider the images of the lines $m_x$ and $m_y$
under powers of $R_xR_y$. Using subscript notation, these can be
written as $m_x$, $m_{xy\bar{x}}$, $m_{xyx\bar{y}\bar{x}},\,\ldots,\,
m_{\bar{y}xy}$, $m_y$. It is not hard to see that there are
exactly $q$ of these lines, all lying in the pencil of complex lines
through $p_{xy}$. (If $q$ is odd then $m_y$ is the image of $m_x$
under $(R_xR_y)^{(q-1)/2}$. If $q$ is even $m_x$ and $m_y$ are in distinct orbits,
each orbit comprising $q/2$ lines.) These $q$ lines form our starting point
when constructing a side of the polyhedron. 

Consider now a third complex reflection $R_z$ with mirror $m_z$ and
angle $\psi=2\pi/p$. For simplicity, we suppose $m_z$ intersects each
of $m_x,\,m_{xy\bar{x}},\,\ldots,\, m_{\bar{y}xy},\, m_y$ in a point
of ${\bf H}^2_{\mathbb C}$, denoted by $p_{xz}$, $p_{xy\bar{x}z}$, and
so on. The $q$ points $p_{xz},\,p_{xy\bar{x}z},\,\ldots,\,p_{yz}$ form
the vertices of a $q$-gon in the line $m_z$. Together with the point
$p_{xy}$, they form the vertices of a pyramid with apex $p_{xz}$ and
base the $q$-gon in $m_z$. The combinatorial model for the side is
obtained by adding edges as follows. We include each of the $q$ edges
of the $q$-gon in the base, together with an edge from each of the
base vertices to $p_{xy}$ (see Figure~\ref{fig:algo3-4}).

In order for our polyhedron to have a well-defined side pairing, given
a pyramid as above, it is natural to construct an opposite isometric
pyramid. Note that the reflection $R_z$ fixes the base of the above
pyramid pointwise. We find a second pyramid by applying $R_z$ to the
first one. In this case $x$ is replaced with $zx\bar{z}$ and $y$ is
replaced with $zy\bar{z}$ and so on.

The above procedure is illustrated in
Figure~\ref{fig:algo3-4}. Part~(a) shows a square pyramid with $x=2$,
$y=3$, $z=1$, and part~(b) shows its opposite pyramid. Part~(c) shows
a triangular pyramid with $x=1$, $y=\bar{3}23$ and $z=23\bar{2}$, and
part~(d) shows (the image under $P^{-2}$ of) its opposite pyramid;
here $P^{-2}z^{-1}xzP^2P^{-2}(\bar{3}\bar{2}3)1(\bar{3}23)P^2=1$, 
and similarly $P^{-2}z^{-1}yzP^2=23\bar{2}$ and
$P^{-2}zP^2=1\bar{3}23\bar{1}$.

The number of sides in these examples come from the fact that $R_2$
and $R_3$ braid with order 4, and $R_1$ and $R_{\bar 32 3}$ braid with
order 3 (see section~\ref{sec:braiding}). Note also that we use
relations in the group to write
$\bar{y}xyz=\bar{3}\bar{2}31\bar{3}2323\bar{2}=\bar{3}\bar{2}3123$ and
$yz=\bar{3}2323\bar{2}=23$.

Schwartz uses a similar procedure in~\cite{schwartz4447}.  Since his
generating reflections have order $p=2$, in his case $x=\bar{x}$ and
$y=\bar{y}$. Consider Fig 4.2 of~\cite{schwartz4447} which depicts
what Schwartz calls an odd $B$ piece. This is made up of two square
pyramids. In the front one, $x=1$, $y=3$ and $z=2$. The apex is
$p_{xy}=[13]$ and the vertices around the base are $p_{xz}=[12]$,
$p_{xy\bar{x}z}=[1312]$, $p_{\bar{y}xyz}=[3132]$ and
$p_{yz}=p_{zy}=[23]$. The rear pyramid is the image of the front one
under $R_z=I_2$ and has $x=212$, $y=232$ and $z=2$.

\newlength{\mywidth} \setlength{\mywidth}{0.2\textwidth}
\begin{figure}[htbp]
\centering
  \subfigure[$r_1$]{\epsfig{figure=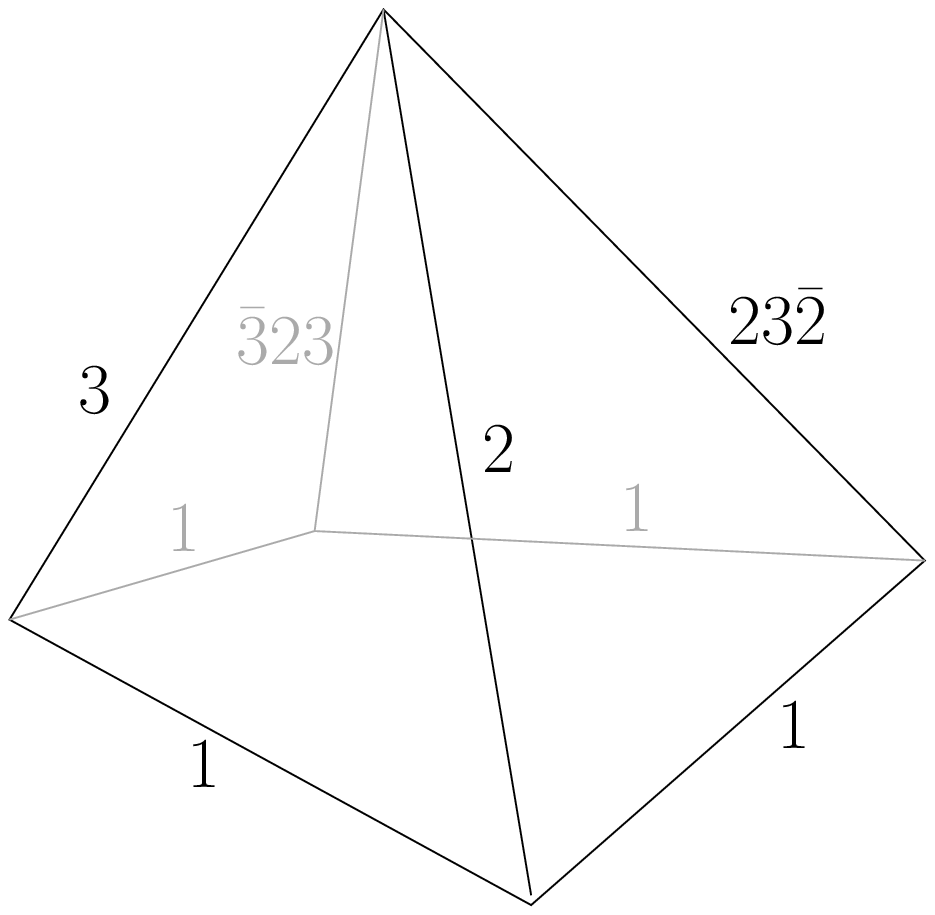, height = \mywidth}}\hfill
  \subfigure[$r_1^-$]{\epsfig{figure=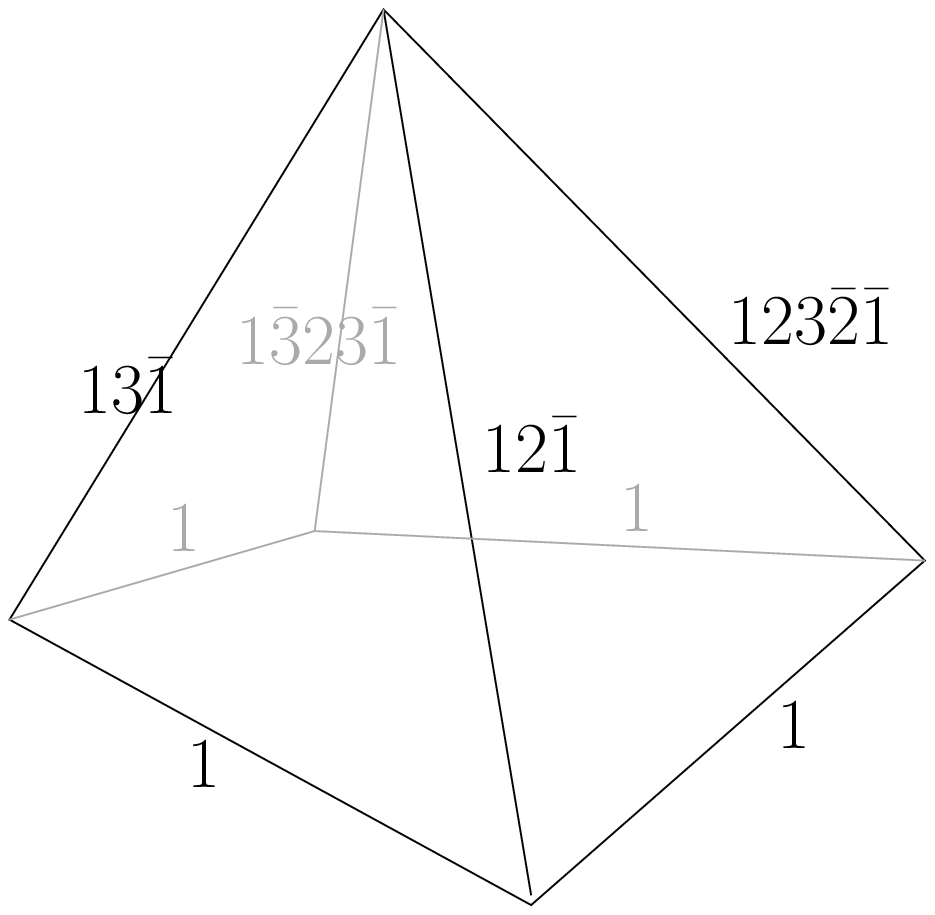, height = \mywidth}}\hfill
  \subfigure[$s_1$]{\epsfig{figure=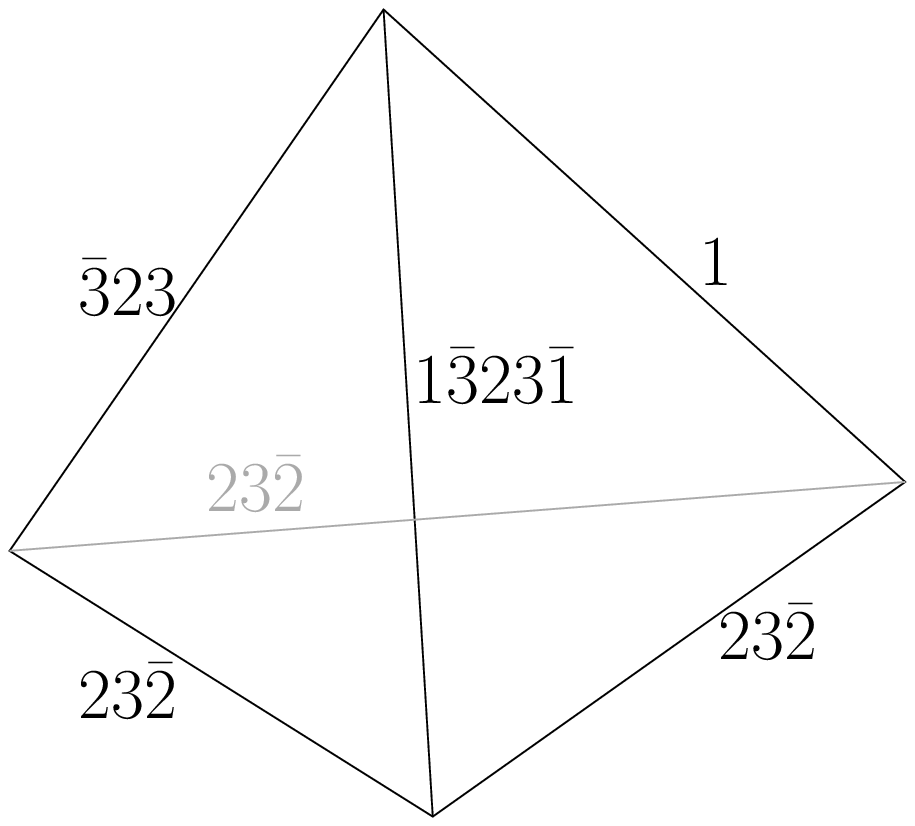, height = \mywidth}}\hfill
  \subfigure[$s_1^-$]{\epsfig{figure=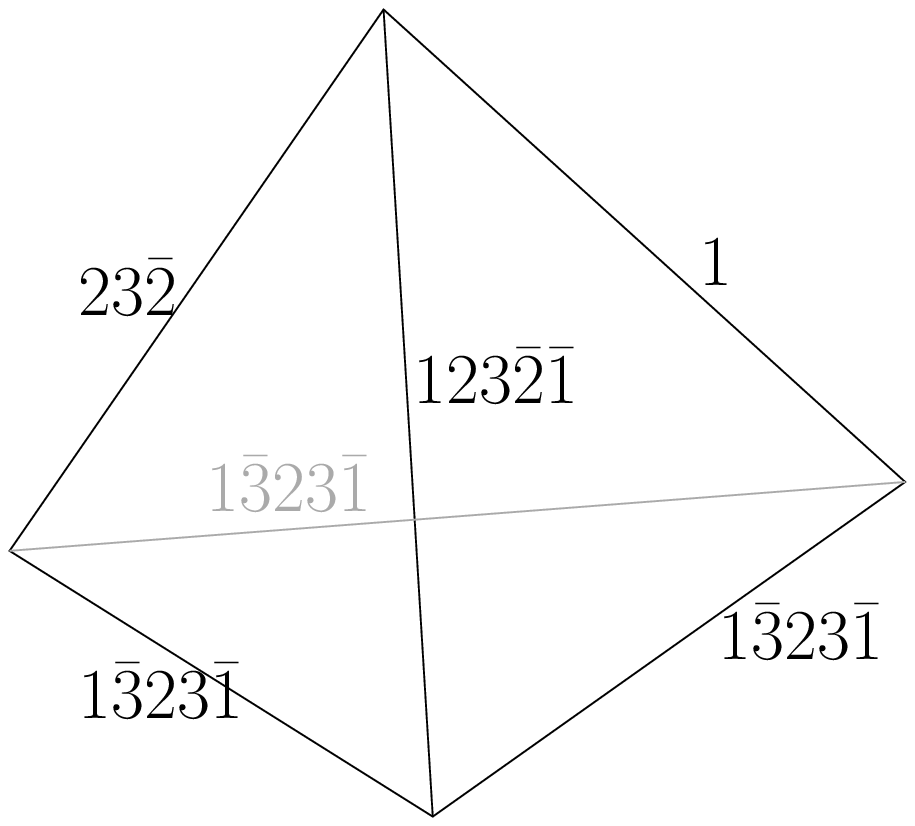, height = \mywidth}}\\
  \caption{Combinatorial structure of the faces, in the case $p=3$ or
    $4$. Since all vertices are inside $\ch 2$, no truncation is needed.}
  \label{fig:algo3-4} 
\end{figure}
\begin{figure}[htbp]
\centering
  \subfigure[$r_1$]{\epsfig{figure=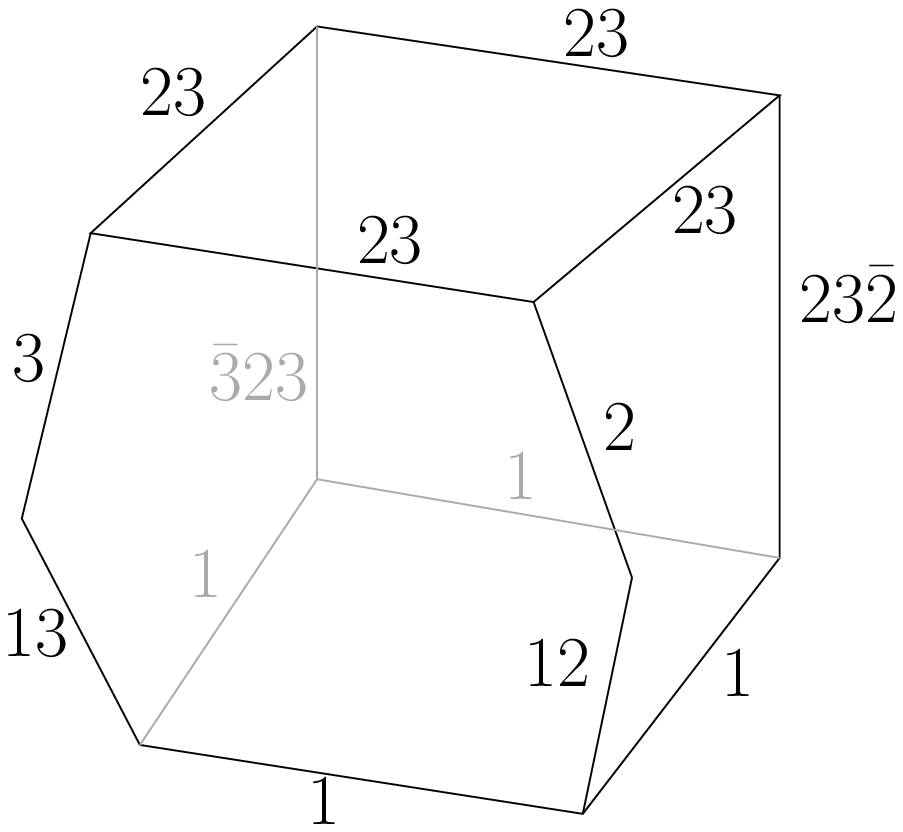, height = \mywidth}}\hfill
  \subfigure[$r_1^-$]{\epsfig{figure=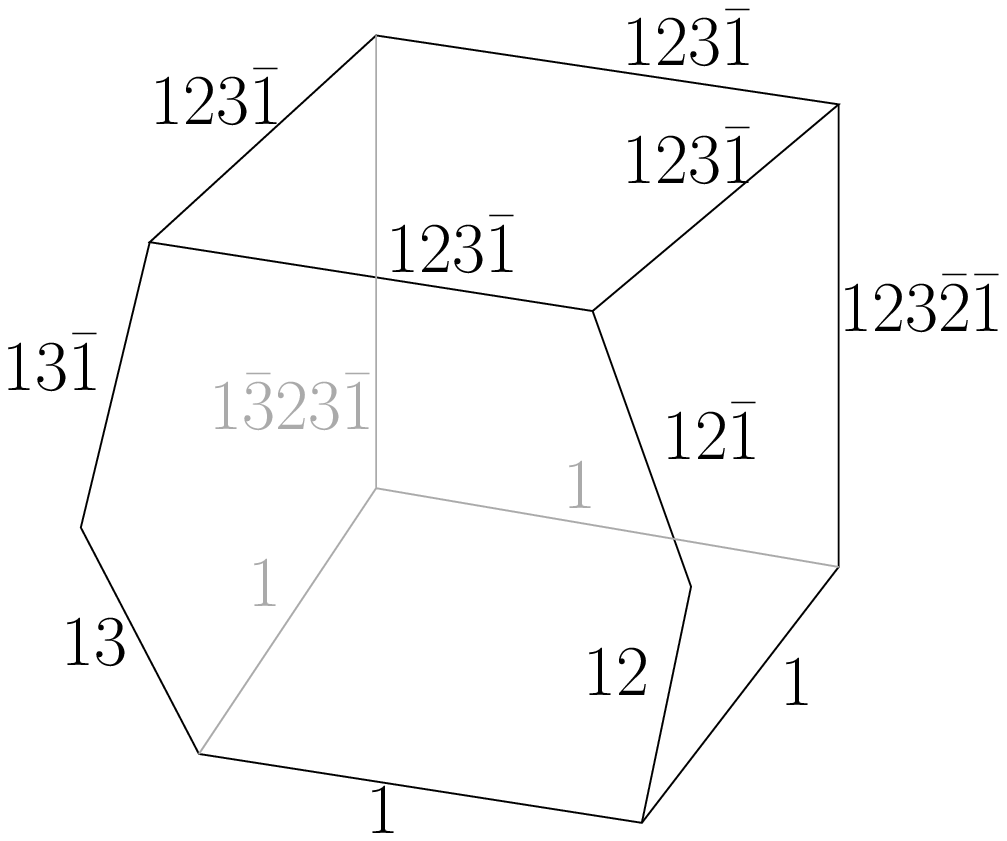, height = \mywidth}}\hfill
  \subfigure[$s_1$]{\epsfig{figure=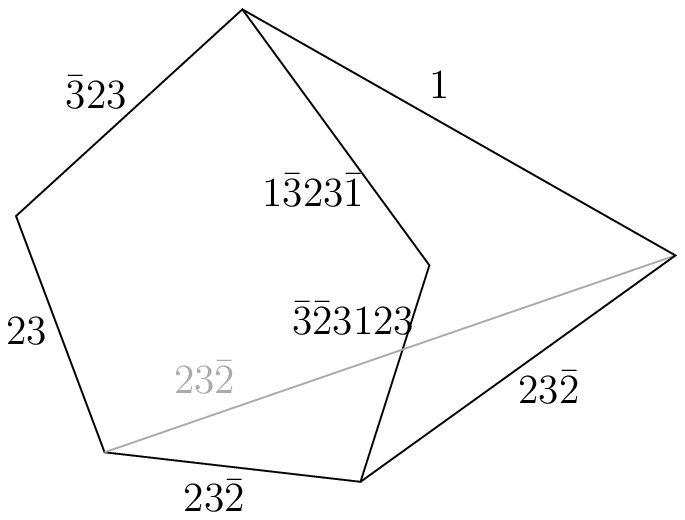, height = 0.9\mywidth}}\hfill
  \subfigure[$s_1^-$]{\epsfig{figure=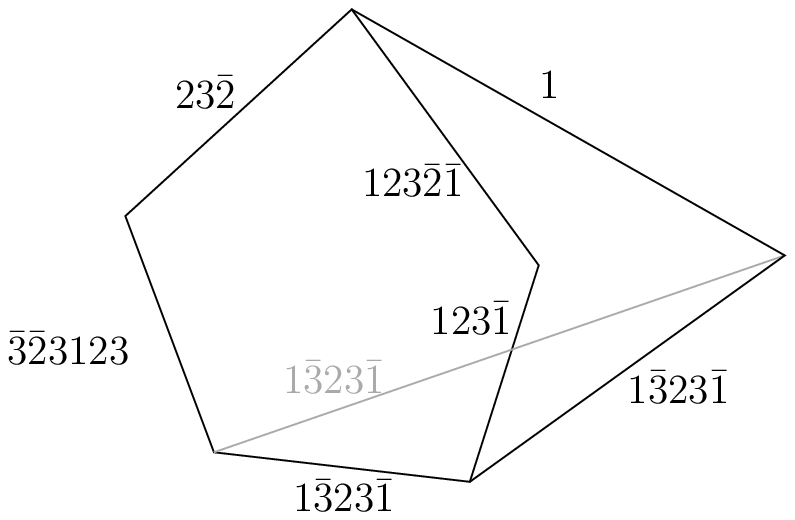, height = 0.9\mywidth}}\\
  \caption{Combinatorial structure of the faces, in the case $p=5$ or
    $6$. A first truncation has occurred since the mirrors of $R_1$, $R_2$ and
    $R_3$ are ultraparallel.}
  \label{fig:algo5-6}  
\end{figure}
\begin{figure}[htbp]
\centering
  \subfigure[$r_1$]{\epsfig{figure=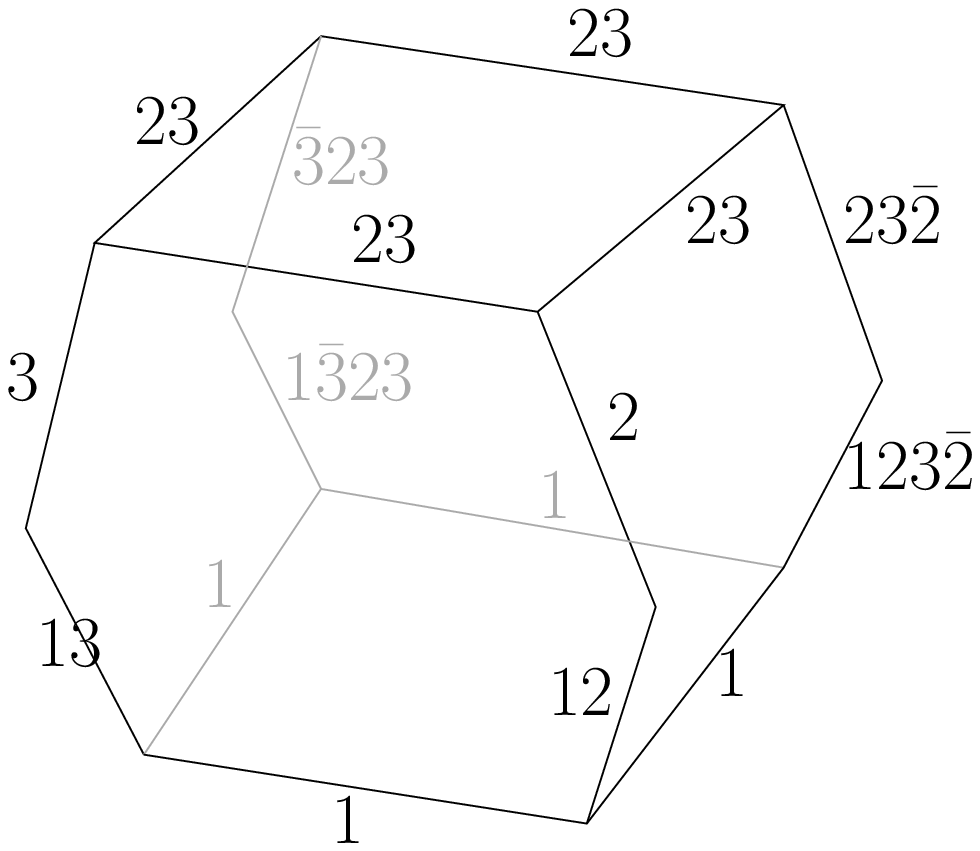, height = \mywidth}}\hfill
  \subfigure[$r_1^-$]{\epsfig{figure=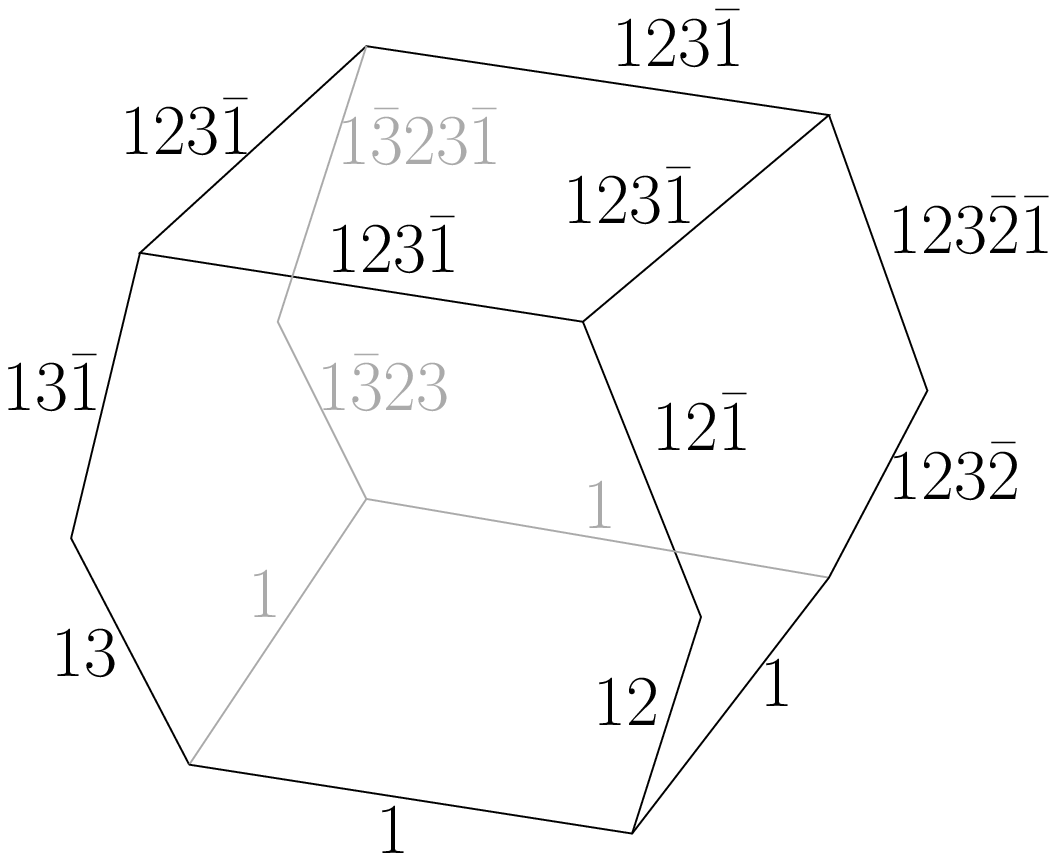, height = \mywidth}}\hfill
  \subfigure[$s_1$]{\epsfig{figure=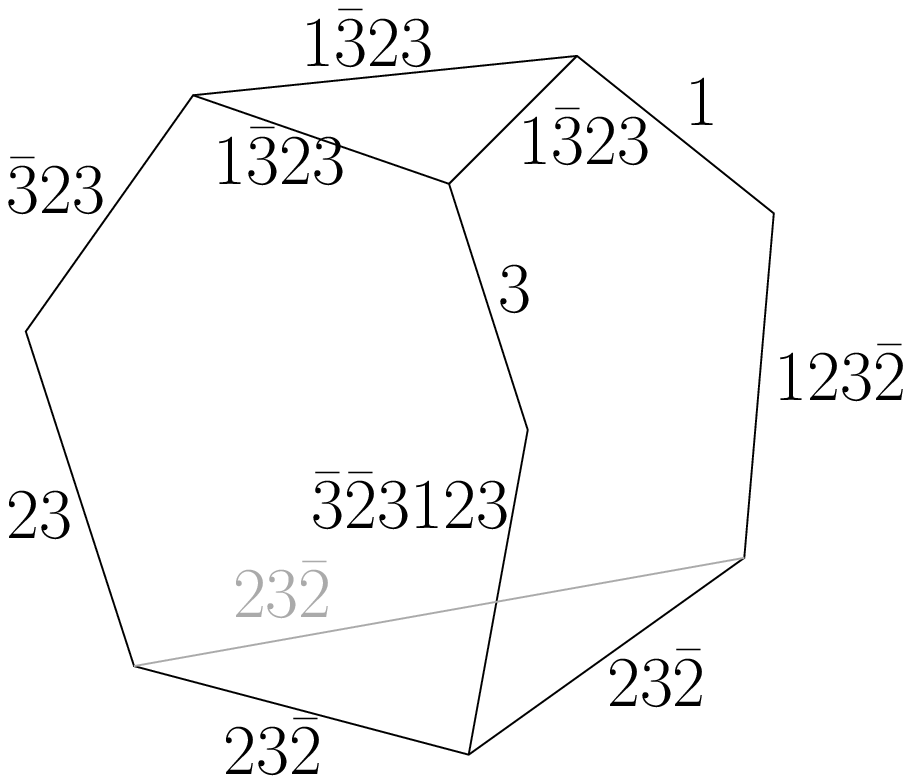, height = \mywidth}}\hfill
  \subfigure[$s_1^-$]{\epsfig{figure=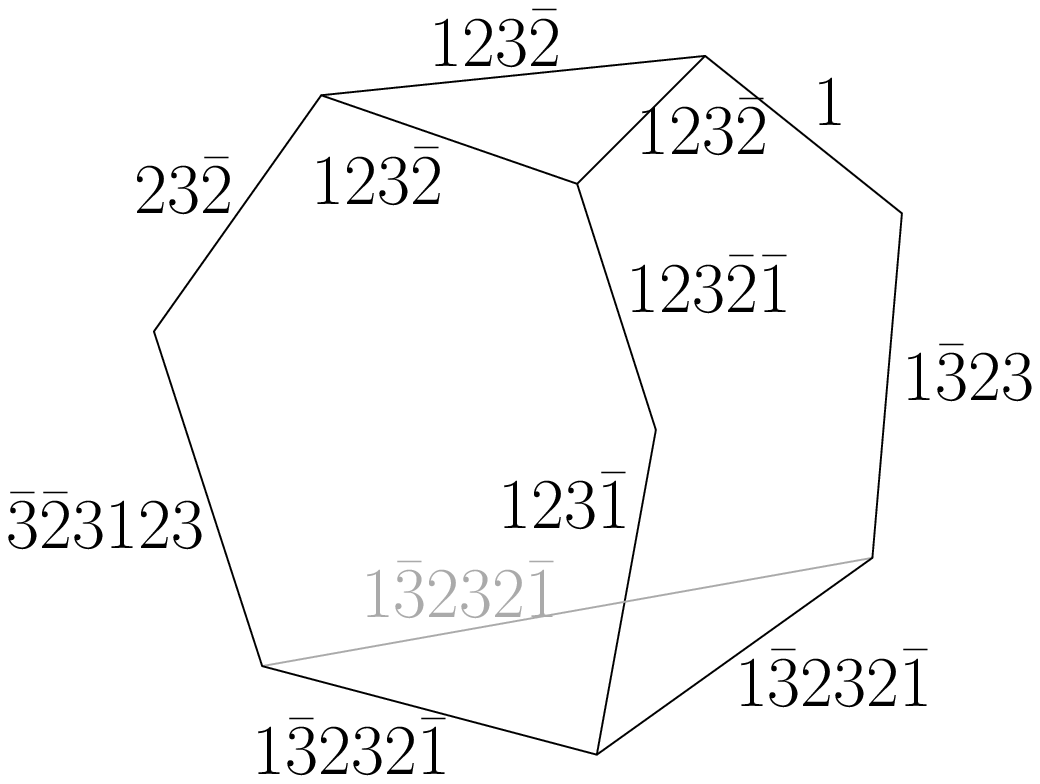, height = \mywidth}}\\
  \caption{Combinatorial structure of the faces, in the case $p=8$ or
    $12$. Further truncation has occurred since the mirrors of $R_1$
    and $R_3^{-1}R_2R_3$ are ultraparallel.}
  \label{fig:algo8-12}  
\end{figure}

Given one pyramid as above, the ridges containing the apex are
all triangles, for example the ridge determined by the three complex lines
$m_x$, $m_y$, $m_z$ with vertices $p_{xy}$, $p_{xz}$, $p_{yz}$. We can
repeat the whole of the above construction starting with a different
pair of these complex lines, for example $m_z$ and $m_y$. This will
yield a new pyramid, whose base lies in $m_x$ and which shares a
common triangular side with our initial pyramid. In order to construct 
a combinatorial model for our polyhedron, we repeat this process.
For each side of each pyramid, there should be exactly one other pyramid
with the same triangular side, and the whole collection should close
up and be invariant under powers of $P$. This game should be compared
with \cite{schwartz4447} or \cite{thompsonthesis} where a similar game 
is played. The main difference is that Schwartz is not interested in 
triangles if all three vertices lie outside complex hyperbolic space.

\subsubsection{Truncation}

It may be that some, or all, of the vertices of the pyramid lie
outside complex hyperbolic space. In this case, it turns out the intersection point
is the polar vector of a mirror of a complex reflection in the group 
(provided $p\ge 3$). With the notation above, suppose $m_x$ and $m_y$
are ultraparallel. As we have seen, $(R_xR_y)^{q/2}$ or
$(R_xR_y)^q$, for $q$ even or odd repectively, is a complex reflection
and commutes with $R_x$ and $R_y$. Therefore the common
orthogonal to $m_x$ and $m_y$ is the mirror $m_{xy}$ of this
power of $R_xR_y$.

The simplest case
of this are the even $A$ piece of Schwartz 
depicted in Figure 4.3 of \cite{schwartz4447}. In this case $q=7$ (and
he still has $p=2$) and he has two pyramids with a heptagonal base.
For the front pyramid $x=131$, $y=2$ and $z=1$.  Thus the apex of the
pyramid is $p_{xy}=[1312]$. Of the base vertices, four are points of
complex hyperbolic space, namely: $p_{xz}=[1311]=[13]$,
$p_{xy\bar{x}z}=[13121311]=[131213]$,
$p_{\bar{y}xyz}p_{z\bar{y}xy}=[121312]$ and $p_{yz}=p_{zy}=[12]$.  The
other three base vertices lie outside complex hyperbolic space.  The
rear pyramid is the image of the first under $R_z=I_1$, and so has
$x=3$, $y=121$ and $z=1$.  Schwartz's picture does not illustrate the
truncation arising from vertices lying outside complex hyperbolic
space.

We now explain the implications of this truncation
process for our polyhedra. The simplest feature is when $m_x$ and
$m_y$ are ultraparallel. We denote their common orthogonal by
$m_{xy}$. The $q$ lines $m_x$, $m_{xy\bar{x}}$,
$m_{xyx\bar{y}\bar{x}},\,\ldots,\,m_{\bar{y}xy}$, $m_y$ are all
orthogonal to $m_{xy}$ and their intersection points with $m_{xy}$
(denoted $p^x_{xy}=m_x\cap m_{xy}$ etc.) now form the vertices of a
$q$-gon in $m_{xy}$; compare Figure~\ref{fig:algo3-4}(a) with
Figure~\ref{fig:algo5-6}(a) and
Figure~\ref{fig:algo8-12}(a), where the apex $p_{23}$ has been
replaced with a quadrilateral.

Next, we explain how the truncation affects base vertices.
Suppose the mirrors $m_x$ and $m_z$ of $R_x$ and $R_z$ are ultraparallel. 
Then consider their common perpendicular line $m_{xz}$. The mirrors $m_x$ 
and $m_z$ intersect $m_{xz}$ in distinct points, which we label $p^x_{xz}$
and $p^z_{xz}$ respectively. The point $p_{xz}$ is then
replaced with a line in $m_{xz}$ between these two points.
Compare Figure~\ref{fig:algo3-4}(c)
with Figure~\ref{fig:algo5-6}(c). 
In particular, the vertex $p_{23}$ has been replaced with
an arc from $p^{\bar{3}23}_{23}$ to $p^{23\bar{2}}_{23}$. 
We can still build up a polyhedron from these truncated polyhedra, the
main difference being that the sides that were formerly triangles may
be quadrilaterals, pentagons or hexagons.

\subsection{Definition of the polyhedron $E$} \label{sec:defE}

In this section we define the fundamental polyhedron $E$. We will show that
$E$ is a piecewise smooth polyhedron (as defined in Section~\ref{sec:poin}), and in parallel we determine its precise combinatorics. This allows us to show that $E$ is a fundamental domain
for the action of the cosets of $\Upsilon={\rm Stab}_\Gamma(E)$ in
$\Gamma$.

\subsubsection{Bounding bisectors}\label{sec:boundingbisectors}

The basic fact we will use is that the pyramids constructed in
section~\ref{sec:algo} have a natural geometric realization inside a
suitable bisector. Indeed, for a given pyramid, there is a unique
bisector that contains the base of the pyramid in a complex slice, and
whose extended real spine contains its apex. We call this bisector the
\emph{supporting bisector} of the pyramid.
Concretely, if the base is given by ${\bf m}^\perp$ and the apex by
${\bf a}$, the real spine of the supporting bisector is given by the
real span of ${\bf a}$ and  ${\bf a}-\langle {\bf a},{\bf m}\rangle{\bf m}/\langle {\bf
    m},{\bf m}\rangle$ (a lift of the orthogonal projection of ${\bf a}$ onto the base). The following gives a geometric justification for considering such bisectors:
    
\begin{proposition}
  The supporting bisector of a pyramid contains all real geodesic
  segments between its vertices.
\end{proposition}
This is obvious for edges in the base, as the base is contained in a complex slice of the supporting bisector.
For edges through the apex, it follows from the fact that a bisector contains all geodesics intersecting the real spine. 
For edges in truncation faces, it follows from the fact that the complex line containing the truncation
face is a complex slice of the supporting bisector.

The basic building blocks of $E$ will be the following four bisectors, which are the
supporting bisectors of the four faces described in
section~\ref{sec:algo}.
 
\begin{dfn}\label{def:bounding}
  \begin{enumerate} \item $\Rc$ is the bisector whose extended real spine contains $\n_1$
      and $\n_{23}$;
    \item $\Rc^-$ is the image of $\Rc$ under $R_1$;
    \item $\Sc$ is the bisector whose extended real spine contains $\n_{23\bar 2}$
      and $\n_{1\bar 323}$;
    \item $\Sc^-$ is the image of $\Sc$ under $S_1$.
  \end{enumerate}
\end{dfn}
Recall from section~\ref{sec:notation} that we denote by ${\bf p}$ the
isolated fixed point of $P$.
\begin{dfn}
Let $F$ be the intersection of the 28 half-spaces containing ${\bf p}$
bounded by the bisectors $P^k\Rc^\pm$ and $P^k\Sc^\pm$ ($k=0,\pm 1,\pm 2,\pm 3$).
The region $E$ is the connected component of $F$ containing ${\bf p}$.
\end{dfn}
We will refer to the 28 bisectors of the form $P^k\Rc^\pm$ and
$P^k\Sc^\pm$ as the {\bf bounding bisectors}.

\begin{rmk} (1) The domain $E$ has much simpler combinatorics than the
Dirichlet domains that were used in~\cite{dpp1} (when
$p\geqslant 4$). As a slight drawback
of not using Dirichlet domains, it is unclear whether $E$ is
star-shaped with respect to ${\bf p}$, and hence it is not convenient
to use geodesic cone arguments; rather we use fundamental domains for
coset decompositions.

(2)  It turns out to be complicated to determine the precise combinatorics of the polyhedron $F$ bounded by the bounding bisectors; in fact it is not even obvious that $F$ is connected. This is why we define $E$ as a component of $F$. In fact one could prove that $E=F$ (which would then give a description of the combinatorial structure of $F$) but this is not needed in order to prove our main theorem. The main reason for working with $E$ rather than $F$ is that it significantly reduces the number of computer calculations.
 \end{rmk}

It will be useful to distinguish between sides of $E$ and the bisector
that contains them, so we will use the following notation:
\begin{dfn}\label{def:sides}
  The intersection of $E$ with a given bounding bisector is called a
  {\bf bounding side}. The $28$ bounding sides will be denoted by
  $$
  P^kr_1^\pm = E\cap P^k\Rc^\pm;\quad P^ks_1^\pm = E\cap P^k\Sc^\pm,
  $$   
  for $k=0, \pm 1, \pm 2, \pm 3$.
\end{dfn}

The combinatorial structure of the bounding sides $P^kr_1^\pm$ and
$P^ks_1^\pm$ is described in detail in Section~\ref{sec:model}. Note
that for a polyhedron $E$ bounded by bisectors in $\ch 2$, it is by no
means obvious how to determine the combinatorics of $E$ from the knowledge
of its vertices, which explains the length of the sections stating and
proving the combinatorics.

\subsubsection{Symmetry}\label{sec:symmetry}

The polyhedron $E$ is by definition $P$-invariant. Our notation for
bisectors is set up so that the $P$-orbits can be conveniently read
off their labels. For vertices (and mirrors of reflections), it is a
bit more tedious but completely straighforward to study these
orbits. For concreteness we treat a couple of examples in detail.

The basic point is that $P=R_1J$, and $JR_k J^{-1}=R_{k+1}$,
i.e. conjugation by $J$ raises indices by one (mod 3). The $P$-image
of the mirror of $R_1$, which corresponds to the complex line polar to
$m_1$, is polar to $Pm_1=m_{12\bar1}$, since
$$
  P1P^{-1}=1J1J^{-1}\bar 1=12\bar 1.
$$ 

The $P$-image of $p_{23}$ is simply $p_{13}$, since it is fixed by
$P23P^{-1}=1J23J^{-1}\bar 1=131\bar1=13$. Similarly,
$Pp_{13}=p_{12}$. Applying $P$ often yields slightly more complicated
results, for instance the $P$-image of $p_{12}$ is given by
$p_{P12P^{-1}}=p_{123\bar 1}$, since
$$
  P12P^{-1}=1J12J^{-1}\bar1=123\bar 1.
$$

It will be important in the sequel to observe that $E$ also has an
anitholomorphic involutive symmetry, which we now explain. The
antiholomorphic map
\begin{equation}\label{eq:iota23}
\iota_{23}:\left[\begin{matrix} z_1 \\ z_2 \\z_3 \end{matrix}\right]
\longmapsto \left[\begin{matrix} \overline{z}_1 \\
\overline{z}_3 \\ \overline{z}_2 \end{matrix}\right].
\end{equation}
clearly induces an isometry of $\ch 2$, since
$
  \langle \iota_{23}{\bf v},\iota_{23}{\bf w}\rangle = \overline{\langle {\bf v},{\bf w}\rangle}
$.
It is easy to see that $\iota_{23}$ conjugates $J$ into $J^{-1}$, $R_1$ into $R_1^{-1}$, $R_2$ into $R_3^{-1}$ and $R_3$ into $R_2^{-1}$. From this it follows that the antiholomorphic map $\iota$ given by 
\begin{equation}\label{eq:sigma}
  \iota=R_1\iota_{23}
\end{equation}
is an involution as well, and this will be useful several times
throughout the paper.

\begin{proposition}\label{prop:symmetry}
The involution $\iota$ preserves $E$, sending $r_1$ to $r_1^-$ and $s_1$ 
to $s_1^-$.
\end{proposition}
\Pf Since $P=R_1J$, it
follows easily that $\iota$ conjugates $R_1$ into
$R_1^{-1}$, and $P$ into $P^{-1}$ (beware that $\iota_{23}$
does not conjugate $P$ into $P^{-1}$). In particular, $\iota$ fixes
${\bf p}$.
  
The fact that $\iota$ conjugates $S_1$ into $S_1^{-1}$ follows 
from the fact that $S_1$ can be expressed as a palindromic word in $R_1$ 
and $P$, namely $S_1=P^2R_1P^{-2}R_1P^2$.

Hence $\iota$ sends the half-space containing ${\bf p}$ bounded by
$P^k\Rc^\pm$ (respectively $P^k\Sc^\pm$) to the half-space containing
${\bf p}$ bounded by $P^{-k}\Rc^\mp$ (respectively $P^{-k}\Sc^\mp$).
This show that $\iota$ preserves $F$, and this implies
that it preserves $E$ (which is the connected component of $F$
containing ${\bf p})$.
\EPf

Note that $\iota$ conjugates an isometry given by a word $w$ (see Section~\ref{sec:wordnotation}) into the isometry given by the word $1w'\bar{1}$, where $w'$ is obtained from $w$ by replacing each occurrence of  $1,2,3$ by  $\bar{1},\bar{3},\bar{2}$ respectively. One can then read the action of $\iota$ on the vertices $p_w$, $q_w$ by applying this rule to the label $w$. For instance, we have
\begin{equation}
\begin{array}{l}
\iota(p_{12})=p_{1\bar{1}\bar{3}\bar{1}}=p_{\bar{3}\bar{1}}=p_{13}\\ 
\iota(p_{23})=p_{1\bar{3}\bar{2}\bar{1}}=p_{123\bar1},
\end{array}
\end{equation}
and
\begin{equation}
  \iota(q_{123\bar2})=q_{\bar3\bar2 3 \bar 1}=q_{1\bar323}
\end{equation}
Observe that $Pp_{13}=p_{12}$ implies
$P^{3}p_{12}=P^{-3}p_{13}$. Since
$p_{\bar3\bar23123}=R_3^{-1}R_2^{-1}R_1^{-1}p_{13}=P^{-3}p_{13}$ and
$\iota(p_{\bar3\bar23123})=R_1(p_{23\bar2\bar1\bar3\bar2})=R_1R_2R_3p_{12}=P^3p_{12}$,
we get
\begin{equation}
  \iota(p_{\bar3\bar23123})=p_{\bar3\bar23123}.
\end{equation}

\subsubsection{Alternative descriptions} \label{sec:alternative}

Since the bounding bisectors $\Rc^\pm$, $\Sc^\pm$ are central to the
construction of the fundamental domain, it will be useful to have
alternative descriptions. In Table~\ref{tab:spines} we write the
bounding bisectors as equidistant from a pair of points and we also
give some information that is useful for computational purposes; see
also Section~\ref{sec:ridgeconsistency}).
We label vectors that appear in
the table following the conventions from
Section~\ref{sec:wordnotation}; $y_0$ denotes the intersection point
of the complex spines of $\Rc$ and $\Sc$ (see below for explicit
computation).

\begin{table}[htbp]
$$
  \begin{array}{|r|l|l|l|l|}
    \hline
    \hbox{Bisector} & \hbox{Equidistant} & 
    \hbox{$\R$-spine} & \hbox{Polar} & \hbox{Reflections}  \\ 
    \hline
    \Rc & \B(y_0, R_1^{-1}y_0)  
             & {\rm Span}_\R(\n_1, \,\n_{23}) 
             & {\bf f}_1=[0,\overline{\tau},\overline{a}\tau]
             & 1,\ (23)^2 \\ 
    \Rc^- & \B(y_0, R_1y_0)  
             & {\rm Span}_\R(a^2 \n_1,\, \overline{a} \n_{123\bar{1}}) 
             & {\bf f}_2={\bf f}_1
             & 1,\  1(23)^2\bar 1 \\ 
    \Sc & \B( y_0, S_1^{-1} y_0) 
             & {\rm Span}_\R(a \n_{23\bar{2}},\, \n_{1\bar{3}23} )
             & {\bf f}_3=[a^2\overline{\tau},-a,-\overline{\tau}] 
             & 23\bar 2, \ (1\bar 323)^3 \\ 
    \Sc^- & \B(y_0, S_1 y_0) 
             & {\rm Span}_\R(\overline{a} \n_{1\bar 3 23\bar 1},\, \n_{123\bar 2} )
             & {\bf f}_4=[a^2\overline{\tau},a\tau,1] 
             & 1\bar 323\bar 1,\  (123\bar 2)^3\\ 
             \hline
  \end{array}
$$
\caption{Descriptions of the bisectors $\Rc^{\pm}$ and $\Sc^{\pm}$ 
as equidistant hypersurfaces, and in terms of vectors spanning their 
extended real spines. We also give polar vectors for their complex spines
and reflections whose mirrors are contained in each bisector
  (the mirror is either a complex line or a point, depending on the
  value of $p$).}\label{tab:spines}
\end{table}
In the table, we list two vectors ${\bf w},{\bf v}$ with real inner
products, so that the extended real spine is obtained by taking real
linear combinations of these two vectors (the real spine corresponds
to such vectors with negative square norm).
For convenience, we also give a polar vector to the complex spine,
i.e. a nonzero vector ${\bf f}$ with $\langle {\bf w},{\bf
  f}\rangle=\langle {\bf v},{\bf f}\rangle=0$, as well as specific
reflections whose mirrors are slices of these bisectors.

\begin{proposition}\label{prop:coequid1}
  The bisectors $\Rc$, $\Rc^-$, $\Sc$ and $\Sc^-$ are all
  coequidistant from a point in $\ch 2$, which we denote by $y_0$.
  More specifically, $\Rc^{\pm}=\B(y_0,R_1^{\mp 1}y_0)$ and
  $\Sc^{\pm}=\B(y_0,S_1^{\mp}y_0)$.
\end{proposition}
\Pf 
From a pair of vectors ${\bf w},{\bf v}$ spanning the real spine of a 
bisector,
one easily finds a vector orthogonal to both of them, since
this amounts to finding a vector ${\bf f}$ that satisfies 
$$
{\bf w}^* H {\bf f} = {\bf v}^* H {\bf f} = 0
$$ 
which is a linear system in ${\bf f}$. This allows us to get
formulae for the four vectors ${\bf f}_j$. Note that $\Rc$ and $\Rc^-$
are cospinal, i.e. they have the same complex spine.

From the ${\bf f}_j$, one can easily find the intersection between two
of the complex spines, again by solving a linear system. The
vector ${\bf y}_0$ can obtained in this way, being the unique
intersection point of the complex spines of $\Rc$ and $\Sc$. For
future reference, we give explicit coordinates:
\begin{equation}\label{eq:x}
  {\bf y}_0 = \left[\begin{matrix} 
     -a^3(\overline{a}^2\overline{\tau}+a)^2 \\
     (a^2\overline{\tau}+\overline{a}\tau)
     (a^2\overline{\tau}-\overline{a}) \\
     (a^2\overline{\tau}+\overline{a}\tau)
     (\overline{a}^2\tau-a) 
  \end{matrix}\right],
\end{equation}
One easily checks that $\langle {\bf y}_0,{\bf f}_j\rangle=0$ 
for all $j=1,2,3,4$, so that all four bisectors $\Rc^\pm$,
$\Sc^\pm$ are equidistant from ${\bf y}_0$. A further simple check
shows that $\langle {\bf y}_0,{\bf y}_0\rangle<0$ and so ${\bf y}_0$
corresponds to a point $y_0\in \ch 2$.

We now show that $\Rc=\B(y_0,R_1^{-1}y_0)$. From Table~\ref{tab:spines}
we see that the spine of $\Rc$ is the real span of ${\bf n}_1$ and ${\bf n}_{23}$.
These vectors are both fixed by the involution $\iota_{23}$ defined in
\eqref{eq:iota23}, and so it also fixes any vector on the real spine of $\Rc$.
Applying $\iota_{23}$ to ${\bf y}_0$ and simplifying, we see that
$\iota_{23}\,{\bf y}_0 = \bar{a}^2R_1^{-1}{\bf y}_0$. Therefore, if ${\bf s}$ is
any vector in the real spine of $\Rc$ we have:
$$
  \langle{\bf y}_0,{\bf s}\rangle
  = \overline{\langle\iota_{23}\,{\bf y}_0,\iota_{23}\,{\bf s}\rangle}
  = a^2\overline{\langle R_1^{-1}{\bf y}_0,{\bf s}\rangle}.
$$
Therefore, any point of $\Rc$ is equidistant from $y_0$ and $R_1^{-1}y_0$.
Applying $R_1$ we see that $\Rc^-$ is equidistant from $R_1y_0$
and $y_0$.

Similarly, the involution $\bar{3}\bar{2}3\bar{1}\bar{3}\iota_{12}$
fixes $a{\bf n}_{23\bar{2}}$ and ${\bf n}_{1\bar{3}23}$ and so
fixes the spine of $\Sc$ pointwise. We calculate that
$\bar{3}\bar{2}3\bar{1}\bar{3}\iota_{12}{\bf
  y}_0=\bar{a}^2S_1^{-1}{\bf y}_0$.  Hence
$\Sc=\B(y_0,S_1^{-1}y_0)$. Applying $S_1$, we find that
$\Sc^-=\B(S_1y_0,y_0)$.  
\EPf

Later, in Table~\ref{tab:bis-xyz} we describe $\Rc^\pm$, $\Sc^\pm$ and
their images under powers of $P$ as equidistant hypersurfaces with
respect to some other points.

\subsection{Combinatorics of $E$} \label{sec:combinatorics}

In this section we describe the combinatorics of $E$ and show that it is
a topological ball with piecewise smooth boundary, each piece being
contained in one of the bounding bisectors.
The proof of the precise combinatorics of $E$ is difficult, but necessary to 
apply the Poincar\'e polyhedron theorem. In fact it is delicate even to give 
a detailed description of the combinatorics (this will be done in
Section~\ref{sec:statement}). Indeed, for a general
polyhedron bounded by bisectors in $\ch 2$, it is not sufficient to
give the vertices, as there can be for instance distinct $1$-cells joining two given
vertices (see~\cite{dfp}).

In order to prove the combinatorics, we first construct a
combinatorial model $\widehat{E}$ for $E$. We then describe its
geometric realization, which is a map $\phi$ from $\widehat{E}$ into
$\chb 2$, sending each $3$-cell into one of the bounding bisectors
(see Definition~\ref{def:bounding}), and we prove that
$\phi(\widehat{E})=E$. This simultaneously proves the statement of the
combinatorics of $E$, and shows that $E$ is a polyhedron in the sense
of Section~\ref{sec:poin}.

Our combinatorial model is somewhat similar to the one used by
Schwartz in Section~6 of~\cite{schwartz4447}. However, our geometric
realization of $\widehat{E}$ uses bisectors whereas Schwartz uses
affine cells.

\subsubsection{Statement of the combinatorics of $E$} \label{sec:statement}

We state the combinatorics in the form of pictures of the $3$-cells,
with incidence information along each piece of the skeleton. Recall
from Section~\ref{sec:domain} that $E$ only has two isometry classes
of faces  (applying $\iota$ and powers of $P$), so we
need only draw pictures for $r_1=\Rc \cap E$ and $s_1=\Sc \cap E$.

\begin{figure}[htbp]
\centering
  \subfigure[ridges of $r_1$]{\epsfig{figure=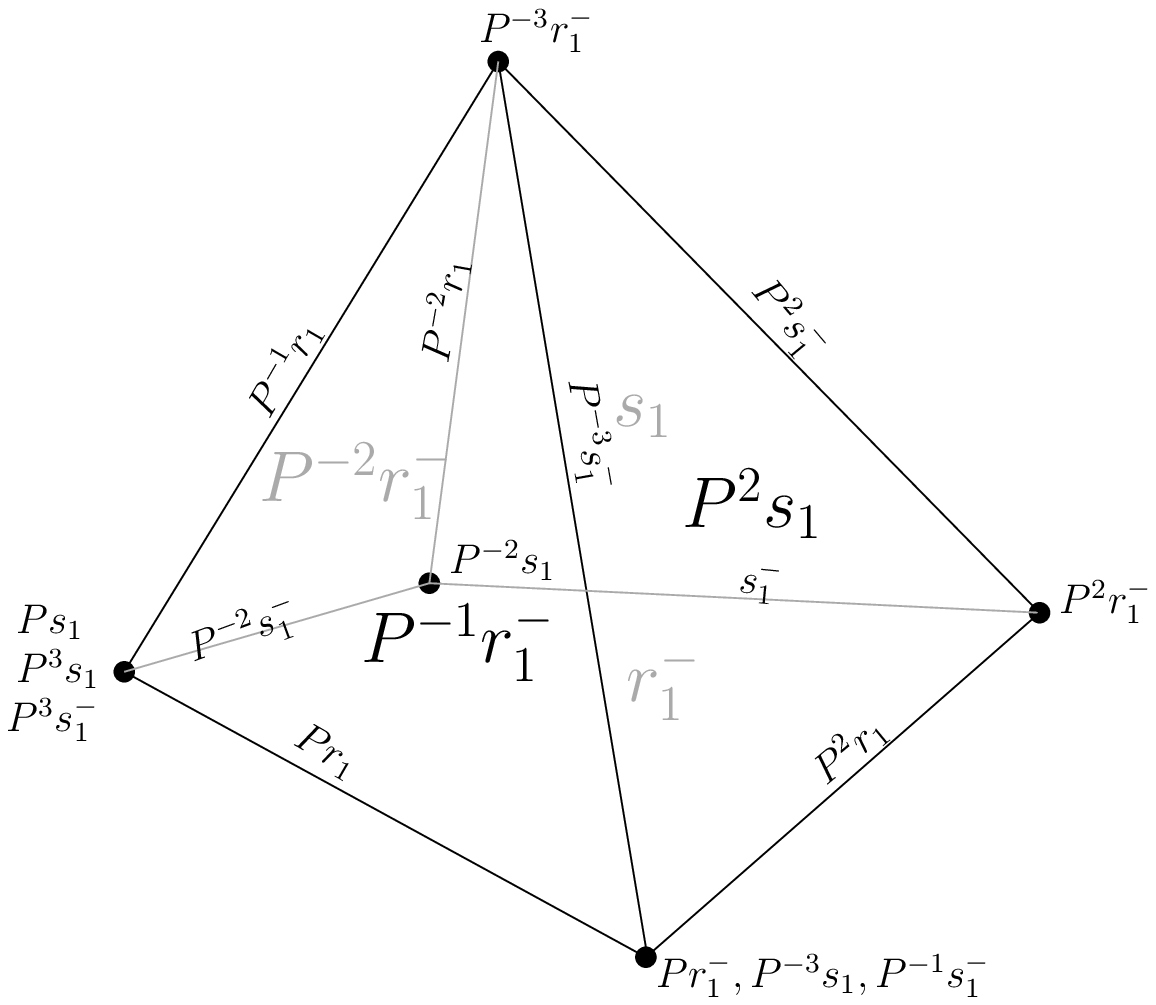, height = 0.36\textwidth}}\hspace{0.1\textwidth}
  \subfigure[ridges of $s_1$]{\epsfig{figure=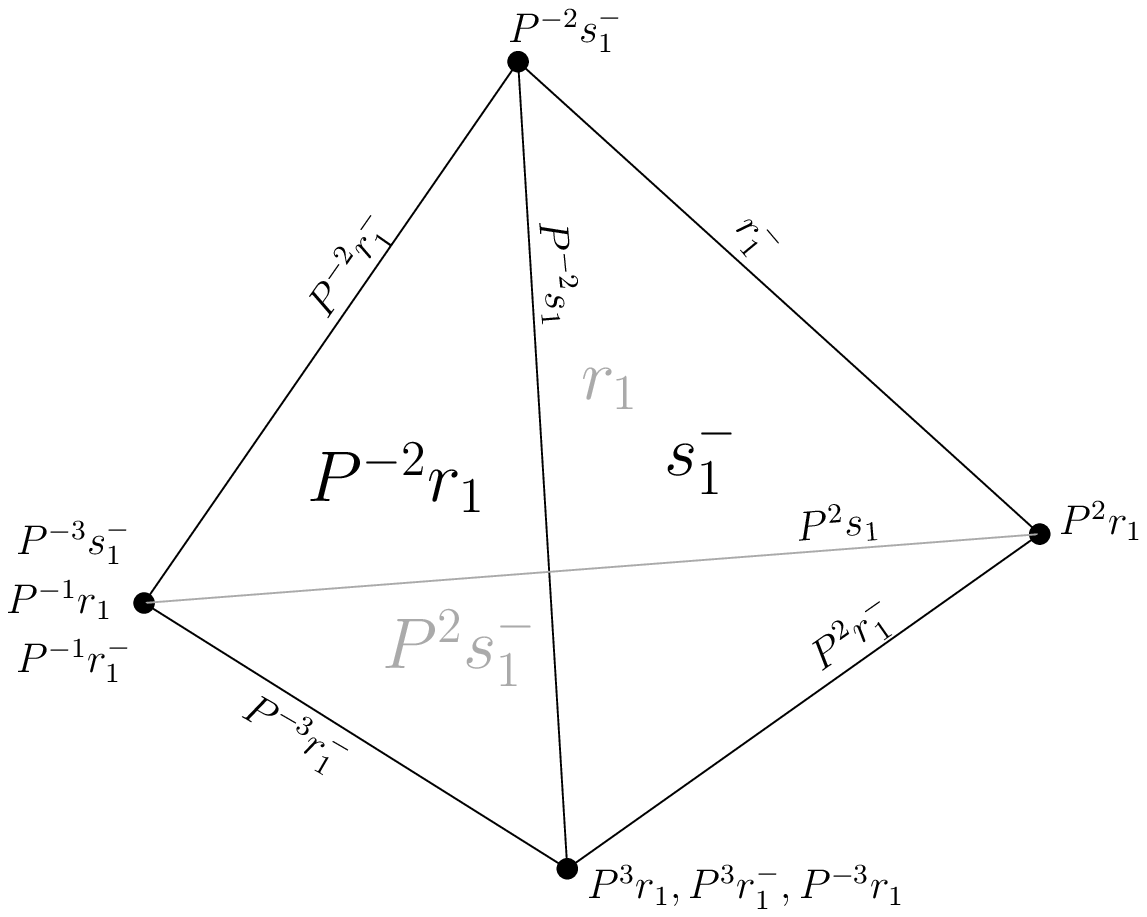, height = 0.36\textwidth}}\\
  \subfigure[vertices of $r_1$]{\epsfig{figure=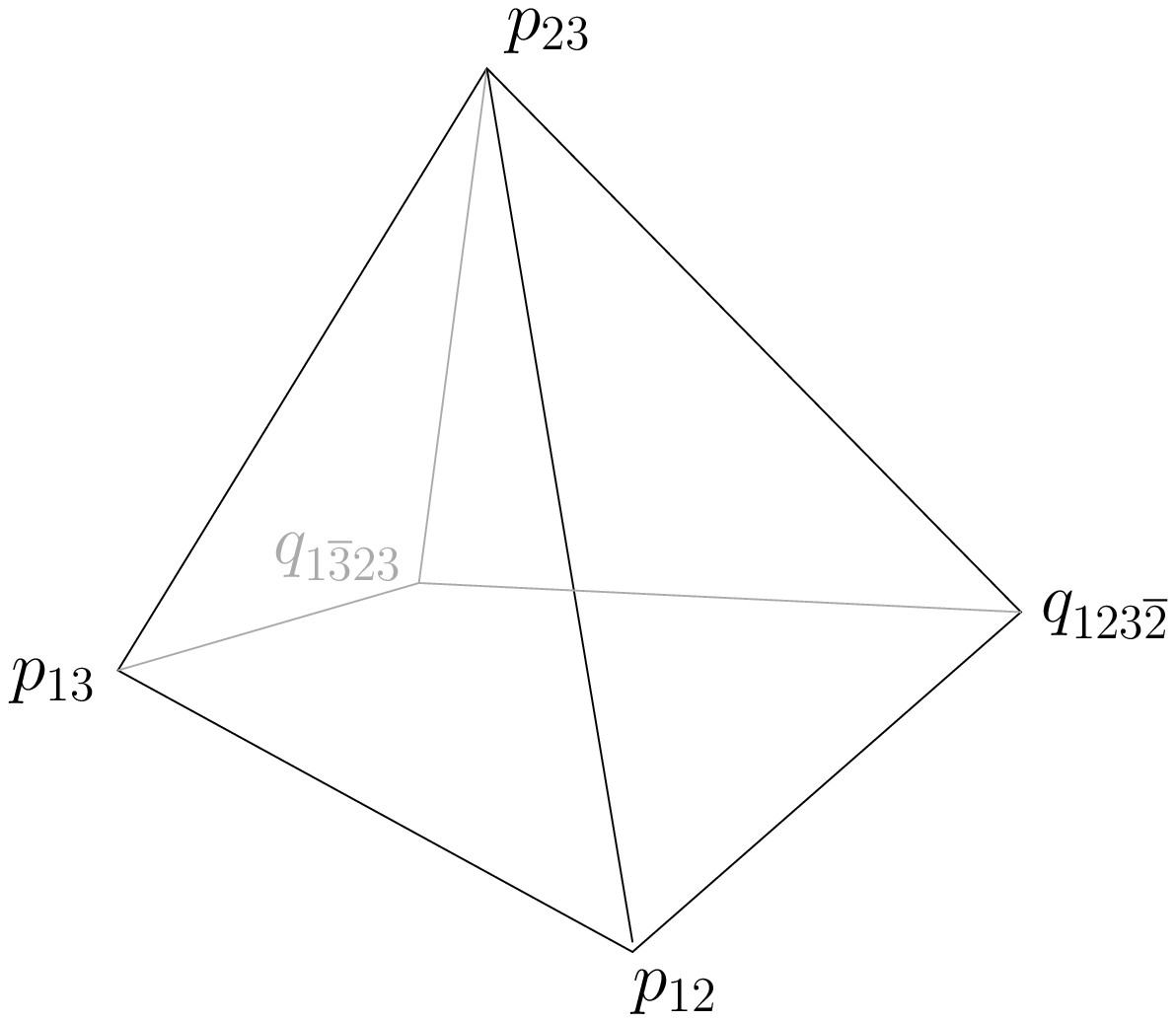, height = 0.35\textwidth}}\hspace{0.1\textwidth}
  \subfigure[vertices of $s_1$]{\epsfig{figure=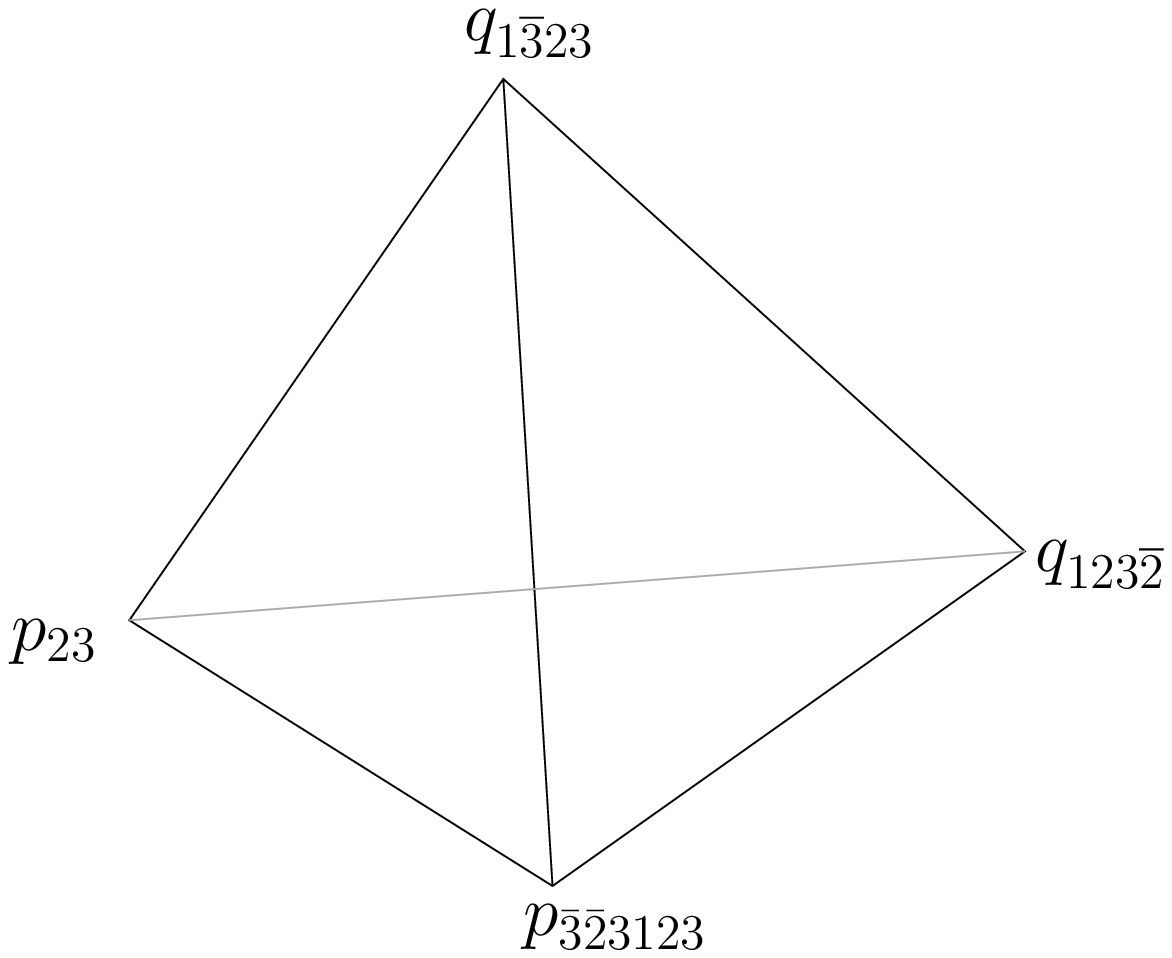, height = 0.35\textwidth}}
  \caption{Adjacency relations of $r_1$ and $s_1$ with other faces of
    $\widehat{E}$, for $p=3,4$. Larger font is used to label ridges,
    and smaller fonts to label lower-dimensional facets.}
  \label{fig:adjacency3-4}  
\end{figure}

\begin{figure}[htbp]
\centering
  \subfigure[ridges of $r_1$]{\epsfig{figure=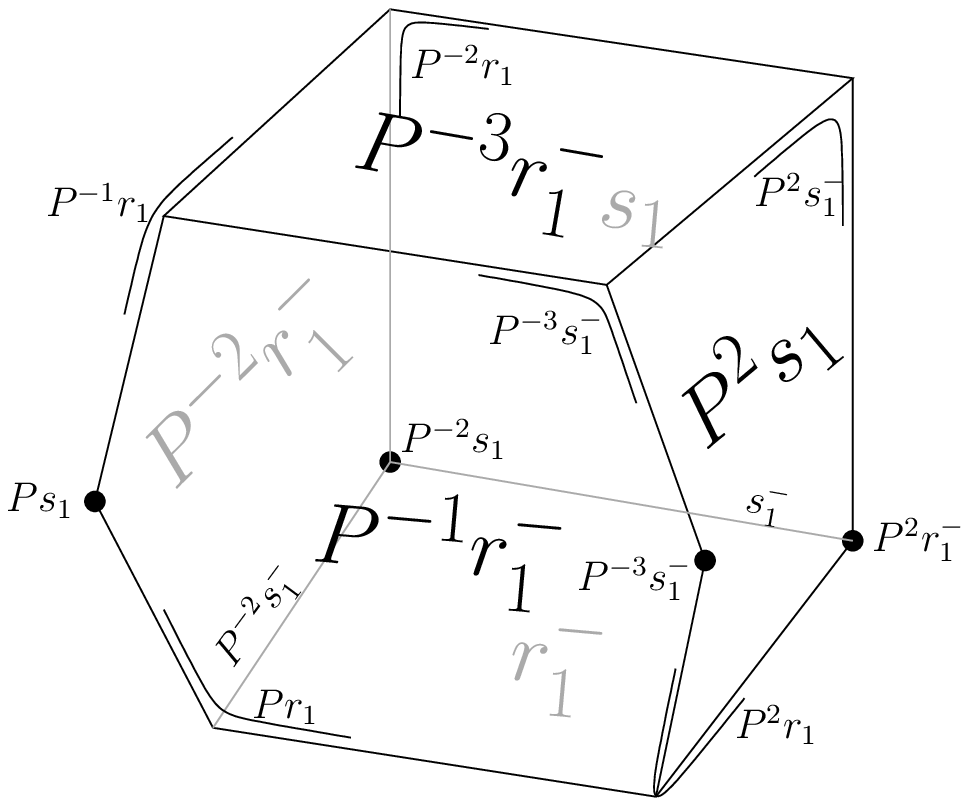, height = 0.33\textwidth}}\hspace{0.1\textwidth}
  \subfigure[ridges of $s_1$]{\epsfig{figure=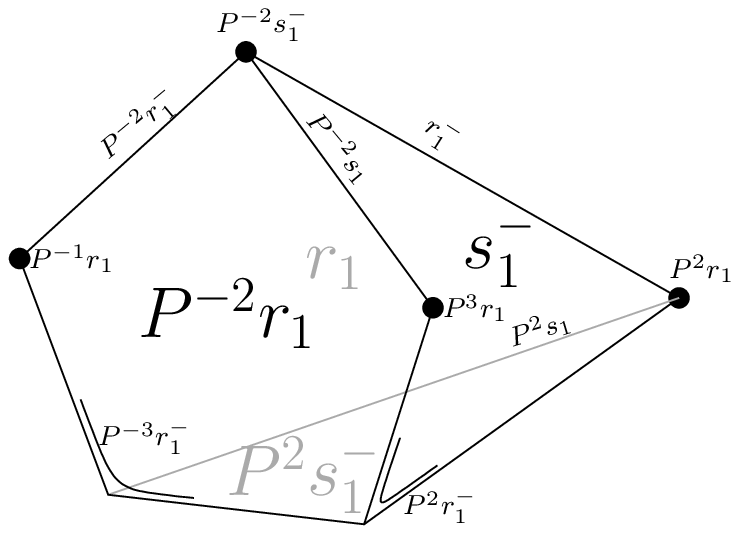, height = 0.3\textwidth}}\\
  \subfigure[vertices of $r_1$]{\epsfig{figure=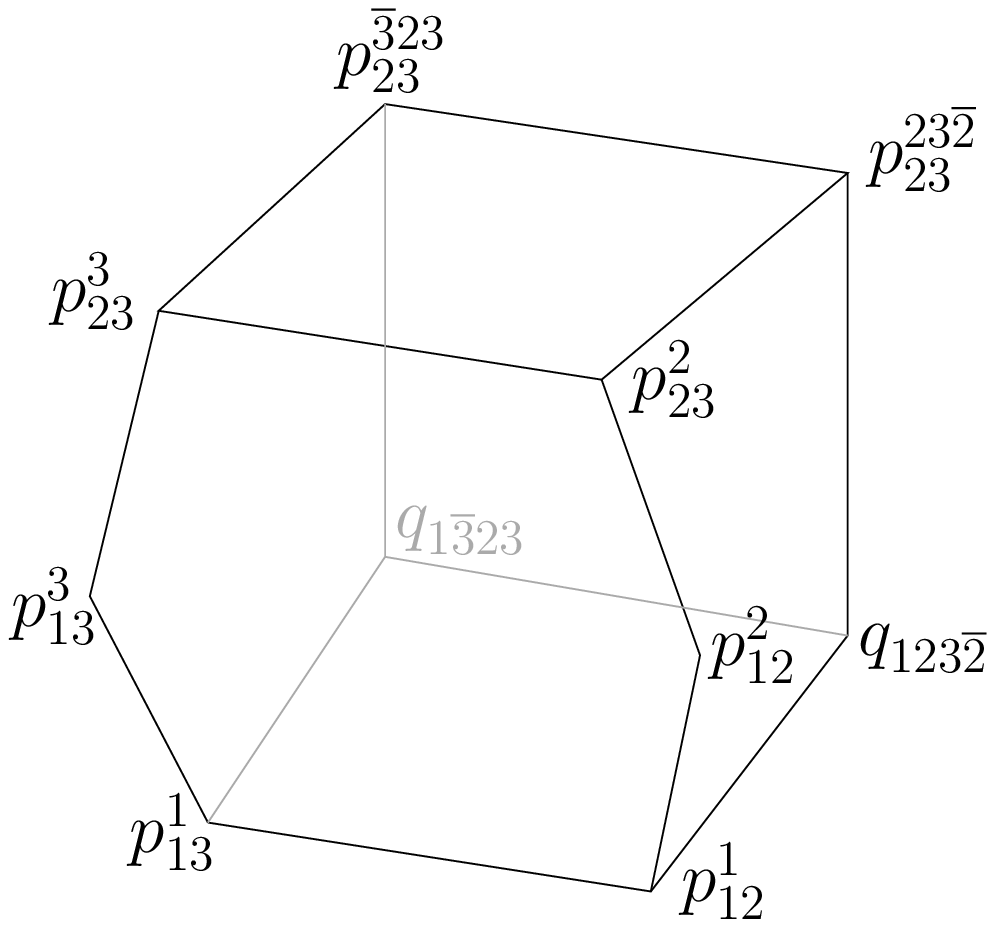, height = 0.33\textwidth}}\hspace{0.1\textwidth}
  \subfigure[vertices of $s_1$]{\epsfig{figure=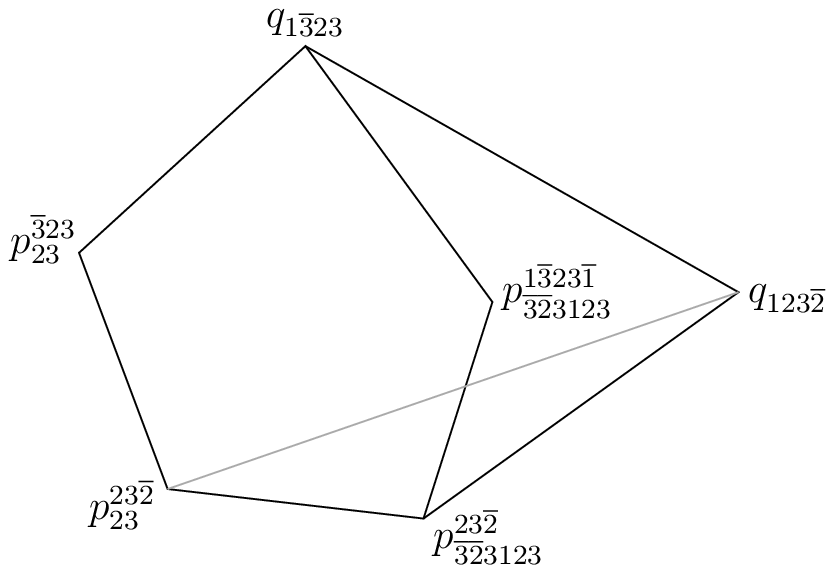, height = 0.33\textwidth}}
  \caption{Adjacency relations of $r_1$ and $s_1$ with other faces of
    $\widehat{E}$, for $p=5,6$. We draw curved lines when two
    neighboring edges lie on the same intersection of three
    bisectors.}
  \label{fig:adjacency5-6}  
\end{figure}

\begin{figure}[htbp]
\centering
  \subfigure[$r_1$]{\epsfig{figure=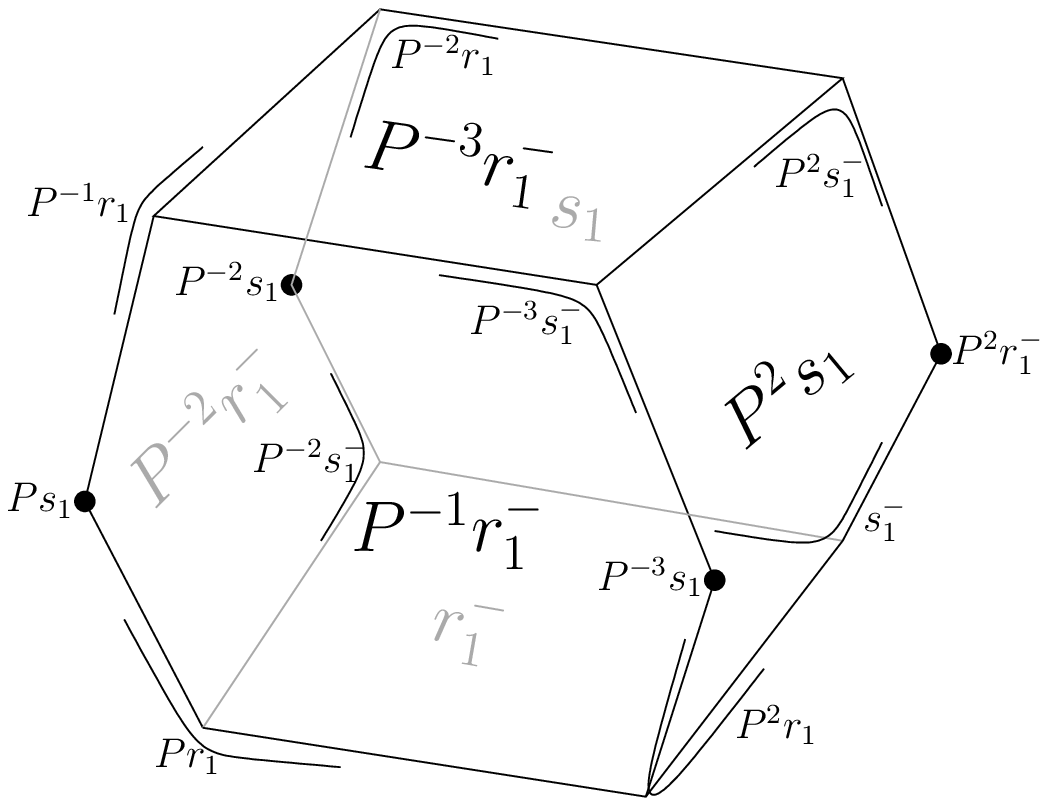, height = 0.33\textwidth}}\hspace{0.1\textwidth}
  \subfigure[$s_1$]{\epsfig{figure=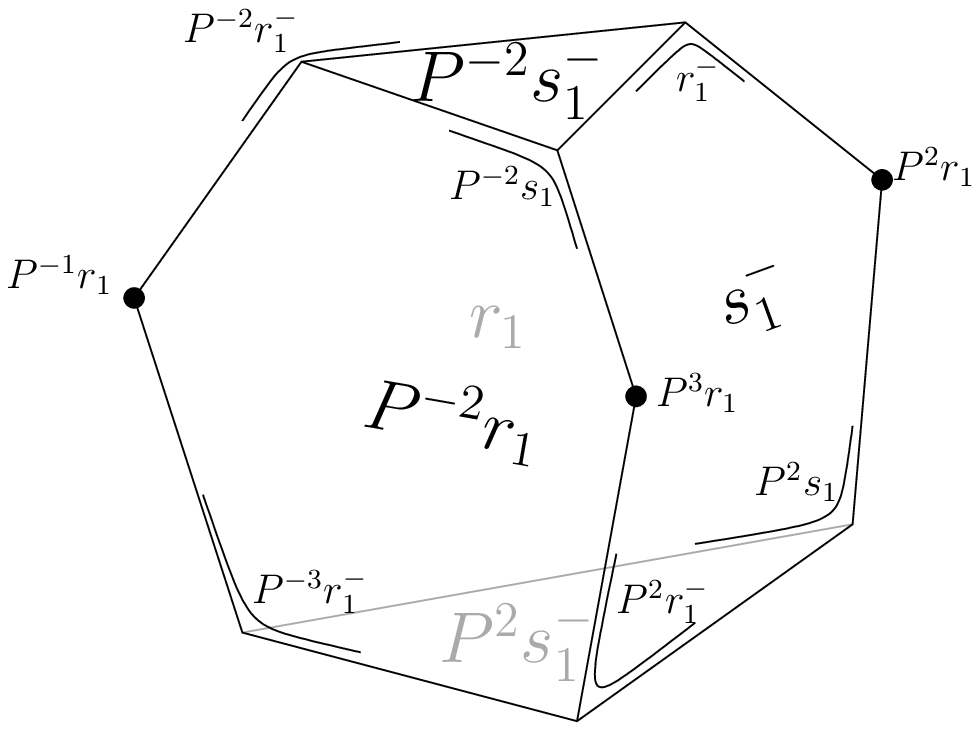, height = 0.33\textwidth}}\\
  \subfigure[$r_1$]{\epsfig{figure=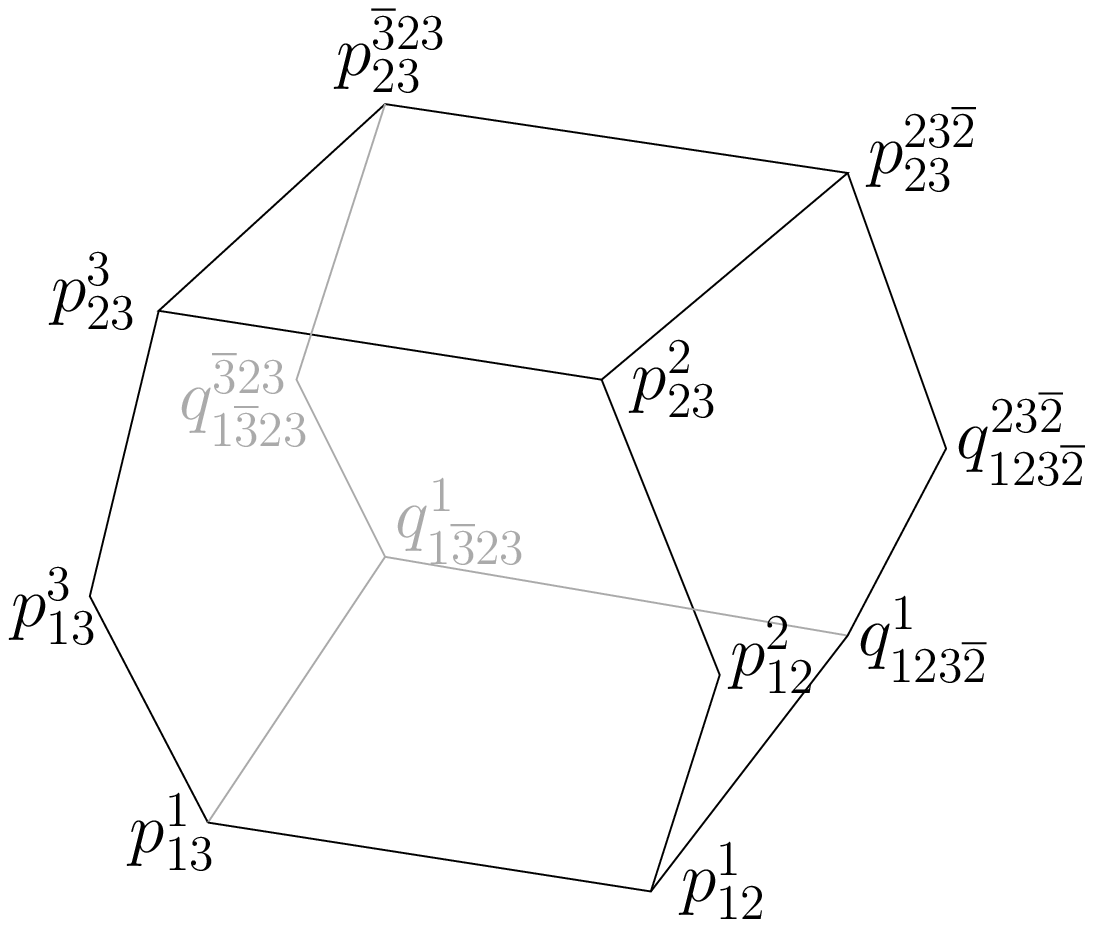, height = 0.33\textwidth}}\hspace{0.1\textwidth}
  \subfigure[$s_1$]{\epsfig{figure=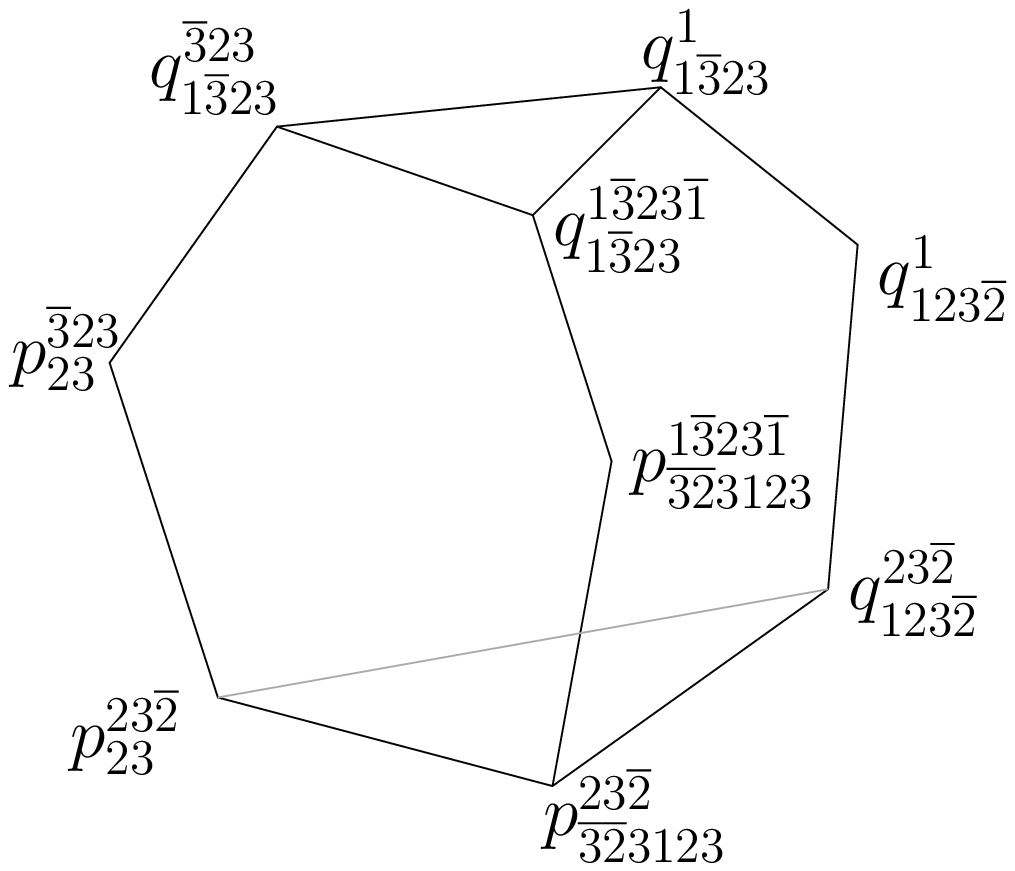, height = 0.33\textwidth}}
  \caption{Adjacency relations of $r_1$ and $s_1$ with other faces of
    $\widehat{E}$, for $p=8,12$.}
  \label{fig:adjacency8-12}  
\end{figure}

Figures~\ref{fig:adjacency3-4} to~\ref{fig:adjacency8-12} give a concise,
efficient description of the combinatorics. For each $3$-face, we give
two different pictures (top and bottom), one listing bisectors
containing each facet, one listing vertices. At first glance, the
pictures may seem a little cryptic, so we briefly explain how the
incidence relations can be read off.

Each cell in the picture is contained in the union of the cells whose names are on facets incident to that cell. As an example, we obtain the list of all bounding sides containing the
vertex labelled $p_{12}$ on part~(c) of
Figure~\ref{fig:adjacency3-4}. Obviously it is contained in $r_1$. The
picture in~(a) suggests that it is on three $2$-cells of $r_1$, whose
labels indicate that $p_{12}$ is also on $r_1^-$, $P^{-1}r_1^-$, and
$P^2s_1$.  Similarly, from the $1$-cells of $r_1$ containing $p_{12}$,
we see that it is also on $Pr_1$, $P^{-3}s_1^-$ and $P^2r_1$. Finally,
the labels near that vertex in part~(a) of the Figure indicate that it
is on $Pr_1^-$, $P^{-3}s_1$ and $P^{-1}s_1^-$. (See Remark 6.1 below for the vertex notation $p_u^v$, $q_u^v$).

Note that the pictures for $P^kr_1$ and $P^ks_1$ are obtained simply
by applying $P^k$, whereas those for $P^{k}r_1^-$ and $P^{k}r_1^-$ one
also needs to apply the antiholomorphic involution $\iota$ from
Section~\ref{sec:symmetry}. 
Recall that $\iota$ conjugates $R_1$ into $R_1^{-1}$, $S_1$ into
$S_1^{-1}$, and $P$ into $P^{-1}$, so the image of the figures in
part~(a) and~(b) under $\iota$ are obtained simply by changing
signs. For the figures in parts~(c) and~(d), one needs to use the
description of the action of $\iota$ on vertices according to their
labels (see Section~\ref{sec:symmetry}).

\begin{rmk}
  \begin{itemize}\item
  The vector ${\bf n}_{23}$ is a negative vector when $p=3$ and a null
  vector when $p=4$. The corresponding point $p_{23}$ is therefore
  in $\ch 2$ and on $\partial\ch 2$ respectively. Hence, for $p=3$ and $p=4$, 
  the combinatorics of the polyhedra $E$ are only
  the same when viewed in $\ch 2\cup \partial \ch 2$, some vertices
  being on the ideal boundary for $p=4$ (see also the discussion
  around Figure~\ref{fig:mirrors}). A similar remark is in order for
  $p=5$ and $p=6$, since for $p=6$ the polyhedron $E$ has vertices on
  the ideal boundary $\partial \ch 2$.
\item The vector ${\bf n}_{23}$ is a positive vector when $p\geqslant 5$. Therefore
the vertex $p_{23}$ of $r_1$ is replaced with a ridge in $m_{23}$, the
complex line polar to ${\bf n}_{23}$. We refer to this process as
{\bf truncation}. In $s_1$, the vertex $p_{23}$ is only replaced 
with an edge in $m_{23}$ since it only intersects two complex lines orthogonal 
to $m_{23}$, namely $m_{23\bar{2}}$ and $m_{\bar{3}23}$.
Similarly, the vertices $p_{12}$, $p_{13}$ of $r_1$ and the vertex
$p_{\bar{3}\bar{2}3123}$ of $s_1$ are replaced with edges in the corresponding
complex lines. A second truncation process occurs for the $q_*$ vertices
for $p\geqslant 8$. We use notation of the from $p_u^v$, $q_u^v$ for the new vertices appearing after truncation: $p_u^v$ or $q_u ^v$ denotes the intersection of the complex lines polar to ${\bf n}_u$ and ${\bf n}_v$.
  \end{itemize}
\end{rmk}

\subsubsection{Outline of the geometrization of the combinatorial model $\widehat{E}$}\label{sec:model}

In this section, we build a combinatorial model $\widehat{E}$ for
$E$.  To keep the notation reasonable, we do not introduce new symbols
for the facets of $\widehat{E}$, and simply label them the same as the
corresponding facets of $E$.

In order to define $\widehat{E}$, we first define its 3-skeleton
$\widehat{E}_3$. This is made up of 28 cells, attached according to
the adjacency relations indicated in Figures~\ref{fig:adjacency3-4}
to~\ref{fig:adjacency8-12}.

Note that each $2$-cell is on precisely two $3$-cells, which makes it
obvious which combinatorial $2$-cells get identified. The gluing is
then uniquely determined by matching labels of the vertices.  The
identifications for $0$ and $1$-cells are most conveniently read off
parts (c), (d) of the figures. Indeed the action of $\iota$ and $P$ on
the vertices can easily be read off the labels; see
Section~\ref{sec:statement} and also Tables~\ref{tab:verts3-4} to
\ref{tab:verts8-12} in the appendix (Section~\ref{sec:comb-data}).

It is easy to see that all $k$-cells of the $3$-skeleton
$\widehat{E}_3$ are homeomorphic to embedded closed balls, and in
particular for each $k$-cell $f$ of $\widehat{E}_3$, $\partial f$ is
an embedded $(k-1)$-sphere (see Figures~\ref{fig:adjacency3-4}
through~\ref{fig:adjacency8-12}).

We will also prove the following, which is much less obvious (see
Section~\ref{sec:sphere}).
\begin{theorem} \label{thm:3-sphere}
  $\widehat{E}_3$ is homeomorphic to $S^3$.  
\end{theorem} 
This allows us to define $\widehat{E}$, which is obtained by attaching a
single 4-cell to $\widehat{E}_3$ in the obvious manner.

We now describe a piecewise smooth geometric realization
$\phi:\widehat{E}\rightarrow \chb 2$, i.e. we describe a specific
$k$-ball in $\chb 2$ for each $k$-cell of $\widehat{E}$. Naturally,
the realization is uniquely determined by the requirement that each
3-cell be mapped into the appropriate bisector (or rather its closure
in $\chb 2$, since for $p=4$ or $6$ some vertices of $E$ lie on the
ideal boundary).

We define $\phi$ inductively on dimension, starting from vertices,
then extending it successively to the $1$-skeleton, then to the
$2$-skeleton, and so forth.  When realizing the $k$-faces, we will
check the following conditions on the restriction of $\phi$ to the
$k$-skeleton $\widehat{E}_k$:
\begin{itemize}
    \item (Consistency) For each $k$-face $f$, $\phi(f)$ is contained
      in all the bisectors indicated by the labels in
      Figures~\ref{fig:adjacency3-4} through~\ref{fig:adjacency8-12};
    \item (Correct side) $\phi(\widehat{E}_k)$ is entirely contained
      in $F$;
    \item (Embeddedness) $\phi$ gives an embedding of $\widehat{E}_k$.
\end{itemize}
The consistency condition is obvious for $2$ and $3$-faces, but it
requires some calculations for vertices and edges (note that these lie
on more bisectors than their codimension).

The fact that $\phi(\widehat{E}_k)$ lies inside $F$ amounts to showing
that for each bounding bisector $\B$, $\phi(\widehat{E}_k)$ is on the
same side of $\B$ as ${\bf p}$.  We will do this by analyzing, for
each $k$-face $f$, the intersection $f\cap\B$ with every bounding
bisector $\B$. Because the intersection patterns of $F$ are governed
by the intersection of bisectors, the ``correct side'' condition
actually implies the embeddedness.

Note that $\phi(\widehat{E})$ is connected, so
$\phi(\widehat{E})\subset F$ actually implies
$\phi(\widehat{E})\subset E$ (recall that $E$ is the connected
component of $F$ containing ${\bf p}$), provided we can connect ${\bf
  p}$ to a single point in $\phi(\widehat{E})$ within $F$. This is
contained in the following result, which is easily proved by giving
explicit parametrizations for the relevant geodesic segments (see
Section~\ref{sec:edges}), and checking that ${\bf p}$ is on the
correct side of the 28 bounding bisectors.
\begin{lemma}\label{lem:correctcomponent}
  For $p=3,4$ (respectively $p=5,6,8,12$) the geodesic segment $[{\bf
      p},p_{12}]$ (resp. $[{\bf p},p_{12}^1]$) intersects the 28
  bounding bisectors at most at the endpoint $p_{12}$
  (resp. $p_{12}^1$). It is contained in $F$, whereas the
  complementary geodesic ray, from $p_{12}$ (resp. from $p_{12}^1$) to
  $\partial\ch 2$ is entirely outside $F$.
\end{lemma}

We now sketch the general scheme that allows us to construct the
geometric realization $\phi$.  For vertices, the labels in the bottom
part of Figures~\ref{fig:adjacency3-4} through~\ref{fig:adjacency8-12}
determine an embedding of the $0$-skeleton of $\widehat{E}$ (each
vertex is the fixed point of a specific isometry in the group). The
embedding part of this statement corresponds simply to the fact that
all vertices in the list are distinct. The consistency condition will
be checked in Section~\ref{sec:vertices}, and the correct side
condition then amounts to checking a small set of (strict)
inequalities, which can easily be done with a computer.

Now let $k\geqslant 1$, and suppose $\phi$ is already known to give an
embedding of the $(k-1)$-skeleton.  In order to extend $\phi$ across a
$k$-cell $f$ of $\widehat{E}$, we exhibit a specific closed
$k$-dimensional ball $D^k$ in $\chb 2$ containing the $(k-1)$-skeleton
of $f$, which is in fact forced by the consistency condition.  More
specifically, for $k=1$, each ball is simply the (unique) real
geodesic segment between the two vertices realizing the $0$-skeleton
of $f$. For $k=2$ or $3$, the $k$-ball is given by the intersection of
$4-k$ (closures of) bisectors. The specific set of bisectors is
obtained from the labels in Figures~\ref{fig:adjacency3-4}
to~\ref{fig:adjacency8-12}. Finally, for $k=4$, the ball is of course
the whole of $\chb 2$.

Since we assume $\phi$ induces an embedding of the $(k-1)$ skeleton,
$\phi(\partial f)$ is a piecewise smoothly embedded $(k-1)$-sphere in
$D^k$ (with finitely many points on $\partial \ch2$). This gives a
well-defined $k$-ball that realizes $f$.  Implicit in the above
description is the assumption that the relevant pairs of bisectors
intersect in disks, which will be proved in Section~\ref{sec:ridgeconsistency}.
We postpone the proof of consistency and embeddedness of the $0$, $1$,
$2$ and $3$-skeletons, which will be discussed in
Sections~\ref{sec:vertices} through~\ref{sec:3-cells} (for clarity, we
collect all results about each dimension in a separate section). Among
all the consistency, correct side and embeddedness conditions, the
most difficult one is correct side condition for the
$2$-skeleton. This relies on difficult computations, explained in
detail in Section~\ref{sec:ridgeembedding}.

Assuming these results, we obtain a realization of all $\widehat{E}$,
by mapping the single $4$-cell to the 4-ball bounded by
$\phi(\widehat{E}_3)$ (this is a piecewise smoothly embedded copy of
$S^3$ so it bounds a well-defined closed ball).  Note that this ball
component contains the fixed point ${\bf p}$ of $P$, by
Lemma~\ref{lem:correctcomponent}.  As a summary, we get the following:
\begin{theorem}\label{thm:embedding}
  $\phi$ defines an embedding of $\widehat{E}$, with image
  contained in $F$.
\end{theorem}
The following result then follows from elementary topology.
\begin{corollary}\label{prop:comb-geom}
  The realization of $\widehat{E}$ is equal to $E$. In particular, $E$
  is a polyhedron in the sense of Section~\ref{sec:poin}.
\end{corollary}

\Pf (of Corollary~\ref{prop:comb-geom}) 
By Theorem~\ref{thm:embedding}, $\phi(\widehat{E})\subset E$.  We now
show $E\subset\phi(\widehat{E})$. No point of $E^\circ$ is contained
in any of the 28 bounding bisectors, hence $E^\circ$ has empty
intersection with $\phi(\widehat{E}_3)$. Since $E^\circ$ is connected
and contains ${\bf p}$, which is contained in $\phi(\widehat{E})$, we
get that $E\subset\phi(\widehat{E})$ as required.  
\EPf

\subsection{$\widehat{E}_3$ is homeomorphic to $S^3$}\label{sec:sphere}

        \begin{figure}[htbp]
          \centering
          \epsfig{figure=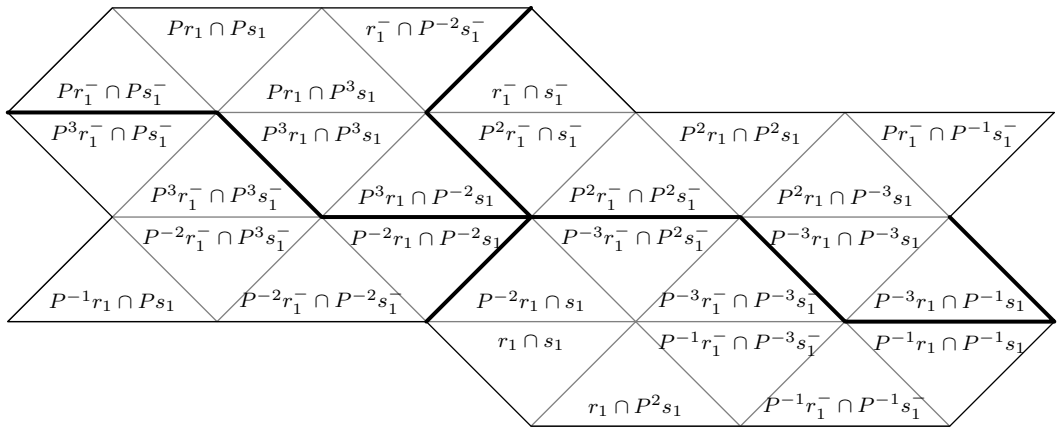, width=0.9\textwidth}
          \caption{The combinatorics of the intersection $T^r\cap
            T^s$, for $p=3,4$. This intersection is
            a torus, and accordingly the top and bottom sides
            (resp. left and right sides) are identified by
            translation.  The bold lines describe 2 splittings of the
            torus into annuli, corresponding to splittings of the solid
            tori into pairs of solid cylinders.}
           \label{fig:toriMeridians3-4}
        \end{figure}

        \begin{figure}[htbp]
          \centering
          \epsfig{figure=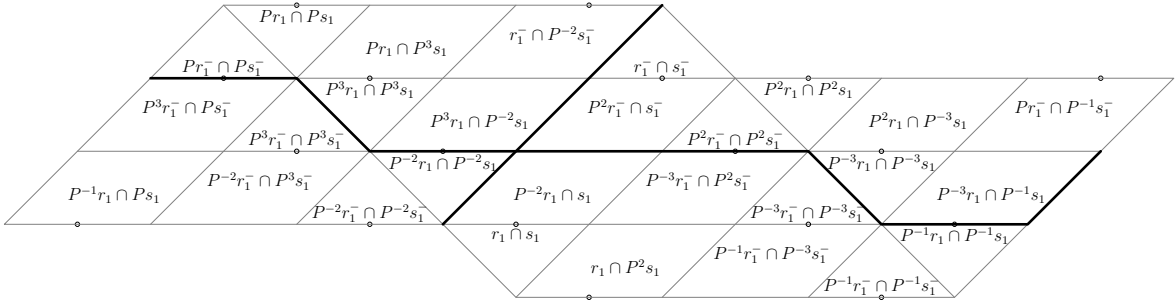, width=\textwidth}
          \caption{The combinatorics of the intersection $T^r\cap
            T^s$, for $p=5,6$.}
           \label{fig:toriMeridians5-6}
        \end{figure}

        \begin{figure}[htbp]
          \centering
          \epsfig{figure=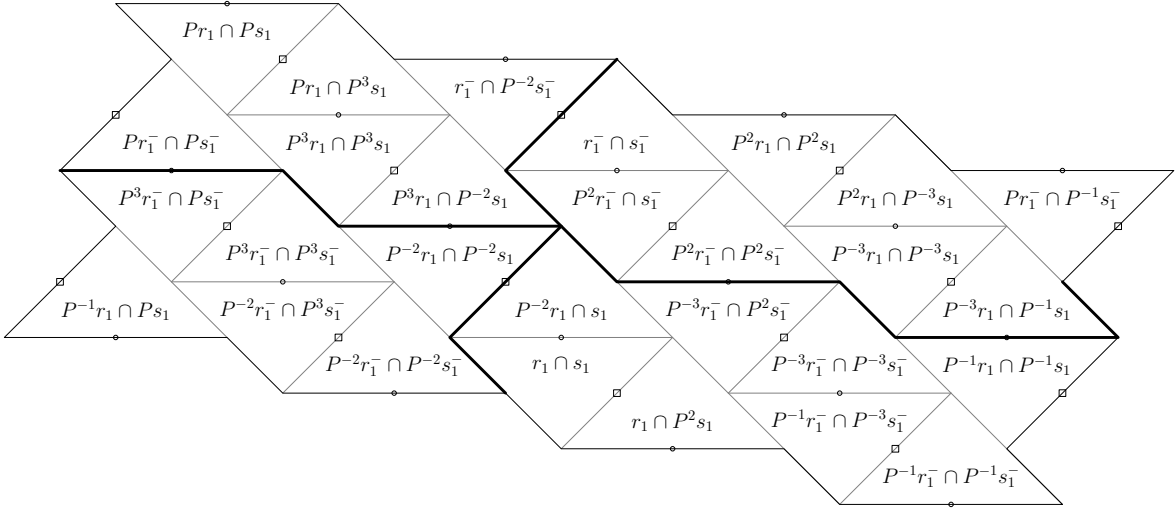, width=\textwidth}
          \caption{The combinatorics of the intersection $T^r\cap
            T^s$, for $p=8,12$. The torus is now not embedded, the
            white squares correspond to pinch points, i.e. they come
            in pairs that get identified in $\ch 2$.}
           \label{fig:toriMeridians8-12}
        \end{figure}

This section contains a proof of Theorem~\ref{thm:3-sphere}.  We will
exhibit $\widehat{E}_3$ as the union of two solid tori with common
boundary, and check that the boundary circles of a meridian in each intersect at a single point.

In order to do this, consider the following two unions of $3$-cells:
$$
  T^r=\bigcup_{k=-3}^{3}P^k(r_1\cup r_1^-)
    \quad\hbox{ and }\quad
  T^s=\bigcup_{k=-3}^{3}P^k(s_1\cup s_1^-).
$$ 
Note that these are clearly $P$-invariant (hence their intersection
is $P$-invariant as well). 

For $p\leqslant 6$, these will turn out to give a solid torus
decomposition. For $p=8,12$, the solid tori are obtained from these by
an adequate local surgery, to be discussed later in the proof (see
Lemma~\ref{lem:solidtori8-12}).

The combinatorial pattern of the intersection of $T^r$ and $T^s$ is
depicted in Figures~\ref{fig:toriMeridians3-4},
~\ref{fig:toriMeridians5-6} and~\ref{fig:toriMeridians8-12} for
various values of $p$ (compare this to Figure~6.2
of~\cite{schwartz4447}).  These pictures can be obtained by somewhat
painful bookkeeping from Figures~\ref{fig:adjacency3-4}
through~\ref{fig:adjacency8-12} and their obvious variations,
i.e. images under powers of $P$ and/or the antiholomorphic symmetry
$\iota$.

A few remarks are in order for the pictures to be read properly. First
note that we cannot draw a Euclidean plane figure, since hexagonal
faces of $T^r\cap T^s$ are often glued along two consecutive sides. For
instance, we represent hexagons by triangles, thinking of the midpoint
of the edge of the triangle as a vertex of the hexagon (see
Figure~\ref{fig:toriMeridians8-12}).

Secondly, the pictures are embedded in $\widehat{E}_3$ only for
$p=3,4,5$ and $6$; for $p=8$ or $12$, the intersection $T^r\cap T^s$
is obtained from the torus in Figure~\ref{fig:toriMeridians8-12} by
identifying $7$ pairs of points, indicated in the figure by square
vertices.  We will refer to these $7$ points as {\bf pinch points}.

We first treat the case $p\leqslant 6$. 
\begin{lemma}\label{lem:solidtori}
If $p\leqslant 6$ then $T^r$ and $T^s$ are both solid tori, with union
homeomorphic to $S^3$. 
\end{lemma}
Note that it would suffice to prove that $
\pi_1(T^r)=\pi_1(T^s)=\mathbb{Z}.  $ It would then follow from
Seifert-Van Kampen that $\pi_1(\widehat{E}_3)=\pi_1(T^r\cup T^s)=1$,
which implies that $\widehat{E}_3$ is homeomorphic to $S^3$ (this
argument requires checking that $\widehat{E}_3$ is a manifold, which
follows from the detailed study of the links of its vertices). Using
the solution of the Poincar\'e conjecture seems like overkill, so we
now give a bare hands proof.

We prove Lemma~\ref{lem:solidtori} by constructing disjoint
closed disks $D_1^{r}$ and $D_2^{r}$ in $T^r$ whose complements
are balls (and similarly $D_1^{s}$ and $D_2^{s}$ in $T^s$), and
check that $D_1^{r}$ and $D_1^{s}$ intersect in a single point.
The boundary of each of these disks is depicted in
Figures~\ref{fig:toriMeridians3-4} and~\ref{fig:toriMeridians5-6}; the
bold horizontal line is $\partial D_1^r$, the top (or bottom)
horizontal line is $\partial D_2^r$, and similarly for vertical lines
and $\partial D_j^s$.

\begin{figure}
  \subfigure[$D_1^r$]{\epsfig{figure=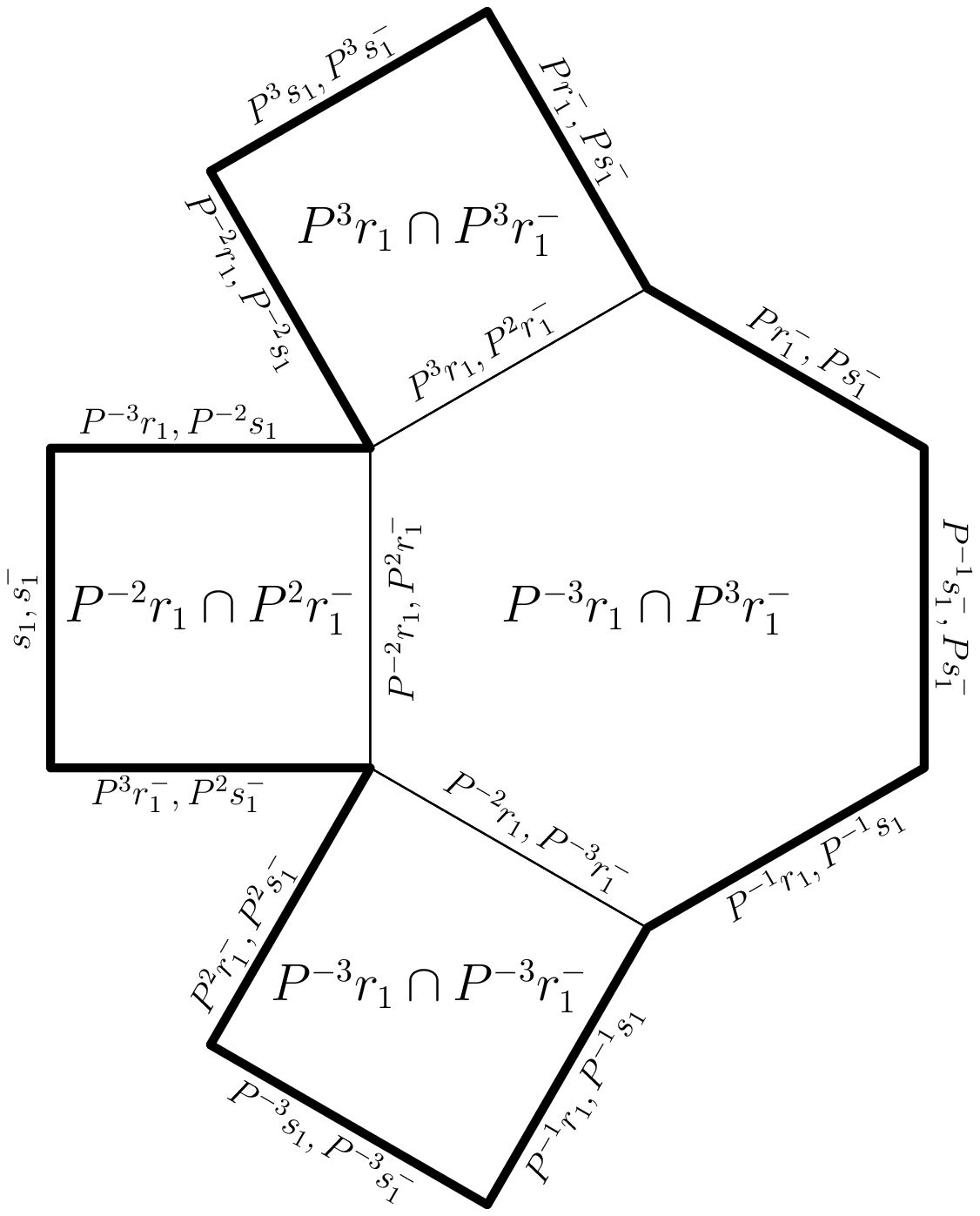,width=0.48\textwidth}}\hfill
  \subfigure[$D_2^r$]{\epsfig{figure=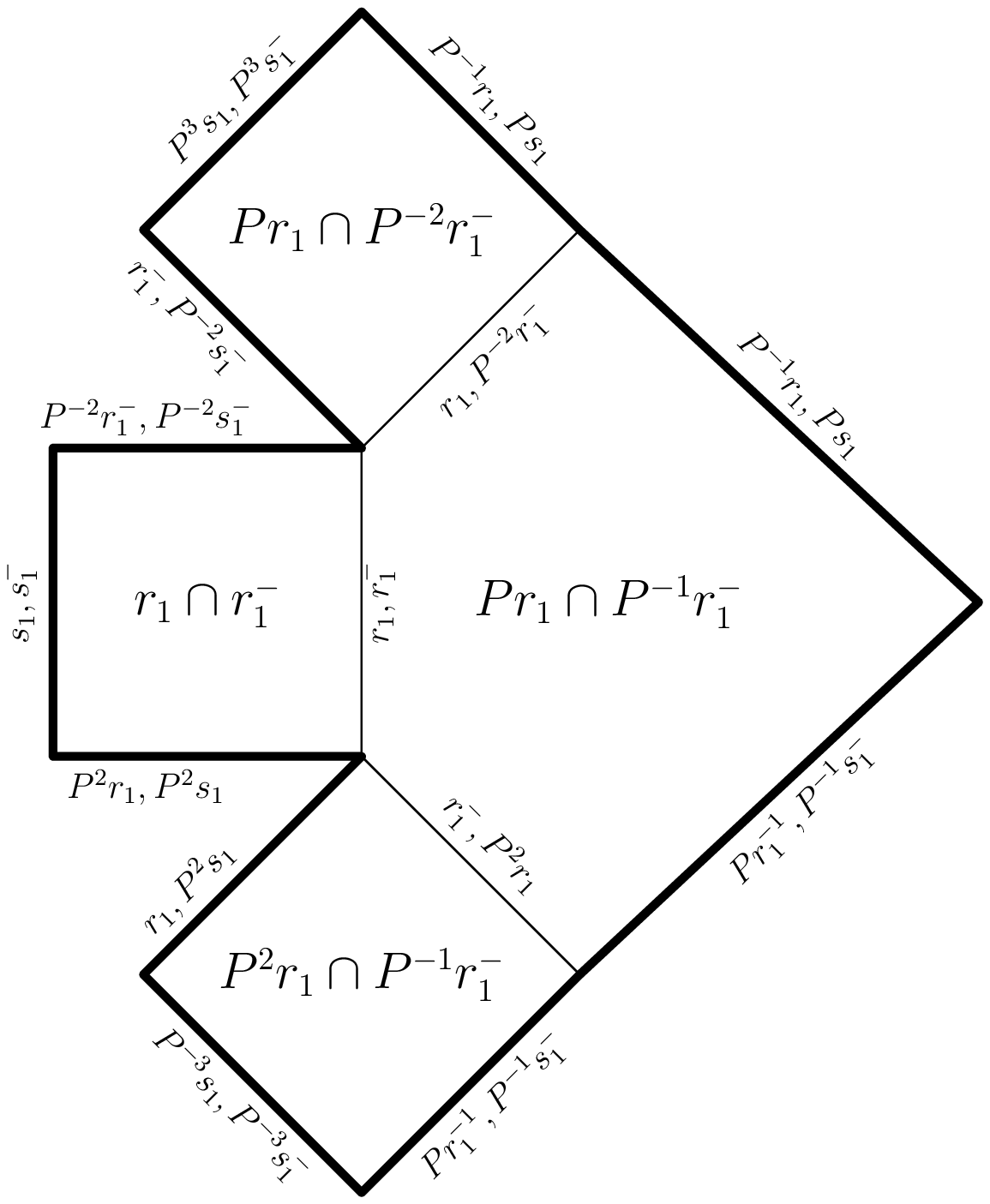,width=0.48\textwidth}}
  \caption{The combinatorics of the disks $D_1^r$ and $D_2^r$, for 
  $p=5,\,6$.}\label{fig:disks}
\end{figure}

Specifically,
\begin{equation}\label{eq:meridians}
\begin{array}{l}
  D_1^s=P^{-2}s_1\cap s_1^-,\\[0.1cm]
  D_2^s=P^{-1}s_1\cap P s_1^-,\\[0.1cm]
  D_1^r=(P^3r_1\cap P^3 r_1^-)\cup(P^{-2}r_1\cap P^2r_1^-)\cup(P^{-3}r_1\cap P^{-3}r_1^-)\cup(P^{-3}r_1\cap P^3r_1^-),\\[0.1cm]
  D_2^r= (Pr_1\cap P^{-1}r_1^-)\cup(Pr_1\cap P^{-2}r_1^-)\cup(r_1\cap r_1^-)\cup(P^2r_1\cap P^{-1}r_1^-).
\end{array}
\end{equation}
The fact that these are indeed disks is readily checked by using the
description of the combinatorics and adjacency relations between the
$3$-cells given in Figures~\ref{fig:adjacency3-4}
to~\ref{fig:adjacency5-6}. The following is clear from those pictures:
\begin{proposition}\label{prop:intersectionnumber}
  The circles $\partial D_1^r$ and $\partial D_1^s$ intersect in a single point.
\end{proposition}
For completeness, we give the combinatorial structure of the disks
$D_1^r$ and $D_2^r$ in Figure~\ref{fig:disks}.

We write $T^r=U^r\cup V^r$ (resp. $T^s=U^s\cup V^s$) for the
decompositions corresponding to splitting along these disks, which are
given in terms of $3$-cells of $\widehat{E}_3$ by the following:
\begin{equation}\label{eq:cylinders}
\begin{array}{l}
  U^s = Ps_1\cup Ps_1^- \cup P^3 s_1\cup P^3s_1^- \cup  P^{-2} s_1\cup P^{-2}s_1^-,\\[0.1cm] 
  V^s = s_1\cup s_1^- \cup P^2 s_1\cup P^2s_1^- \cup P^{-3} s_1\cup P^{-3}s_1^- \cup P^{-1} s_1\cup P^{-1}s_1^-,\\[0.1cm]
  U^r = r_1 \cup P^{-1}r_1\cup P^{-1} r_1^- \cup 
  P^{-2}r_1\cup P^{-2} r_1^- \cup  
  P^3r_1^-\cup P^{- 3} r_1^-,\\[0.1cm]
  V^r = r_1^-\cup 
  Pr_1\cup P r_1^- \cup 
  P^2r_1\cup P^2 r_1^- \cup  
  P^3 r_1\cup P^{- 3} r_1.
\end{array}
\end{equation}
In order to prove Lemma~\ref{lem:solidtori}, it is enough to prove
  the following.
\begin{lemma}\label{lem:cylinders}
  $U^s$ and $V^s$ are homeomorphic to balls, and intersect precisely
  in the two disks $D_1^s$, $D_2^s$.  Likewise, $U^r$ and $V^r$ are
  homeomorphic to balls, and intersect precisely in the two disks
  $D_1^r$, $D_2^r$.
\end{lemma}
\Pf
This follows from a careful study of the gluings, using the fact that
each $3$-cell of $\widehat{E}_3$ is a $3$-ball. This is easiest for $U^s$
and $V^s$, since $F=s_1\cup s_1^-$ is a ball, $F\cap P^{\pm 2}F$
are both disks and $F\cap P^4 F$ is empty; note that
$$
  V^s=F\cup P^2 F \cup P^4 F\cup P^6 F\quad\hbox{ and }\quad
  U^s=P F\cup P^3 F \cup P^5 F.
$$

For the $r$-splitting, we write each of $U^r$ and $V^r$ as an
increasing union, gluing a ball along a single closed disk at each
stage. More specifically, given a subcomplex $Z$ (which is either
$U^r$ or $V^r$), we give an explicit sequence $C_n$ of subcomplexes
with $C_0=F_0$ a $3$-cell,
$$ 
  C_{n}=C_{n-1}\cup F_{n}
$$ 
for a single $3$-cell $F_n$ of $\widehat{E}_3$, and which terminates
with $C_N=Z$. The point is to choose the $3$-cells $F_n$ so that $F_n\cap
C_{n-1}$ is homeomorphic to a disk.

To be specific, we give an explicit such sequence for $V^r$ that works
for all values of $p$:
$$
  P^{-3} r_1,\  P^2r_1^-,\ P^3r_1,\  Pr_1^-,\  P^2r_1,\  r_1^-,\  P r_1.
$$
At each stage, we check that $C_{n-1}\cap F_n$ is indeed a
disk by using the description of the combinatorics given in
Figures~\ref{fig:adjacency3-4} to~\ref{fig:adjacency8-12}
(and obvious variations, obtained by applying the suitable
power of $P$ and/or the symmetry $\iota$).

The corresponding sequence for $U^r$ is easily deduced from the one
for $V^r$ by applying the antiholomorphic involution$~\iota$. This finishes the proof of
Lemma~\ref{lem:cylinders}, hence also the proof of
Lemma~\ref{lem:solidtori}.
\EPf

We now consider the case $p\geqslant 8$. What remains true from the
statement of Lemma~\ref{lem:solidtori} is the following:
\begin{lemma}\label{lem:solidtori8-12}
  For $p=8$ and $12$, $T^r\cup T^s$ is homeomorphic to $S^3$.
\end{lemma}
However, $T^r$ and $T^s$ are now singular handlebodies, with
complementary singularities as we now explain.

A lot of the description of $T^r$ and $T^s$ for $p\leqslant 6$ goes
through. In particular, we use the same definition for $D_1^r$,
$D_2^r$, $D_1^s$, $D_2^s$, see equation~\eqref{eq:meridians}, and for
$U^r$, $V^r$, $U^s$, $V^s$, see equation~\eqref{eq:cylinders}.  The
picture in Figure~\ref{fig:toriMeridians8-12} makes it clear that
$U^r$ and $V^r$ are cylinders, whereas $U^s$ and $V^s$ are both
singular. We state this more precisely in Lemmas~\ref{lem:singular-r}
and~\ref{lem:singular-s}.

\begin{lemma}\label{lem:singular-r}
  Let $p=8$ or $12$. Then $D_1^r$ and $D_2^r$ are disjoint embedded
  closed disks, $U^r$ and $V^r$ are both homeomorphic to
  $3$-balls. The intersection $U^r\cap V^r$ is given by the disjoint
  union of $D_1^r$, $D_2^r$ and seven isolated points. The interior
  of $T^r$ is an open solid torus.
\end{lemma}

The seven isolated points are given by $q_{123\bar{2}}^{23\bar{2}}$
(which is the single intersection point of the 3-cells $r_1$ and
$P^{2}r_1^-$), and the six other points in its $P$-orbit. These are
indicated by square vertices in Figure~\ref{fig:toriMeridians8-12}. As
mentioned earlier, we refer to these as pinch points.

The corresponding result for $T^s$ is the following.
By a {\bf bowtie}, we mean the result of gluing two triangles along a single
vertex, which we call the {\bf apex} of the bowtie.
\begin{lemma}\label{lem:singular-s}
  Let $p=8$ or $12$. Then $D_1^s$, $D_2^s$ are disjoint bowties, with
  apex at a pinch point. $U^s$ contains $3$ pinch points, $V^s$
  contains $4$ pinch points, and the interior of $T^s$ is an open
  handlebody of genus $8$.
\end{lemma}

We now analyze the local structure around any of the seven pinch
points, say $q_{123\bar{2}}^{23\bar{2}}$. It lies on
precisely six 3-cells of $\widehat{E}_3$ (see Figure~\ref{fig:adjacency8-12}),
namely $r_1$, $P^2r_1^-$, $s_1$, $P^2s_1$, $s_1^-$ and $P^2s_1^-$, and
appears as the intersection of two different pairs of hexagons on the
boundary torus,
$$
r_1\cap s_1,\ r_1\cap P^2s_1
$$
on the one hand, and 
$$
P^2r_1^-\cap s_1^-,\ P^2r_1^-\cap P^2s_1,
$$
see Figure~\ref{fig:toriMeridians8-12}.

\begin{figure}[htbp]
  \centering
  \epsfig{figure=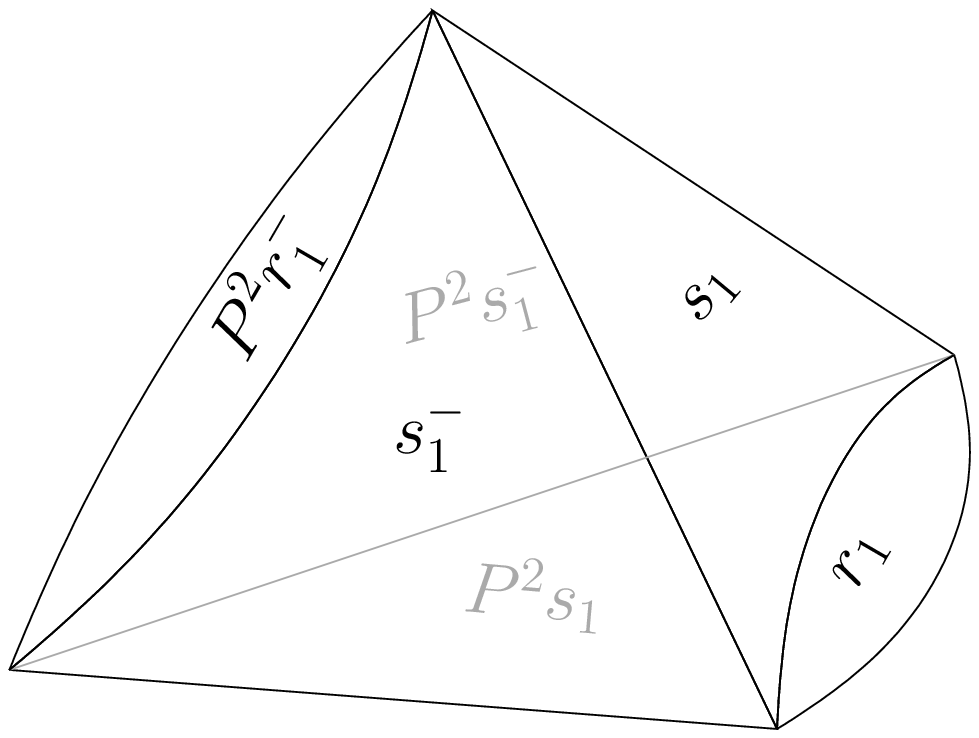, height=0.2\textwidth}
  \caption{Link of $q_{123\bar{2}}^{23\bar{2}}$, for $p=8$ or $12$.}
  \label{fig:link8-12}
\end{figure}

The link of that vertex is shown in Figure~\ref{fig:link8-12}. Note
that it is homeomorphic to a sphere, so even though both pieces $T^r$
and $T^s$ are singular, the union $T^r\cup T^s$ is a $3$-manifold.

We construct a non-singular solid torus $\widetilde{T}^r$
(resp. $\widetilde{T}^s$) from $T^r$ (resp. $T^s$) by performing an
obvious surgery in a small ball near each pinch point, using the
picture of the link given in Figure~\ref{fig:link8-12} as a guide.
\begin{figure}[htbp]
  \centering
  \epsfig{figure=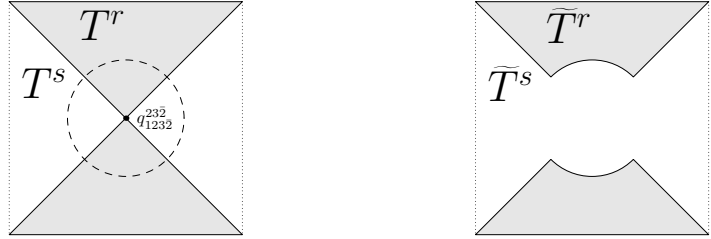, width=0.6\textwidth}
  \caption{We draw vertical sections of a rotationally symmetric local
    model (think of $T^r$ near the pinch point as a circular
    cone). After surgery near the pinch points, the white bowties become
    disks.}
  \label{fig:coneLink}
\end{figure}

After surgery, the two bowties (shown in white in
Figure~\ref{fig:coneLink}) become disjoint disks, $\widetilde{D}_1^s$
and $\widetilde{D}_2^s$, which exhibits $\widetilde{T}^s$ as a solid
torus (union of two cylinders).  The disk $D_1^r$ is unaffected by the
surgery, and it clearly intersects $\widetilde{D}_1^s$ in a single
point, since the intersection point of the bold lines in
Figure~\ref{fig:toriMeridians8-12} occurs away from the pinch
points.\EPf

\subsection{Geometric realization} \label{sec:realization}

\subsubsection{Realization of the vertices of $\widehat{E}$} \label{sec:vertices}

In this section, we describe in detail the $0$-skeleton of
$\widehat{E}$, depending on the parameter $p$. We describe the
geometric realization $\phi$ on the level of the $0$-skeleton, check
the consistency of this realization with the labels in
Figures~\ref{fig:adjacency3-4} through~\ref{fig:adjacency8-12}, as
well as the embeddedness of the realization (see
Section~\ref{sec:model} for the explanation of this terminology).

Consistency is checked by finding the complete list of bounding
bisectors (see Section~\ref{sec:domain}) that contain the realization
of each vertex. This will be done case by case in
Propositions~\ref{prop:vertbisectors3-4},~\ref{prop:vertbisectors5-6}
and~\ref{prop:vertbisectors8-12}.  Checking embeddedness then requires
no further verifications, since the vertices turn out to be uniquely
determined by the set of bounding bisectors that contain them. The
correct side condition for vertices amounts to proving a small number
of strict inequalities, which is easily done with a computer
(the point is that the relevant vertices and bisectors are
  $\k$-rational, where $\k=\Q(a,\tau)$, see the discussion in
  Section~\ref{sec:numberfields}).

As mentioned in Section~\ref{sec:model}, the realization of the
vertices of $\widehat{E}$ is determined by
Figures~\ref{fig:adjacency3-4} through~\ref{fig:adjacency8-12}, where
vertices are described as fixed point sets of the isometry given in
word notation by their labels (see
Section~\ref{sec:wordnotation}). The figures only list the vertices of
$r_1$ and $s_1$, but the other ones are obtained from those by
applying $\iota$ and/or a power of $P$. For completeness, we list the
$P$-orbits of these vertices, depending on the value of $p$. There are
two types of vertices, which we call $p_*$-vertices and $q_*$-vertices
(see Section~\ref{sec:wordnotation}).

\smallskip
\noindent
{\bf Case $p=3$ or $4$:} In this case there are two $P$-orbits of
vertices, so that $\widehat{E}$ has $14$ vertices. There is one
$P$-orbit of $p_*$-vertices and one $P$-orbit of $q_*$-vertices.
\begin{enumerate}
\item The realization of the vertex $p_w$ is the unique fixed point in
  $\chb 2$ of the element of $\Gamma$ represented by the word $w$. The
  seven $p_*$-vertices are realized by ideal vertices when $p=4$.
\item The realization of the vertex $q_w$ is the unique fixed point in
  $\ch 2$ of the map $w\in\Gamma$.
\end{enumerate}
A representative of the $P$-orbit of $p_*$-vertices is $p_{12}$,
realized as the unique fixed point of $R_1R_2$. A representative of
the $q_*$-vertices is $q_{123\bar{2}}$, realized by the unique fixed
point of $R_1R_2R_3R_2^{-1}$.  The $P$-orbits are easy to compute from
the fact that $P=R_1J$ (and $P^7=\textrm{Id}$), for instance
\begin{eqnarray*}
  P 1 \bar P & = &  1 J 1 \bar J \bar 1 \  = \  1 2 \bar 1, \\
  P 23\bar 2 \bar P & =  & 1 J 23\bar 2 \bar J \bar 1
  \  = \  1 31\bar 3 \bar 1\  =\  \bar 3 1 3,
\end{eqnarray*}
where the last equality follows from the higher braid relation
$(31)^2=(13)^2$. These two equalities give $P q_{123\bar 2} =
q_{12\bar 1 \bar 3 13}$. The vertices are given in
Tables~\ref{tab:verts3-4} and \ref{tab:verts3-6} in the appendix
(Section~\ref{sec:comb-data}).  From the preceding discussion, we
obtain a realization of the $0$-skeleton of $\widehat{E}$ into $\chb
2$, i.e. a definition of $\phi$ on the level of the $0$-skeleton; we
now check the consistency of this embedding with the labels of
Figures~\ref{fig:adjacency3-4} to~\ref{fig:adjacency8-12}.

First, let us show that $\phi$ sends the vertices of $r_1$ and $s_1$
to points in the bisectors $\Rc$ and $\Sc$ respectively.
\begin{lemma}\label{lem:phi-verts3-4}
\begin{enumerate}
\item The bisector $\Rc$ contains $p_{12}$, $p_{13}$, $p_{23}$,
  $q_{123\bar{2}}$ and $q_{1\bar{3}23}$.
\item The bisector $\Sc$ contains $p_{23}$, $p_{\bar{3}\bar{2}3123}$,
  $q_{123\bar{2}}$ and $q_{1\bar{3}23}$.
\end{enumerate}
\end{lemma}
\Pf 
By construction, ${\bf n}_{23}$ lies on the real spine of
$\Rc$. This vector projects to $p_{23}$, the unique fixed point in
$\chb 2$ of $23$.  Likewise, the complex line $m_1$ fixed by $R_1$ is
a slice of $\Rc$.  The fixed points of $12$, $13$, $123\bar{2}$ and
$1\bar{3}23$ all lie on $m_1$ and so are contained in $\Rc$ (see notation from section 4.3). This
proves part (1).  For part (2), begin by observing that the braid
relation (see Proposition~\ref{prop:3-braiding}) implies
$\bar{3}\bar{2}3123=1\bar{3}23\bar{1}23\bar{2}$ and so its fixed point
lies on $m_{23\bar{2}}$ and $m_{1\bar{3}23\bar{1}}$.  The rest of the
proof follows as in part (1).  
\EPf

To finish the consistency verification, we list the bounding bisectors
containing $p_{12}$ and $q_{123\bar 2}$ (the list of bisectors
containing \emph{any} vertex is then obtained by applying suitable powers of
$P$).

\begin{proposition} \label{prop:vertbisectors3-4}
Let $p=3$ or $4$. Then among the 28 bounding bisectors:
\begin{enumerate}
\item The vertex $p_{12}$ is on 
  $\Rc$, $\Rc^-$,
  $P\Rc$, $P^{-1}\Rc^-$, $P^2\Rc$, $P^2\Sc$, $P^{-3}\Sc$, 
  $P^{-3}\Sc^-$, $P\Rc^-$, $P^{-1}\Sc^-$;
\item The vertex $q_{123\bar 2}$ is on
  $\Rc$, $\Rc^-$,
  $\Sc$, $\Sc^-$, $P^2\Rc$, $P^2\Rc^-$, $P^2\Sc$, $P^2\Sc^-$.
\end{enumerate}
\end{proposition}
\Pf 
This follows from Lemma~\ref{lem:phi-verts3-4} and
Tables~\ref{tab:verts3-4},~\ref{tab:verts3-6} in
Section~\ref{sec:comb-data} (see also the action of $\iota$ on the
vertices described in Section~\ref{sec:statement}). We only treat a couple of
cases. For instance, $p_{12}\in P^{-1}\Rc^-$ is equivalent to
$Pp_{12}=p_{123\bar1}\in\Rc^-$. Applying $\iota$, this is equivalent
to $p_{23}\in\Rc$, which was proved in Lemma~\ref{lem:phi-verts3-4}.

The fact that $p_{12}$ is on $P^{-3}\Sc^-$ translates into
$P^3p_{12}=p_{\bar3\bar23123}\in \Sc^-$, or equivalently by applying
$\iota$, $P^{-3}p_{13}\in\Sc$. Since $Pp_{13}=p_{12}$ and $P^7=id$, we
have $P^{-3}p_{13}=P^3p_{12}$. 
\EPf

\smallskip
\noindent
{\bf Case $p=5$ or $6$:} For $p\geqslant 5$ the word $(R_1R_2)^2$ is a
complex reflection fixing the complex line $m_{12}$. In that case, we
replace the vertex $p_{12}$ with four vertices $p^1_{12}$,
$p^{2}_{12}$, $p_{12}^{\bar 212}$ and $p_{12}^{12\bar 1}$, that are
given by the intersections of $m_{12}$ with $m_1$, $m_2$, $m_{\bar{2}12}$, and
$m_{12\bar{1}}$ respectively. The $q_*$-vertices are unchanged.  Thus,
there are five $P$-orbits of vertices, so that $\widehat{E}$ has $35$
vertices.
\begin{enumerate}
\item The realization of the vertex $p^v_w$ is the intersection of the
  complex lines $m_v$ and $m_w$, fixed by $v$ and $w^2$ respectively.
\item The realization of the vertex $q_w$ is the unique fixed point of
  $w$.  The seven $q_*$ vertices map to ideal points when $p=6$.
\end{enumerate}
The cycles of $q_*$-vertices are the same as for $p=3,\,4$ (see the second
row of Table~\ref{tab:verts3-6}). 
The four orbits of $p_*$-vertices are listed in Table~\ref{tab:verts5-12}. 
\begin{lemma}\label{lem:phi-verts5-6}
\begin{enumerate} 
\item The bisector $\Rc$ contains the vertices $p^2_{12}$, $p^1_{12}$,
  $p^1_{13}$, $p^3_{13}$, $p^3_{23}$, $p^2_{23}$, $p^{23\bar 2}_{23}$, $p^{\bar 323}_{23}$,
  $q_{123\bar{2}}$ and $q_{1\bar{3}23}$.
\item The bisector $\Sc$ contains the vertices $p^{\bar{3}23}_{23}$,
  $p^{23\bar{2}}_{23}$, $p^{23\bar{2}}_{\bar{3}\bar{2}3123}$,
  $p^{1\bar{3}23\bar{1}}_{\bar{3}\bar{2}3123}$, $q_{123\bar{2}}$ and
  $q_{1\bar{3}23}$.
\end{enumerate}
\end{lemma}

\Pf 
For $\Rc$, the statement is obvious for eight of the ten vertices,
since vertices of the form $p^1_*$ and $p^*_{23}$ lie on $m_1$ and
$m_{23}$ respectively, and these are complex slices of $\Rc$.  The
last two are handled using Lemma~\ref{lem:orthslice}. We treat the
case of $p^2_{12}$ in detail, $p^3_{13}$ is entirely similar.

Note that $m_{12}$ is orthogonal to $m_1$, since $R_1$ and
$(R_1R_2)^2$ commute. Similary, $m_{2}$ is orthogonal to
$m_{23}$. Lemma~\ref{lem:orthslice} implies that $\Rc$ intersects each
of the complex lines $m_2$ and $m_{12}$ in a real geodesic, and by
definition $p_{12}^2$ is the intersection point of these two
geodesics.
This proves part (1). Part (2) is proved similarly.  
\EPf

The analogue of Proposition~\ref{prop:vertbisectors3-4} is the following:
\begin{proposition} \label{prop:vertbisectors5-6}
Let $p=5$ or $6$. Then among the 28 bounding bisectors:
\begin{enumerate}
\item the vertex $p_{12}^1$ is on $\Rc$, $\Rc^-$, $P\Rc$, $P^{2}\Rc$,
  $P^2\Sc$, $P^{-1}\Rc^-$;
\item the vertex $p_{12}^2$ is on $\Rc$, $P^2\Rc$, $P^2\Sc$,
  $P^{-3}\Sc$, $P^{-3}\Sc^-$, $P^{-1}\Rc^-$;
\item the vertex $p_{23}^2$ is on $\Rc$, $P^2\Sc$, $P^2\Sc^-$,
  $P^{-3}\Rc^-$, $P^{-3}\Sc^-$, $P^{-1}\Rc^-$;
\item the vertex $p_{23}^3$ is on $\Rc$, $P^{-3}\Rc^-$, $P^{-3}\Sc^-$,
  $P^{-2}\Rc^-$, $P^{-1}\Rc$, $P^{-1}\Rc^-$;
\item the vertex $q_{123\bar 2}$ is on $\Rc$, $\Rc^-$, $\Sc$, $\Sc^-$,
  $P^2\Rc$, $P^2\Rc^-$, $P^2\Sc$ and $P^2\Sc^-$.
\end{enumerate}
\end{proposition}

\smallskip
\noindent
{\bf Case $p=8$ or $12$:} In this case there are seven $P$-orbits of
vertices (so that $\widehat{E}$ has $49$ vertices). The cycles of
$p_*$-vertices are the same as for $p=5,6$ (there are four of them,
see Table~\ref{tab:verts5-12}), and there are three orbits of
$q_*$-vertices, listed in Table~\ref{tab:verts8-12}. Since
$(R_1R_2R_3R_2^{-1})^3$ is a complex reflection fixing the complex
line $m_{123\bar{2}}$, we replace the vertex $q_{123\bar{2}}$ with
three vertices lying in this complex line.
\begin{enumerate}
\item The realization of the vertex $p^v_w$ is the intersection of the
  complex lines $m_v$ and $m_w$, fixed by $v$ and $w^2$ respectively.
\item The realization of the vertex $q^v_w$ is the intersection of the complex
lines $m_v$ and $m_w$, fixed by $v$ and $w^3$ respectively.
\end{enumerate}

The following lemma is proved as before.
\begin{lemma}\label{lem:phi-verts8-12}
\begin{enumerate} 
\item The bisector $\Rc$ contains the vertices $p^2_{12}$, $p^1_{12}$,
  $p^1_{13}$, $p^3_{13}$, $p^3_{23}$, $p^2_{23}$, $p^{23\bar 2}_{23}$,
  $p^{\bar 323}_{23}$, $q^{23\bar{2}}_{123\bar{2}}$,
  $q^1_{123\bar{2}}$, $q^1_{\bar{3}23}$ and
  $q^{\bar{3}23}_{1\bar{3}23}$.
\item The bisector $\Sc$ contains the vertices $p^{\bar{3}23}_{23}$,
  $p^{23\bar{2}}_{23}$, $p^{23\bar{2}}_{\bar{3}\bar{2}3123}$,
  $p^{1\bar{3}23\bar{1}}_{\bar{3}\bar{2}3123}$,
  $q^{23\bar{2}}_{123\bar{2}}$, $q^1_{123\bar{2}}$,
  $q^1_{1\bar{3}23}$, $q^{\bar{3}23}_{1\bar{3}23}$ and
  $q^{1\bar{3}23\bar{1}}_{1\bar{3}23}$.
\end{enumerate}
\end{lemma}

\begin{proposition} \label{prop:vertbisectors8-12}
  Let $p=8$ or $12$. The $p_*$-vertices are on the same bounding
  bisectors as in the case $p=5,6$, and the $q_*$-vertices are on the
  following bounding bisectors:
  \begin{enumerate}
  \item the vertex $q^1_{1\bar323}$ is on $\Rc$, $\Rc^-$, $\Sc$,
    $\Sc^-$, $P^{-2}\Rc^-$ and $P^{-2}\Sc^-$;
  \item the vertex $q^{\bar 323}_{1\bar323}$ is on $\Rc$, $\Sc$,
    $P^{-2}\Rc$, $P^{-2}\Rc^-$, $P^{-2}\Sc$ and $P^{-2}\Sc^-$;
  \item the vertex $q^{23\bar2}_{123\bar2}$ is on $\Rc$, $\Sc$,
    $\Sc^-$, $P^{2}\Rc^-$, $P^{2}\Sc$ and $P^{2}\Sc^-$.
  \end{enumerate}
\end{proposition}

\subsubsection{Realizing edges of $\widehat{E}$} \label{sec:edges}

In this section, we check the consistency and the embeddedness
condition for the realization of the $1$-skeleton.  Each $1$-cell is
realized as a geodesic segment joining the realization of its two
endpoints.  According to the description of $\widehat{E}$, each
$1$-cell is on four $3$-cells, so we need to check that the geodesic
segment realizing each $1$-cell is contained in the appropriate set of
four bisectors (the set of bisectors is indicated by the labels of
Figures~\ref{fig:adjacency3-4} to~\ref{fig:adjacency8-12}).  Recall
from Section~\ref{sec:bisectorsandgeodsics} that in order to check
that a geodesic segment $[v_0,v_1]$ is on a bisector $\B$, it is not
enough to check that its endpoints $v_0$ and $v_1$ are on $\B$.

\begin{proposition}\label{prop:1consistency}
  The realization of the $1$-skeleton of $\widehat{E}$ is consistent.
\end{proposition}

\Pf 
This follows from repeated use of Lemmas~\ref{lem:realspine}
or~\ref{lem:orthslice}.  We illustrate this in a couple of specific
cases.  First, suppose $p=3$ or $4$, and consider the geodesic through
$p_{12}$ and $p_{13}$. According to
Proposition~\ref{prop:vertbisectors3-4}, it should be contained in
$$
  \Rc,\quad \Rc^-,\quad P\Rc,\quad  P^{-1}\Rc^-.
$$ 
By definition of the points $p_{12}$ and $p_{13}$, they are both in
the mirror of $R_1$, which is a common complex slice of $\Rc$ and
$\Rc^-$.  Moreover, $p_{13}$ (resp. $p_{12}$) is on the real spine of
$P\Rc$ (resp. $P^{-1}\Rc^-$), so the geodesic is contained in $P\Rc$
(resp. $P^{-1}\Rc^-$).

Now suppose $p\geqslant 5$, and consider the geodesic $\alpha$ through
$p^1_{13}$ and $p^1_{12}$.  These two points are on $m_1$, which is a
common slice of $\Rc$ and $\Rc^-$. Moreover, $m_1$ is by construction
orthogonal to the complex line $m_{12}$, which is a complex slice of
$P^{-1}\Rc^-$, so Lemma~\ref{lem:orthslice} implies that $\alpha$ is
contained in $P^{-1}\Rc^-$. Similarly, $m_1$ is orthogonal to the
slice of $P\Rc$ given by $\n_{13}^\perp$, so $\alpha$ is contained in
$P\Rc$.

In some cases, one needs to use both Lemmas~\ref{lem:realspine}
and~\ref{lem:orthslice}.  For instance, for $p=5$ or $6$, $p^1_{12}$
and $q_{123\bar 2}$ are contained in $m_{1}$, so clearly the geodesic
joining them is in $\Rc$ and $\Rc^-$. The fact that it is contained in
$P^2\Rc$ follows from Lemma~\ref{lem:orthslice} (and the fact that
$\langle \n_1,\n_{12}\rangle = 0$), whereas the fact that it is
contained in $P^2\Sc$ follows from Lemma~\ref{lem:realspine} (and the
fact that $q_{123\bar 2}$ is on the real spine of $P^2\Sc$).

In some cases, the application of Lemma~\ref{lem:orthslice} is not
completely obvious, since the orthogonality check can be a bit more
involved. For instance, consider the geodesic through $p^{1\bar
  323\bar 1}_{\bar 3\bar 2 3123}$ and $q_{1\bar 323}$. Clearly these
two points are on $m_{1\bar 323\bar 1}$, so the geodesic is contained
in $\Sc^-$ and $P^{-2}\Sc$. It is contained in $\Sc$ simply because
$q_{1\bar323}$ is a point of its real spine, but in order to check
that it is contained in $P^{-2}\Rc$, one needs to check that $\langle
\n_{1\bar 323\bar 1}, \n_{\bar 3\bar 23123}\rangle = 0$.  Since
$\n_{1\bar 323\bar 1}=P^3\n_2$ and $\n_{\bar 3\bar 23123}=P^3\n_{12}$
the result follows from $\langle \n_2,\n_{12}\rangle=0$.  
\EPf

\begin{rmk}
  Lemma~\ref{lem:orthslice} shows a bit more than the relevant
  inclusions of geodesics in bisectors. It also implies that the
  geodesic is in a real slice of certain bisectors - when two adjacent
  edges are contained in real slices of the same bisector (and their
  intersection is not on the real spine of that bisector), the two
  real slices must actually coincide.
\end{rmk}

\begin{proposition}
  The realization of the $1$-skeleton of $\widehat{E}$ is embedded
  and contained in $F$.
\end{proposition}

\Pf Let $f$ be a $1$-cell of $\widehat{E}$. We have just proved that
its realization $\phi(f)$ is on (at least) four bounding
bisectors. For each other bounding bisector $\B$, we check that
$\phi(f)$ is on the correct half-space bounded by $\B$. Since the endpoints of $\phi(f)$ are $k$-rational (with
  $\k=\Q(a,\tau)$), this can be done using arithmetic in $\k$, see
  Section~\ref{sec:numberfields}.

Note that the situation where some 1-cells are tangent to some of
the bounding bisectors, alluded to in Section~\ref{sec:numberfields},
really does occur. Specifically, it occurs when an endpoint of the
1-face is a cusp, which happens at certain specific 
vertices when $p=4$ and $p=6$.

Note that because of the symmetry given by $P$ and $\iota$, it is
enough to study the tangencies that appear in (the intersection of
$\phi(\widehat{E})$ with)
\begin{equation}\label{eq:ridgesneeded}
  \Rc\cap \Rc^-,\ \Rc\cap \Sc,\ \Rc\cap P^{-3}\Rc^-,\ \Rc\cap
  P^{-2}\Rc^-,\ \Rc\cap P^{-1}\cap\Rc^-,\ \Sc\cap P^2\Sc^-,\ \Sc\cap
  P^{-2}\Sc^-.
\end{equation}
Although not technically needed for the proof, it can be instructive
to draw figures of these Giraud disks and their intersection with all
bounding bisectors. This is done in Figure~\ref{fig:1-3} for $\Rc\cap
\Sc$, where it should be apparent that tangencies only occur for $p=4$
and $p=6$ (see parts~(b) and~(d) of the figure).
In order for the figure to be read properly, notice that the labels on
curves in the graph correspond to the index of the bisectors given in
Table~\ref{tab:numbering}.
\begin{table}[htbp]
  \centering
  \begin{tabular}{l|c|c|c|c}
    Bisector & $P^k\Rc$ & $P^k\Rc^-$ & $P^k\Sc$ & $P^k\Sc^-$ \\ \hline
    Index notation & $\Bn_{1+4k}$ & $\Bn_{2+4k}$ & $\Bn_{3+4k}$ & $\Bn_{4+4k}$ 
  \end{tabular}
  \caption{Numbering of the bounding bisectors according to $P$-symmetry,
    for $k=0,\dots,6$.}\label{tab:numbering}
\end{table}

We treat these two tangencies in
Proposition~\ref{prop:onlyonepoint}. The other ones are similar.
\begin{proposition}\label{prop:onlyonepoint}
  \begin{enumerate}
  \item For $p=4$, the real geodesic through $p_{23}$ and $q_{1\bar
    323}$ intersects $P^{-1}\Rc$ only at the ideal vertex
    $p_{23}$.
  \item For $p=6$, the real geodesic through $p^{23\bar2}_{23}$ and
    $q_{123\bar2}$ intersects $P^2\Rc^-$ only at the ideal vertex
    $q_{123\bar2}$.
 \item For $p=6$, the real geodesic through $p^{\bar 323}_{23}$ and
    $q_{1\bar 323}$ intersects $P^{-2}\Sc$ only at the ideal vertex
    $q_{1\bar 323}$.
 \item For $p=6$, the real geodesic through $q_{123\bar2}$ and
    $q_{1\bar323}$ intersects $P^2\Rc$ only at the ideal vertex
    $q_{123\bar2}$.
  \end{enumerate}
\end{proposition}
\Pf 
First let $p=4$. The geodesic arc from $p_{23}$ to $q_{1\bar{3}23}$
is a side of $r_1\cap s_1$. We check that it only intersects $P^{-1}\Rc$
at $p_{23}$. Apart from $q_{1\bar 323}$, the extended real
geodesic is parametrized by vectors of the form ${\bf a}+t{\bf b}$
with
$$
  {\bf a} = \n_{23},\ {\bf b} = \langle \n_{23},\n_{1\bar 323}\rangle \n_{1\bar 323},
$$ 
and $t\in \R$ (there are further restrictions on $t$, but this will
be irrelevant in what follows). The corresponding point is on
$P^{-1}\Rc$ (or rather its extension to projective space) if and only if
$$ 
  | \langle {\bf a}+t{\bf b}, P^{-1} {\bf y}_0\rangle |^2 = | \langle
  {\bf a}+t{\bf b}, P^{-1}R_1^{-1}{\bf y}_0\rangle |^2.
$$ 
Using $\tau=-(1+i\sqrt{7})/2$, then $a=e^{2\pi i/12}=(\sqrt{3}+i)/2$, we compute
$$ 
  | \langle {\bf a}+t{\bf b}, P^{-1} {\bf y}_0\rangle |^2 - | \langle
     {\bf a}+t{\bf b}, P^{-1}R_1^{-1}{\bf y}_0\rangle |^2 = - (196156+74140\sqrt{7})\ t^2,
$$
the latter vanishing precisely at $t=0$.

For $p=6$ and ${\bf a}=q_{123\bar2}$, ${\bf b}=p_{23}^{23\bar2}$, one gets
$$ 
  | \langle {\bf a}+t{\bf b}, P^{2} {\bf y}_0\rangle |^2 - | \langle
     {\bf a}+t{\bf b}, P^{2}R_1{\bf y}_0\rangle |^2 = - (2525 + 551\sqrt{7}\sqrt{3})\  t^2.
$$
For $p=6$ and ${\bf a}=q_{123\bar2}$, ${\bf b}=q_{1\bar 323}$, one gets
$$ 
  | \langle {\bf a}+t{\bf b}, P^{2} {\bf y}_0\rangle |^2 - | \langle
  {\bf a}+t{\bf b}, P^{2}R_1^{-1}{\bf y}_0\rangle |^2 = -
  (161797+35307\sqrt{7}\sqrt{3}) t^2/2.
$$
For $p=6$ and ${\bf a}=q_{1\bar323}$, ${\bf b}=p_{23}^{\bar323}$, one gets
$$ 
  | \langle {\bf a}+t{\bf b}, P^{-2} {\bf y}_0\rangle |^2 - | \langle
     {\bf a}+t{\bf b}, P^{-2}S_1^{-1}{\bf y}_0\rangle |^2 = - (3524+769\sqrt{7}\sqrt{3})\  t^2.
$$
\EPf

\subsubsection{Realizing ridges of $\widehat{E}$ - consistency} \label{sec:ridgeconsistency}

The goal of this section is to show the consistency result for the
2-skeleton of $\widehat{E}$. Recall that we want to realize each
2-face $f$ by using the realization of $\partial f$, which is a
piecewise smooth circle (possibly with some ideal vertices) in the
intersection of two of the bounding bisectors.  We now need to show
that the intersection of these two bisectors is a disk. We split the
discussion in two cases, depending on the geometry of the intersection
(see Section~\ref{sec:bisectors}).

\noindent\textbf{Complex ridges}

We first list the ridges of $\widehat{E}$ that get realized in complex
lines, calling them simply ``complex ridges''. The detailed list is
given in the appendix (see Table~\ref{tab:mirrors}).

For $p=3,4$ there are 14 complex lines containing complex ridges,
namely the $P$-orbits of $m_1$ and $m_2$. For $p\geqslant 5$ the point
$p_{12}$ is replaced with a ridge in the complex line $m_{12}$ and so
there are 21 complex ridges. The $P$-orbits of $m_{12}$ may be found
from Table~\ref{tab:verts3-4} by replacing $p_w$ with $m_w$ (these are
closely related, in fact they are polar to each other). For $p=8,12$
we similarly replace $q_{123\bar{2}}$ with a ridge in
$m_{123\bar{2}}$, and its $P$-orbit may be read of from the same table
by replacing $q_w$ with $m_w$ in Table~\ref{tab:verts3-6}; this gives
a total of 28 complex ridges.

We only analyze the ridges that appear on $r_1$ and $s_1$, see the
bottom and top 2-cells of Figures~\ref{fig:adjacency3-4}
to~\ref{fig:adjacency8-12} (the top 2-cell reduces to a point for
small values of $p$). We have already shown that, for each such 2-face
$f$, the vertices and edges in $\phi(\partial f)$ are contained in a
complex line, which is the mirror of a certain complex reflection in
the group (see Section~\ref{sec:edges}). Moreover, the adjacency
pictures suggest that we should realize the 2-face in the intersection
of two specific bisectors.

The consistency verification amounts to the verification that these
two specific bisectors intersect precisely along the complex line
mentioned above. We state this in the following.

\begin{proposition}\label{prop:cridges}
  \begin{enumerate}
  \item $\Rc\cap \Rc^-$ is a complex line, which
    is $m_1$.
  \item For $p\geqslant 5$, $\Rc\cap P^{-3}\Rc^-$ is a
    complex line, which is $m_{23}$.
  \item $\Sc\cap P^2\Sc^-$ is a complex line,
    which is $m_{23\bar 2}$.
  \item For $p\geqslant 8$, $\Sc\cap P^{-2}\Sc^-$ is a
    complex line, which is $m_{1\bar{3}23}$. 
  \end{enumerate}
\end{proposition}
\Pf 
Parts~1 and~3 are relatively easy, since the corresponding pairs
of bisectors are not only cotranchal, they are in fact cospinal.
This is obvious for $\Rc\cap\Rc^-$, since the complex spine of $\Rc$
is by construction orthogonal to the mirror of $R_1$, so $\Rc$ and
$\Rc^-=R_1(\Rc)$ have the same complex spine.

The case of $\Sc\cap P^2\Sc^-$ is just a little bit more complicated. By
construction $P^2S_1$ maps $\Sc$ to $P^2\Sc^-$, and
Proposition~\ref{prop:betterS1} implies that $P^2S_1$ is a complex
reflection fixing $\n_{23\bar{2}}$, hence it preserves the
complex spine of $\Sc$.

For parts~2 and~4, we first prove that the corresponding pair of
bisectors is cotranchal, i.e. that the mirror of the relevant complex
reflection is indeed contained in the bisector intersection. Then we
prove that the intersection really consists only of that common complex
slice, by using Proposition~\ref{prop:cotranchal}.

\smallskip
\noindent {\bf Part~2:} 
In this case we need to verify the criterion from
Proposition~\ref{prop:cotranchal}.  That is, we must find the
intersection of the real spine $\sigma(\Rc)$ of $\Rc$ with $m_{23}$
and lift it to a point ${\bf v}_1$ so that $\sigma(\Rc)$ is given by
real linear combinations of ${\bf n}_{23}$ and ${\bf v}_1$. We know
that $\sigma(\Rc)$ is the real span of ${\bf n}_1$ and ${\bf n}_{23}$
and that $\langle{\bf n}_1,{\bf n}_{23}\rangle$ is real. Hence we take
$$
  {\bf v}_1={\bf n}_1
  -\frac{\langle{\bf n}_1,{\bf n}_{23}\rangle}
  {\langle{\bf n}_{23},{\bf n}_{23}\rangle}{\bf n}_{23}
  =\frac{1}{a^3+\bar{a}^3}\left[\begin{matrix} 0 \\ -a^2\bar{\tau}+\bar{a} \\ -\bar{a}^2\tau+a
    \end{matrix}\right].
$$
Now $P^{-3}\Rc^-=(R_2R_3)^{-1}\Rc$ and 
$\bar{a}^2(R_2R_3)^{-1}$ fixes ${\bf n}_{23}$ and so we take
$$
  {\bf v}_2=\bar{a}^2(R_2R_3)^{-1}{\bf v}_1
  =\frac{1}{a^3+\bar{a}^3}\left[\begin{matrix} 0 \\ \bar{a}^4\bar{\tau}+\bar{a} \\
      -\bar{a}^2\tau-\bar{a}^5\end{matrix}\right].
$$
Then
\begin{eqnarray*}
  \langle{\bf v}_1,{\bf v}_1\rangle \ = \ \langle{\bf v}_2,{\bf v}_2\rangle
  & = & 
  (a^3-1)\bigl((1-\tau)-\bar{a}^3(1-\bar{\tau})\bigr)/(a^3+\bar{a}^3), \\
  \langle{\bf v}_2,{\bf v}_1\rangle & = & 
  -(a^3-1)\bigl((1-\bar{\tau})\bar{a}^3+(1-\tau)\bar{a}^6\bigr)/(a^3+\bar{a}^3).
\end{eqnarray*}
Thus
\begin{eqnarray*}
  \lefteqn{ 
    \frac{(a^3+\bar{a})^2}{(2-a^3-\bar{a}^3)} 
    \Bigl(4\langle{\bf v}_1,{\bf v}_1\rangle\langle{\bf v}_2,{\bf v}_2\rangle
    -\bigl(\langle{\bf v}_1,{\bf v}_2\rangle+\langle{\bf v}_2,{\bf v}_1\rangle\bigr)^2\Bigr)}\\
  & = & 
  (a^3+\bar{a}^3)\bigl(2|1-\tau|^2-a^3(1-\bar{\tau})^2-\bar{a}^3(1-\tau)^2\bigr)
  +2(a^3-\bar{a}^3)\bigl(a^3(1-\bar{\tau})^2-\bar{a}^3(1-\tau)^2\bigr).
\end{eqnarray*}
Both terms are positive for all $p>4$.

\smallskip
\noindent {\bf Part~4:} 
We argue as in Part~2. The real spine $\sigma(\Sc)$ of $\Sc$ is the real 
span of $a{\bf n}_{23\bar{2}}$ and ${\bf n}_{1\bar{3}23}$. Its image under 
$a^2 1\bar{3}23$ is the real span of $-a{\bf n}_{13\bar{1}}$ and 
${\bf n}_{1\bar{3}23}$. Therefore we construct:
\begin{eqnarray*}
  {\bf v}_1 & = & a{\bf n}_{23\bar{2}}
  -\frac{\langle a{\bf n}_{23\bar{2}},{\bf n}_{1\bar{3}23}\rangle}
  {\langle{\bf n}_{1\bar{3}23},{\bf n}_{1\bar{3}23}\rangle}{\bf n}_{1\bar{3}23}
  =\frac{1}{a^3+\bar{a}^3-1}\left[\begin{matrix}
      -a^2-\bar{a}\bar{\tau} \\ a\bar{\tau}+\bar{a}^2\tau \\ -1+\tau+2\bar{a}^3
    \end{matrix}\right], \\
  {\bf v}_2 & = & -a{\bf n}_{13\bar{1}}
  -\frac{\langle -a{\bf n}_{13\bar{1}},{\bf n}_{1\bar{3}23}\rangle}
  {\langle{\bf n}_{1\bar{3}23},{\bf n}_{1\bar{3}23}\rangle}{\bf n}_{1\bar{3}23}
  =\frac{1}{a^3+\bar{a}^3-1}\left[\begin{matrix}
      a^5\bar{\tau}+a^2\bar{\tau} \\ -a^4\tau-a \\ -2a^3-\bar{\tau}
    \end{matrix}\right]. \\
\end{eqnarray*}
We have
\begin{eqnarray*}
  \langle{\bf v}_1,{\bf v}_1\rangle \ =\ \langle{\bf v}_2,{\bf v}_2\rangle 
  & = &
  (a^3-1)\bigl((2-\tau)-\bar{a}^3(2-\bar{\tau})\bigr)/(a^3+\bar{a}^3-1), \\
  \langle{\bf v}_2,{\bf v}_1\rangle & = &
  (a^3-1)(2-\bar{\tau})a^3/(a^3+\bar{a}^3-1).
\end{eqnarray*}
Then
\begin{eqnarray*}
  \lefteqn{
    \frac{(a^3+\bar{a}^3-1)^2}{(2-a^3-\bar{a}^3)}
    \Bigl(4\langle{\bf v}_1,{\bf v}_1\rangle\langle{\bf v}_2,{\bf v}_2\rangle
    -\bigl(\langle{\bf v}_1,{\bf v}_2\rangle+\langle{\bf v}_2,{\bf v}_1\rangle\bigr)^2\Bigr)} \\
  & = & 6|2-\tau|^2-(4a^3+\bar{a}^6)(2-\tau)^2-(4\bar{a}^3+a^6)(2-\bar{\tau})^2.
\end{eqnarray*}
This is positive for all $p>6$.

This completes the proof of Proposition~\ref{prop:cridges}. \EPf

\noindent\textbf{Generic ridges} 

In the first part of this section, we checked the consistency of the
embedding of complex ridges. We now prove the analogous statement for
all the other ridges, which will be realized inside Giraud disks.
Recall that given a 2-cell $f$ of $\widehat{E}$,
Figures~\ref{fig:adjacency3-4} to~\ref{fig:adjacency8-12} suggest a pair
$\B,\B'$ of bounding bisectors that should contain the realization
$\phi(f)$. We have already defined the realization on
the level of the $1$-skeleton, and shown that $\phi(\partial f)$ lies
inside $\B\cap \B'$ (or possibly its union with finitely many ideal
vertices).

The corresponding consistency verification is almost automatic, we
simply verify that all the corresponding intersections $\B\cap \B'$
are indeed disks, hence $\phi(\partial f)$ bounds a well defined disk,
which we use as the realization $\phi(f)$.
In order to get these disks, we simply prove that the corresponding
pairs of bisectors are coequidistant (see Section~\ref{sec:bisectors}).

Let ${\bf y}_0$, ${\bf y}_1$, ${\bf y}_2$ be given by the following
formulae:
  \begin{equation}\label{eq:xyz}
  {\bf y}_0 = \left[\begin{matrix} 
     -a^3(\bar{a}^2\bar{\tau}+a)^2 \\
     (a^2\bar{\tau}+\bar{a}\tau)
     (a^2\bar{\tau}-\bar{a}) \\
     (a^2\bar{\tau}+\bar{a}\tau)
     (\bar{a}^2\tau-a) 
   \end{matrix}\right],
  \ 
    {\bf y}_1 = \left[\begin{matrix} 
     a(1-a^3\bar{\tau})^2 \\
     (1-a^3\bar{\tau})(a^3-\tau)\\
     \bar{a}(a^3-\tau)^2
   \end{matrix}\right]
  \ 
    {\bf y}_2 = \left[\begin{matrix} 
     (a^3-\tau)(\bar{a}^3+\tau) \\
     \bar{a}(1-a^3\bar{\tau})(1+\bar{a}^3\bar{\tau})\\
     a(1-\bar{a}^3\tau)(1+\bar{a}^3\bar{\tau})
   \end{matrix}\right].
  \end{equation}
Note that $\langle {\bf y}_j,{\bf y}_j\rangle<0$ and so ${\bf y}_j$
corresponds to a point $y_j$ of $\ch 2$. Indeed, the following lemma,
which is verified by direct computation, implies that all three points are in $E$.

\begin{lemma}\label{lem:xyz}
  Let $y_0$, $y_1$, $y_2$ be the points in $\ch 2$ corresponding to
  the vectors given in \eqref{eq:xyz}. If $\B$ is any of the 28
  bounding bisectors, then the half-space determined by $\B$ that
  contains the fixed point of $P$ also contains $y_0$, $y_1$ and
  $y_2$.
\end{lemma}

We have already seen in Proposition~\ref{prop:coequid1} that $\Rc$,
$\Rc^-$, $\Sc$ and $\Sc^-$ are coequidistant from ${\bf y}_0$.  Since
$y_0$ does not lie on any of their real spines, their pairwise
intersections are Giraud disks (apart from $\Rc\cap \Rc^-$, since
these two bisectors have the same complex spine).  Using the points
$y_1$ and $y_2$ we can describe some of the other intersections
between bounding bisectors in a similar way.

\begin{table}
$$
  \begin{array}{|rl|rl|}
    \hline
    \hbox{Bisector} & \hbox{Equidistant} & \hbox{Bisector} & \hbox{Equidistant} \\
    \hline
    \Rc & \B(y_0,R_1^{-1}y_0) & \Rc^- & \B(y_0,R_1y_0) \\
    \Sc & \B(y_0,S_1^{-1}y_0) & \Sc^- & \B(y_0,S_1y_0) \\
    \hline
    \Rc & \B(y_1,R_1^{-1}Py_1) & P^{-1}\Rc^- & \B(y_1,P^{-1}R_1y_1) \\
    \hline
    \Rc & \B(y_2,R_1^{-1}P^2y_2) & P^{-2}\Rc^- & \B(y_2,P^{-2}R_1y_2) \\
    P^2\Sc & \B(y_2,P^2S_1^{-1}P^{-3}y_2) & P^3\Sc^- & \B(y_2,P^3S_1P^{-2}y_2) \\
    \hline
  \end{array}
$$
\caption{Description of bounding bisectors in terms of $y_0$, $y_1$ and $y_2$}
\label{tab:bis-xyz}
\end{table}

\begin{table}
$$
  \begin{array}{|r l|}
    \hline
    \hbox{Intersection} & \hbox{Giraud disk} \\
    \hline
    \Rc\cap\Sc & \G(y_0,R_1^{-1}y_0,S_1^{-1}y_0) \\
    \Rc^-\cap\Sc^- &  \G(y_0,R_1y_0,S_1y_0) \\
    \Sc\cap\Sc^- &  \G(y_0,S_1^{-1}y_0,S_1y_0) \\
    \hline
    \Rc\cap P^{-1}\Rc^-
    & \G(y_1,R_1^{-1}Py_1,P^{-1}R_1y_1) \\
    \hline
    \Rc\cap P^2\Sc &  
    \G(y_2,R_1^{-1}P^2y_2,P^2S_1^{-1}P^{-3}y_2) \\
    \Rc\cap P^{-2}\Rc^- & 
    \G(y_2,R_1^{-1}P^2y_2,P^{-2}R_1y_2) \\
    P^{-2}\Rc^-\cap P^3\Sc^- & 
    \G(y_2,P^{-2}R_1y_2,P^3S_1P^{-2}y_2) \\
    \hline
\end{array}
$$
\caption{Giraud intersections among bounding bisectors.}
\label{tab:giraud1}
\end{table}

\begin{proposition}\label{prop:giraud-disk}
  Let $y_0$, $y_1$ and $y_2$ be the points in $\ch 2$ corresponding to
  the vectors ${\bf y}_0$, ${\bf y}_1$, ${\bf y}_2$ given in
  \eqref{eq:xyz}. The intersections among bounding bisectors listed in
  Table~\ref{tab:giraud1} are contained in Giraud disks defined by
  $y_0$, $y_1$ and $y_2$.
\end{proposition}

\Pf 
We must first prove that certain bounding bisectors are
coequidistant from $y_0$, $y_1$ and $y_2$ listed in
Table~\ref{tab:bis-xyz}. This was proved for $y_0$ in
Proposition~\ref{prop:coequid1}.  We use a similar method to show
equidistance from $y_1$ and $y_2$.  Namely, consider two of the
bisectors $\B$ and $\B'$.  We first show that the complex spines of
$\B$ and $\B'$ intersect in ${\bf y}_j$, for $j=1$ or $2$. We can find
vectors spanning the real spine of $\B$ by applying a suitable power
of $P$ to the vectors listed in Table~\ref{tab:spines}. Using
$\iota$,we find an antiholomorphic isometry fixing the spines of
$\B$. From this we can then find an image of $y_j$ so that $\B$ is
coequidistant from $y_j$ and this second point.

To be specific, the real spine of $\Rc$ is invariant under the
involution $\iota_{23}$ defined in \eqref{eq:iota23}. We see that
$\iota_{23}\,{\bf y}_1=\bar{a}^6 J\,{\bf y}_1$ and
$\iota_{23}\,{\bf y}_2=a^2R_1^{-1}P^2\,{\bf y}_2$. Therefore,
using $J=R_1^{-1}P$ we see that $\Rc=\B(y_1,R_1^{-1}P\,y_1)
=\B(y_2,R_1^{-1}P^2\,y_2)$.  Similarly, the spine of $P^2\Sc$ is the
real span of ${\bf n}_2=-P^2\,{\bf n}_{23\bar{2}}$ and ${\bf
  n}_{123\bar{2}}=P^2\bar{a}\,{\bf n}_{1\bar{3}23}$. These two vectors
are fixed by $\bar{a}^2R_2\,\iota_{13}$, where $\iota_{13}$ is defined
analogously to $\iota_{23}$. We find that
$\bar{a}^2R_2\iota_{13}\,{\bf y}_2=a^2P^2S_1^{-1}P^{-3}{\bf y}_2$.
Therefore $P^2\Sc=\B(y_2,P^2S_1^{-1}P^{-3}y_2)$.

The claims in Table~\ref{tab:giraud1} follow at once from those in
Table~\ref{tab:bis-xyz}.
\EPf

\subsubsection{Realizing ridges of $\widehat{E}$ - embedding} 
\label{sec:ridgeembedding}

From the results in Section~\ref{sec:ridgeconsistency}, we get that
$\phi$ gives a well defined realization of each 2-face $f$ of
$\widehat{E}$, inside the intersection of two specific bounding
bisectors.  

For every other bounding bisector $\B$, we now prove that $\phi(f)$ is
on the correct side of $\B$. In particular, we will show that
$\phi(f^\circ)$ does not intersect $\B$, and $\B$ intersects $\phi(f)$
precisely as predicted from the labels in
Figures~\ref{fig:adjacency3-4} to~\ref{fig:adjacency8-12}. From this,
it follows that $\phi$ is an embedding on the level of
$\widehat{E}_2$.

First, observe that when $\phi(f)$ is in a complex line, no further
verifications are needed, due to the following.
\begin{lemma}\label{lem:openmapping}
  Let $Q$ be a polygon contained in a complex line $C$, and let $\B$
  be a bisector. If $\partial Q$ is entirely on one side of $\B$, then
  so is $Q$.
\end{lemma}
\Pf 
The restriction to $C$ of the orthogonal projection onto the
complex spine of $\B$ is an open map (because it is holomorphic).
Recall that a set is entirely on one side of $\B$ if and only if its
orthogonal projection onto the complex spine is in the corresponding
half disk bounded by the real spine.\EPf

Orthogonal projections restricted to Giraud disks are, in general, not
open.  Hence the cases where $\phi(f)$ lies in a Giraud disk are
considerably more complicated. They will occupy us for the remainder
of this section.

Let us start with a Giraud polygon $T$ in a Giraud disk $\G$ (i.e. a
disk bounded by a piecewise geodesic simple closed curve) and a
bisector $\B$ not containing $\G$. We describe a general method to
prove that $T$ is entirely on one side of $\B$.

Let $\B=\B(x_0,x_1)$, and assume that ${\bf x}_0$, ${\bf x}_1$ are
lifts of $x_0$, $x_1$ with $\langle {\bf x}_0, {\bf x}_0 \rangle =
\langle {\bf x}_1,{\bf x}_1 \rangle$.  Using the parametrization given
in Section~\ref{sec:coequidistant}~\eqref{paramgiraud}, we may
write an equation for the trace of $\B$ on $\G$ as
\begin{equation}\label{eq:thirdbisector}
  g(t_1,t_2)=0, \textrm{ where }
  g(t_1,t_2)=|\langle V(t_1,t_2),{\bf x}_0 \rangle |^2 - |\langle
  V(t_1,t_2),{\bf x}_1 \rangle |^2.
\end{equation}
Since $V(t_1,t_2)$ is affine in each variable,
equation~(\ref{eq:thirdbisector}) is quadratic in each variable.

If we claim that $\B\cap T$ is empty, we need to prove that $g$ does
not vanish on $T$, for instance that $g>0$ (switch $x_0$ and $x_1$ if
needed). In order to do this, we
\begin{enumerate}
  \item check that $g>0$ on the boundary of $T$ (see
    Section~\ref{sec:numberfields}).
  \item find the critical points of $g$. If they are in $T$, check the
    sign of the corresponding values of $g$.
\end{enumerate}

This strategy works well in most cases, because of the following.
\begin{proposition}\label{prop:nondegenerate}
  The coefficients of $g$ lie in the field $\Q(a,\tau)$. The fact that $g$ only has
    non-degenerate critical points can be proved using computations in
    that number field, and the critical points can then be computed
    with arbitrary precision.
\end{proposition}

Since the key to these computations is to be able to compute critical
points of $g$, we expand a little on how this can be done in
practice. We write $g$ as
$$ 
  g(t_1,t_2) =
  a_0+a_1t_1+a_2t_2
   +a_{11}t_1^2+a_{12}t_1t_2+a_{22}t_2^2
   +a_{112}t_1^2t_2+a_{122}t_1t_2^2
   +a_{1122}t_1^2t_2^2,
$$ 
where all coefficients are real, and depend explicitly on the vectors
${\bf v}_j,{\bf w}_j$ (and are known in exact form, see
Proposition~\ref{prop:nondegenerate}).  The partial derivatives of
this function are given by
$$
  \left\lbrace\begin{array}{l}
  \frac{\partial g}{\partial t_1}(t_1,t_2) = 
             (a_1+a_{12}t_2+a_{122}t_2^2)+2t_1(a_{11}+a_{112}t_2+a_{1122}t_2^2)
             =P_2(t_2)+t_1Q_2(t_2)\\
  \frac{\partial g}{\partial t_2}(t_1,t_2) = 
             (a_2+a_{12}t_1+a_{112}t_1^2)+2t_2(a_{22}+a_{122}t_1+a_{1122}t_1^2)
             =P_1(t_1)+t_2Q_1(t_1)
  \end{array}\right.
$$ 
where the polynomials $P_j,Q_j$ have degree at most two.  Coordinates
of the critical points must be solutions of a polynomial of degree
five, obtained as the resultant of $\frac{\partial g}{\partial t_1}$
and $\frac{\partial g}{\partial t_2}$ (with respect to, say $t_1$).

Explicitly, if $(t_1,t_2)$ is a critical point, then $t_2$ is a root
of $\displaystyle\sum_{k=0}^5 b_k x^k$ where
\begin{eqnarray*}
  b_0 & = & a_1^2 a_{112} + 4 a_2 a_{11}^2 - 2 a_{12} a_{11}a_1\\
  b_1 & = & -4 a_{122} a_{11}a_1 + 2 a_1^2 a_{1122} + 8 a_{22} a_{11}^2 - 2 a_{12}^2 a_{11}+8 a_2a_{11}a_{112}\\
  b_2 & = & 8 a_2a_{11}a_{1122} + 2a_{1}a_{1122}a_{12} - 2a_{1}a_{112}a_{122} - 6a_{12}a_{11}a_{122} \\
&&\quad + 4a_{2}a_{112}^2 - a_{12}^2a_{112} + 16a_{22}a_{11}a_{112}\\
  b_3 & = & 8 a_{2} a_{112} a_{1122} + 8 a_{22} a_{112}^2 - 4 a_{12} a_{112} a_{122} - 4 a_{122}^2 a_{11} + 16 a_{22} a_{11} a_{1122}\\
  b_4 & = & -3 a_{122}^2 a_{112} - 2 a_{122} a_{1122} a_{12} + 4 a_{2} a_{1122}^2 + 16 a_{22} a_{112} a_{1122}\\
  b_5 & = & -2 a_{122}^2 a_{1122} + 8 a_{22} a_{1122}^2
\end{eqnarray*}
The real roots of this polynomial can be computed with arbitrary
precision, most easily when it only has simple roots (in that case the
Sturm algorithm suffices). In particular, we can count critical points
and locate them.

There are situations where the resultant has multiple roots, but there
is an easy way to work around these, which turns out to work for all
verifications needed in our paper. The point is that there is a simple
reason why the resultant can have multiple roots, namely that the
intersection $\G\cap \B$ contains a horizontal or a vertical segment.

Each time this happens, we can factor out a linear expression in one
of the two variables (or in some cases one factor for each variable),
and in fact we will be able to do this exactly (i.e. no numerical
approximation is involved).

More specifically, suppose the line $t_1=\alpha$ is contained in the
graph (the situation is obviously entirely similar for horizontal
lines). Then we write the function in the form
$$
g(t_1,t_2)=p_0(t_1)+p_1(t_1)t_2+p_2(t_1)t_2^2,
$$ 
where the $p_j$ are polynomials of degree at most two. Then the
coefficients $p_0$, $p_1$ and $p_2$ must have $\alpha$ as a common root,
and these roots are easily computed.

\begin{lemma}
  The vertical lines in the graph of $g(t_1,t_2)=0$ are given by the equations 
  $t_1=\alpha$ where $\alpha$ is a common root of the polynomials $p_j$,
  $j=0,1,2$.
\end{lemma}

It can easily be guessed using numerical computations whether or not
these three polynomials have a common root; when they do not, these
numerical computations constitute a proof (possibly after adjusting
the precision). When they seem to have a common root, we need to
verify this by exact calculations. We will work out the details of
these verifications for some explicit examples below.

Note that when $\{t_1=\alpha\}$ is indeed contained in the graph, it
is of course easy to extract a factor $(t_1-\alpha)$ from
$g(t_1,t_2)$. Explicitly, we write
\begin{equation}
  g(t_1,t_2) =q_0(t_2)+q_1(t_2)t_1+q_2(t_2)t_1^2 =
  (t_1-\alpha)h(t_1,t_2)\label{eq:t1powers}
\end{equation}
and equate the coefficients of degree $0$ and $2$ in $t_1$, to get
$$
h(t_1,t_2)=q_2(t_2)t_1+q_1(t_2)+\alpha q_2(t_2).
$$ 

\begin{rmk}
  When we do this factorization, it is convenient to adjust the sign
  of $h$ so that it has the same sign as $g$ on the open ridge we are
  studying. Depending on whether the side lies to the left or right of
  the line $t_1=\alpha$, we factor either $(t_1-\alpha)$ or
  $(\alpha-t_1)$.
\end{rmk}

The point of this remark is that in order to verify $g\geqslant 0$ on a
Giraud polygon $T$, it is enough to verify that $h\geqslant 0$ on $T$, and
we will do this by computing the critical points of $h$. Note that the
critical points of $h$ are of course not the same as those of $g$, but
they are easier to handle computationally (the relevant resultants
will have simple roots, see the above discussion).

\begin{figure}[htbp]
  \centering
  \subfigure[$p=3$]{\epsfig{figure=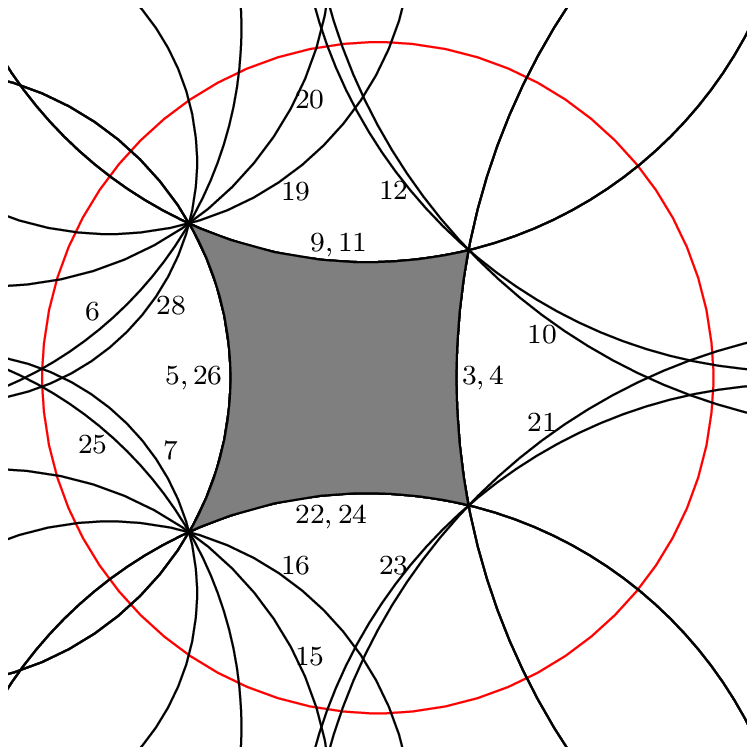, width=0.3\textwidth}}\hfill
  \subfigure[$p=4$]{\epsfig{figure=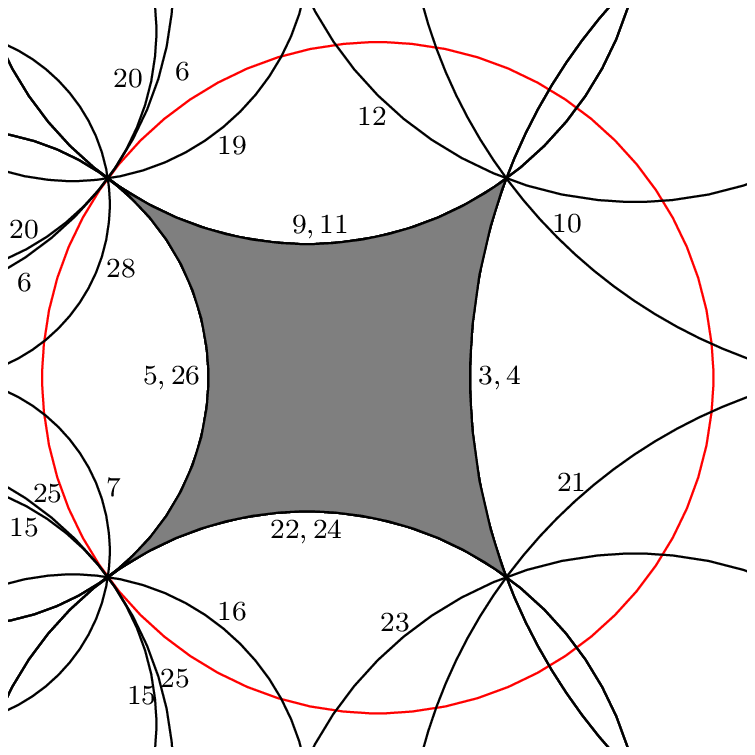, width=0.3\textwidth}}\hfill
  \subfigure[$p=5$]{\epsfig{figure=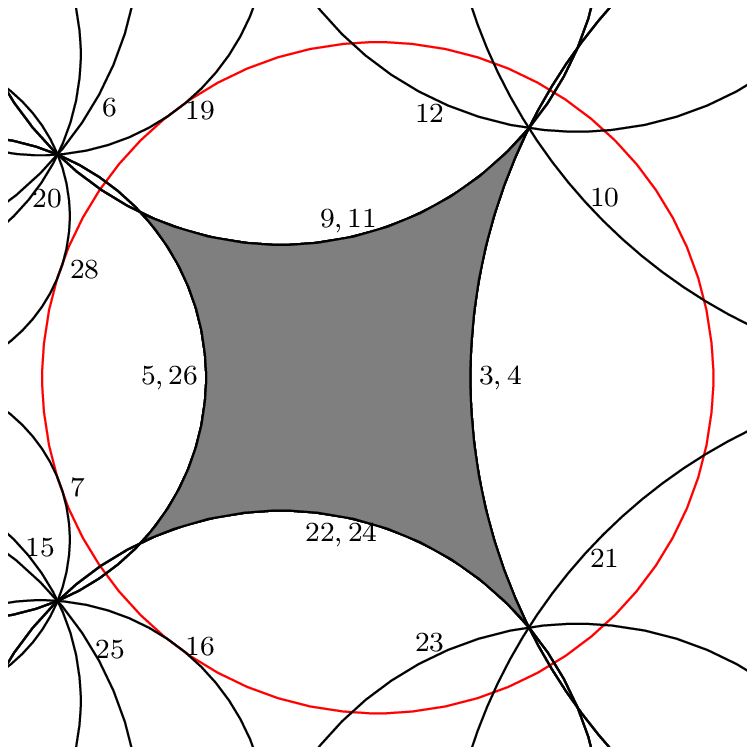, width=0.3\textwidth}}\\
  \subfigure[$p=6$]{\epsfig{figure=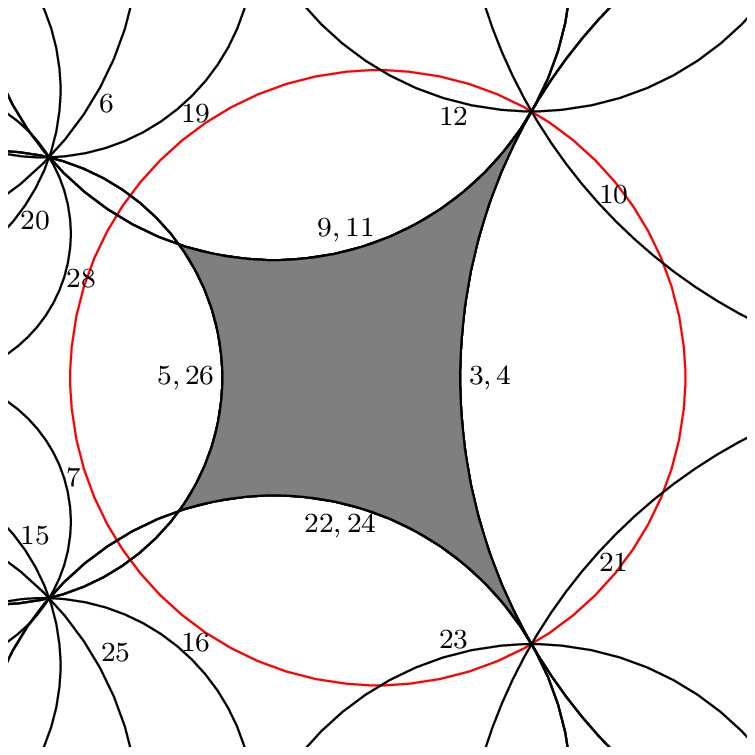, width=0.3\textwidth}}\hfill
  \subfigure[$p=8$]{\epsfig{figure=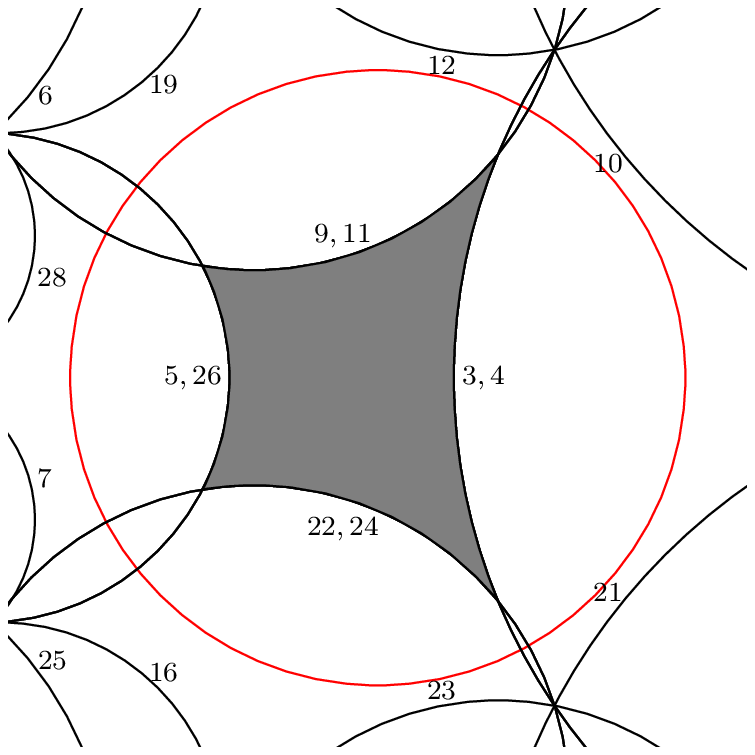, width=0.3\textwidth}}\hfill
  \subfigure[$p=12$]{\epsfig{figure=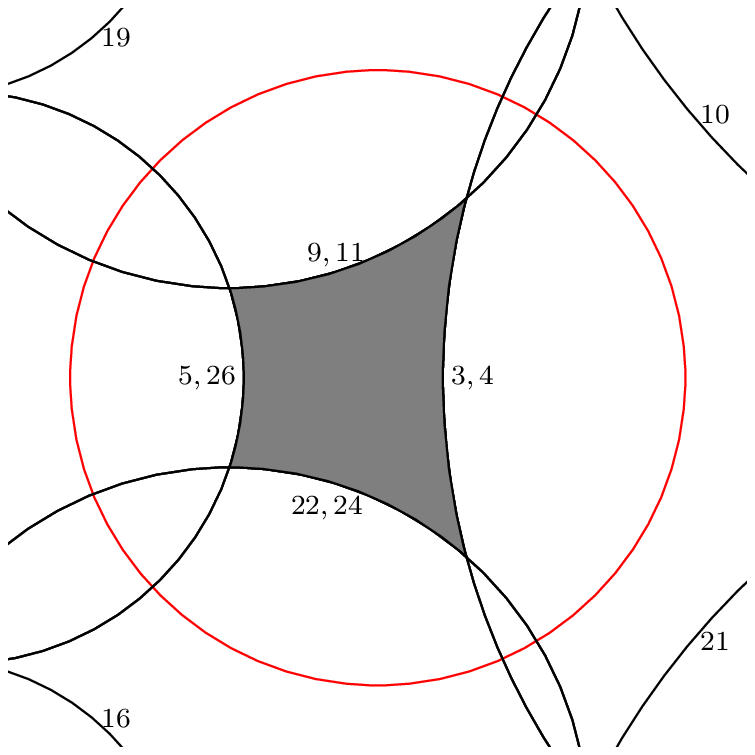, width=0.3\textwidth}}
  \caption{Intersection pattern of all the bounding sides with the complex line $\Rc\cap \Rc^-$, according to the values of $p$.}\label{fig:1-2}
\end{figure}

\begin{figure}[htbp]
  \centering
  \subfigure[$p=3$]{\epsfig{figure=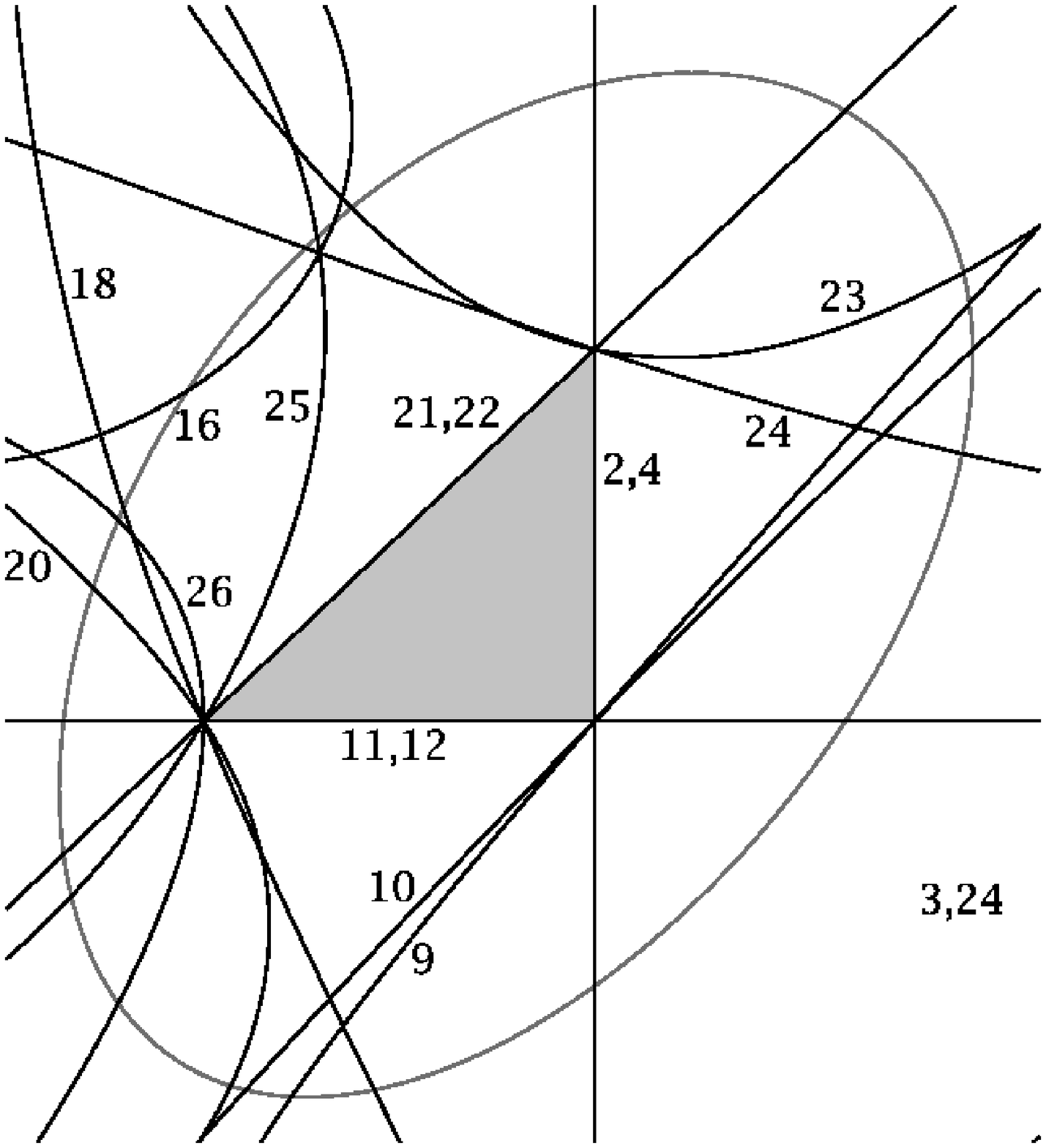, width=0.3\textwidth}}\hfill
  \subfigure[$p=4$]{\epsfig{figure=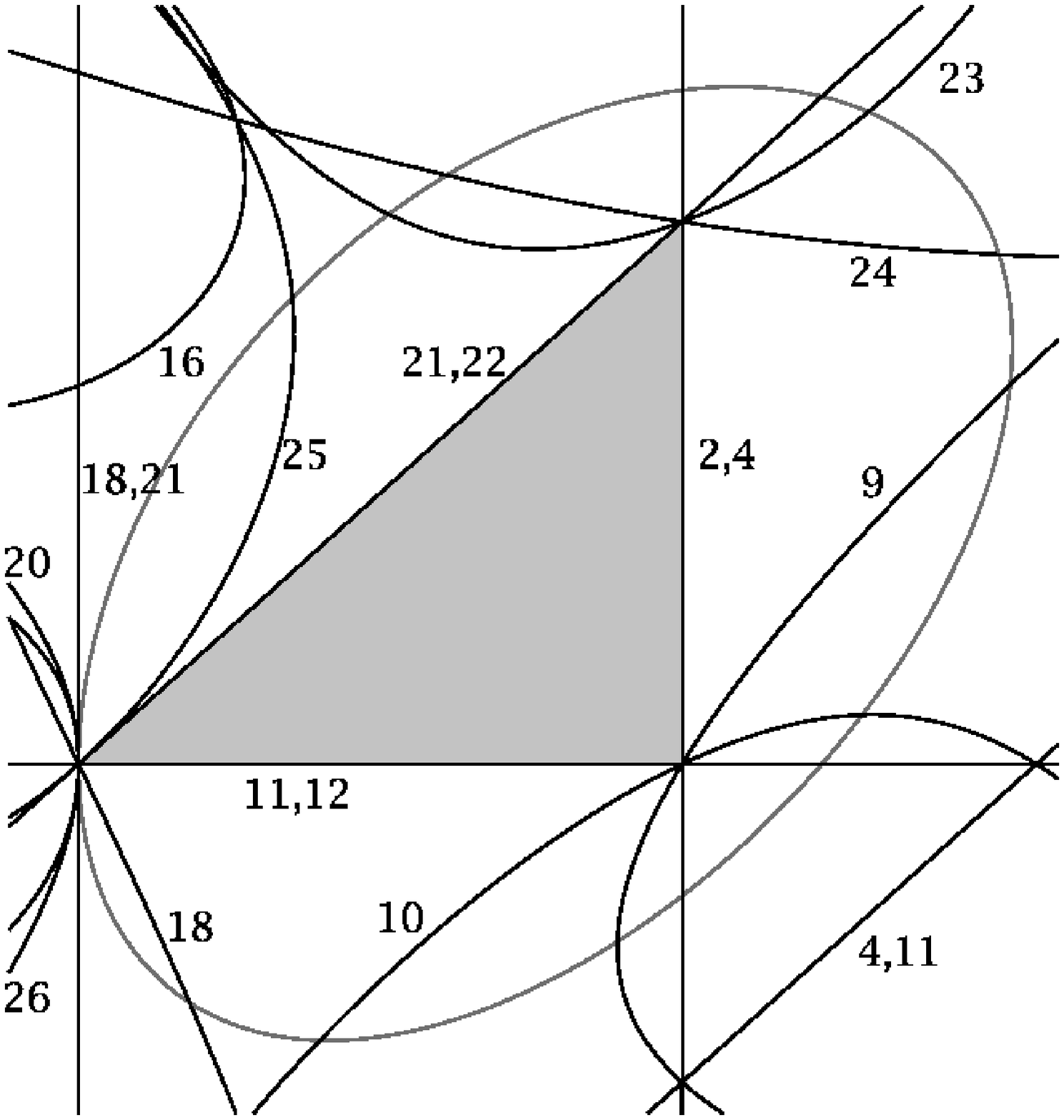, width=0.3\textwidth}}\hfill
  \subfigure[$p=5$]{\epsfig{figure=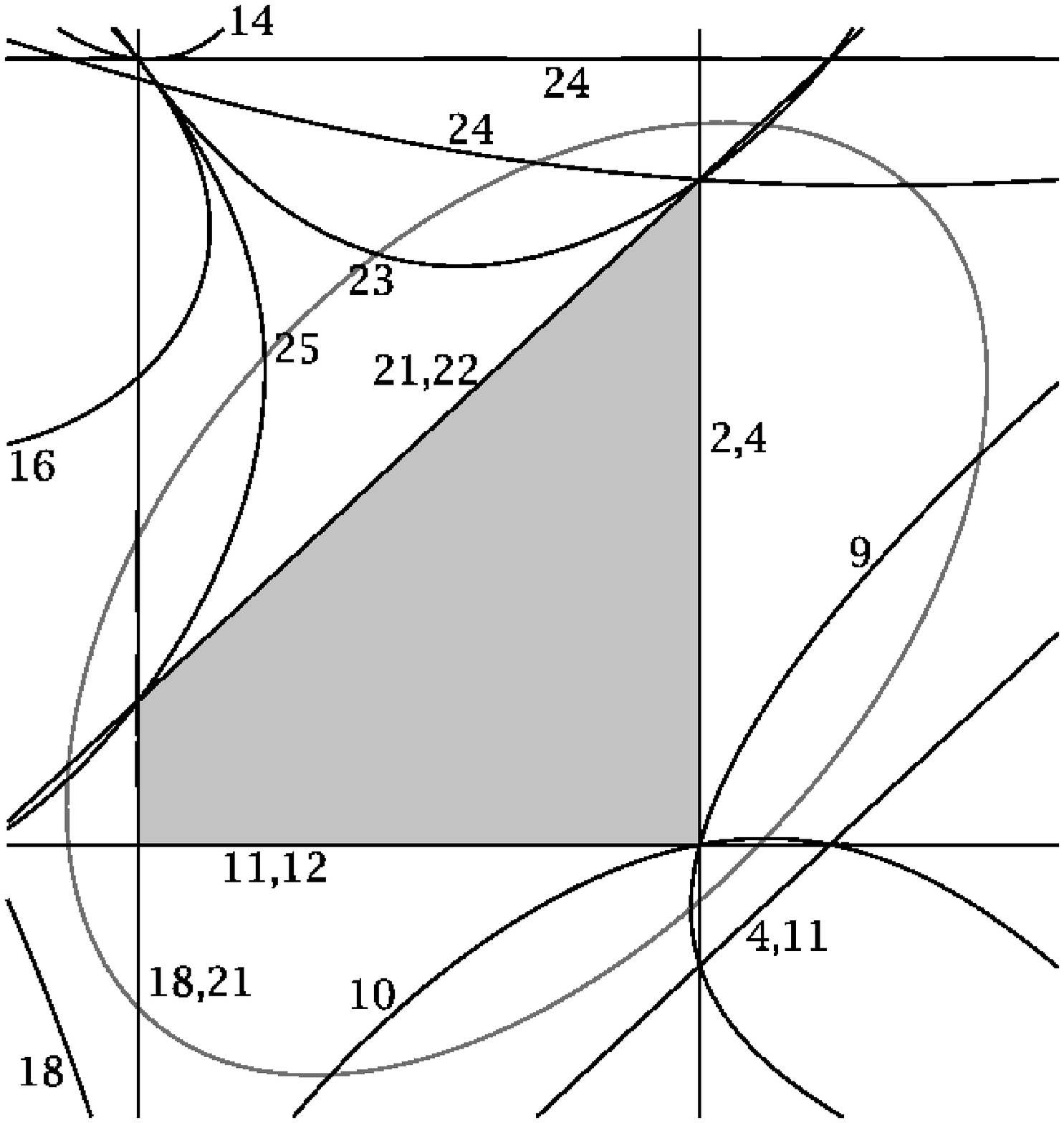, width=0.3\textwidth}}\\
  \subfigure[$p=6$]{\epsfig{figure=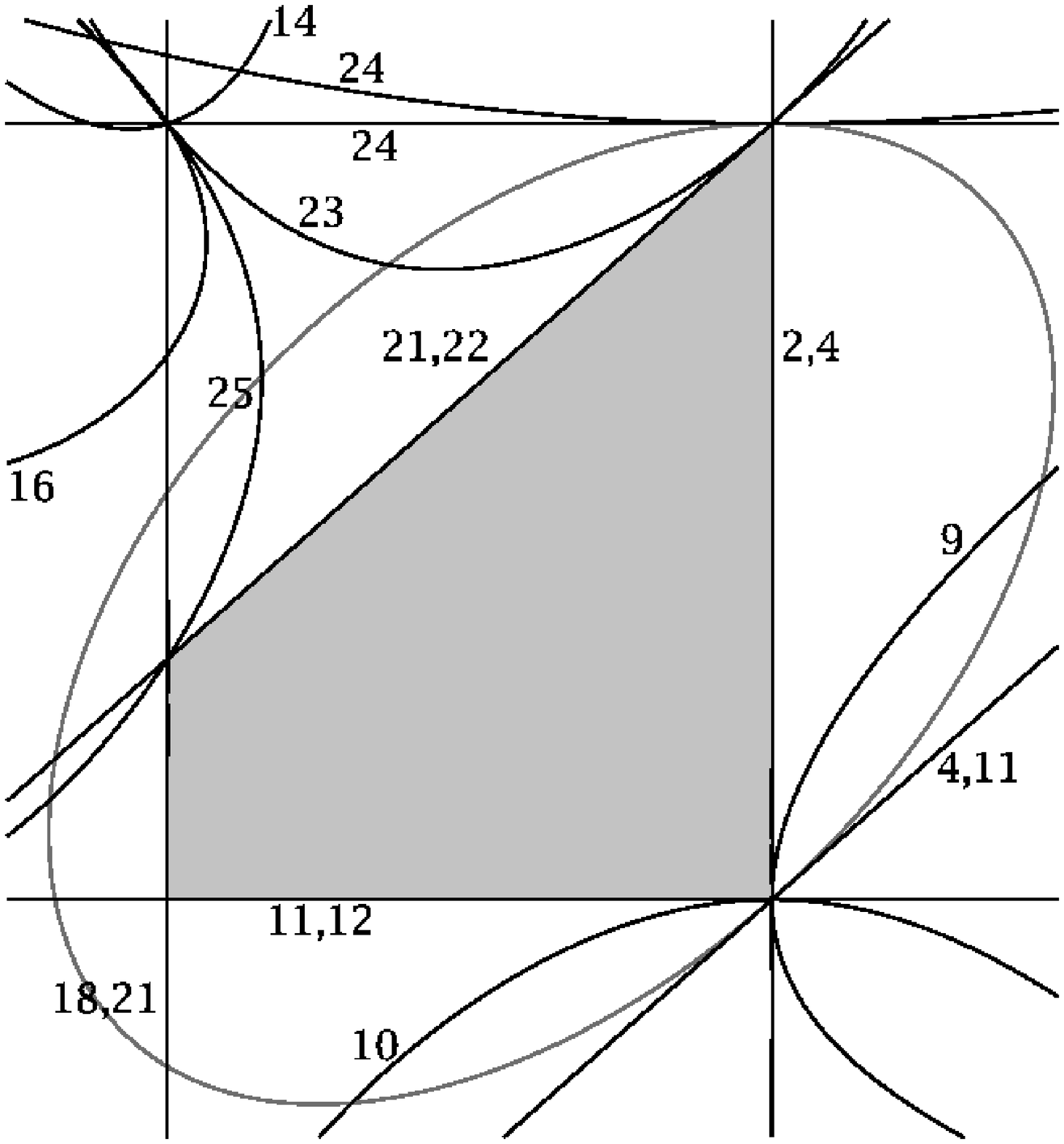, width=0.3\textwidth}}\hfill
  \subfigure[$p=8$]{\epsfig{figure=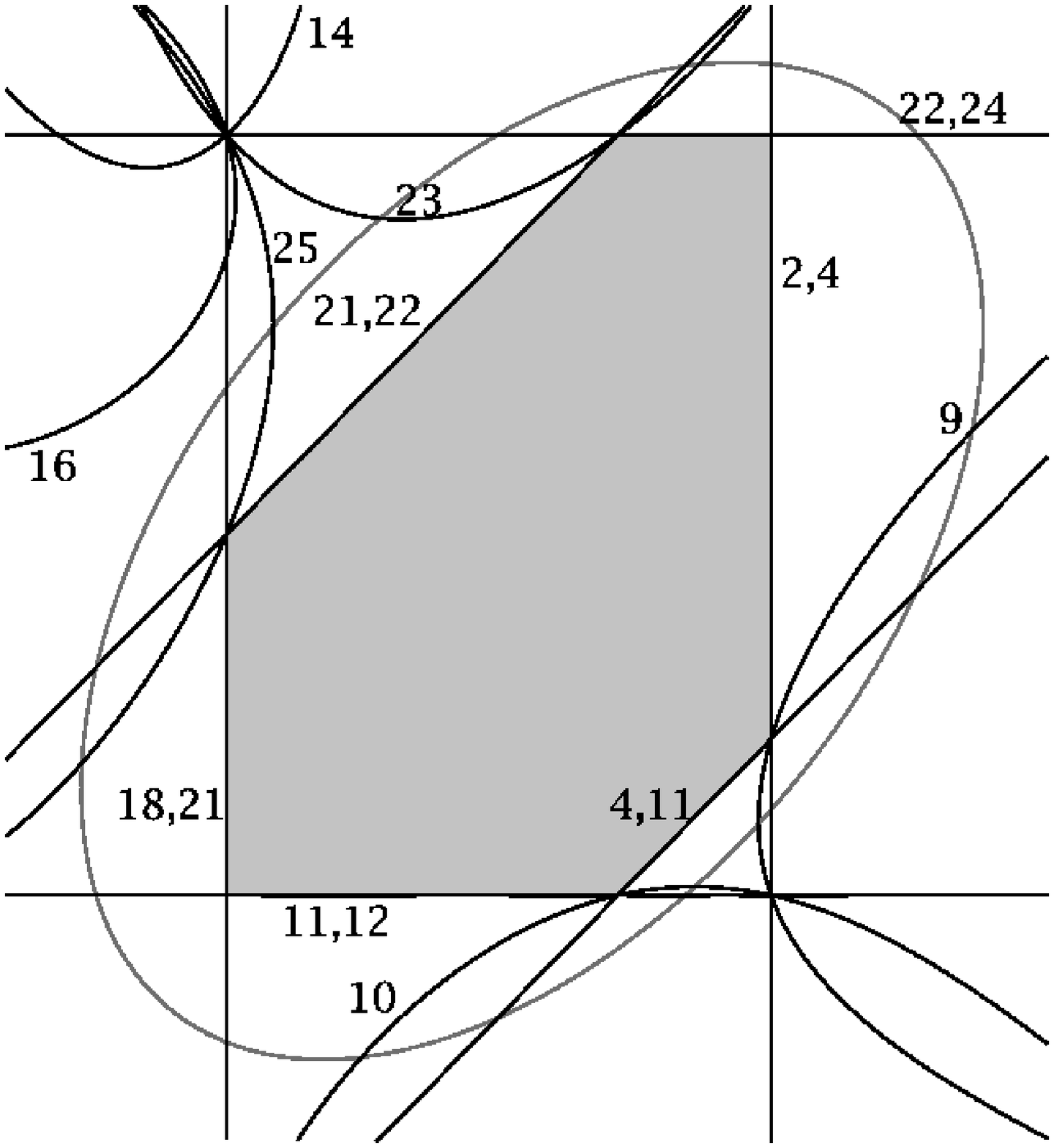, width=0.3\textwidth}}\hfill
  \subfigure[$p=12$]{\epsfig{figure=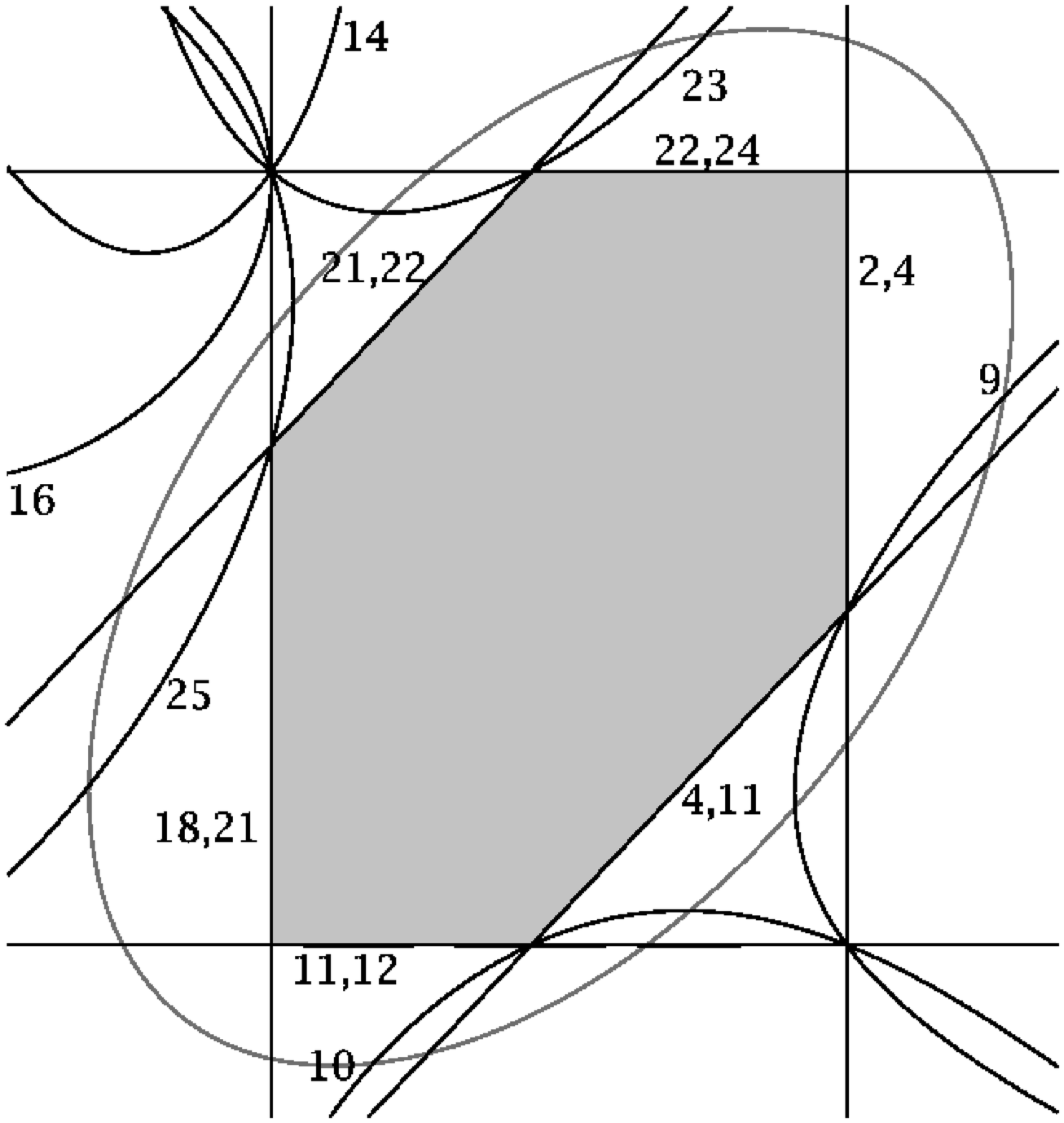, width=0.3\textwidth}}
  \caption{Intersection pattern of all the bounding sides with the Giraud
    disk $\Rc\cap \Sc$, according to the values of $p$.}\label{fig:1-3}
\end{figure}

We now work out the details of the above computations in some specific
cases. In order to allow for the reader to verify the numerical
figures given below, we specify explicit coordinates for the relevant
Giraud disks, giving specific choices of vectors ${\bf v}_j$ and ${\bf
  w}_j$ defining spinal coordinates as in~\eqref{paramgiraud}.

It is convenient to adjust the Giraud disk to be bounded in these
coordinates, hence we choose the ${\bf w}_j$ to be negative
vectors. For instance, for $\Rc$ one can choose ${\bf v} = \n_1$ and
$$
  {\bf w} = \langle \n_1,\n_1\rangle \n_{23} - \langle \n_{23},\n_1 \rangle \n_1,
$$ 
which corresponds to the intersection of the real spine of $\Rc$ with
the mirror of $R_1$.
For $\Sc$, one can choose ${\bf v} = \n_{23\bar 2}$ and
$$
  {\bf w} = \langle \n_{23\bar2},\n_{1\bar323}\rangle \{ \langle \n_{23\bar2},\n_{23\bar2}\rangle  \n_{1\bar323} 
                  - \langle \n_{1\bar323},\n_{23\bar2}\rangle \n_{23\bar2} \}.
$$
There is then a well-defined parametrization for every Giraud
intersection of bounding bisectors, obtained from the above vectors by
applying $R_1$, $S_1$ and/or powers of $P$.

We give the verifications for one situation, namely $\G=\Rc\cap
\Sc=\Bn_1\cap \Bn_3$ (see Table~\ref{tab:numbering} on
page~\pageref{tab:numbering} for an explanation of the labelling of
the sides). For concreteness, we start by assuming $p=3$ (for other
values of $p$ and other ridges, the computations are only longer, but
not more complicated). In this case $T$ is a triangle.

The intersection $\G\cap \Bn_2$ is particularly easy, since
$\Rc\cap\Rc^-$ is a complex line, namely the mirror of $R_1$ (see
Section~\ref{sec:ridgeconsistency}). If we take ${\bf x}_0={\bf y}_0$,
${\bf x}_1=R_1{\bf y}_0$ 
and write an equation for $\G\cap \Bn_2$, we get
$$
  g(t_1,t_2)=81(582+127\sqrt{21})(-14+3\sqrt{21}+3t_2-18t_2^2)t_1,
$$ 
where the degree $2$ polynomial in $t_2$ never vanishes, so the
intersection is given by the $t_2$-axis.

When $\B=\Sc^-=\Bn_4$, we take ${\bf x}_0={\bf y}_0$, ${\bf x}_1=S_1{\bf
  x}$. In that case, we already know that the $t_2$-axis is part of
$\G\cap \Sc^-$ (see Section~\ref{sec:edges}), so we must have
$a_0=a_2=a_{22}=0$, but in fact no other coefficient vanishes. The
reduced equation $g=f/t_1$ is given by
$$
  h(t_1,t_2)=81(6741+1471\sqrt{21}) 
              (36-8\sqrt{21}-15t_1+3t_1\sqrt{21}-3t_2\sqrt{21}+9t_2+18t_1t_2)/4.
$$

Here the partial derivatives are both functions of one variable only,
which immediately gives a critical point with coordinates
$$
  \left((-3+\sqrt{21})/6,(5-\sqrt{21})/6\right)
  \approx(0.26376261\dots,0.06957071\dots).
$$
It is easy to verify that this is outside $T$.

\begin{table}[htbp]
\centering
\begin{tabular}{|r|c|c|c|c|c|}\hline
$p$ &  $\Bn_2$  & $\Bn_4$     & $\Bn_{11}$     & $\Bn_{12}$  & $\Bn_{18}$                      \\\hline
$3$ & $t_1=0$ & $t_1=0$               & $t_1=(21+\sqrt{21})/42$ &       &$t_1=(-3+\sqrt{21})/6$     \\
    &         & $t_2=(7-\sqrt{21})/21$& $t_2=0$                 &$t_2=0$&                           \\\hline
$4$ & $t_1=0$ & $t_1=0$               & $t_1=(-7+5\sqrt{7})/14$ &       &$t_1=(3-\sqrt{7})/2$       \\
    &         & $t_2=1-5\sqrt{7}/14$  & $t_2=0$                 &$t_2=0$&                           \\\hline
$5$ & $t_1=0$ & $t_1=0$               & $t_1=0.44786394\dots$   &       &$t_1=0.18371174\dots$      \\
    &         & $t_2=0.05801393\dots$ & $t_2=0$                 &$t_2=0$&                           \\\hline
$6$ & $t_1=0$ & $t_1=0$               & $t_1=(-7+3\sqrt{21})/14$&       &$t_1=(5-\sqrt{21})/2$      \\
    &         & $t_2=4-6\sqrt{21}/7$  & $t_2=0$                 &$t_2=0$&                           \\\hline
$8$ & $t_1=0$ & $t_1=0$               & $t_1=0.57827900\dots$   &       &$t_1=0.27949078\dots$      \\
    &         & $t_2=0.11986937\dots$ & $t_2=0$                 &$t_2=0$&                           \\\hline
$12$& $t_1=0$ & $t_1=0$               & $t_1=0.80219658\dots$   &       &$t_1=(\sqrt{3}-\sqrt{7})/2$\\
    &         & $t_2=0.28759793\dots$ & $t_2=0$                 &$t_2=0$&                           \\\hline 
\end{tabular}

\vspace{.3cm}

\begin{tabular}{|r|c|c|c|c|}\hline
$p$ &  $\Bn_{21}$   & $\Bn_{22}$     & $\Bn_{24}$      & $\Bn_{26}$                \\\hline
$3$ & $t_1=(-3+\sqrt{21})/6$ & $t_1=1/2+\sqrt{21}/42$  &                         & $t_1=1/2+\sqrt{21}/42$  \\
    & $t_2=(7-\sqrt{21})/21$ & $t_2=(5-\sqrt{21})/6$   & $t_2=(5-\sqrt{21})/6$   &                         \\\hline
$4$ & $t_1=(3-\sqrt{7})/2$   & $t_1=-1/2+5\sqrt{7}/14$ &                         & $t_1=-1/2+5\sqrt{7}/14$ \\
    & $t_2=1-5\sqrt{7}/14$   & $t_2=4-3\sqrt{7}/2$     & $t_2=4-3\sqrt{7}/2$     &                         \\\hline
$5$ & $t_1=0.18371174\dots$  & $t_1=0.44786394\dots$   &                         & $t_1=0.44786394\dots$   \\
    & $t_2=0.05801393\dots$  & $t_2=0.03375000\dots$   & $t_2=0.03375000\dots$   &                         \\\hline
$6$ & $t_1=(5-\sqrt{21})/2$  & $t_1=-1/2+3\sqrt{21}/14$&                         & $t_1=-1/2+3\sqrt{21}/14$\\
    & $t_2=4-6\sqrt{21}/7$   & $t_2=(23-5\sqrt{21})/2$ & $t_2=(23-5\sqrt{21})/2$ &                         \\\hline
$8$ & $t_1=0.27949078\dots$  & $t_1=0.57827900\dots$   &                         & $t_1=0.57827900\dots$   \\
    & $t_2=0.11986937\dots$  & $t_2=0.07811509\dots$   & $t_2=0.07811509\dots$   &                         \\\hline
$12$& $t_1=0.45685025\dots$  & $t_1=0.80219658\dots$   &                         & $t_1=0.80219658\dots$   \\
    & $t_2=0.28759793\dots$  & $t_2=0.20871215\dots$   & $t_2=(5-\sqrt{21})/2$   &                         \\\hline
\end{tabular}
\caption{Lines parallel to coordinate axes contained in the graphs,
  for $\Bn_1\cap \Bn_3$.}\label{tab:lines}
\end{table}

The case $\B=P\Rc=\Bn_5$ is in a sense generic; there are no horizontal
or vertical lines, and the resultant of $\partial f/\partial t_1$ and
$\partial f/\partial t_2$ with respect to $t_2$ is given by
\begin{eqnarray*}
  &&-223205130784873566-48707352729720306\sqrt{21}\\
  &&+(1289075326223152230+281299298043708180\sqrt{21})t_1\\
  &&-(1960742086219144494+427869001351675404\sqrt{21})t_1^2\\
  &&+(2269952366448531240+495344216342598840\sqrt{21})t_1^3\\
  &&-(1166845194031834926+254626496473636566\sqrt{21})t_1^4\\
  &&+(259717950970322292+56675103317162172\sqrt{21})t_1^5.
\end{eqnarray*}
This polynomial has precisely one real root, which is a simple root,
given approximately by $t_1=0.23860554\dots$. Substituting this value
into the partial derivatives yields a polynomial of degree $1$ and a
polynomial of degree $2$; solving the degree one equation gives one
value for $t_2$, and we get precisely one critical point with
coordinates $(0.23860554\dots,0.01603880\dots)$. This point is outside
$T$.

\begin{table}[htbp]
\centering
\begin{tabular}{|r|c|c|c|c|c|c|c|c|c|} \hline
$p$ &  Side & Coordinates                                                            & Function value   \\\hline
  8 &  $\Bn_9$  & (0.22347669\dots, 0.06214532\dots)  & 11.44128654\dots \\\hline
 12 &  $\Bn_9$  & (0.38008133\dots, 0.18112900\dots)  & 0.23597323\dots  \\
    &  $\Bn_{10}$ & (0.38823985\dots, 0.13852442\dots)  & 0.18097266\dots  \\\hline
\end{tabular}  
\caption{Critical points inside the ridge for $\Bn_1\cap \Bn_3$.}\label{tab:criticalpointsinside}
\end{table}

From now on, we do not write down all the calculations, but simply
mention some of the possible phenomena. 

Consider the intersection of $\G$ with $P\Sc^-=\Bn_{11}$, for
instance. In that case $g$ is given by
$$
  27(6741+1471\sqrt{21})(-15t_1+3t_1\sqrt{21}+36-8\sqrt{21}+18t_1t_2-3t_2\sqrt{21}+9t_2)(21-42t_1+\sqrt{21})t_2/56.
$$ 
Clearly there are both a horizontal and a vertical line in this intersection, in which case
we factor out two linear factors, one linear in $t_1$ and the other
linear in $t_2$. The reduced equation is given by the function
$$
  81(6741+1471\sqrt{21})(-15t_1+3t_1\sqrt{21}+36-8\sqrt{21}+18t_1t_2-3t_2\sqrt{21}+9t_2)/4.
$$
The latter has one critical point given by
$$
  \left((-3+\sqrt{21})/6, (5-\sqrt{21})/6\right)
  \approx(0.26376261\dots,0.06957071\dots)
$$
and this point is outside $T$.

In fact for $p\leqslant 6$, all critical points turn out to be outside
the Giraud polygon, whereas this is not quite true for $p=8$ or $12$
(see Table~\ref{tab:criticalpointsinside}).  For higher values of $p$,
the computations are a bit longer and more tedious than for $p=3$,
since the coefficients of the equations live in slightly larger number
fields, but they are not more difficult in any essential way.

For the convenience of the skeptic reader (who may be willing to check
our claims), we list all horizontal and vertical lines for curves on
$\Bn_1\cap\Bn_3$, see Table~\ref{tab:lines}. These are needed to reduce to
the case of resultants with only simple root. The numbers given as
decimal expansions are written in this manner simply to avoid writing
huge formulae, one can easily get explicit algebraic
expressions instead.

In Table~\ref{tab:criticalpointsinside}, we list the critical points
that lie inside the ridge $\Bn_1\cap \Bn_3$ (there are critical points
inside only for $p=8$ and $p=12$). For each of them, we need to check
that the corresponding value of the function $h$ is strictly positive
(recall that $h$ is equal to the equation of the intersection in most
cases, it is obtained by dividing by appropriate linear factors in one
variable).

For future reference we mention the following consequence of the above arguments.
\begin{proposition} \label{prop:ridgenbd}
  Let $\rho$ be a ridge of $\widehat{E}$, and let $D$ be the (Giraud
  or complex) disk containing $\phi(\rho)$. Then there is an open
  neighborhood $U$ of $\phi(\rho)$ in $D$ such that $U\cap
  E=\phi(\rho)$.
\end{proposition}
In particular, even though we did not prove that $\phi(\rho)=D\cap F$, we
know that $\phi(\rho)$ is a connected component of $D\cap F$.

\subsubsection{Realizing sides of $\widehat{E}$} \label{sec:3-cells}

As mentioned early in Section~\ref{sec:combinatorics}, no consistency
verifications are needed for $3$-cells. We do need to check
embeddedness. More specifically, we check the following.
\begin{proposition}\label{prop:3faces}
  For each $3$-cell $f$ of $\widehat{E}$ and every bounding bisector
  $\B$, $\phi(f)$ is on the correct side of $\B$. In particular,
  $\phi(f)$ is contained in $F$.
\end{proposition}

This will follow from the following result, which is of general
independent interest. 
\begin{lemma}\label{lem:3faces}
  Let $f$ be one of the bounding sides and let $\B$ be any
  bisector. If $\B$ intersects the interior $f^\circ$, then it
  intersects its boundary $\partial f$ as well.
\end{lemma}

\Pf (of Lemma~\ref{lem:3faces}) 
We may assume that $\B$ does not
contain $f$, and we denote by $\B_f$ the bisector containing $f$. By
Lemma~9.1.5 in~\cite{goldman}, $\B\cap \B_f$ is either empty or a
union of properly embedded disks.  If $\B$ intersects the interior of
$f$, then by properness it must also intersect its boundary (in $\chb
2$).
\EPf

\Pf (of Proposition~\ref{prop:3faces}) 
Recall that for each 3-cell $f$
of $\widehat{E}$, $\partial f$ is a 2-sphere, and we have checked in
Section~\ref{sec:ridgeembedding} that $\phi$ embeds that 2-sphere
inside one of the bounding bisectors $\B_f$ (as well as inside $E$).
Now if $\B$ is any other bounding bisector, $\B$ cannot intersect
$\phi(f^\circ)$. For this we argue by contradiction; if it did, then
by Lemma~\ref{lem:3faces}, it would also intersect $\phi(\partial
f)$. But this cannot happen, since the closed set $\phi(f)\cap \B\cap
\B_f$ has a neighborhood in $\B\cap \B_f$ which does not intersect
$\phi(f^\circ)$ (see Proposition~\ref{prop:ridgenbd}). 
\EPf

For completeness, we mention the following consequence of
Proposition~\ref{prop:3faces}.
\begin{proposition}\label{prop:sideonly}
  Let $\B$ be one of the bounding bisectors, and $f$ the side of
  $\widehat{E}$ such that $\phi(f) \subset \B$. Then $\B\cap E=\phi(f)$.
\end{proposition}
Note that we have not proved that $\phi(f)$ can also be described as
$\B\cap F$, but once again this is not needed in order to apply the
Poincar\'e polyhedron theorem (see the discussion at the beginning of
Section~\ref{sec:domain}).

\section{Applying the Poincar\'e polyhedron theorem}\label{sec:applyingpoincare}

We proved in Section~\ref{sec:combinatorics} that $E$ is a
polyhedron. In this section we verify the remaining hypotheses of the
Poincar\'e polyhedron theorem (Sections~\ref{sec:tess}
and~\ref{sec:stabilisers}). Using the theorem, we obtain a
presentation for $\Gamma$ and we calculate the orbifold Euler
characteristic of the quotient orbifold.

Along with the verification of the hypotheses, we will determine local
stabilizers of every facet $f$ of $E$. By local stabilizer, we mean
the group obtained from the side-pairing maps by following all
possible cycles involving $f$. The description of local stabilizers
for cusps will imply the existence of a consistent system of
horoballs.

Note that the local stabilizer of a facet $f$ is always a subgroup of
the full stabilizer of $f$ in $\Gamma$. It will follow from the
Poincar\'e polyhedron that local stabilizers are actually equal to the
corresponding full stabilizer.  We will determine local stabilizers of
ridges in Section~\ref{sec:tess}, then for vertices and edges in
Section~\ref{sec:stabilisers}.

\subsection{Hypotheses - side pairing}\label{sec:pairing}

In this section, we describe a side pairing for $E$ (in the sense of
Section~\ref{sec:poin}). For each side $s$, we define an isometry
$S=\sigma(s)$ that pairs two sides $s$ and $s^-$ of $E$. We also show
that $E\cap S(E)$ is given precisely by $s^-$.

\begin{lemma}\label{lem:side-pairing}
  The map $R_1$ sends $r_1=\Rc\cap E$ to $r_1^-=\Rc^-\cap E$
  preserving the cell structure. The map $S_1$ sends $s_1=\Sc\cap E$
  to $s_1^-=\Sc^-\cap E$ preserving the cell structure.
\end{lemma}

\Pf 
Using the description of the bisectors given in
Table~\ref{tab:spines} it is clear that $R_1$ maps $\Rc$ to $\Rc^-$.
Moreover, $R_1$ maps the complex lines containing ridges of $r_1$ to
the complex lines containing ridges of $r_1^-$. Specifically, $R_1$
sends $m_1=\Rc\cap\Rc^-$ to itself and maps $m_{23}=\Rc\cap
P^{-3}\Rc^-$ to $m_{123\bar{1}}=\Rc^-\cap P^3\Rc$.

We claim that the Giraud disks containing ridges of $\Rc$ are mapped
to Giraud disks containing ridges of $\Rc^-$. For a list of these
Giraud disks, see Table~\ref{tab:giraud1}. We give the details for one
such Giraud disk, namely $\Rc\cap P^2\Sc$, the other arguments are
entirely similar.  Recall that the two corresponding bisectors have
the following description in terms of the point $y_2$:
$$
  \Rc = \B(y_2,R_1^{-1}P^2y_2), \quad P^2\Sc=\B(y_2,P^2S_1^{-1}P^{-3}y_2).
$$
Applying $R_1$ to these bisectors, we obtain
$$
  \Rc^-=R_1\Rc=\B(R_1y_2,P^2y_2),\quad 
  R_1P^2\Sc=\B(R_1y_2,P^2R_1^{-1}P^2y_2).
$$
Here we have used $S_1=P^2R_1P^{-2}R_1P^2$ and $P^{-5}=P^2$ to write
$R_1P^2S_1^{-1}P^{-3}=P^2R_1^{-1}P^2$. The third bisector containing this
Giraud disk is
$$
  \B(P^2y_2,P^2R_1^{-1}P^2y_2)=P^2\B(y_2,R_1^{-1}P^2y_2)=P^2\Rc.
$$
Therefore the Giraud disk $\Rc\cap P^2\Sc$ is sent by $R_1$ to 
the Giraud disk $\Rc^-\cap P^2\Rc$. 

This shows that the complex lines and Giraud disks containing the
ridges of $r_1$ are sent to complex lines and Giraud disks containing
ridges of $r_1^-$. The corresponding statement about vertices
(resp. edges) follows from the one about ridges, since the vertices
(resp. edges) can be described as intersections of ridges. This shows
that $R_1$ maps $r_1=\Rc\cap E$ to $r_1^-=\Rc^-\cap E$, preserving
the cell structure.  
\EPf

We now define the side pairing used in our application of the
Poincar\'e polyhedron theorem.

\begin{proposition}\label{prop:side-pairing}
  The following map $\sigma$ defines a side pairing on the sides of
  $E$:
  $$
    \sigma(P^kr_1^\pm)=P^kR_1^{\pm 1}P^{-k},\quad 
    \sigma(P^ks_1^\pm)=P^kS_1^{\pm 1}P^{-k},
  $$ 
  where $k=-3,\,\ldots,\,3$.  Moreover, this side pairing is
  compatible with $\Upsilon=\langle P\rangle$.
\end{proposition}

\Pf 
We only give the argument for $r_1=\Rc\cap E$. Applying powers of
$P$ and the symmetry $\iota$ gives the result for the other faces
$P^kr_1^\pm$, and the faces $P^ks_1^\pm$ are similar.

In Lemma~\ref{lem:side-pairing} we showed that $R_1$ sends $r_1$ to
$r^-_1$ preserving the cell structure. Moreover, the interior of $E$
is contained in the half-space closer to $y_0$ than to $R_1y_0$, where
$y_0$ is given in Lemma~\eqref{lem:xyz}. Hence, the interior of
$R_1^{-1}E$ is contained in the half-space closer to $R_1^{-1}y_0$ than
to $y_0$. Thus $E$ and $R_1^{-1}E$ intersect in $\Rc\cap E=r_1$ (see
Proposition~\ref{prop:sideonly}) and their interiors are
disjoint. Furthermore, any point in the interior of $r_1$ has an open
neighborhood contained in $E\cup R_1^{-1}E$.  \EPf

\subsection{Hypotheses - local tessellation}\label{sec:tess}

We now study the ridge cycles of $E$, as explained in
Section~\ref{sec:poin}.

\begin{table}
$$
  \begin{array}{|r|l|llll|}
    \hline
    R_1 & r_1\cap r_1^-
    & p_{12} & p_{13} & q_{1\bar{3}23} & q_{123\bar{2}} \\
    & r_1\cap r_1^-
    & p_{12} & p_{13} & q_{1\bar{3}23} & q_{123\bar{2}} \\
    \hline
    R_1 & r_1\cap P^{-1}r_1^-
    & p_{12} & p_{13} & p_{23} & \\
    P^{-1} & Pr_1\cap r_1^-
    & p_{12} & p_{13} & p_{123\bar{1}} & \\
    & r_1\cap P^{-1}r_1
    & p_{13} & p_{23} & p_{12} & \\
    \hline
    R_1 & r_1\cap P^2s_1          
    & p_{23} & p_{12} & q_{123\bar{2}} & \\
    P^2R_1P^{-2} & P^2r_1\cap r_1^- 
    & p_{123\bar{1}} & p_{12} & q_{123\bar{2}} & \\
    S_1^{-1} & s_1^-\cap P^2r_1^-   
    & p_{123\bar{1}} & p_{\bar{3}\,\bar{2}3123}  & q_{123\bar{2}} & \\
    P^2 & P^{-2}r_1\cap s_1        
    & p_{\bar{3}\,\bar{2}3123} & p_{23} & q_{1\bar{3}23} & \\  
    & r_1\cap P^2s_1 & p_{23} & p_{12} & q_{123\bar{2}} & \\
    \hline
    R_1 & r_1\cap s_1 
    & p_{23} & q_{1\bar{3}23} & q_{123\bar{2}} & \\
    S_1^{-1} & s_1^-\cap r_1^-
    & p_{123\bar{1}} & q_{1\bar{3}23} & q_{123\bar{2}} & \\
    S_1^{-1} & s_1^-\cap s_1
    & p_{\bar{3}\,\bar{2}3123} & q_{123\bar{2}} & q_{1\bar{3}23} & \\
    & r_1\cap s_1
    & p_{23} & q_{1\bar{3}23} & q_{123\bar{2}} & \\
    \hline
    S_1 & s_1\cap P^2s_1^-
    & p_{23} & p_{\bar{3}\,\bar{2}3123} & q_{123\bar{2}} & \\
    P^2 & P^{-2}s_1\cap s_1^-
    & p_{\bar{3}\,\bar{2}123} & p_{123\bar{1}} & q_{1\bar{3}23} & \\
    & s_1 \cap P^2s_1^-
    & p_{23} & p_{\bar{3}\,\bar{2}3123} & q_{123\bar{2}} & \\
    \hline
  \end{array}
$$
\caption{The ridge cycles for $p=3$, $4$. For higher values of $p$ the cycles
are the same, but additional vertices appear in the right and column.}
\label{tab:cycles3-4}
\end{table}
For $p=3$ and $p=4$ we obtain the ridge cycles given in
Table~\ref{tab:cycles3-4}. The second column lists the successive
ridges in the cycle, given as the intersection of two sides, whereas
in the third column we list the corresponding vertices (with ordering
coherent with the successive side pairing maps). Applying the side
pairing map in the first column (which corresponds to the first side
in the second column) gives the line below. The product of these side
pairing maps gives the cycle transformation.

This cycle transformation may not be the identity on $\ch 2$, or even
on the ridge. However, we can find a power that is the identity, which
gives the cycle relation. For example, $R_1$ maps $r_1\cap r_1^-$ to
itself and acts as the identity on this ridge (which can be seen from
the vertices in the third column). The cycle relation is
$R_1^p=id$. On the other hand, for $r_1\cap P^{-1}r_1^-$ we obtain the
cycle transformation $P^{-1}R_1$. This maps $r_1\cap P^{-1}r_1^-$ to
itself but cyclically permutes the vertices. Therefore the cycle
relation is $(P^{-1}R_1)^3=id$.

As a last example, we work out the cycle for the ridge $r_1\cap
P^2s_1$. It is mapped by $R_1$ to $P^2r_1\cap r_1^-$ (see the proof of
Proposition~\ref{prop:side-pairing}), which is then mapped by
$P^2R_1P^{-2}$ to $s_1^-\cap P^2r_1^-$, which is then mapped by
$S_1^{-1}$ to $P^{-2}r_1\cap s_1$. This is clearly in the image of the
original ridge under a power of $P$, and the corresponding cycle
transformation is
$$
P^2\cdot S_1^{-1}\cdot P^2R_1P^{-2} \cdot R_1.
$$ 
This isometry is in fact the identity (see its expression in terms
of the vertices in the last column of Table~\ref{tab:cycles3-4}),
which corresponds to our definition of $S_1$ (see
Section~\ref{sec:notation}).

For $p=5$, $6$, $8$ and $12$, the ridge cycles from $p=3$, $4$
persist, even though the corresponding polygons have extra vertices
corresponding to the truncations.  We also have an extra ridge cycle
associated to $r_1\cap P^{-3}r_1^-$ given in
Table~\ref{tab:cycles-higher}.
\begin{table}
$$
  \begin{array}{|l|r|l|llll|}
    \hline
    p\geqslant 5 & R_1 & r_1\cap P^{-3}r_1^- 
    & p_{23}^2 & p_{23}^{23\bar{2}} & p_{23}^{\bar{3}23} & p_{23}^3 \\
    & P^{-3} & P^3r_1\cap r_1^-
    & p_{123\bar{1}}^{12\bar{1}} & p_{123\bar{1}}^{123\bar{2}\,\bar{1}}
    & p_{123\bar{1}}^{1\bar{3}23\bar{1}} & p_{23}^{13\bar{1}} \\
    & & r_1\cap P^{-3}r_1^-
    & p_{23}^{\bar{3}23} & p_{23}^3 & p_{23}^2 & p_{23}^{23\bar{2}} \\
    \hline
    p\geqslant 8 & S_1 & s_1\cap P^{-2}s_1^- 
    & q_{1\bar{3}23}^1 & q_{1\bar{3}23}^{\bar{3}23} &
    q_{1\bar{3}23}^{1\bar{3}23\bar{1}} & \\
    & P^{-2} & P^2 \cap s_1^-
    & q^1_{123\bar{2}} & q^{23\bar{2}}_{123\bar{2}}
    & q_{123\bar{2}}^{123\bar{2}\,\bar{1}} & \\
    & & s_1\cap P^{-2}s_1^-
    & q_{1\bar{3}23}^{\bar{3}23} 
    & q_{1\bar{3}23}^{1\bar{3}23\bar{1}} & q_{1\bar{3}23}^1 & \\
    \hline
  \end{array}
$$ 
\caption{Additional ridge cycles for higher values of $p$.}
\label{tab:cycles-higher}
\end{table}
Note that $P^{-3}R_1$ maps $r_1\cap P^{-3}r_1^-$ to itself but does not act as
the identity on this ridge; $(P^{-3}R_1)^2$ acts as the identity on  
$r_1\cap P^{-3}r_1^-$ but not on $\ch 2$ and $(P^{-3}R_1)^{4p/(p-4)}=id$.
When $p=8$ and $p=12$ we also have the ridge cycle associated to 
$s_1\cap P^{-2}s_1^-$ given in Table~\ref{tab:cycles-higher}.

As a summary, thanks to the action of powers of $P$ and the freedom to
choose the initial ridge inside a given ridge cycle (which does not
affect the cycle relation, see Section~\ref{sec:poin}), it is enough
to consider the ridges $r_1\cap r_1^-$, $r_1\cap
P^{-1}r_1^-$, $r_1\cap P^2s_1$, $r_1\cap s_1$, $s_1\cap P^2s_1^-$,
$r_1\cap P^{-3}r_1^-$ and $s_1\cap P^{-2}s_1^-$. In the penultimate
case we only need consider $p=5$, $6$, $8$ and $12$ and for the last
case, we only need consider $p=8$ and $12$. We now show local tiling
around each of these ridges.

Recall that the ridges of $E$ are of two very different types,
depending on whether they are contained in complex lines or in Giraud
disks (see Section~\ref{sec:ridgeconsistency}).

We first consider the complex ridges (see
Proposition~\ref{prop:cridges} and Section~\ref{sec:ridgeembedding}),
and show that the local images of $E$ tessellate around those.

\begin{lemma}\label{lem:cx-line-tess}
  The images of $E$ tessellate a neighborhood of $r_1\cap r_1^-$,
  $r_1\cap P^{-3}r_1^-$, $s_1\cap P^2s_1^-$ and $s_1\cap P^{-2}s_1^-$
  for the appropriate values of $p$.  Specifically:
\begin{enumerate}
\item $R_1^kE$ for $k=0,\,1\,\ldots,\, p-1$ cover a neighborhood of
  the interior of $r_1\cap r_1^-$;
\item if $p=5$, $6$, $8$, $12$ then $(R_2R_3)^kE$ for
  $k=0,\,1\,\ldots,\, 2c-1=4p/(p-4)-1$ cover a neighborhood of the
  interior of $r_1\cap P^{-3}r_1^-$;
\item $(R_2R_3R_2^{-1})^kE$ for $k=0,\,1\,\ldots,\, p-1$ cover a
  neighborhood of the interior of $s_1\cap P^2s_1^-$;
\item if $p=8$, $12$ then $(R_1R_3^{-1}R_2R_3)^kE$ for
  $k=0,\,1\,\ldots,\, 3d-1=6p/(p-6)-1$ cover a neighborhood of the
  interior of $s_1\cap P^{-2}s_1^-$.
\end{enumerate}
\end{lemma}

\Pf 
Using Lemma~\ref{lem:cx-tess} it suffices to consider a fixed
point $o$ in each of these ridges and to show that on the orthogonal
complex line $C^\perp_o$ the cycle transformation acts as a rotation
through angle $2\pi/\ell$.

For Parts~1 and~3 this is straightforward since $R_1$ and
$R_2R_3R_2^{-1}$ are complex reflections with angle $2\pi/p$ fixing
$r_1\cap r_1^-$ and $s_1\cap P^2s_1^-$ respectively.  For Parts~2
and~4 we find the rotation angle by finding eigenvalues associated to
the complex line and the fixed point.

For $p\geqslant 5$, the complex line $m_{23}$ and the fixed point
$o_{23}$ of $R_2R_3$ contained in $m_{23}$ correspond to the
eigenvectors ${\bf n}_{23}$ and ${\bf o}_{23}$ of $R_2R_3$ given by:
$$
  {\bf n}_{23}=\left[\begin{matrix} 
      a^3+\overline{a}^3\\
      a^2\overline{\tau}-\overline{a}\\
      \overline{a}^2\tau-a
    \end{matrix}\right],\quad
  {\bf o}_{23}=\left[\begin{matrix}
      0 \\ a\tau \\ 1-i \end{matrix}\right].
$$
(Note that $\langle {\bf o}_{23},{\bf o}_{23}\rangle<0$ for $p\ge 5$.)
One easily checks that the eigenvalues are $\bar{a}^2$ and $-ia$ respectively.
Hence the rotation angle is $(-4\pi/3)-(2\pi/3p-\pi/2)=(p-4)\pi/2p=2\pi/2c$. 

For $p\geqslant 8$, the complex line $m_{1\bar{3}23}$ and the fixed point 
$o_{1\bar{3}23}$ of $R_1R_3^{-1}R_2R_3$ correspond to eigenvectors
${\bf n}_{1\bar{3}23}$ and ${\bf o}_{1\bar{3}23}$ given by:
$$
  {\bf n}_{1\bar{3}23}=\left[\begin{matrix}
      a^2+\bar{a}\bar{\tau} \\ a^4\tau+a \\ a^3-\bar{a}^3-\tau
    \end{matrix}\right],\quad
  {\bf o}_{1\bar{3}23}=\left[\begin{matrix}
      -a^2\bar{\omega} \\ a \\ \bar{\tau}\end{matrix}\right]
$$
where $\omega=e^{2\pi i/3}$. One easily checks that the eigenvalues
are $\bar{a}^2$ and $-a\omega$ respectively and so the rotation angle is
$(-4\pi/3)-(2\pi/3p-\pi/3)=(p-6)\pi/3p=2\pi/3d$.
\EPf

We now consider the Giraud ridges (see Table~\ref{tab:giraud1} for a
list of the relevant Giraud intersections, and
Section~\ref{sec:ridgeembedding} for the fact that these are indeed
ridges of $E$).

\begin{lemma}\label{lem:Giraud-tess}
  The images of $E$ tessellate a neighborhood of $r_1\cap
  P^{-1}r_1^-$, $r_1\cap P^2s_1$ and $r_1\cap s_1$. Specifically:
  \begin{enumerate}
  \item $E$, $R_1^{-1}E$ and $P^{-1}R_1E$ cover a neighborhood of the
    interior of $r_1\cap P^{-1}r_1^-$;
  \item $E$, $R_1^{-1}E$ and $P^2S_1^{-1}E$ cover a neighborhood of
    the interior of $r_1\cap P^2s_1$;
  \item $E$, $R_1^{-1}E$ and $S_1^{-1}E$ cover a neighborhood of the
    interior of $r_1\cap s_1$.
  \end{enumerate}
\end{lemma}

\Pf 
Using Lemma~\ref{lem:giraud-tess}, this follows from the fact that
these three ridges are contained in the Giraud disks given in
Table~\ref{tab:giraud1}.  The bounding bisectors defining these Giraud
disks are coequidistant from one of $y_0$, $y_1$ or $y_2$, as
indicated in Table~\ref{tab:bis-xyz}. These bisectors divide $\ch 2$
into three regions as in Lemma~\ref{lem:giraud-tess}. We know by
Lemma~\ref{lem:xyz} that $E$ is contained in the region containing
$y_j$. By applying the side pairing maps in the ridge cycle we obtain
copies of $E$ contained in the other two regions.
\EPf

This proves that the side pairings of the polyhedron $E$ satisfy the
cycle conditions of the Poincar\'e polyhedron theorem.

\subsection{Hypotheses - consistent horoballs}\label{sec:stabilisers}

In order to prove the existence of a consistent system of horoballs
(see Corollary~\ref{cor:4-horo} and Corollary~\ref{cor:6-horo}), we
need to study local stabilizers of cusps. 

In fact we will determine local stabilizers of all vertices (cusps or
not), as well as edges. The techniques are similar to those in
Section~18.2 of~\cite{mostowpacific}. The order of these stabilizers
will be used when we calculate the orbifold Euler characteristic, see
Section~\ref{sec:euler}.

\begin{proposition}\label{prop:stab-p-34}
When $p=3$ or $4$ there is a single-orbit of $p_*$-vertices. A representative
is $p_{12}$ and its stabilizer is $\langle R_1,R_2\rangle$. Moreover, there is
a single orbit of $(p_*,p_*)$ edges. A representative is $(p_{13}, p_{12})$ and
its stabilizer is $\langle R_1\rangle$.
\end{proposition}

\Pf
Since $p_{23}=P^{-2}p_{12}$, $p_{13}=P^{-1}p_{12}$,
$p_{123\bar{1}}=Pp_{12}$ and $p_{\bar{3}\,\bar{2}3123}=P^3p_{12}$ it
is clear that all $p_*$-vertices lie in the same orbit.

The action of the side pairing maps on the $p_*$-vertices may be
expressed in terms of $p_{12}$ as follows. We can write such a vertex
as $p_a=P^ip_{12}$.  The side pairing $R$ maps $p_a$ to $p_b$ and we
can write $p_b=P^jp_{12}$.  This corresponds to an element
$P^jRP^{-i}$ of the stabilizer of $p_{12}$.  We list the side pairing,
the vertex, its image and the corresponding word in the stabilizer of
$p_{12}$ in terms of $R_1$ and $R_2$. Taking the inverse of a side
pairing map (resp. conjugating it by a power of $P$) gives the inverse
word (resp. a conjugate word).
$$
  \begin{array}{|l|l|l|l|}
    \hline
    & \hbox{Vertex} & \hbox{Image} & \hbox{Word} \\
    \hline
    R_1 & p_{12} & p_{12} & R_1 \\
    R_1 & p_{13}=P^{-1}p_{12} & p_{13}=P^{-1}p_{12} & PR_1P^{-1}=R_1R_2R_1^{-1} \\
    R_1 & p_{23}=P^{-2}p_{12} & p_{123\bar{1}}=Pp_{12} & P^{-1}R_1P^{-2}=R_2^{-1}R_1^{-1} \\
    S_1 & p_{23}=P^{-2}p_{12} & p_{\bar{3}\,\bar{2}3123}=P^3p_{12} 
    & P^{-3}S_1P^{-2}=R_2^{-1} \\
    S_1 & p_{\bar{3}\,\bar{2}3123}=P^3p_{12} & p_{123\bar{1}}=Pp_{12} 
    & P^{-1}S_1P^3=R_2^{-1}R_1^{-1}R_2 \\
    \hline
  \end{array}
$$ 
To show that the stabilizer of $(p_{13},p_{12})$ is generated by $R_1$
we argue in the same way as we did for vertex stabilizers.  All edges
joining two $p_*$-vertices are the images under a power of $P$ of
either the edge $(p_{13},p_{12})$ or the edge
$(p_{12},p_{23})=R_1^{-1}P(p_{13},p_{12})$. We describe each edge
$(p_*,p_*)$ and its image under the side pairing in this way and so
obtain an element of the stabilizer of $(p_{13},p_{12})$.  We do this
for the three $(p_*,p_*)$ edges in $r_1$ and the single $(p_*,p_*)$
edge in $s_1$.  The others follow by applying powers of $P$.
$$
\begin{array}{|l|l|l|l|}
  \hline
  & \hbox{Edge} & \hbox{Image} & \hbox{Word}\\
  \hline
  R_1 & (p_{13},p_{12}) & (p_{13},p_{12}) & R_1 \\
  R_1 & (p_{12},p_{23})=R_1^{-1}P(p_{13},p_{12}) 
  & (p_{12},p_{123\bar{1}})= P(p_{13},p_{12}) & id \\
  R_1 & (p_{23},p_{12})=P^{-1}(p_{13},p_{12}) 
  & (p_{123\bar{1}},p_{13})=PR_1^{-1}P(p_{13},p_{12}) & R_1^{-1} \\
  S_1 & (p_{23},p_{\bar{3}\bar{2}3123})=P^{-2}R_1^{-1}P(p_{13},p_{12})
  & (p_{\bar{3}\bar{2}3123},p_{123\bar{1}}) =P^3R_1^{-1}P(p_{13},p_{12})) & R_1^{-1}\\
  \hline
\end{array}
$$
In fact no group element interchanges $p_{12}$ and $p_{13}$, and so the 
stabilizer of the edge $(p_{13},p_{12})$ is the intersection of the vertex 
stabilizers, namely $\langle R_1\rangle$. 
\EPf

In Proposition~\ref{prop:stab-12} we found the structure of $\langle
R_1,R_2\rangle$. In particular, for $p=3$ it is finite.  When $p=4$,
the point $p_{12}$ lies on $\partial\ch 2$ and we have just verified
the conjectural stabilizer given in \cite{dpp1}. Since there is a single orbit
of cusps, a consistent system of horoballs comprises a single
horoball at $p_{12}$ and all its images under the local stabilizer.
Since the local stabilizer of $p_{12}$ is generated by complex reflections, 
any cycle of side pairing maps or elements of $\Upsilon$ that maps
$p_{12}$ to itself also maps any horoball at $p_{12}$ to itself.
Thus, we immediately have the following corollary:

\begin{corollary}\label{cor:4-horo}
  When $p=4$, the (local) stabilizer in $\Gamma$ of
  $p_{12}\in\partial{\bf H}^2_\C$ contains no loxodromic maps. In
  particular, there is a consistent system of horoballs.
\end{corollary}

\begin{proposition}\label{prop:stab-p56812}
  When $p\geqslant 5$ there are two orbits of
  $p_*$-vertices. Representatives are $p^1_{12}$ and $p^1_{13}$ with
  stabilizers $\langle R_1,(R_1R_2)^2\rangle$ and $\langle
  R_1,(R_1R_3)^2\rangle$ respectively. There are two orbits of
  $(p_*,p_*)$ edges. Representatives are $(p^1_{13},p^1_{12})$ and
  $(p^1_{12},p^2_{12})$ with stabilizers $\langle R_1\rangle$ and
  $\langle (R_1R_2)^2\rangle$ respectively.
\end{proposition}

\Pf
The $\Gamma$-orbit of $p^1_{12}$
in $E$ comprises $P$-orbits of $p^1_{12}$ and $p^{\bar{2}12}_{12}=R_2^{-1}p^1_{12}$. 
$$
  \begin{array}{|l|l|l|l|}
    \hline
    & \hbox{Vertex} & \hbox{Image} & \hbox{Word} \\
    \hline
    R_1 & p^1_{12} & p^1_{12} & R_1 \\
    R_1 & p^3_{13}=P^{-1}p^1_{12} & p^{13\bar{1}}_{13}=P^{-1}R_2^{-1}p^1_{12} 
    & R_2PR_1P^{-1}=(R_1R_2)^2R_1^{-2} \\
    R_1 & p^2_{23}=P^{-2}R_2^{-1}p^1_{12} & p^{12\bar{1}}_{123\bar{1}}=Pp^1_{12} 
    & P^{-1}R_1P^{-2}R_2^{-1}=R_1(R_1R_2)^{-2} \\
    R_1 & p^{\bar{3}23}_{23}=P^{-2}p^1_{12} 
    & p^{1\bar{3}23\bar{1}}_{123\bar{1}}=PR_2^{-1}p^1_{12} 
    & R_2P^{-1}R_1P^{-2}=R_1^{-1} \\
    S_1 & p^{\bar{3}23}_{23}=P^{-2}p^1_{12} 
    & p^{23\bar{2}}_{\bar{3}\,\bar{2}3123}=P^3R_2^{-1}p^1_{12} 
    & R_2P^{-3}S_1P^{-2}=id \\
    S_1 & p^{23\bar{2}}_{\bar{3}\,\bar{2}3123}=P^3R_2^{-1}p^1_{12} 
    & p^{1\bar{3}23\bar{1}}_{123\bar{1}}=PR_2^{-1}p^1_{12} & R_2P^{-1}S_1P^3R_2^{-1}
    =R_1^{-1}\\
    \hline
  \end{array}
$$
Finding the stabilizer of $p^1_{13}$ is done in a similar way.

The orbit and stabilizer of the edge $(p^1_{13},p^1_{12})$ follow as before. 
A similar argument shows that no element of the group interchanges 
$p^1_{12}$ and $p^2_{12}$. Therefore the stabilizer of $(p^1_{12},p^2_{12})$ is 
the intersection of the vertex stabilizers. 
\EPf

Note that, since the stabilizer of $p^1_{12}$ is generated by the 
commuting complex reflections $R_1$ and $(R_1R_2)^3$, its order is 
simply $2p^2/(p-4)$ which is the product of the orders of these reflections.

\begin{proposition}\label{prop:stab-q3456}
  When $p\leqslant 6$ there is a single orbit of $q_*$-vertices. A
  representative is $q_{123\bar{2}}$ with stabilizer $\langle
  R_1,R_2R_3R_2^{-1}\rangle$. There is a single orbit of $(q_*,q_*)$
  edges. A representative is $(q_{1\bar{3}23},q_{123\bar{2}})$ with
  stabilizer $\langle S_1\rangle$.
\end{proposition}

\Pf
This is similar to the constructions of the stabilizers of the $p_*$-vertices.
$$
  \begin{array}{|l|l|l|l|}
    \hline
    & \hbox{Vertex} & \hbox{Image} & \hbox{Word} \\
    \hline
    R_1 & q_{123\bar{2}} & q_{123\bar{2}} & R_1 \\
    R_1 & q_{1\bar{3}23}=P^{-2}q_{123\bar{2}} 
    & q_{1\bar{3}23}=P^{-2}q_{123\bar{2}} 
    & P^2R_1P^{-2}=R_1(R_2R_3R_2^{-1})R_1^{-1} \\
    S_1 & q_{123\bar{2}} & q_{1\bar{3}23}=P^{-2}q_{123\bar{2}} 
    & P^2S_1=(R_2R_3R_2^{-1})^{-1} \\
    S_1 & q_{1\bar{3}23}=P^{-2}q_{123\bar{2}} & q_{123\bar{2}}
    & S_1P^{-2}=R_1(R_2R_3R_2^{-1}) \\
    \hline
  \end{array}
$$
There is one $P$-orbit of $(q_*,q_*)$ edges. Since $S_1$ interchanges 
$q_{1\bar{3}23}$ and $q_{123\bar{2}}$ we see that the intersection of 
the vertex stabilizers, namely $\langle R_1\rangle$, has index two in the 
stabilizer of the edge $(q_{1\bar{3}23},q_{123\bar{2}})$, which is 
$\langle S_1\rangle$. This is clear since $R_1=S_1^2$.
\EPf

When $p\leqslant 5$, the stabilizer of $q_{123\bar{2}}$ has order
$24p^2/(6-p)$ by Proposition~\ref{prop:stab-1232i}. When $p=6$
the point $q_{123\bar{2}}$ lies on $\partial\ch 2$. Once again
its stabilizer agrees with that conjectured in \cite{dpp1} and we have:

\begin{corollary}\label{cor:6-horo}
  When $p=6$, the (local) stabilizer in $\Gamma$ of
  $q_{123\bar{2}}\in\partial{\bf H}^2_\C$ contains no loxodromic maps.
  In particular, there is a consistent system of horoballs.
\end{corollary}

A similar argument gives the following result.

\begin{proposition}\label{prop:stab-q812}
  When $p=8$, $12$ there is a single orbit of $q_*$-vertices. A
  representative is $q^1_{123\bar{2}}$ with stabilizer $\langle
  R_1,(R_1R_2R_3R_2^{-1})^3\rangle$.  There are two orbits of
  $(q_*,q_*)$ edges. Representatives are
  $(q^1_{1\bar{3}23},q^1_{123\bar{2}})$ and
  $(q^1_{123\bar{2}},q^{23\bar{2}}_{123\bar{2}})$ with stabilizers
  $\langle S_1\rangle$ and $\langle R_1R_2R_3R_2^{-1}R_1\rangle$
  respectively.
\end{proposition}

Observe that the stabilizer of $q^1_{123\bar{2}}$ is generated by the
commuting complex reflections $R_1$ and $(R_1R_2R_3R_2^{-1})^3$.
Its order is $2p^2/(p-6)$, the product of the orders
of these reflections.

\begin{proposition}\label{prop:stab-pq}
  For all values of $p$ there are two orbits of $(p_*,q_*)$
  edges. Representatives are $(p^1_{12},q^1_{123\bar{2}})$ and
  $(p^1_{13},q^1_{1\bar{3}23})$, both with stabilizer $\langle
  R_1\rangle$. (For low values of $p$ omit the superscript $1$ as
  appropriate.)
\end{proposition}

\Pf
Since the vertices of such an edge are in different orbits, the stabilizer of an
the edge must be the intersection of the vertex stabilizers. The only thing to
check is that these two edges are in distinct orbits. This may be checked 
easily as before.
\EPf

\begin{figure}[htbp]
  \centering
  \subfigure[Link of $p_{12}$ ($p=4$)]{
    \epsfig{figure=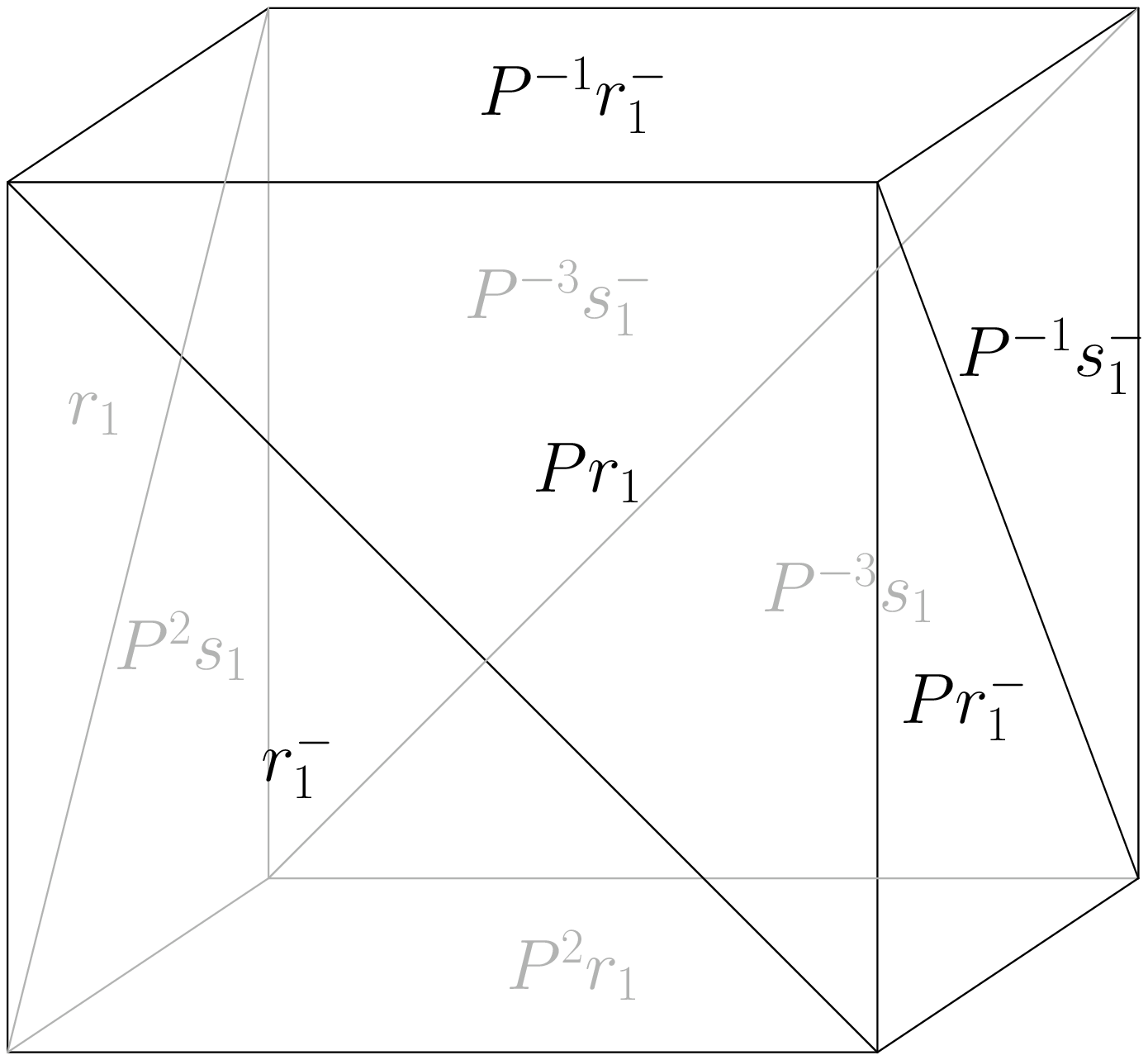, width=0.3\textwidth}}
  \hspace{0.1\textwidth}
  \subfigure[Link $q_{123\bar2}$ ($p=6$)]{
    \epsfig{figure=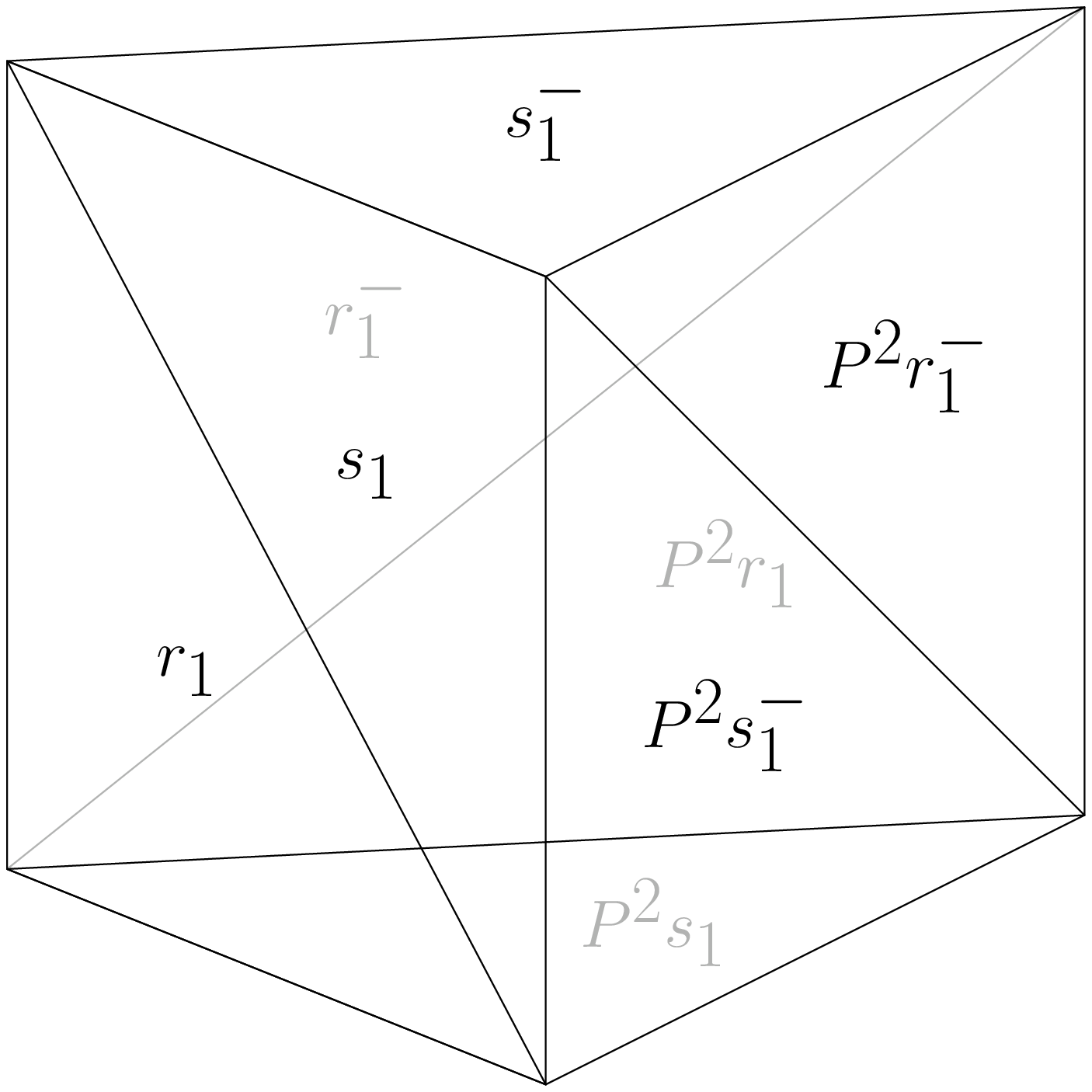, width=0.3\textwidth}}
  \hfill 
  \caption{Combinatorics of the links of ideal vertices.}\label{fig:cusplinks}
\end{figure}

For completeness, we end this section with a description of the links
of ideal vertices, which describe the structure of cusp
neighborhoods (see Figure~\ref{fig:cusplinks}).  Note that these are compact, so even for $p=4$ and
$p=6$, our polyhedron $E$ has finite volume (see
Section~\ref{sec:poin}).

\subsection{Conclusion - presentation for $\Gamma$}\label{sec:pres}

In Section~\ref{sec:tess} we verified that images of $E$ under 
the cosets of $\Upsilon$ in $\Gamma$ satisfied the local
tessellation hypotheses of the Poincar\'e polyhedron theorem.
In Corollaries~\ref{cor:4-horo} and~\ref{cor:6-horo} above, we showed
that when $p=4$ and $p=6$ there are consistent horoballs at the cusps. 
Note that $E$ has no cusps for the other values of $p$. Thus we have verified
the hypotheses of the Poincar\'e polyhedron theorem. Hence we have proved that
$\Gamma$ is discrete and $E$ is a fundamental polyhedron for the cosets
of $\Upsilon$ in $\Gamma$. 

We may read off a presentation for $\Gamma$ directly from the Poincar\'e 
polyhedron theorem; compare Section~20 of~\cite{mostowpacific}.
The generators of $\Gamma$ are the side pairing maps $R_1$ and $S_1$ of $E$ 
together with the generator $P$ of $\Upsilon$. 
The relations come from the cycle transformations
associated to the ridges and the relation from $\Upsilon$ 
(there are no reflection relations in this case). 
For all values of $p$ we obtain the cycle transformations from Table~\ref{tab:cycles3-4}:
$$
  id=R_1^p=(P^{-1}R_1)^3=P^2S_1^{-1}P^2R_1P^{-2}R_1=S_1^{-2}R_1=(P^2S_1)^p.
$$
In addition, for $p=5$, $6$, $8$ and $12$, from Table~\ref{tab:cycles-higher}, we obtain
$id=(P^{-3}R_1)^{4p/(p-4)}$
and for $p=8$, $12$ we obtain
$id=(P^{-2}S_1)^{6p/(p-6)}$.
Finally, we obtain the relation $id=P^7$ from the stabilizer 
of $\langle P\rangle$ of $E$. 
Thus, the Poincar\'e theorem gives the presentation:
\begin{equation}\label{eq:pres-poin}
  \left\langle R_1,\,S_1,\,P\ \Big\vert\ 
  \begin{array}{c}
    id=R_1^p=P^7=(P^{-1}R_1)^3=P^2S_1^{-1}P^2R_1P^{-2}R_1\\
    id =S_1^{-2}R_1=(P^2S_1)^p
    =(P^{-3}R_1)^{4p/(p-4)}=(P^{-2}S_1)^{6p/(p-6)}
  \end{array}
  \right\rangle
\end{equation}
where we omit the last two relations when the exponent is infinite or
negative.  Note that for the values of $p$ where these exponents are
negative the two relations still hold and follow from the other
relations (compare Corollary~5.6 of \cite{parkerlivne} and Section~2.2 of
\cite{mostowpacific}). 

We now simplify the relations and change generators in order to
recover the presentation of Theorem~\ref{thm:pres}. In particular, we
recover the relations found in Section~\ref{sec:braiding} only using
the cycle relations.

The relation $id=P^2S_1^{-1}P^2R_1P^{-2}R_1$ clearly implies that 
$S_1=P^2R_1P^{-2}R_1P^2$; see \eqref{eq:S1}.
Also, Proposition~\ref{prop:betterS1} implies that 
$P^2S_1$ is conjugate to $R_1^{-1}$. Hence the relation $id=(P^2S_1)^p$ 
is equivalent to $id=R_1^p$, and so may be omitted. We now show that the 
relation $S_1^2=R_1$, which also implies $S_1R_1=R_1S_1$, recovers both 
the generalized braid relation $(R_1R_2)^2=(R_2R_1)^2$ and
the braid relation $R_1(R_2R_3R_2^{-1})R_1=(R_2R_3R_2^{-1})R_1(R_2R_3R_2^{-1})$,
which were discussed in Section~\ref{sec:braiding}.

\begin{lemma}\label{lem:poin-braid}
  Assume that $P$ has order $7$ and that $J=P^{-1}R_1$ has order 3.
  Suppose that $S_1$, $R_2$ and $R_3$ are defined by
  $$
    S_1=P^2R_1P^{-2}R_1P^2,\quad 
    R_2=JR_1J^{-1}=R_1^{-1}PR_1P^{-1}R_1, \quad R_3=J^{-1}R_1J=P^{-1}R_1P.
  $$ 
  Then, the relation $S_1^2=R_1$ is equivalent to
  $(R_1R_2)^2=(R_2R_1)^2$.  Also, the relation $S_1R_1=R_1S_1$ is
  equivalent to
  $R_1(R_2R_3R_2^{-1})R_1=(R_2R_3R_2^{-1})R_1(R_2R_3R_2^{-1})$.
\end{lemma}

\Pf
Using $P^4=P^{-3}$, $R_1P^{-1}R_1=PR_1^{-1}P$ and $PR_1P^{-1}=R_1R_2R_1^{-1}$, we have:
\begin{eqnarray*}
  S_1^2R_1^{-1} 
  & = & P^2R_1P^{-1}(P^{-1}R_1P^{-1})P^{-1}(P^{-1}R_1P^{-1})
  P^{-1}R_1P(PR_1^{-1}) \\
  & = & P(PR_1P^{-1})
  (R_1^{-1}PR_1^{-1}P^{-1})(R_1^{-1}PR_1^{-1}P^{-1})
  R_1(PR_1P^{-1})R_1P^{-1} \\
  & = & P\bigl(R_1R_2(R_2R_1)^{-2}R_1R_2\bigr)P^{-1}.
\end{eqnarray*}
Secondly, one easily checks $R_2R_3R_2^{-1}=R_1^{-1}P^2R_1P^{-2}R_1$. Therefore
\begin{eqnarray*}
  S_1R_1^{-1}S_1^{-1}R_1
  & = & P^2R_1P^{-2}R_1P^2R_1^{-1}P^{-2}R_1^{-1}P^2R_1^{-1}P^{-2}R_1 \\
  & = & R_1(R_2R_3R_2^{-1})R_1(R_2R_3^{-1}R_2^{-1})R_1^{-1}(R_2R_3^{-1}R_2^{-1}).
\end{eqnarray*}
\EPf

Substituting for $P^{-1}R_1=J$, $R_2=JR_1J^{-1}=R_1^{-1}PR_1P^{-1}R_1$ 
and $R_3=P^{-1}R_1P=J^{-1}R_1J$, we obtain the presentation given in 
Theorem~\ref{thm:pres}, which agrees with the conjectural presentation 
given in \cite{dpp1}.

\subsection{Conclusion - orbifold Euler characteristics}\label{sec:euler}

In this section we prove Theorem~\ref{thm:euler} which gives the orbifold Euler 
characteristic of the quotient of $\ch 2$ by $\Gamma$. 
We write $\chi(2\pi/p,\overline{\sigma}_4)$ for the Euler characteristic
$\chi\bigl(\Gamma(2\pi/p,\overline{\sigma}_4)\backslash{\bf H}^2_\C\bigr)$
where $\tau=\overline{\sigma}_4=-(1+i\sqrt{7})/2$.

We begin by calculating the orbifold Euler characteristic when $p=3$
and $p=4$.  We choose one representative from each $\Gamma$-orbit of
facets and give its stabilizer.  The stabilizers of vertices and
edges, and their orders, were found in
Section~\ref{sec:stabilisers}. The orders of some vertex stabilizers
had also been given in Section~\ref{sec:braiding}.  The ridge cycles
described in Section~\ref{sec:tess} give the stabilizers of the
ridges. A representative of each orbit of facets and its stabilizer,
together with the order of the stabilizer, is given in
Table~\ref{tab:euler3-4}.
\begin{table}
  $$
  \begin{array}{|l|l|l|l|}
    \hline
    & \hbox{Facet} & \hbox{Stabilizer} & \hbox{Order}
    \\
    \hline
    \hbox{Vertices} 
    & p_{12} & \langle R_1,R_2 \rangle & 8p^2/(4-p)^2
    \\
    & q_{123\bar{2}} & \langle R_1,R_2R_3R_2^{-1} \rangle & 24p^2/(6-p)^2
    \\
    \hline
    \hbox{Edges} 
    & (p_{12},p_{13}) & \langle R_1\rangle & p \\
    & (p_{12},q_{231\bar{2}}) & \langle R_1\rangle & p \\
    & (p_{13},q_{1\bar{3}23}) & \langle R_1\rangle & p \\
    & (q_{123\bar{2}},q_{1\bar{3}23}) & \langle S_1\rangle & 2p \\
    \hline
    \hbox{Ridges} 
    & r_1\cap r_1^- & \langle R_1\rangle & p \\
    & r_1\cap P^{-1}r_1^- & \langle P^{-1}R_1 \rangle & 3 \\
    & r_1\cap P^2s_1 & id & 1 \\
    & r_1\cap s_1 & id & 1 \\
    & s_1\cap P^2s_1^- & \langle R_2R_3R_2^{-1} \rangle & p \\
    \hline
    \hbox{Sides} 
    & r_1 & id & 1 \\
    & s_1 & id & 1 \\
    \hline
    \hbox{Polyhedron} & E & \langle P\rangle & 7 \\
    \hline
\end{array}
$$
\caption{The Euler characteristic calculation for $p=3,\,4$.}
\label{tab:euler3-4}
\end{table}
Using these values, the orbifold Euler characteristic is
\begin{eqnarray*}
  \chi(2\pi/p,\overline{\sigma}_4) 
  & = & \left(\frac{(4-p)^2}{8p^2}+\frac{(6-p)^2}{24p^2}\right)
  -\left(\frac{3}{p}+\frac{1}{2p}\right)
  +\left(\frac{2}{p}+\frac{1}{3}+2\right)-2+\frac{1}{7} \\
  & = &  \frac{49-42p+9p^2}{14p^2} 
\end{eqnarray*}
Putting in $p=3$ and $p=4$ we obtain $\chi(2\pi/3,\overline{\sigma}_4)=2/63$ 
and $\chi(2\pi/4,\overline{\sigma}_4)=25/224$. 

We now do the same thing for $p=5$ and $p=6$. The main difference is that the
vertex $p_{12}$ has become a complex ridge and we introduce several more vertices.
\begin{table}
$$
  \begin{array}{|l|l|l|l|l|}
    \hline
    & \hbox{Facet} & \hbox{Stabilizer} & \hbox{Order} \\
    \hline
    \hbox{Vertices}
    & p^1_{12} & \langle R_1, (R_1R_2)^2 \rangle & 2p^2/(p-4) \\
    & p^1_{13} & \langle R_1, (R_1R_3)^2 \rangle & 2p^2/(p-4) \\
    & q_{123\bar{1}} & \langle R_1,R_2R_3R_2^{-1}\rangle & 24p^2/(6-p)^2
    \\
    \hline
    \hbox{Edges}
    & (p^1_{12},p^1_{13}) & \langle R_1 \rangle & p \\
    & (p^1_{12},p^2_{12}) & \langle (R_1R_2)^2 \rangle & 2p/(p-4) \\
    & (p^1_{12},q_{123\bar{2}}) &  \langle R_1\rangle & p \\
    & (p^1_{13},q_{1\bar{3}23}) &  \langle R_1\rangle & p \\
    & (q_{123\bar{2}},q_{1\bar{3}23}) & \langle S_1\rangle & 2p \\
    \hline
    \hbox{Ridges}
    & r_1\cap r_1^- & \langle R_1\rangle & p \\
    & r_1\cap P^{-1}r_1^- & \langle P^{-1}R_1 \rangle & 3 \\
    & r_1\cap P^2s_1 & id & 1 \\
    & r_1\cap s_1 & id & 1 \\
    & s_1\cap P^2s_1^- & \langle R_2R_3R_2^{-1}\rangle & p \\ 
    & r_1\cap P^{-3}r_1^- &\langle R_2R_3\rangle & 4p/(p-4) \\
    \hline
    \hbox{Sides}
    & r_1 & id & 1 \\
    & s_1 & id & 1 \\
    \hline
    \hbox{Polyhedron} & E & \langle P\rangle & 7 \\
    \hline
\end{array}
$$
\caption{The Euler characteristic calculation for $p=5,\,6$.}
\label{tab:euler5-6}
\end{table}
The facets, stabilizers and orders in Table~\ref{tab:euler5-6} imply that
the orbifold Euler characteristic is
\begin{eqnarray*}
  \chi(2\pi/p,\overline{\sigma}_4) 
  & = & \left(\frac{2(p-4)}{2p^2}+\frac{(6-p)^2}{24p^2}\right)
  -\left(\frac{3}{p}+\frac{p-4}{2p}+\frac{1}{2p}\right) \\
  && \quad +\left(\frac{2}{p}+\frac{1}{3}+2+\frac{p-4}{4p}\right)-2+\frac{1}{7} \\
  & = & \frac{15p^2-140}{56p^2}.
\end{eqnarray*}
Putting in $p=5$ and $p=6$ we obtain $\chi(2\pi/5,\overline{\sigma}_4)=47/280$ 
and $\chi(2\pi/6,\overline{\sigma}_4)=25/126$. 

Finally, we consider the cases $p=8$ and $p=12$. The facets,
stabilizers and orders are given in Table~\ref{tab:euler8-12},
and show that the orbifold Euler characteristic is
\begin{table}
  $$
  \begin{array}{|l|l|l|l|}
    \hline
    & \hbox{Orbit} & \hbox{Stabilizer} & \hbox{Order} \\
    \hline
    \hbox{Vertices}
    & p^1_{12} & \langle R_1, (R_1R_2)^2 \rangle & 2p^2/(p-4) \\
    & p^1_{13} & \langle R_1, (R_1R_3)^2 \rangle & 2p^2/(p-4) \\
    & q^1_{123\bar{2}} & \langle R_1,(R_1R_2R_3R_2^{-1})^3\rangle & 2p^2/(p-6) \\
    \hline
    \hbox{Edges}
    & (p^1_{12},p^1_{13}) & \langle R_1 \rangle & p \\
    & (p^1_{12},p^2_{12}) & \langle (R_1R_2)^2 \rangle & 2p/(p-4) \\
    & (p^1_{12},q^1_{123\bar{2}}) &  \langle R_1\rangle & p \\
    & (p^1_{13},q^1_{1\bar{3}23}) &  \langle R_1\rangle & p \\
    & (q^1_{123\bar{2}},q^1_{1\bar{3}23}) & \langle S_1\rangle & 2p \\
    & (q^1_{123\bar{2}},q^{23\bar{2}}_{123\bar{2}}) 
    & \langle R_1R_2R_3R_2^{-1}R_1\rangle & 4p/(p-6) \\
    \hline
    \hbox{Ridges}
    & r_1\cap r_1^- & \langle R_1\rangle & p \\
    & r_1\cap P^{-1}r_1^- & \langle P^{-1}R_1 \rangle & 3 \\
    & r_1\cap P^2s_1 & id & 1 \\
    & r_1\cap s_1 & id & 1 \\
    & s_1\cap P^2s_1^- & \langle R_2R_3R_2^{-1}\rangle & p \\ 
    & r_1\cap P^{-3}r_1^- &\langle R_2R_3\rangle & 4p/(p-4) \\
    & s_1\cap P^{-2}s_1^- & \langle R_1R_3^{-1}R_2R_3 \rangle & 6p/(p-6) \\ 
    \hline
    \hbox{Sides}
    & r_1 & id & 1 \\
    & s_1 & id & 1 \\
    \hline
    \hbox{Polyhedron} & E & \langle P\rangle & 7 \\
    \hline
\end{array}
$$
\caption{The Euler characteristic calculation for $p=8,\,12$.}
\label{tab:euler8-12}
\end{table}
\begin{eqnarray*}
  \chi(2\pi/p,\overline{\sigma}_4) 
  & = & \left(\frac{2(p-4)}{2p^2}+\frac{(6-p)^2}{24p^2}\right)
  -\left(\frac{3}{p}+\frac{p-4}{2p}+\frac{1}{2p}+\frac{p-6}{4p}\right) \\
  && \quad 
  +\left(\frac{2}{p}+\frac{1}{3}+2+\frac{p-4}{4p}+\frac{p-6}{6p}\right)
  -2+\frac{1}{7} \\
  & = & \frac{-98+21p+2p^2}{14p^2}.
\end{eqnarray*}
Putting in $p=8$ and $p=12$ gives $\chi(2\pi/8,\overline{\sigma}_4)=99/448$ and
$\chi(2\pi/12,\overline{\sigma}_4)=221/1008$.

\section{Arithmeticity and Galois theory}\label{sec:arithmeticgalois}

\subsection{Trace fields and commensurability}\label{sec:fields}

Recall from~\cite{parkerpaupert} that we can adjust
the generators $R_1$, $J$ and the Hermitian form $H$ 
so that their entries lie in
$$
  \l=\Q(\tau,e^{2\pi i/p}).
$$
Note that after making this adjustment, $R_1$ no longer has unit determinant.

In this section, we consider the adjoint trace field
$\k=\mathbb{Q}(\rm{tr}\,\rm{Ad}\,\Gamma)$ for the relevant
sporadic triangle groups, i.e. the ones with
$\tau=\overline{\sigma}_4=-(1+i\sqrt{7})/2$, and $p=3$, $4$, $5$, $6$
, $8$ or $12$.  Recall from \cite{delignemostow}, \cite{mostowpacific}
and \cite{paupert} that:

\begin{lemma}
  The field $\mathbb{Q}(\rm{tr}\,\rm{Ad}\,\Gamma)$ is a
  commensurability invariant of subgroups $\Gamma < {\rm U}(n,1)$,
  where ${\rm tr}\, {\rm Ad} (\gamma)=|{\rm tr} (\gamma)|^2$ for
  $\gamma \in {\rm U}(n,1)$.
\end{lemma}

Since we have arranged for $R_1$ and $J$ (as well as the matrix of the
Hermitian form $H$) to have entries in $\l=\Q(\tau,e^{2\pi i/p})$, the
adjoint trace field is contained in $\R\cap \Q(\tau,e^{2\pi i/p})$.

In fact it is easy to see by computing the trace of a handful of group
elements that $\k$ contains
$\mathbb{Q}(\cos\frac{2\pi}{p},\sqrt{7}\sin\frac{2\pi}{p})$.  Note
that the adjustment indicated above changes the trace of elements of
$\Gamma$, but the absolute values of these traces remain the same. In
fact the following two traces suffice:
\begin{eqnarray*}
  \bigl|{\rm tr}(R_1)\bigr|^2
  & = & |e^{4i\pi/3p}+2e^{-2i\pi/3p}|^2 \\
  & = & 5+4\cos(2\pi/p), \\
  \bigl|{\rm tr}(JR_1JR_1^{-1})\bigr|^2
  & = & |e^{-2\pi i/3p}-\overline{\tau}e^{4i\pi/3p}|^2 \\
  & = & 3-2\Re(\tau e^{2\pi i/p}) \\
  & = & 3+\cos(2\pi/p)+\sqrt{7}\sin(2\pi/p).
\end{eqnarray*}

Now clearly we can go from $\k$ to $\l$ simply by adjoining
$\sqrt{-7}$, so $\l$ has degree $2$ over $\k$. This shows that $\k$ is
equal to $\R\cap \l$, and
$$ 
  \k=\mathbb{Q}\bigl(\cos(2\pi/p),\sqrt{7}\sin(2\pi/p)\bigr).
$$
More specifically, for the relevant values of $p$, we find the fields
given in Table~\ref{tab:fields}. This shows that the corresponding six
lattices are in different commensurability classes; indeed there is
only one pair having the same adjoint trace field, but
$\Gamma(2\pi/3,\overline{\sigma}_4)$ is cocompact whereas
$\Gamma(2\pi/6,\overline{\sigma}_4)$ is not.
\begin{table}
$$
 \begin{array}{|l|l|l|}
  \hline
    p & \k & \hbox{Degree} \\ 
    \hline
    3 & \Q(\sqrt{21}) & 2\\
    4 & \Q(\sqrt{7}) & 2\\
    5 & \Q(\sqrt{14}\sqrt{5+\sqrt{5}}) & 4\\
    6 & \Q(\sqrt{21}) & 2\\
    8 & \Q(\sqrt{2},\sqrt{7}) & 4\\
    12 & \Q(\sqrt{3},\sqrt{7}) & 4\\
    \hline
  \end{array}
  $$
\caption{The list of the adjoint trace fields of the lattices
  $\Gamma(2\pi/p,\overline{\sigma}_4)$ shows that these $6$ lattices
  lie in different commensurability classes.} \label{tab:fields}
\end{table}
For completeness, we also recall the following result
from~\cite{paupert}.
\begin{theorem}\label{thm:paupert}
  None of the six lattices $\Gamma(2\pi/p,\overline{\sigma}_4)$ is
  commensurable to any Deligne-Mostow lattice.
\end{theorem}
Indeed, the adjoint trace fields of the lattices
in~\cite{delignemostow} are all of the form $\R\cap
\Q(\zeta_n)=\Q(\cos\frac{2\pi}{n})$, and only certain values of $n$
actually appear in that list. More specifically, each Deligne-Mostow
group in ${\rm PU}(2,1)$ is described by a $5$-tuple of rational
numbers, and $n$ is simply the least common denominator of these
rational numbers (see~\cite{delignemostow}, Section~12). Note that if
$n=2k$ where $k$ is odd, then $\Q(\zeta_n)=\Q(\zeta_k)$.

Given the degrees that appear in Table~\ref{tab:fields}, it is enough
to check the cyclotomic fields of degree $4$ or $8$. The relevant
values of $n$ are listed in Table~\ref{tab:fieldsDM}.  Since all
groups except the one with $p=3$ are non-arithmetic, we only look at
the non-arithmetic Deligne-Mostow groups in ${\rm PU}(2,1)$.  We list
the corresponding adjoint trace fields in Table~\ref{tab:fieldsDM}.
This proves Theorem~\ref{thm:paupert}.
\begin{table}[htbp]
 $$
  \begin{array}{|l|l|l|}
  \hline
    n & \Q(\zeta_n)\cap\R & \hbox{Degree}\\ 
    \hline  
    5 \hbox{ or }10 & \Q(\sqrt{5}) & 2\\
    8 & \Q(\sqrt{2}) & 2\\
    12 & \Q(\sqrt{3}) & 2\\ 
    \hline
    15\hbox{ or }30 & \Q(\sqrt{3}\sqrt{\frac{5+\sqrt{5}}{2}}) & 4\\
    16 & \Q(\sqrt{2+\sqrt{2}}) & 4\\
    20 & \Q(\sqrt{\frac{5+\sqrt{5}}{2}}) & 4\\
    24 & \Q(\sqrt{2},\sqrt{3}) & 4\\
    \hline
  \end{array}
  $$
  \caption{The list of cyclotomic fields whose totally real field has
    degree $2$ or $4$.} \label{tab:fieldsDM}
\end{table}

\subsection{Arithmeticity} \label{sec:arithmetic}

In this section we apply the Mostow--Vinberg arithmeticity criterion
(Proposition~\ref{prop:crit}) to study arithmeticity. The result
(Proposition~\ref{prop:arithmetic}) follows from Theorem~4.1 of
\cite{paupert}, but for completeness we provide details for this
special case.  In the coordinates used in the previous section $H$ is
the following
$$
H=\left[
  \begin{array}{ccc}
    \alpha & \beta & \gamma\\
    \overline{\beta} & \alpha & \beta\\
    \overline{\gamma} & \overline{\beta} & \alpha
  \end{array}\right]
$$
where
$\alpha=2-e^{2\pi i/p}-e^{-2\pi i/p}$, $\beta=(e^{-2\pi i/p}-1)\tau$ and 
$\gamma=(1-e^{-2\pi i/p})\overline{\tau}$.
In particular, all entries of $H$ lie in $\l=\Q(\tau,e^{2\pi i/p})$.
Observe also that
$\tau=\overline{\zeta}_7+\overline{\zeta}_7^2+\overline{\zeta}_7^4$
for $\zeta_7=e^{2\pi i/7}$, so that $\l \subset
\Q(\zeta_7,\zeta_p)=\Q(\zeta_{7p})$, where the last equality follows
from the fact that $7$ is coprime to $p$.

Let $k\in\N$. By the Chinese remainder theorem, there exists $l\in\N$
such that $l\equiv 1 (\textrm{mod } 7)$ and $l\equiv k (\textrm{mod }
p)$.  In particular, if $k$ is coprime to $p$,
$\zeta_{7p}\mapsto\zeta_{7p}^l$ induces an automorphism of
$\Q(\zeta_{7p})$ which fixes $\zeta_7$ and sends $\zeta_p$ to
$\zeta_p^k$; hence $\l=\Q(\tau,e^{2\pi i/p})$ has an automorphism
fixing $\tau$ and sending $\zeta_p$ to $\zeta_p^k$. Similarly, there
is an automorphism of $\l$ fixing $\zeta_p$ and sending $\tau$ to
$\overline{\tau}$.

Note that all these automorphisms clearly commute with complex
conjugation, so that for any $\sigma\in \textrm{Gal}(\l/\Q)$,
$H^\sigma$ is a Hermitian form. Although $H$ has signature $(2,1)$,
its Galois conjugates $H^\sigma$ may have various signatures for
different $\sigma$'s. With the current notation,
Proposition~\ref{prop:crit} translates into the following.
\begin{proposition}
  The lattice $\Gamma(2\pi/p,\tau)$ is arithmetic if and only if for
  any $\sigma\in \textrm{Gal}(\l/\Q)$ not inducing the identity on
  $\k=\R\cap \l$, the form $H^\sigma$ is definite.
\end{proposition}

Since the diagonal elements of $H^\sigma$ are all positive, it is
indefinite (i.e it has signature $(1,2)$ or $(2,1)$) precisely when
$\det(H^\sigma)<0$. We list values of $\det(H)$ and all $\sigma$ with
$\det(H^\sigma)<0$ in Table~\ref{tab:conjugates} (for $p=5$ there are
two non-trivial Galois automorphisms for which $\det(H^\sigma)<0$, we
list both of them). It is quite clear from Table~\ref{tab:fields}
that these automorphisms act non-trivially on the adjoint trace field
$\k$.

\begin{table}[htbp]
  $$
  \begin{array}{| l | l | l | l | l |}
    \hline
    p & \tau^\sigma,\zeta_p^\sigma & \Gamma^\sigma & \det(H)                   & \det(H^\sigma)\\
    \hline
    3 &                              &  & -\frac{1}{2}(9+3\sqrt{21})  &      \\[4pt]
    4 &  \overline{\tau}, \ \zeta_4      & \Gamma(2\pi/4,\sigma_4) & -6-2\sqrt{7}            & -6+2\sqrt{7} \\[4pt]
    5 &  \overline{\tau},\  \zeta_5      & \Gamma(2\pi/5,\sigma_4) & 
      -\frac{\sqrt{5}}{4}\bigl(5+\sqrt{5}+\sqrt{14}\sqrt{5-\sqrt{5}}\bigr) &
      -\frac{\sqrt{5}}{4}\bigl(5+\sqrt{5}-\sqrt{14}\sqrt{5-\sqrt{5}}\bigr) \\
      &  \tau,\  \zeta_5^2      & \Gamma(4\pi/5,\overline{\sigma}_4) & 
      & \frac{\sqrt{5}}{4}\bigl(5-\sqrt{5}+\sqrt{14}\sqrt{5+\sqrt{5}}\bigr) \\[4pt]
   6 &  \overline{\tau},\  \zeta_6      & \Gamma(2\pi/6,\sigma_4) &
     -\frac{1}{2}(5+\sqrt{21}) & -\frac{1}{2}(5-\sqrt{21})\\[4pt] 
   8 &  \tau , \  \zeta_8^3  & \Gamma(6\pi/8,\overline{\sigma}_4) &
     -1+\sqrt{7} - \sqrt{14}  & - 1-\sqrt{7} - \sqrt{14} \\[4pt]
  12 & \tau,\ \zeta_{12}^5 & \Gamma(10\pi/12,\overline{\sigma}_4) &
    \frac{3}{2} - \sqrt{3} - \sqrt{7} + \frac{1}{2}\sqrt{21} &
    \frac{3}{2} + \sqrt{3} - \sqrt{7} - \frac{1}{2}\sqrt{21} \\[4pt]
  \hline
\end{array}
$$
\caption{Computation of $\det(H)$ for various values of $p$; when
  there is one, we list an automorphism $\sigma$ such that $H^\sigma$
  is indefinite.}
  \label{tab:conjugates}
\end{table}

\begin{proposition}\label{prop:arithmetic}
  Consider the groups $\Gamma(2\pi/p,\overline{\sigma}_4)$,
  $p=3,4,5,6,8$ and $12$. Then $\Gamma(2\pi/p,\overline{\sigma}_4)$ is
  arithmetic if and only if $p=3$.
\end{proposition}

\section{Appendix: combinatorial data}\label{sec:comb-data}

In this section we gather together the data about the facets of $E$ and their
orbits under $P$.

\begin{table}[htbp]
\centering
\begin{tabular}{|c|}
\hline
$
  \begin{array}{cccccccccccc}
    m_1 & {\buildrel P\over\longmapsto} & 
    m_{12\bar{1}} & {\buildrel P\over\longmapsto} & 
    m_{123\bar{2}\bar{1}} & {\buildrel P\over\longmapsto} & 
    m_{1231\bar{3}\bar{2}\bar{1}} & {\buildrel P\over\longmapsto} &
    m_{\bar{3}\bar{2}123}
    \\
    \smallskip & & & {\buildrel P\over\longmapsto} & 
    m_{\bar{3}23} & {\buildrel P\over\longmapsto} & 
    m_3 & {\buildrel P\over\longmapsto} & m_1
  \end{array}
$
\\ \hline
$
  \begin{array}{cccccccccccc}
    m_2 & {\buildrel P\over\longmapsto} & 
    m_{13\bar{1}} & {\buildrel P\over\longmapsto} & 
    m_{\bar{2}12} & {\buildrel P\over\longmapsto} & 
    m_{1\bar{3}23\bar{1}} & {\buildrel P\over\longmapsto} &
    m_{\bar{3}12\bar{1}3}
    \\ \smallskip & & & {\buildrel P\over\longmapsto} & 
    m_{23\bar{2}} & {\buildrel P\over\longmapsto} & 
    m_{\bar{3}13} & {\buildrel P\over\longmapsto} & m_2
  \end{array}
$
\\ \hline
\end{tabular}
\caption{$P$-orbits of mirrors of $R_1$ and $R_2$.}\label{tab:mirrors}
\end{table}

Suppose that $u$ and $v$ are words so that $uv$ is conjugate to $12$.
If the mirrors of $m_u$ and $m_v$ intersect in $\chb 2$ then their
intersection point is denoted $p_w$. The relationship between $u$, $v$
and $w$ is slightly subtle. The easiest situation is when $uv=w$, such
as $m_1\cap m_2=p_{12}$.  However, any mirror in the group generated
by $u$ and $v$ also passes through this point. So, for example,
$m_1\cap m_{12\bar{1}}=m_{12\bar{1}}\cap m_{\bar{2}12}=p_{12}$.

\begin{table}[htbp]
\centering
\begin{tabular}{|c|}
\hline
$
  \begin{array}{ccccccccc}
    p_{12} & {\buildrel P\over\longmapsto} & 
    p_{123\bar{1}} & {\buildrel P\over\longmapsto} &
    p_{1231\bar{2}\bar{1}} & {\buildrel P\over\longmapsto} & 
    p_{\bar{3}\bar{2}3123} & {\buildrel P\over\longmapsto} & 
    p_{\bar{3}123}  \\
    & & & {\buildrel P\over\longmapsto} & p_{23} & {\buildrel P\over\longmapsto} &
    p_{13} & {\buildrel P\over\longmapsto} & p_{12}
  \end{array}
$ 
\\ \hline
\end{tabular}
\caption{$P$-orbits of $p_*$ vertices, for $p=3$ or $4$.}\label{tab:verts3-4}
\end{table}

If $m_u$ and $m_v$ are ultraparallel then we denote their common
perpendicular by $m_w$ where $u$, $v$ and $w$ are related as above. In
this case the $P$-orbit of $m_w$ may be found using the same
subscripts as appear in the $P$-orbit of $p_w$ in the
Table~\ref{tab:verts3-4}. In this case, $m_u$ and $m_w$ intersect in
$p^u_w$. The $P$-orbits of these points are given in
Table~\ref{tab:verts5-12}.

\begin{table}[htbp]
\centering
\begin{tabular}{|c|}
\hline
$
  \begin{array}{ccccccccc}
    p^1_{12} & {\buildrel P\over\longmapsto} & 
    p^{12\bar{1}}_{123\bar{1}} & {\buildrel P\over\longmapsto} &
    p^{123\bar{2}\,\bar{1}}_{1231\bar{2}\,\bar{1}} & 
    {\buildrel P\over\longmapsto} & 
    p^{1231\bar{3}\,\bar{2}\,\bar{1}}_{\bar{3}\,\bar{2}3123} & 
    {\buildrel P\over\longmapsto} & 
    p^{\bar{3}\,\bar{2}123}_{\bar{3}123} \\ \smallskip
    & & & {\buildrel P\over\longmapsto} & 
    p^{\bar{3}23}_{23} & {\buildrel P\over\longmapsto} &
    p^3_{13} & {\buildrel P\over\longmapsto} & p^1_{12} 
  \end{array}
$
\\ \hline
$
  \begin{array}{ccccccccc}
    p^2_{12} & {\buildrel P\over\longmapsto} & 
    p^{13\bar{1}}_{123\bar{1}} & {\buildrel P\over\longmapsto} &
    p^{\bar{2}12}_{1231\bar{2}\,\bar{1}} & {\buildrel P\over\longmapsto} & 
    p^{1\bar{3}23\bar{1}}_{\bar{3}\,\bar{2}3123} & 
    {\buildrel P\over\longmapsto} & 
    p^{\bar{3}12\bar{1}3}_{\bar{3}123} \\ \smallskip
    & & & {\buildrel P\over\longmapsto} & 
    p^{23\bar{2}}_{23} & {\buildrel P\over\longmapsto} &
    p^{\bar{3}13}_{13} & {\buildrel P\over\longmapsto} & p^2_{12}
  \end{array}
$
\\ \hline
$
  \begin{array}{ccccccccc}
    p^{\bar{2}12}_{12} & {\buildrel P\over\longmapsto} & 
    p^{1\bar{3}23\bar{1}}_{123\bar{1}} & {\buildrel P\over\longmapsto} &
    p^{\bar{3}12\bar{1}3}_{1231\bar{2}\,\bar{1}} & 
    {\buildrel P\over\longmapsto} & 
    p^{23\bar{2}}_{\bar{3}\,\bar{2}3123} & {\buildrel P\over\longmapsto} & 
    p^{\bar{3}13}_{\bar{3}123} 
    \\ \smallskip
    & & & {\buildrel P\over\longmapsto} & p^2_{23} & {\buildrel P\over\longmapsto} &
    p^{13\bar{1}}_{13} & {\buildrel P\over\longmapsto} & 
    p^{\bar{2}12}_{12}, 
  \end{array}
$
\\ \hline
$
  \begin{array}{ccccccccc}
    p^{12\bar{1}}_{12} & {\buildrel P\over\longmapsto} & 
    p^{123\bar{2}\,\bar{1}}_{123\bar{1}} & 
    {\buildrel P\over\longmapsto} &
    p^{1231\bar{3}\,\bar{2}\,\bar{1}}_{1231\bar{2}\,\bar{1}} & {\buildrel P\over\longmapsto} & 
    p^{\bar{3}\,\bar{2}123}_{\bar{3}\,\bar{2}3123} & 
    {\buildrel P\over\longmapsto} & 
    p^{\bar{3}23}_{\bar{3}123} \\ \smallskip
    & & & {\buildrel P\over\longmapsto} & p^3_{23} & {\buildrel P\over\longmapsto} &
    p^1_{13} & {\buildrel P\over\longmapsto} & p^{12\bar{1}}_{12}.
  \end{array}
$
\\ \hline
\end{tabular}
\caption{For $p=5$ or $6$, there are four $P$-orbits of $p_*$
  vertices.}\label{tab:verts5-12}
\end{table}

Similarly, suppose $u$ and $v$ are words in $\Gamma$ so that $uv$ is
conjugate to $123\bar{2}$. If $m_u$ and $m_v$ intersect in $\chb 2$
then we denote their intersection point by $q_w$, where $u$, $v$ and
$w$ where $w$ is conjugate in the group generated by $u$ and $v$ to
$uv$. These vertices are listed in Table~\ref{tab:verts3-6}.

\begin{table}[htbp]
\centering
\begin{tabular}{|c|}
\hline
$
  \begin{array}{ccccccccc}
    q_{123\bar{2}} & {\buildrel P\over\longmapsto} & 
    q_{12\bar{1}\bar{3}13} & {\buildrel P\over\longmapsto} &
    q_{2\bar{3}\bar{2}123} & {\buildrel P\over\longmapsto} & 
    q_{13\bar{1}\bar{3}23} & {\buildrel P\over\longmapsto} & 
    q_{\bar{2}123} \\
    & & & {\buildrel P\over\longmapsto} & q_{1\bar{3}23} 
    & {\buildrel P\over\longmapsto} & 
    q_{12\bar{1}3} & {\buildrel P\over\longmapsto} & q_{123\bar{2}}
  \end{array}
$
\\ \hline
\end{tabular}
\caption{$P$-orbits of $q_*$ vertices, for $3\leqslant p \leqslant 6$.}\label{tab:verts3-6}
\end{table}

If $m_u$ and $m_v$ are ultraparallel then we denote their common
perpendicular by $m_w$, where $u$, $v$, $w$ are related as before. The
intersection point of $m_u$ and $m_w$ is denoted $q^u_w$; see
Table~\ref{tab:verts8-12}.

\begin{table}[htbp]
\centering
\begin{tabular}{|c|}
\hline
$
  \begin{array}{cccccccccccc}
    q^1_{1\bar{3}23} 
    & {\buildrel P\over\longmapsto} & 
    q^{12\bar{1}}_{12\bar{1}3} & {\buildrel P\over\longmapsto} & 
    q^{123\bar{2}\,\bar{1}}_{123\bar{2}} &
    {\buildrel P\over\longmapsto} & 
    q^{1231\bar{3}\,\bar{2}\,\bar{1}}_{12\bar{1}\,\bar{3}13} 
    & {\buildrel P\over\longmapsto} &
    q^{\bar{3}\,\bar{2}123}_{2\bar{3}\,\bar{2}123} \\
    \smallskip
    & & & {\buildrel P\over\longmapsto} & 
    q^{\bar{3}23}_{13\bar{1}\,\bar{3}23} & 
    {\buildrel P\over\longmapsto} & 
    q^3_{\bar{2}123} & {\buildrel P\over\longmapsto} & q^1_{1\bar{3}23} 
  \end{array}
$
\\ \hline
$
  \begin{array}{cccccccccccc}
    q^{\bar{3}23}_{1\bar{3}23} 
    & {\buildrel P\over\longmapsto} & 
    q^3_{12\bar{1}3} & {\buildrel P\over\longmapsto} & 
    q^1_{123\bar{2}} & {\buildrel P\over\longmapsto} & 
    q^{12\bar{1}}_{12\bar{1}\,\bar{3}13} 
    & {\buildrel P\over\longmapsto} &
    q^{123\bar{2}\,\bar{1}}_{2\bar{3}\,\bar{2}123} 
    \\ \smallskip & & &
    {\buildrel P\over\longmapsto} & 
    q^{1231\bar{3}\,\bar{2}\,\bar{1}}_{13\bar{1}\,\bar{3}23} & 
    {\buildrel P\over\longmapsto} & 
    q^{\bar{3}\,\bar{2}123}_{\bar{2}123} &
    {\buildrel P\over\longmapsto} & q^{\bar{3}23}_{1\bar{3}23}
  \end{array}
$
\\ \hline
$
  \begin{array}{cccccccccccc}
    q^{1\bar{3}23\bar{1}}_{1\bar{3}23} 
    & {\buildrel P\over\longmapsto} & 
    q^{\bar{3}12\bar{1}3}_{12\bar{1}3} & 
    {\buildrel P\over\longmapsto} & 
    q^{23\bar{2}}_{123\bar{2}} & {\buildrel P\over\longmapsto} & 
    q^{\bar{3}13}_{12\bar{1}\,\bar{3}13} & 
    {\buildrel P\over\longmapsto} &
    q^2_{2\bar{3}\,\bar{2}123} \\
    \smallskip 
    & & & {\buildrel P\over\longmapsto} & 
    q^{13\bar{1}}_{13\bar{1}\,\bar{3}23} & 
    {\buildrel P\over\longmapsto} & q^{\bar{2}12}_{\bar{2}123} &
    {\buildrel P\over\longmapsto} & 
    q^{1\bar{3}23\bar{1}}_{1\bar{3}23} 
  \end{array}
$
\\ \hline
\end{tabular}
\caption{For $p=8$ or $12$, there are three $P$-orbits of $q_*$ vertices.}
 \label{tab:verts8-12}
\end{table}

\begin{flushleft}
  \textsc{Martin Deraux\\
  Institut Fourier, Universit\'e de Grenoble 1}\\
  \verb|deraux@ujf-grenoble.fr|
\end{flushleft}
\begin{flushleft}
  \textsc{John R. Parker\\
    Department of Mathematical Sciences, Durham University}\\
    \verb|j.r.parker@durham.ac.uk|
\end{flushleft}
\begin{flushleft}
  \textsc{Julien Paupert\\
   School of Mathematical and Statistical Sciences, Arizona State University}\\
       \verb|paupert@asu.edu|
\end{flushleft}

\end{document}